\documentclass[11pt]{amsart}

\usepackage{amsmath,amssymb,amsthm,mathtools,mathrsfs,enumerate,hyperref}
\usepackage{tikz}
\usetikzlibrary{arrows.meta,positioning}

\usepackage{hyperref}
\numberwithin{equation}{section}
\allowdisplaybreaks
\newtheorem{assumption}{Assumption}
\newtheorem{theorem}{Theorem}[section]
\newtheorem{proposition}[theorem]{Proposition}
\newtheorem{lemma}[theorem]{Lemma}
\newtheorem{corollary}[theorem]{Corollary}
\theoremstyle{definition}

\theoremstyle{remark}
\newtheorem{remark}[theorem]{Remark}
\newcommand{\R}{\mathbb R}
\newcommand{\C}{\mathbb C}
\newcommand{\Z}{\mathbb Z}

\newcommand{\na}{\nabla}
\newcommand{\cN}{\mathcal N}
\newcommand{\cP}{\mathcal P}
\newcommand{\cL}{\mathcal L}
\newcommand{\cA}{\mathcal A}
\newcommand{\cB}{\mathcal B}
\newcommand{\cE}{\mathcal E}
\newcommand{\cX}{\mathcal X}
\newcommand{\cH}{\mathcal H}

\newcommand{\cS}{\mathcal S}

\newcommand{\cR}{\mathcal R}

\newcommand{\bvec}{\mathbf b}

\newcommand{\app}{\mathrm{app}}
\newcommand{\spec}{\operatorname{spec}}
\newcommand{\diag}{\operatorname{diag}}
\newcommand{\diver}{\operatorname{div}}
\newcommand{\sym}{\operatorname{sym}}
\newcommand{\supp}{\operatorname{supp}}
\newcommand{\arsinh}{\operatorname{arsinh}}
\newcommand{\dist}{\operatorname{dist}}

\newcommand{\ee}{\mathrm e}
\newcommand{\av}[1]{\left\langle #1\right\rangle}
\newcommand{\norm}[1]{\left\|#1\right\|}

\newcommand{\Order}{\mathcal O}
\newcommand{\pair}[2]{\begin{pmatrix}#1\\#2\end{pmatrix}}

\hypersetup{pdftitle={Real-valued blow-up for the supercritical defocusing nonlinear wave equation},pdfauthor={Feng Shao, Dongyi Wei, Zhifei Zhang},pdfsubject={Defocusing nonlinear wave equation and analytic Euler profiles}}
\newcommand{\Hol}{\operatorname{Hol}}
\newcommand{\cQ}{\mathcal Q}

\begin{document}
\title[Blow-up for the supercritical defocusing NLW]{On blow-up for the supercritical defocusing nonlinear scalar wave equation}

\author[F. Shao]{Feng Shao}
\address[F. Shao]{Academy of Mathematics \& Systems Science and Hua Loo-Keng Center for Mathematical Sciences, Chinese Academy of Sciences, Beijing 100190, China}
\email{fshao@amss.ac.cn}

\author[D. Wei]{Dongyi Wei}
\address[D. Wei]{School of Mathematical Sciences, Peking University, Beijing 100871,  China}
\email{jnwdyi@pku.edu.cn}

\author[Z. Zhang]{Zhifei Zhang}
\address[Z. Zhang]{School of Mathematical Sciences, Peking University, Beijing 100871, China}
\email{zfzhang@math.pku.edu.cn}

\date{\today}

\begin{abstract}
	We study finite-time singularity formation for the scalar defocusing nonlinear wave equation $-\partial_t^2u+\Delta u=u^p$.  First, for fixed $d\geq4$ and odd $p\geq3$, we establish a conditional blow-up result under the assumption that a smooth self-similar imploding profile exists for the leading-order relativistic Euler system, with associated normalized coefficients satisfying quantitative analytic radius conditions. Second, we verify both the profile and the analytic radius conditions in dimensions four and five for all sufficiently large odd  $p$. Using finite speed of propagation, we obtain a dimension-independent threshold $p_0$ such that, for every $d\geq4$ and every odd $p\geq p_0$, finite-time blow-up occurs for some smooth, compactly supported, real-valued initial data. The proof combines a normalized formulation that separates the leading periodic oscillation from the unknown modulation, growth of the spatial analytic radius induced by an outgoing property, and frequency-dependent matrix multipliers that control the cumulative growth caused by the periodic linearized potential. 

\end{abstract}

\maketitle
\tableofcontents

\section{Introduction}\label{Sec.Introduction}

In this paper, we consider the real-valued defocusing nonlinear wave equation
\begin{equation}\label{Eq.wave_eq}
	\Box u=u^p,\qquad (t,y)\in\R\times\R^d,
\end{equation}
where $u:\R^{1+d}\to\R$ is the unknown field, $p\in2\Z_{\geq1}+1$, and $\Box:=\partial^\mu\partial_\mu=-\partial_t^2+\Delta_y$ is the d'Alembertian on Minkowski spacetime. We raise and lower indices using the metric $\diag(-1,1,\ldots,1)$ and sum repeated Greek indices from $0$ to $d$. Since $p$ is odd, the nonlinearity is smooth and $u^p=|u|^{p-1}u$.

Given $(u_0,u_1)\in C_c^\infty(\R^d;\R)^2$, the local Cauchy theory gives a unique smooth solution on the maximal existence interval $[0,T_*)$; see \cite{Luk,Sogge}. If $T_*<\infty$, the continuation criterion gives
\begin{equation}\label{Eq.continuation}
	\limsup_{t\uparrow T_*}\norm{u(t,\cdot)}_{L^\infty(\R^d)}=\infty.
\end{equation}
The solution has finite speed of propagation. In particular, compact support is preserved on each bounded time interval, and the conserved energy is
\[E[u](t):=\int_{\R^d}\left(\frac12|\partial_tu|^2+\frac12|\na_yu|^2+\frac1{p+1}|u|^{p+1}\right)\,\mathrm dy=E[u](0).\]
All three terms are nonnegative. Nevertheless, in the energy-supercritical regime the conservation of $E[u]$ does not provide an a priori bound for the norms needed in the continuation criterion.

To describe this range, recall that \eqref{Eq.wave_eq} is invariant under the scaling transformation
\[u(t,y)\longmapsto u_\lambda(t,y):=\lambda^{\frac2{p-1}}u(\lambda t,\lambda y),\qquad \lambda>0.\]
The homogeneous Sobolev norm preserved by the transformation $u\mapsto u_\lambda$ is given by 
\[s_c:=\frac d2-\frac2{p-1},\quad \norm{u_\lambda(t,\cdot)}_{\dot H^{s_c}}=\norm{u(\lambda t,\cdot)}_{\dot H^{s_c}}.\]
The equation is energy-subcritical, energy-critical, or energy-supercritical according to $s_c<1$, $s_c=1$, or $s_c>1$. For $d\geq3$, the critical power is $1+4/(d-2)$; for $d\leq2$, every finite power is energy-subcritical.

For the subcritical problem, the global regularity theory in dimension three goes back to J\"orgens \cite{Jorgens}. Ginibre and Velo \cite{Ginibre-Velo1,Ginibre-Velo2} established global well-posedness in the energy space for the corresponding class of nonlinear Klein--Gordon equations in arbitrary dimension. Propagation of higher regularity was studied by Brenner and von Wahl \cite{Brenner-Wahl}, in particular for dimensions at most nine. The energy-critical problem requires an additional argument to rule out concentration of the energy. Struwe \cite{Struwe1988} first proved global regularity for the three-dimensional quintic equation with radial smooth data. The works of Grillakis \cite{Grillakis1990,Grillakis1992} and Shatah and Struwe \cite{Shatah-Struwe1993} removed the symmetry restriction and established global regularity in dimensions at most nine. Shatah and Struwe \cite{Shatah-Struwe1994} subsequently obtained global well-posedness in the energy space in all dimensions. For general powers, the global regularity theorems require the smoothness assumptions on the nonlinearity specified in the cited works. For the asymptotic behavior of global defocusing solutions, we refer to Yang \cite{Yang} and the references therein.

In \cite[p.~37]{Bourgain2000}, Bourgain posed the following problem:
\begin{quote}
	\textsc{Problem.} \emph{Global existence of classical solutions for defocusing supercritical NLW, NLS.}
\end{quote}

For the supercritical problem, global weak solutions exist in the energy class \cite{Strauss}, and the small-data theory gives global solutions and scattering under appropriate regularity and decay assumptions \cite{Lindblad-Sogge}. On the other hand, loss of regularity for supercritical wave equations was studied by Lebeau \cite{Lebeau2005}. Christ, Colliander and Tao \cite{Christ-Colliander-Tao} established ill-posedness below the scaling-critical regularity, and Ibrahim, Majdoub and Masmoudi \cite{Ibrahim-Majdoub-Masmoudi} studied ill-posedness in the supercritical energy space. Another line of work establishes global existence and scattering under an a priori bound for the critical Sobolev norm. For radial solutions in three dimensions, such a result was proved by Kenig and Merle \cite{Kenig-Merle2011}. Further conditional results were obtained in \cite{Bulut2012, Killp-Visan2011Trans,Killip-Visan2011Radial}. 

An important blow-up result for defocusing wave systems was obtained by Tao \cite{Tao2016}. He considered $\Box u=(\na_{\R^m}F)(u)$, with the unknown $u:\R^{1+d}\to\R^m$, where the potential $F$ is smooth, and is positive and homogeneous of degree $p+1$ outside the unit ball. For $d=3$ and every $p>5$, Tao's theorem gives such a potential with $m=40$ components and smooth data such that the corresponding solution develops a finite-time singularity. The data can be made compactly supported by the finite speed of propagation. The freedom to choose the potential and the number of components is essential in \cite{Tao2016}; the construction concerns a system with a chosen potential rather than the prescribed scalar potential $F(u)=|u|^{p+1}/(p+1)$ in \eqref{Eq.wave_eq}. Tao \cite{Tao2018} also constructed finite-time blow-up solutions for a defocusing nonlinear Schr\"odinger system with a suitable potential and a smooth compactly supported forcing term.

For the energy-supercritical defocusing nonlinear Schr\"odinger equation, Merle, Rapha\"el, Rodnianski and Szeftel \cite{Merle-4-Invent} constructed finite-time blow-up solutions, giving a negative answer to Bourgain's global existence problem \cite{Bourgain2000}. Precisely, for $(d,p)\in\{(5,9),(6,5),(8,3),$ $(9,3)\}$, the theorem in \cite{Merle-4-Invent} constructs a finite-codimensional family of smooth, localized, radial initial data whose solutions blow up in finite time, subject to a profile nondegeneracy condition checked numerically in the companion work \cite{Merle-4-Ann1}. The leading profiles come from smooth self-similar imploding solutions of the compressible Euler equations. The compression of the associated fluid flow produces a growing amplitude together with increasingly rapid oscillations. The corresponding singularity formation theory for compressible Euler and Navier-Stokes equations was developed in \cite{Merle-4-Ann2}. See also \cite{BCLGS2025,CLGSSS2025,CLGSSS-2024,SWWZ2025} for further developments. A recent preprint of OpenAI \cite{OpenAI-NLS2026} presents a construction of stable self-similar blow-up for the supercritical defocusing nonlinear Schr\"odinger equation on $\mathbb T^{12}$. More precisely, for arbitrarily large odd exponents $p$, there exist an integer $k>8$ and a nonempty open subset of $H^k(\mathbb T^{12};\C)$ such that every initial datum in the open subset generates a solution that blows up in finite time. 

Motivated by \cite{Merle-4-Invent, Merle-4-Ann1, Merle-4-Ann2}, the connection between a defocusing dispersive equation and an imploding fluid provided the starting point for our study of the nonlinear wave equation. In our preceding work \cite{SWZ-NLW}, we constructed finite-time blow-up solutions of the complex-valued equation $\Box u=|u|^{p-1}u$ for $d=4$ and odd $p\geq29$, and for $d\geq5$ and odd $p\geq17$. The leading amplitude and phase were obtained from the self-similar solutions of the relativistic Euler equations constructed in \cite{Shao-Wei-Zhang}. The blow-up result of \cite{SWZ-NLW} was extended by Buckmaster and Chen \cite{Buckmaster-Chen2024} to a smaller exponent $(d,p)=(4,7)$. To the best of our knowledge, our earlier works \cite{Shao-Wei-Zhang, SWZ-NLW} were the first to identify the relativistic Euler equations as a limiting system for the nonlinear wave equation and to construct blow-up solutions from relativistic Euler profiles. More precisely, the representation $u=w\ee^{\mathrm i\phi}$ separates the equation into real and imaginary parts, and the leading-order terms give the relativistic Euler system. The remaining terms are then treated as perturbations. A complex-valued solution cannot be turned into a solution of \eqref{Eq.wave_eq} merely by taking $\operatorname{Re}u$, because $\operatorname{Re}(|u|^{p-1}u)$ is not determined by $\operatorname{Re}u$ alone.

Recently, Cadiot, Ionescu and Palasek \cite{CIP2026} constructed finite-time blow-up solutions of \eqref{Eq.wave_eq} for $(d,p)=(10,5)$ from smooth compactly supported radial real-valued initial data. The proof in \cite{CIP2026} uses a computer-assisted Newton--Kantorovich argument to construct an exact discretely self-similar solution in a backward light cone. The profile is periodic in the logarithmic self-similar time and analytic across the lateral boundary of the cone. The theorem in \cite{CIP2026} concerns the quintic equation in dimension ten, whereas our unconditional result below (Theorem \ref{Thm.main_thm}) applies to every $d\geq4$ and all sufficiently large odd $p$.

In the present paper, we construct real-valued blow-up solutions from self-similar imploding profiles of the relativistic Euler equations. As suggested in our work \cite{SWZ-NLW}, we replace the complex exponential by the real-valued periodic oscillator
\begin{equation}\label{Eq.oscillator_intro}
	f''+f^p=0,\qquad f(0)=1,\qquad f'(0)=0,
\end{equation}
and start from the ansatz $u=w f(\phi)$. The leading amplitude and phase satisfy the same system as in \cite{SWZ-NLW}. {\bf The main new difficulty is the linearized potential $p f^{p-1}$, which is periodic and nonconstant.} We must control the oscillation of $pf^{p-1}$ together with the spatial derivatives and the characteristic boundary of the backward light cone. We use a periodic matrix to define an analytic energy, followed by an averaging argument for the perturbation pair.

{\bf A further key ingredient} is the outgoing property of the normalized transport field. In our previous work \cite{SWZZ-Euler2026}, an outgoing property yields weighted decay estimates for the transport semigroup in the stability analysis of incompressible Euler profiles. In the present problem, a quantitative outgoing property gives a favorable term proportional to the derivative order in the energy estimate. This contribution permits the spatial analytic radius to increase and leads to compact evolution operators on fixed analytic spaces. For the full wave equation, the outgoing contribution must also dominate the growth caused by the periodic linearized potential, which is exactly the content of our Assumption \ref{Ass.analytic_margin}.

We first state a conditional result for a fixed Euler profile.

\begin{theorem}\label{Thm.conditional}
	Let $d\in\Z_{\geq4}$ and $p\in2\Z_{\geq1}+1$. Suppose that Assumption \ref{Ass.analytic_margin} holds for $(d,p)$. Then there exist real-valued radial data
	\[(u_0,u_1)\in C_c^\infty(\R^d;\R)^2\]
	such that there is no global smooth solution $u:[0,+\infty)\times\R^d\to\R$ to the defocusing nonlinear wave equation \eqref{Eq.wave_eq} with initial data $u(0,\cdot)=u_0$ and $\partial_tu(0,\cdot)=u_1$. In particular, the lifespan of the classical solution $u$ with initial data $(u_0, u_1)$ is $T_*\in(0,\infty)$ and \eqref{Eq.continuation} holds.
\end{theorem}

We explain the roles of the two assumptions in Theorem \ref{Thm.conditional}.

\begin{remark}\label{Rmk.assumptions}
	In the proof of Theorem \ref{Thm.conditional}, we impose two assumptions (Assumptions \ref{Ass.background} and \ref{Ass.analytic_margin}) for Euler profiles. For a given pair $(d,p)$, Assumption \ref{Ass.background} requires the existence of a self-similar Euler profile, consisting of a homogeneity parameter $\beta>1$ and a velocity profile $v=v(Z)$. Assumption \ref{Ass.analytic_margin} requires that the same profile admit a common coefficient analytic radius $R_v$ and constants $\chi,\varrho$ satisfying \eqref{Eq.analytic_margin}. In particular, Assumption \ref{Ass.analytic_margin} implies Assumption \ref{Ass.background}.
\end{remark}

\begin{remark}\label{Rmk.dimension_three}
	The conditional implication in Theorem \ref{Thm.conditional} remains valid with $d=3$, but Assumption \ref{Ass.background} itself cannot hold in $d=3$. Indeed, Lemma \ref{Lem.sonic_geometry} shows that every profile in Assumption \ref{Ass.background} satisfies
	\[d-1>\beta(\ell+\sqrt\ell),\qquad \beta>1,\qquad \ell=1+\frac4{p-1}>1.\]
	Thus $d-1>2$, and $\bigl\{(3,p):p\in2\Z_{\geq1}+1,\ \text{Assumption \ref{Ass.background} holds for }(3,p)\bigr\}=\varnothing$.

	This dimensional restriction is intrinsic to the radial Euler-based front-compression mechanism considered here: the obstruction already appears at the level of the leading profile, before any perturbative estimates are imposed. The related wave-equation constructions in \cite{Buckmaster-Chen2024,SWZ-NLW} are likewise restricted to dimensions $d\geq4$. An analogous low-dimensional obstruction occurs in the pioneering work of Merle, Rapha\"el, Rodnianski and Szeftel \cite{Merle-4-Invent} on the defocusing nonlinear Schr\"odinger equation, where the dominance of the Euler regime requires $d\geq5$.

	Recently, Buck, Carballeira, G\'omez-Serrano and Shi \cite{BCGS2026} bypassed this restriction by directly constructing self-similar profiles for the full Schrödinger equation. Their computer-assisted construction in \cite{BCGS2026} yields asymptotically self-similar finite-time blow-up solutions from smooth compactly supported radial initial data for $(d,p)\in\{(3,27),(4,9),(5,7)\}$. For the wave equation \eqref{Eq.wave_eq}, however, a construction based on a smooth continuously self-similar profile is unavailable in dimension three. Indeed, the monotonicity argument of Tao \cite[Proposition~2.2]{Tao2016} shows that, for $d=3$ and $p>5$, a solution of the form
	\[u(t,y)=(T-t)^{-\frac{2}{p-1}} U\!\left(\frac{y}{T-t}\right), \qquad U\in C^\infty\left(\overline{B_1^3};\R\right),\]
	must be identically zero. Here smoothness up to the lateral boundary of the light cone is essential.
\end{remark}

Our second result gives a common exponent threshold for all dimensions $d\geq4$. We verify Assumption \ref{Ass.analytic_margin} in dimensions four and five and use the finite speed of propagation for the remaining dimensions.

\begin{theorem}\label{Thm.main_thm}
	There exists an odd integer $p_0\geq15$, independent of $d$, such that, for every $d\in\Z_{\geq4}$ and every odd integer $p\geq p_0$, there exist real-valued data $(u_0,u_1)\in C_c^\infty(\R^d;\R)^2$ such that there is no global smooth solution $u:[0,+\infty)\times\R^d\to\R$ to the defocusing nonlinear wave equation \eqref{Eq.wave_eq} with initial data $u(0,\cdot)=u_0$ and $\partial_tu(0,\cdot)=u_1$.
	For $d\in\{4,5\}$, both $u_0$ and $u_1$ can be chosen radial.
\end{theorem}

Theorem \ref{Thm.main_thm} gives a negative answer to Bourgain's global existence problem \cite{Bourgain2000} for the real-valued scalar nonlinear wave equation \eqref{Eq.wave_eq}: in every dimension $d\geq4$ and for all sufficiently large odd $p$, there are smooth compactly supported real-valued initial data whose solutions blow up in finite time.

The construction gives quantitative bounds for the solution inside the backward light cone.

\begin{remark}[The local singularity]\label{Rmk.cone}
	Fix $d\in\Z_{\geq4}$, an odd $p\geq3$, and a choice $(\beta,v,R_v,\chi,\varrho_0)$ in Assumption \ref{Ass.analytic_margin}. The direct construction gives numbers $T>t_0$ and a real-valued radial solution $u^{\mathrm c}$ on
	\[\mathcal C_{t_0,T} :=\{(t,y):t_0\leq t<T,\ |y|\leq T-t\}.\]
	The solution $u^{\mathrm c}$ is smooth up to the lateral boundary away from $(T,0)$. There are $c,C>0$ and $t_1\in[t_0,T)$ such that
	\begin{equation}\label{Eq.cone_rate}
		c(T-t)^{-\frac{2\beta}{p-1}} \leq\norm{u^{\mathrm c}(t,\cdot)}_{L^\infty(|y|\leq T-t)} \leq C(T-t)^{-\frac{2\beta}{p-1}}, \qquad t_1\leq t<T.
	\end{equation}
	Moreover, there are $K_*>0$ and points $y_\pm(t)$ for every $t\in[t_1,T)$ such that
	\begin{equation}\label{Eq.peaks}
		|y_\pm(t)|\leq K_*(T-t)^{\frac{\beta+1}{2}},\qquad \lim_{t\uparrow T}(T-t)^{\frac{2\beta}{p-1}} u^{\mathrm c}(t,y_\pm(t))=\pm1.
	\end{equation}
	For the family in Proposition \ref{Prop.quantitative_background}, the parameters satisfy, along odd integers $p\to\infty$,
	\begin{equation}\label{Eq.beta_limits}
		\lim_{p\to\infty}\beta_{4,p}=\frac76, \qquad \lim_{p\to\infty}\beta_{5,p}=\frac53.
	\end{equation}
	The estimates \eqref{Eq.cone_rate}--\eqref{Eq.peaks} concern $u^{\mathrm c}$ inside the cone. A compact extension of the initial data of $u^{\mathrm c}$ has maximal time $T_*\leq T-t_0$ after translating $t_0$ to zero; an exterior singularity could occur before the cone vertex is reached.
\end{remark}

We also explain why the exponent threshold in Theorem \ref{Thm.main_thm} is independent of the final dimension.

\begin{remark}\label{Rmk.dimensions}
	For Theorem \ref{Thm.main_thm}, Corollary \ref{Cor.quantitative_margin} verifies Assumption \ref{Ass.analytic_margin} in dimensions four and five with one exponent threshold. In each dimension $d>5$, we extend a four-dimensional cone solution constantly in the additional variables and cut off the initial data outside the initial cone section. Finite speed of propagation gives the required solution without imposing a further restriction on $p$. Thus $p_0$ is independent of the final dimension. The higher-dimensional extension is generally nonradial in $\R^d$.
\end{remark}

\smallskip

\noindent  {\bf  Key ingredients of the proof.}\smallskip

The argument separates the construction for a fixed Euler profile from the quantitative estimates needed to choose the profile. See Figure \ref{Fig.proof_scheme} for a more refined summary. See Section \ref{Sec.Proof} for a more detailed roadmap of the proof.
\begin{enumerate}[(i)]
	\item \emph{Euler profiles and normalization.} Following our previous work \cite{SWZ-NLW}, we choose a self-similar imploding solution of the potential-flow form of the relativistic Euler equations
	\[\rho_0^{p-1}+\partial^\mu\phi_0\partial_\mu\phi_0=0,\qquad \partial^\mu(\rho_0^2\partial_\mu\phi_0)=0.\]
	Assumption \ref{Ass.background} requires the corresponding smooth velocity profile. Writing $u=\rho_0z$, taking $\sigma=\phi_0$, and changing the radial coordinate as in \eqref{Eq.sigma_z}, we obtain the exact equation on a fixed ball
	\[\cN(z)=D^2z+\sigma^{-2}(A_0-V)z+z^p=0,\qquad D=\partial_\sigma+\sigma^{-1}B.\]
	Here $B=\bvec\cdot\nabla$, $A_0=-J_0^{-1}\diver(J_0\nabla)$, and the coefficients are determined by the Euler profile through \eqref{Eq.normalized_coefficients}.
	
	\item \emph{The phase correction.} Solvability of the antiperiodic correction equations requires orthogonality to $f'$. We introduce $\theta_\sigma=\sigma+\sigma^{-1}\Psi(\log\sigma,x)$, where $\Psi$ is polynomial in $\tau=\log\sigma$ and solves
	\[\bigl((\partial_\tau+B-2)(\partial_\tau+B-1)+c_a^2A_0\bigr)\Psi=-\av{f^2}(B-2)V,\qquad c_a^2=\frac{p-1}{p+3}.\]
	The averaged evolution in step (iii) makes the polynomial construction possible. The remaining (easier) antiperiodic corrections give $z_{\app}$ with $\|\cN(z_{\app})\|_{\mathcal B_{\varrho_1}}\lesssim\sigma^{-5}(1+\log\sigma)^M$, for some $M\geq0$, by Proposition \ref{Prop.approximation}; $\mathcal B_r$ denotes the spatial analytic space in \eqref{Eq.analytic_norm}.
	
	\item \emph{Modulation and compactness.} For $\zeta=z-z_{\app}$, we use $P_0$ from \eqref{Eq.Pzero} and define
	\[W=(\psi,\pi)^T=\diag(1,\sigma)P_0(\theta_\sigma)^{-1}\binom\zeta{\lambda^{-1}D\zeta},\qquad \lambda=D\theta_\sigma.\]
	By \eqref{Eq.state_operator_form}, the exact equation is 
	\begin{equation}\label{Eq.W-eq-intro}
		\partial_\tau W=(\mathscr A+\cB(\sigma)+\mathscr E_\sigma)W+\mathbf g_\sigma h,
	\end{equation}
	where $h=-\sigma^2(\cN(z_{\app})+\cQ_\sigma(\zeta))$ and $\cQ_\sigma$ is the nonlinear remainder in \eqref{Eq.remainder_equation}. The operator $\cB$ is periodic with zero average, and the coefficients of $\mathscr E_\sigma$ tend to zero. The homogeneous averaged equation is
	\[(\partial_\tau+B)\psi=\pi,\qquad (\partial_\tau+B)\pi=\pi-c_a^2A_0\psi.\]
	The outgoing bounds $\sym\mathrm D_x\bvec\geq\gamma_vI_d$ and $\bvec\cdot n=1>c_a$ in \eqref{Eq.acoustic_outgoing_property} permit an increase of the analytic radius at sufficiently small radii. Compact embedding at a smaller radius gives compact evolution operators and compact resolvent of $\cA_r$ (Proposition \ref{Prop.macro_spectrum}). The phase equation reduces to an evolution generated by $\cA_{\varrho_2}+I$; any resonance at zero is finite-dimensional and is handled by polynomial dependence on $\tau$.
	
	\item \emph{Outgoing transport versus periodic growth.} The full linearized equation retains the rapidly oscillating potential $\lambda^2pf(\theta_\sigma)^{p-1}$. For $\chi>\chi_p$, Proposition \ref{Prop.matrix} constructs a family of matrix multipliers $S$ with
	\begin{align*}
		\partial_\theta S+\mathsf P_\xi^TS+S\mathsf P_\xi&\leq(2\chi\xi+C\delta)S,\qquad \mathsf P_\xi=
		\begin{pmatrix}
			0&1\\
			-pf(\theta)^{p-1}-\xi^2&0
		\end{pmatrix}.
	\end{align*}
	Using $S$ in the energy \eqref{Eq.energy_weights}, with $\xi=N/(\varrho\sigma\lambda)$ and $\delta=\sigma^{-1}$, bounds the oscillatory contribution at derivative order $N$ by $2\chi N/(\varrho\sigma)+C/\sigma$. The outgoing contribution is $-2\gamma_vN/\sigma$, and the transport commutator loss is measured by $\mathcal L_v(\varrho)$ in \eqref{Eq.geometry_budget}. Thus the energy estimate \eqref{Eq.full_energy_differential} permits an increase of the analytic radius when
	\[\gamma_v-\mathcal L_v(\varrho)-\frac\chi\varrho-\frac{\varrho'}{\varrho}>0.\]
	Assumption \ref{Ass.analytic_margin} imposes the quantitative margin needed for this inequality and the auxiliary radii. This is the most difficult part in our proof.
	
	\item \emph{Averaging and a decaying solution.} The full energy estimate, the periodic corrector \eqref{Eq.averaging_corrector}, and the analytic-radius gaps give
	\[\|\mathcal P_\Sigma(1)-\cS(1)\|_{X\to X}\longrightarrow0\quad\text{as }\Sigma\to\infty,\qquad X=\cX_{\varrho_0}.\]
	Here $\mathcal P_\Sigma(1)$ and $\cS(1)$ are the full homogeneous and averaged solution operators over one unit of logarithmic time, respectively. Proposition \ref{Prop.compact_decay}, applied to the nonlinear endpoint maps near the compact operator $\ee^2\cS(1)$, selects a solution with $\|W(\sigma)\|_{\cX_{\varrho_-}}\lesssim\sigma^{-a}$ for some $2<a<3$. This reduces the full equation \eqref{Eq.W-eq-intro} into a fixed-point formulation. Reconstruction and finite speed of propagation prove Theorem \ref{Thm.conditional}.
	
	\item \emph{Uniform estimates for large $p$.} Assumption \ref{Ass.background} alone does not provide the quantitative margin in Assumption \ref{Ass.analytic_margin}. We revisit the Euler construction in \cite{Shao-Wei-Zhang} for $d=4,5$, retaining bounds uniform in $p$. Propositions \ref{Prop.quantitative_background} and \ref{Prop.scattering_growth} give $\gamma_v\geq\gamma_*>0$, an analytic radius $R_p\asymp p^{-3/2}$ and $\chi_p\lesssim p^{-2}$. Choosing $\varrho_p$ to be a sufficiently small fixed fraction of $R_p$ makes $\mathcal L_v(4\varrho_p)$ small compared with $\gamma_*$. Choosing $\chi>\chi_p$ of order $p^{-2}$ then gives $\chi/\varrho_p=O(p^{-1/2})$. Corollary \ref{Cor.quantitative_margin} verifies Assumption \ref{Ass.analytic_margin} for all sufficiently large odd $p$. Constant extension of the four-dimensional solution and finite propagation give Theorem \ref{Thm.main_thm} in every $d\geq4$.
\end{enumerate}
\smallskip

\noindent\textbf{Organization of the paper.} \smallskip

Section \ref{Sec.Proof} gives a roadmap of the proof. We derive the normalized equation, construct the finite approximation, and formulate the perturbation equation. We then state the analytic energy and averaging estimates and apply an abstract compact-operator proposition to the nonlinear endpoint maps. Section \ref{Sec.approximation_proofs} proves the averaged energy estimate, constructs the polynomial phase and the finite antiperiodic corrections, and establishes the perturbation identities. Section \ref{Sec.matrix_proof} constructs the periodic matrices at every frequency. Section \ref{Sec.energy_proof} proves the energy estimate for the full equation. Section \ref{Sec.state_proof} establishes quantitative averaging and constructs the nonlinear perturbation. Section \ref{Sec.background_proof} proves the uniform scattering and background estimates. The linear and nonlinear arguments are formulated for a fixed background under Assumptions \ref{Ass.background} and \ref{Ass.analytic_margin}; Corollary \ref{Cor.quantitative_margin} verifies Assumptions \ref{Ass.background} and \ref{Ass.analytic_margin} for the exponents in Theorem \ref{Thm.main_thm}. Appendix \ref{App.analyticity} establishes the regularity and expansion properties of the Euler profiles under Assumption \ref{Ass.background}, and Subsection \ref{App.normalization} proves the normalization lemma. Appendix \ref{App.analytic_spaces} collects the properties of the analytic function spaces.

\smallskip
\noindent\textbf{Notations and conventions.}

\begin{itemize}
	\item We use standard multi-index notation and write $\Z_+:=\{1,2,\ldots\}$. The natural logarithm is denoted by $\log$, and the natural logarithmic base by $\ee=\exp(1)$. The letters $c,C$ denote positive constants that may depend on the fixed parameters and change from line to line. Any uniformity in parameters is stated explicitly. For nonnegative quantities, $a\lesssim b$ means $a\leq Cb$, and $a\asymp b$ means $a\lesssim b\lesssim a$.

	\item The physical variables are $(t,y)\in\R\times\R^d$, with $r=|y|$. In normalized coordinates $(\sigma,x)$, we write $X=|x|$, $\na=\na_x$ and $\Delta=\Delta_x$. The letter $\varrho$ denotes an analytic radius. We write $B_r^n(a)$ for the $n$-dimensional Euclidean ball of radius $r$ centered at $a$, and $B_r^n:=B_r^n(0)$, omitting $n$ when clear.

	\item Partial derivatives are written with $\partial$, whereas primes denote ordinary derivatives of one-variable functions. Subscripts on functions may denote parameters or components; in particular, the prescribed phase $\theta_\sigma$ is not a derivative. For $I=(i_1,\ldots,i_j)\in\{1,\ldots,n\}^j$, we define $\partial_I:=\partial_{x_{i_1}}\cdots\partial_{x_{i_j}}$. We use the full ordered derivative tensor $\na^j h=(\partial_Ih)_I$, with
	\[|\na^j h|^2=\sum_{I\in\{1,\ldots,n\}^j}|\partial_Ih|^2=\sum_{|\alpha|=j}\frac{j!}{\alpha!}|\partial^\alpha h|^2,\qquad \na^0h=h.\]
	Here $|\cdot|$ is the Euclidean or Hermitian norm over all components of the unknowns. The ordered-tensor convention also applies to weak derivatives in $L^q$. For matrices, $\|\cdot\|$ denotes the Euclidean operator norm and $\|\cdot\|_{\rm HS}$ the Hilbert--Schmidt norm. Matrix-valued coefficients are differentiated using
	\[\|\na^jM\|^2:=\sum_{I\in\{1,\ldots,n\}^j}\|\partial_IM\|^2.\]
	For a vector field $h$, $\mathrm D_xh$ is the Jacobian of $h$, so that $|\na h|=\|\mathrm D_xh\|_{\rm HS}$, whereas $\|\mathrm D_xh\|$ is the operator norm. For square matrices with real entries, we set $\sym M:=(M+M^T)/2$.

	\item For a finite-dimensional space $E$, a compact set $K\subset\R^n$ and $s\in\{\infty,\omega\}$, $C^s(K;E)$ consists of restrictions of $E$-valued $C^s$ functions defined on an open neighborhood of $K$; $\omega$ denotes real analyticity. On other non-open sets, the same convention applies locally at included boundary points. In particular, $C^\omega([a,b];E)$ requires analytic extension across both endpoints. Subscripts $\mathrm e$ and $\mathrm o$ denote even and odd functions. On $[0,a]$, with $a>0$, the subscripts require a $C^s$ even or odd extension to $[-a,a]$; equivalently,
	\begin{equation}\label{Eq.even_odd_spaces}
		\begin{aligned}
			C^s_{\rm e}([0,a];E)&=\{Z\mapsto g(Z^2):g\in C^s([0,a^2];E)\},\\
			C^s_{\rm o}([0,a];E)&=\{Z\mapsto Zg(Z^2):g\in C^s([0,a^2];E)\}.
		\end{aligned}
	\end{equation}
	We use the same notation for a function and the prescribed even or odd extension, and omit $E$ when $E=\R$.

	\item We write $\Hol(\mathcal U;E)$ for holomorphic functions on an open set $\mathcal U\subset\C^n$ with values in a complex Banach space $E$, and $\Hol(\mathcal U):=\Hol(\mathcal U;\C)$. For compact $K\subset\R^n$, bounded $\Omega\subset\R^n$ and $\delta>0$, set
	\[K^\C_\delta:=\{\zeta\in\C^n:\dist(\zeta,K)<\delta\},\qquad\Omega^\C_\delta:=(\overline\Omega)^\C_\delta,\]
	where the distance is Euclidean. We also write $\mathbb D_r:=\{\zeta\in\C:|\zeta|<r\}$ and $\mathbb S_r:=\{\zeta\in\C:|\operatorname{Im}\zeta|<r\}$. For spectral and semigroup arguments, spaces and operators are complexified without changing notation, using the same norm formulas and preserving all symmetry and compatibility conditions.
\end{itemize}

\section{A roadmap of the proof}\label{Sec.Proof}

In Subsection \ref{Subsec.leading}, we reduce Theorem \ref{Thm.conditional} to an existence result for the normalized equation. Subsection \ref{Subsec.unconditional} verifies the hypotheses for a family of profiles and proves Theorem \ref{Thm.main_thm}. After constructing an approximate solution, we formulate the perturbation equation and state the principal linear estimates. Growth of the analytic radius yields compact endpoint maps. An abstract spectral proposition then selects a decaying nonlinear trajectory.

\begin{figure}[tbp]
	\centering

	\begin{tikzpicture}
		[
		box/.style={draw,rectangle,align=center,text width=6.35cm,
		minimum height=1.18cm,inner xsep=4pt,inner ysep=5pt,
		font=\small},
		arrow/.style={-{Stealth[length=2mm]},line width=0.45pt},
		head/.style={font=\small\bfseries,align=center},
		node distance=0.48cm
		]
		\node[head] at (0,0) {The conditional construction};
		\node[head] at (7.75,0) {The uniform exponent threshold};
		\node[box] (l1) at (0,-1.05)
		{A fixed $d\geq4$ and odd $p$\\
		An Euler profile satisfying\\Assumption \ref{Ass.background}};
		\node[box,below=of l1] (l2)
		{Fourth-order approximate solution\\
		(Proposition \ref{Prop.approximation})\\
		Residual $O(\sigma^{-5}(1+\log\sigma)^M)$};
		\node[box,below=of l2] (l3)
		{Full energy estimate and analytic margin\\
		(Proposition \ref{Prop.energy}, Assumption \ref{Ass.analytic_margin})\\
		Then averaging (Proposition \ref{Prop.averaging})};
		\node[box,below=of l3] (l4)
		{Compact-operator decay and nonlinear endpoints\\
		(Proposition \ref{Prop.compact_decay}, Theorem \ref{Thm.state})};
		\node[box,below=of l4] (l5)
		{Radial blow-up in the fixed dimension $d$\\
		(Theorem \ref{Thm.conditional})};
		\node[box] (r1) at (7.75,-1.05)
		{$d\in\{4,5\}$, odd $p$ sufficiently large\\
		Quantitative Euler construction\\
		(Proposition \ref{Prop.quantitative_background})};
		\node[box,below=of r1] (r2)
		{Assumption \ref{Ass.analytic_margin} for $d=4,5$\\
		and every odd $p\geq p_0$\\
		One common $p_0$};
		\node[box,below=of r2] (r3)
		{Radial blow-up in $d\in\{4,5\}$ (Theorem \ref{Thm.conditional})};
		\node[box,below=of r3] (r4)
		{For $d>5$: constant extension from $d=4$\\
		Cutoff of the data and finite propagation};
		\node[box,below=of r4] (r5)
		{Blow-up in every dimension $d\geq4$\\
		The same $p_0$; Theorem \ref{Thm.main_thm}};
		\foreach \a/\b in {l1/l2,l2/l3,l3/l4,l4/l5,r1/r2,r2/r3,r3/r4,r4/r5}
		\draw[arrow] (\a.south)--(\b.north);
	\end{tikzpicture}

	\caption{The conditional construction and the proof of the dimension-independent exponent threshold. The bounds in the right column are uniform for $d\in\{4,5\}$.}
	\label{Fig.proof_scheme}
\end{figure}
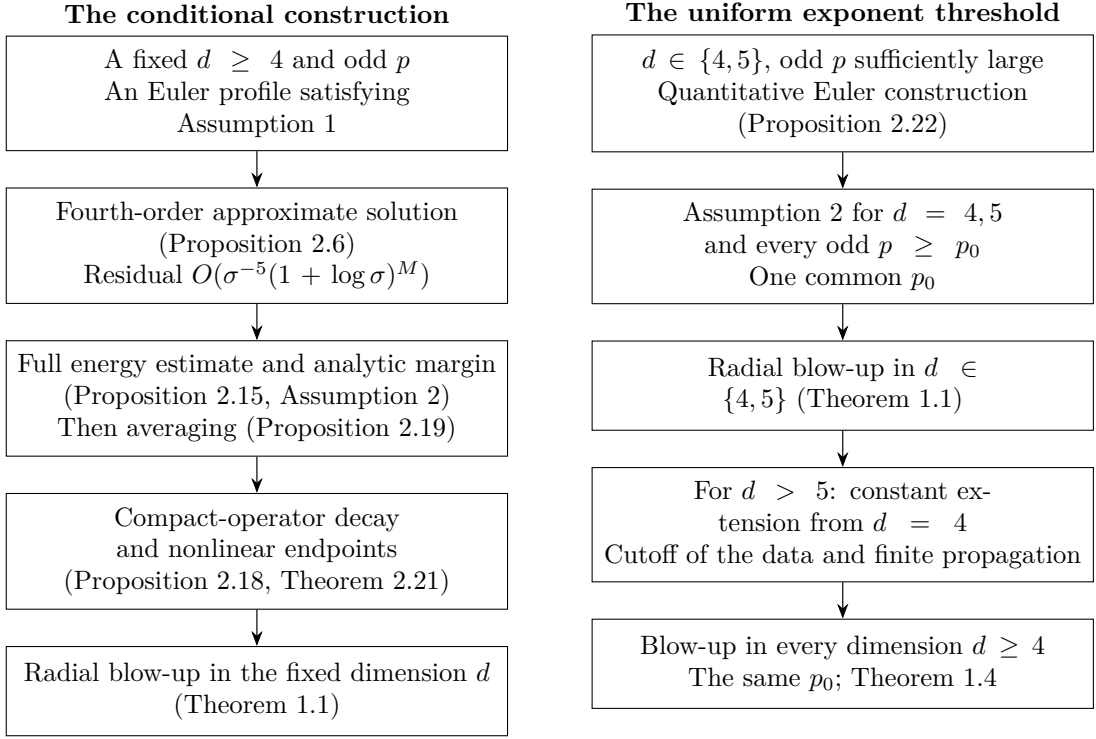

Figure \ref{Fig.proof_scheme} summarizes the conditional construction and the quantitative profile estimates. The left column leads from an Euler profile to a decaying perturbation of the fourth-order approximation. The right column combines the background estimates and the oscillator bound to obtain an exponent threshold independent of the final dimension.

\subsection{The leading profile and the normalized scalar equation}
\label{Subsec.leading}

Fix the solution $f$ of \eqref{Eq.oscillator_intro}. The oscillator energy identity is
\begin{equation}\label{Eq.oscillator_energy}
	\frac12(f')^2+\frac1{p+1}f^{p+1}=\frac1{p+1}.
\end{equation}
Consequently, $f$ is defined on $\R$ and oscillates between $1$ and $-1$. The ODE in \eqref{Eq.oscillator_intro} and uniqueness give $f\in C^\omega_{\rm e}(\R)$. We denote the period of $f$ by $2L$, where
\[L:=2\sqrt{\frac{p+1}{2}}\int_0^1\frac{\mathrm da}{\sqrt{1-a^{p+1}}}.\]
Uniqueness for the oscillator equation gives $f(\theta+L)=-f(\theta)$. In particular, $p f^{p-1}\in C^\omega_{\rm e}(\R)$ has period $L$. We use the normalization $f(0)=1$, $f'(0)=0$ throughout the paper.

Let $w=w(t,y)>0$ and let $\phi=\phi(t,y)$ be real-valued. The chain rule gives the exact identity
\begin{equation}\label{Eq.scalar_ansatz}
	\begin{aligned}
		\Box\bigl(w f(\phi)\bigr)-\bigl(w f(\phi)\bigr)^p ={}&f(\phi)\Box w +f'(\phi)\bigl(2\partial^\mu w\partial_\mu\phi+w\Box\phi\bigr)\\
		&\qquad-w f(\phi)^p\bigl(w^{p-1}+\partial^\mu\phi\partial_\mu\phi\bigr).
	\end{aligned}
\end{equation}
Thus, for a proposed pair $(w,\phi)$, the scalar equation \eqref{Eq.wave_eq} is the vanishing of the right-hand side of \eqref{Eq.scalar_ansatz}. Unlike the complex modulus-phase decomposition, this identity does not give two independent equations by separating real and imaginary parts. The leading-order approximation is obtained by imposing
\begin{equation}\label{Eq.leading_system}
	\rho_0^{p-1}+\partial^\mu\phi_0\partial_\mu\phi_0=0, \qquad \partial^\mu(\rho_0^2\partial_\mu\phi_0)=0.
\end{equation}
For $u_0=\rho_0f(\phi_0)$ satisfying \eqref{Eq.leading_system}, the residual is
\begin{equation}\label{Eq.leading_error}
	\Box\bigl(\rho_0f(\phi_0)\bigr) -\bigl(\rho_0f(\phi_0)\bigr)^p =f(\phi_0)\Box\rho_0.
\end{equation}
We use \eqref{Eq.leading_system} to determine $\rho_0,\phi_0$ and correct the residual $f(\phi_0)\Box\rho_0$ below.

We use the same self-similar notation as in \cite{SWZ-NLW}:
\begin{equation}\label{Eq.selfsimilar}
	\phi_0(t,r):=(T-t)^{1-\beta}\widehat\phi_0(Z),\quad \rho_0(t,r):=(T-t)^{-\frac{2\beta}{p-1}}\widehat\rho_0(Z), \quad\ Z:=\frac r{T-t},\quad r:=|y|.
\end{equation}
Here
\[k:=d-1,\qquad \ell:=1+\frac4{p-1},\qquad \beta>1.\]
The normalization of $\beta$ is that of \cite{SWZ-NLW}; the homogeneity parameter denoted by $\beta$ in \cite{Shao-Wei-Zhang} is $(\ell+1)\beta$ in the present notation. The leading system \eqref{Eq.leading_system} is the potential-flow formulation used in the relativistic Euler construction of \cite{Shao-Wei-Zhang,SWZ-NLW}. With $v(Z):=\partial_r\phi_0/\partial_t\phi_0$, the radial profile equation of \eqref{Eq.leading_system} is
\begin{equation}\label{Eq.velocity}
	\begin{aligned}
		\Delta_Z(Z,v)v'&=\Delta_v(Z,v),\\
		\Delta_Z(Z,v)&:=Z\bigl((1-Zv)^2-\ell(v-Z)^2\bigr),\\
		\Delta_v(Z,v)&:=(1-v^2)\bigl(\beta\ell(1-v^2)Z-kv(1-Zv)\bigr).
	\end{aligned}
\end{equation}
See \cite[Subsection A.1]{SWZ-NLW} for the detailed derivation of \eqref{Eq.velocity} from \eqref{Eq.leading_system}. We impose the following smooth-profile assumption on the interval corresponding to the closed backward light cone, as in \cite{SWZ-NLW}.

\begin{assumption}\label{Ass.background}
	For the given $d\in\Z_{\geq2}$ and $p\in2\Z_{\geq1}+1$, we define $k:=d-1$ and $\ell:=1+4/(p-1)$. There exist $\beta>1$ and $v\in C^\infty_{\rm o}([0,1])$ such that
	\[|v(Z)|<1,\qquad \Delta_Z(Z,v(Z))v'(Z)=\Delta_v(Z,v(Z))\qquad(0\leq Z\leq1).\]
	The polynomials $\Delta_Z$ and $\Delta_v$ are defined by \eqref{Eq.velocity} for the fixed values of $k,\ell,\beta$.
\end{assumption}

Under Assumption \ref{Ass.background}, the profile extends analytically across the origin and the sonic point, as proved in Proposition \ref{Prop.automatic_analyticity}. We use the analytic extension to define the normalized coefficients.

\begin{remark}\label{Rmk.background_analyticity}
	For every pair $(\beta,v)$ satisfying Assumption \ref{Ass.background}, the odd extension of $v$ satisfies
	\[k>\beta(\ell+\sqrt\ell),\qquad v\in C^\omega_{\rm o}([-1,1]).\]
	More precisely, there exist $\eta>0$ and an extension $\widetilde v$ such that
	\begin{equation}\label{Eq.profile_tube}
		\mathcal U_\eta:=\{\zeta\in\C:\dist(\zeta,[-1,1])<\eta\},\qquad \widetilde v\in\Hol(\mathcal U_\eta),\qquad \widetilde v|_{[-1,1]}=v.
	\end{equation}
	The parameter inequality and the nondegeneracy of the sonic point are established in \cite[Lemma 2.1 and Lemmas A.2--A.3]{SWZ-NLW}. Proposition \ref{Prop.automatic_analyticity} proves \eqref{Eq.profile_tube}, including the case of an integer eigenvalue ratio at the sonic point. The radius $\eta$ may depend on the chosen profile.
\end{remark}

Fix a pair $(\beta,v)$ satisfying Assumption \ref{Ass.background}; by Remark \ref{Rmk.background_analyticity}, $v\in C^\omega_{\rm o}([0,1])$. Since $1-Zv(Z)>0$ on $[0,1]$, define
\begin{equation}\label{Eq.profiles}
	\begin{aligned}
		\widehat\phi_0(Z)&:=\frac1{\beta-1} \exp\left((\beta-1)\int_0^Z\frac{v(a)}{1-av(a)}\,\mathrm da\right),\\
		\widehat\rho_0(Z) &:=\frac{(\beta-1)^{\frac2{p-1}}\widehat\phi_0(Z)^{\frac2{p-1}}(1-v(Z)^2)^{\frac1{p-1}}}{(1-Zv(Z))^{\frac2{p-1}}}.
	\end{aligned}
\end{equation}
Equation \eqref{Eq.profiles} gives $\widehat\rho_0(0)=1$ and $\widehat\phi_0(0)=(\beta-1)^{-1}$. Moreover, we also have $\widehat{\phi}_0,\widehat\rho_0\in C_{\operatorname{e}}^\infty([0,1])$, see \cite[Lemma A.6]{SWZ-NLW}. Consequently, \eqref{Eq.selfsimilar} gives a smooth imploding solution to \eqref{Eq.leading_system} through \eqref{Eq.selfsimilar}.

The error term in \eqref{Eq.leading_error} suggests seeking an exact solution to \eqref{Eq.wave_eq} of the form $u=\rho_0z$, with $z$ close to $f(\phi_0)$. We use $\sigma=\phi_0$ as the time variable, so that the leading oscillation becomes the fixed periodic function $f(\sigma)$. The self-similar radius $Z$ already maps the backward light cone to $0\leq Z\leq1$. We make a further change of the radial variable to simplify the two characteristic speeds of the wave operator.

Indeed, differentiation of \eqref{Eq.profiles} and \eqref{Eq.selfsimilar} gives
\begin{equation}\label{Eq.phase_derivatives}
	\partial_t\phi_0=\frac{(\beta-1)\phi_0}{(T-t)(1-Zv)},\qquad \partial_r\phi_0=\frac{(\beta-1)\phi_0v}{(T-t)(1-Zv)}.
\end{equation}
Along the radial characteristic curves $\mathrm dr/\mathrm dt=\pm1$, we obtain
\[\sigma\frac{\mathrm dZ}{\mathrm d\sigma} =\frac{(Z\pm1)(1-Zv)}{(\beta-1)(1\pm v)} =\frac{(1-Zv)^2}{(\beta-1)(1-v^2)} \left(\frac{Z-v}{1-Zv}\pm1\right).\]
Thus, the common positive factor depending on $Z$ can be removed by taking
\begin{equation}\label{Eq.x_b}
	X(Z):=(\beta-1)\int_0^Z\frac{1-v(a)^2}{(1-av(a))^2}\,\mathrm da,\quad b(X(Z)):=\frac{Z-v(Z)}{1-Zv(Z)},\qquad X_1:=X(1).
\end{equation}
Then $X'(Z)>0$ for $Z\in [0,1]$, and the two characteristic speeds in the $(\sigma,X)$ variables are $(b(X)\pm1)/\sigma$. In particular, the difference between the two speeds is $2/\sigma$, and the lateral boundary of the cone is the fixed sphere $X=X_1$, where $b(X_1)=1$. The coordinates $(\sigma,X)$ give the spatial coefficient $\sigma^{-2}$ in \eqref{Eq.N}. We write $Z=Z(X)$ for the inverse of $X=X(Z)$.

We now seek a solution of \eqref{Eq.wave_eq} in the form
\begin{equation}\label{Eq.sigma_z}
	u(t,y)=\rho_0(t,r)z(\sigma,x),\qquad \sigma:=\phi_0(t,r),\qquad x:=
	\begin{cases}
		X(Z)\frac{y}{r},&r>0,\\
		0,&r=0.
	\end{cases}
\end{equation}
Here $z(\sigma,\cdot)$ is a real-valued radial function on the normalized ball $\Omega:=B_{X_1}(0)\subset\R^d$. We also define the vector field
\[\bvec(x):=\begin{cases} b(|x|)\frac{x}{|x|},&x\ne0,\\ 0,&x=0. \end{cases}\]
We write $\mathrm D_x\bvec:=(\partial_jb_i)_{i,j=1}^d$. The equation for $z=z(\sigma, x)$ is given by the following lemma.

Define the operators
\begin{equation}\label{Eq.normalized_operators}
	B:=\bvec\cdot\na,\quad D:=\partial_\sigma+\sigma^{-1}B,\quad A_0:=-J_0^{-1}\diver(J_0\na),\quad \cP:=D^2+\sigma^{-2}A_0,
\end{equation}
where
\[c_0:=b'(0)=\frac{d-\beta\ell}{d(\beta-1)}>0,\]
the coefficients are defined by
\begin{equation}\label{Eq.normalized_coefficients}
	\begin{aligned}
		J_0(x)&:=\exp\left(\int_0^X\frac{dc_0-b'(a)-kb(a)/a}{b(a)}\,\mathrm da\right),\\
		V(x)&:=\left.\widehat\phi_0(Z)^2\widehat\rho_0(Z)^{-p} \left[(1-Z^2)\widehat\rho_0''(Z) +\left(\frac{k}{Z}-2\left(1+\frac{2\beta}{p-1}\right)Z\right)\widehat\rho_0'(Z) \right.\right.\\[-2pt] &\hspace{42mm}\left.\left. -\frac{2\beta}{p-1}\left(1+\frac{2\beta}{p-1}\right)\widehat\rho_0(Z) \right]\right|_{Z=Z(X)},\qquad X=|x|. \end{aligned} \end{equation} In \eqref{Eq.normalized_coefficients}, the values at zero are \[\left.\frac{b(a)}a\right|_{a=0}=c_0,\qquad \left.\frac{dc_0-b'(a)-kb(a)/a}{b(a)}\right|_{a=0}=0,\qquad \left.\frac{\widehat\rho_0'(Z)}Z\right|_{Z=0}=\widehat\rho_0''(0).\]
One has $J_0>0$, $J_0(0)=1$, and
\begin{equation}\label{Eq.density_identity}
	\diver(J_0\bvec)=dc_0J_0,\qquad b(X_1)=1.
\end{equation}
\begin{lemma}[The normalized equation]\label{Lem.normalized_equation}
	Under Assumption \ref{Ass.background}, the transformation \eqref{Eq.sigma_z} makes the wave equation equivalent to
	\begin{equation}\label{Eq.N}
		\cN(z):=D^2z+\sigma^{-2}(A_0-V)z+z^p=0.
	\end{equation}
	More precisely, let $I\subset(0,\infty)$ be an open interval and let $z\in C^\infty(I\times\overline\Omega;\R)$ be radial in $x$. Set $z^\sharp(t,y):=z(\sigma(t,y),x(t,y))$ on $\{(t,y):t<T,\ |y|\leq T-t,\ \sigma(t,y)\in I\}$. Then
	\begin{equation}\label{Eq.physical_P}
		(\cP z)^\sharp=-\rho_0^{-p-1}\partial^\mu(\rho_0^2\partial_\mu z^\sharp),
	\end{equation}
	and
	\begin{equation}\label{Eq.exact_reconstruction}
		\begin{aligned}
			\Box(\rho_0z^\sharp)-(\rho_0z^\sharp)^p&=-\rho_0^p\bigl(\cN(z)\bigr)^\sharp,\\
			\cN(z)&=\partial_\sigma^2z+2\sigma^{-1}B\partial_\sigma z +\sigma^{-2}(B^2-B+A_0-V)z+z^p.
		\end{aligned}
	\end{equation}
\end{lemma}

See Subsection \ref{App.normalization} for the detailed proof of Lemma \ref{Lem.normalized_equation}, which is a brute force computation. For later use, define the coefficient map
\begin{equation}\label{Eq.background_coefficient_map}
	\mathbf c(x):=\bigl(\bvec(x),J_0(x),J_0(x)^{-1},V(x), \widehat\rho_0(Z(|x|)),\widehat\phi_0(Z(|x|))\bigr).
\end{equation}
By Proposition \ref{Prop.automatic_geometry}, $\mathbf c\in C^\omega(\overline\Omega;\R^{d+5})$, and there are $\gamma_v,R_v,M_v>0$ and
\begin{equation}\label{Eq.qualitative_coefficient_extension}
	\widetilde{\mathbf c}\in\Hol(\Omega_{R_v}^{\C};\C^{d+5}),\qquad \widetilde{\mathbf c}|_{\overline\Omega}=\mathbf c,
\end{equation}
depending on the fixed profile, such that
\begin{equation}\label{Eq.vecb-analytic-bounds}
	\sym\mathrm D_x\bvec\geq\gamma_v I_d,\qquad \norm{\na^j\bvec}_{L^\infty(\Omega)}\leq M_v(j-1)!R_v^{1-j},\quad \forall\ j\in\Z_{+}.
\end{equation}
We require the same radius $R_v$ to satisfy, for each scalar component $a$ of $\mathbf c$,
\begin{equation}\label{Eq.common_coefficient_radius}
	\|\na^ja\|_{L^\infty(\Omega)}\leq C_a j!R_v^{-j} \quad \forall\ j\in\Z_{\geq 0}
\end{equation}
for some $C_a<\infty$. Proposition \ref{Prop.automatic_geometry} gives a common analytic radius $R_v>0$ for \eqref{Eq.qualitative_coefficient_extension}, \eqref{Eq.common_coefficient_radius} and the ordered derivative bounds for $\bvec$ in \eqref{Eq.vecb-analytic-bounds}.

In particular, Lemma \ref{Lem.normalized_equation} and the oscillator equation give
\begin{equation}\label{Eq.Nf}
	\cN(f(\sigma))=-\sigma^{-2}Vf(\sigma).
\end{equation}
We shall construct a solution of \eqref{Eq.N} approaching $f(\sigma)$ as $\sigma\to\infty$.

\begin{theorem}[Existence of the normalized scalar solution]\label{Thm.normalized}
	Let $d\in\Z_{\geq4}$ and let $p\geq3$ be odd. Suppose that Assumption \ref{Ass.analytic_margin} holds for $(d,p)$, and fix $(\beta,v,R_v,\chi,\varrho_0)$. There exist $\sigma_0\geq2$, $C>0$, an integer $M\geq0$, and a real-valued radial function $z\in C^\infty([\sigma_0,\infty)\times\overline\Omega)$
	such that \eqref{Eq.N} holds and
	\[\sup_{x\in\overline\Omega}|z(\sigma,x)-f(\sigma)| \leq C\sigma^{-1}(1+\log\sigma)^M, \qquad \sigma\geq\sigma_0.\]
\end{theorem}

We prove Theorem \ref{Thm.conditional} by using Theorem \ref{Thm.normalized}. The additional estimates in Remark \ref{Rmk.cone} are proved at the end of Subsection \ref{Subsec.cone_rates}. See Subsection \ref{Subsec.fixedpoint} for the proof of Theorem \ref{Thm.normalized}.

\begin{proof}[Proof of Theorem \ref{Thm.conditional}]
	Fix $d\in\Z_{\geq4}$, an odd integer $p\geq3$, and a choice $(\beta,v,R_v,\chi,\varrho_0)$ in Assumption \ref{Ass.analytic_margin} for $(d,p)$. Set $T:=1$ and let $z$ be given by Theorem \ref{Thm.normalized} for $(\beta,v,R_v,\chi,\varrho_0)$. Set
	\[m_\phi:=\min_{0\leq Z\leq1}\widehat\phi_0(Z)>0,\qquad h_0:=\min\left\{\frac12,\left(\frac{m_\phi}{2\sigma_0}\right)^{1/(\beta-1)}\right\},\qquad t_0:=1-h_0.\]
	For $t_0\leq t<1$ and $|y|\leq1-t$, one has $\sigma=\phi_0(t,|y|)\geq m_\phi(1-t)^{1-\beta} \geq m_\phi h_0^{1-\beta}\geq2\sigma_0$.
	Consequently the function
	\begin{equation}\label{Eq.direct_cone_solution}
		u^{\mathrm c}_{d}(t,y):=\rho_0(t,|y|)z\bigl(\phi_0(t,|y|),x(t,y)\bigr), \qquad t_0\leq t<1,\quad |y|\leq1-t,
	\end{equation}
	is well defined. Lemma \ref{Lem.normalized_equation} shows that $u_d^{\mathrm c}$ solves \eqref{Eq.wave_eq}.

	We verify the Cartesian smoothness of $u_d^{\mathrm c}$. By \eqref{Eq.even_odd_spaces}, there exist $C^\infty([0,1])$ functions $\Phi_*,R_*,X_*$ such that
	\[\widehat\phi_0(Z)=\Phi_*(Z^2),\qquad \widehat\rho_0(Z)=R_*(Z^2),\qquad X(Z)=ZX_*(Z^2),\qquad X_*(0)=\beta-1.\]
	With $s:=|y|^2/(1-t)^2$, the arguments in \eqref{Eq.direct_cone_solution} are
	\begin{equation}\label{Eq.Cartesian_reconstruction}
		\rho_0(t,|y|)=(1-t)^{-2\beta/(p-1)}R_*(s),\qquad \sigma=(1-t)^{1-\beta}\Phi_*(s),\qquad x=(1-t)^{-1}X_*(s)y.
	\end{equation}
	For $t_0\leq t<1$, $\sigma(t,y)$, $\rho_0(t,|y|)$ and $x(t,y)$ are smooth, including at $y=0$. The reconstruction $u_d^{\mathrm c}=\rho_0z(\sigma,x)$ is therefore smooth in the interior of the cone. At a point of the lateral boundary with $t<1$, the functions $\Phi_*,R_*,X_*$ extend to a neighborhood of $s=1$. Moreover, the corresponding normalized point has $\sigma\geq2\sigma_0$ and $x\in\partial\Omega$, and $z$ has a local smooth extension there, by $z\in C^\infty([\sigma_0,\infty)\times\overline\Omega)$. Composing the local extensions of $\Phi_*,R_*,X_*$ and $z$ in \eqref{Eq.Cartesian_reconstruction} gives a smooth extension of $u^{\mathrm c}_{d}$ across the lateral boundary point. At $t=t_0$ we have $\sigma>\sigma_0$, so the reconstruction extends smoothly across the initial slice as well. Thus $u^{\mathrm c}_{d}$ is smooth up to the boundary of every compact truncated cone. Radiality follows directly from rotational invariance.

	Translate the initial time to zero and define
	\[U_{d}(t,y)=u^{\mathrm c}_{d}(t+t_0,y),\qquad 0\leq t<h_0,\quad |y|\leq h_0-t.\]
	For $j\in\{0,1\}$, we have $\partial_t^jU_{d}(0,\cdot)\in C^\infty(\overline B_{h_0}^{d})$. We extend $\partial_t^jU_d(0,\cdot)$ to compactly supported radial functions on $\mathbb R^{d}$: for each $j\in\{0,1\}$, pick $u_{j,d}\in C_c^\infty(\R^{d};\R)$ such that
	\[u_{j,d}\ \text{is radial},\quad \supp u_{j,d}\subset B_{2h_0}^{d},\qquad u_{j,d}(y)=\partial_t^jU_{d}(0,y)\quad\text{for }|y|\leq h_0.\]
	Let $u_{d}$ be the unique smooth solution with initial data $(u_{0,d},u_{1,d})$, and let $T_*>0$ be the maximal existence time of $u_d$. The local existence theorem and the continuation criterion \eqref{Eq.continuation} apply, since the data are smooth and compactly supported; see \cite{Luk,Sogge}. The finite speed of propagation implies that $u_{d}$ agrees with $U_{d}$ in the backward cone, i.e.,
	\begin{equation}\label{Eq.cone_agreement}
		u_{d}(t,y)=U_{d}(t,y),\qquad \forall\ 0\leq t<\min\{T_*,h_0\},\ |y|\leq h_0-t.
	\end{equation}
    
	We next prove $T_*\leq h_0$. At the origin of the direct-dimensional cone, \eqref{Eq.profiles} gives
	\[\phi_0(t+t_0,0)=\frac{(h_0-t)^{1-\beta}}{\beta-1},\qquad \rho_0(t+t_0,0)=(h_0-t)^{-2\beta/(p-1)}.\]
	Set
	\[n_0:=1+\left\lfloor\frac{h_0^{1-\beta}}{2L(\beta-1)}\right\rfloor, \qquad h_n:=\bigl(2Ln(\beta-1)\bigr)^{-1/(\beta-1)},\qquad t_n:=h_0-h_n\quad(n\geq n_0).\]
	Then $0<t_n<h_0$, $t_n\uparrow h_0$, and $\phi_0(t_n+t_0,0)=2Ln$. Since $f(2Ln)=1$, Theorem \ref{Thm.normalized} gives
	\begin{equation}\label{Eq.origin_blowup_sequence}
		\left|h_n^{2\beta/(p-1)}U_{d}(t_n,0)-1\right| =|z(2Ln,0)-f(2Ln)| \leq\frac{C(1+\log(2Ln))^M}{2Ln}\longrightarrow0.
	\end{equation}
	In particular, $|U_{d}(t_n,0)|\to\infty$.

	Suppose on the contrary that $T_*>h_0$. Then $u_{d}$ is continuous on the compact set $[0,h_0]\times\overline B_{h_0}^{d}$, and consequently
	\[\max_{[0,h_0]\times\overline B_{h_0}^{d}}|u_{d}|<\infty.\]
	On the other hand, \eqref{Eq.cone_agreement} and \eqref{Eq.origin_blowup_sequence} give $|u_{d}(t_n,0)|\to\infty$, a contradiction. Therefore $0<T_*\leq h_0<\infty$, and \eqref{Eq.continuation} applies. The constructed data are real-valued and radial. Theorem \ref{Thm.conditional} follows.
\end{proof}

\subsection{The analytic spaces and the approximate solution}\label{Subsec.approximation}

Throughout the conditional construction, $d\geq4$, $p\geq3$ is odd, and the background is fixed. We define
\begin{equation}\label{Eq.dimension_indices}
	s_d:=\left\lfloor\frac d2\right\rfloor+1,\qquad s_*:=s_d+1.
\end{equation}
For a finite-dimensional Euclidean space $E$ or a complex Euclidean space, define
\begin{equation}\label{Eq.analytic_norm}
	\begin{aligned}
		\|h\|_{\mathcal B_\varrho}^2 &:=\sum_{N=0}^\infty(1+N)^{2s_*}\frac{\varrho^{2N}}{(N!)^2} \|\nabla^Nh\|_{L^2(\Omega;E)}^2,\\
		\mathcal B_\varrho(\Omega;E) &:=\left\{h\in\bigcap_{N\geq0}H^N(\Omega;E):\|h\|_{\mathcal B_\varrho}<\infty\right\}.
	\end{aligned}
\end{equation}
We write $\cX_\varrho:=\mathcal B_\varrho\times\mathcal B_\varrho$. Scalar unknowns and sources are radial and real-valued unless a complexification is explicitly used. Since $s_*>s_d+1/2$, we have the multiplication estimate
\begin{equation}\label{Eq.fixed_algebra}
	\|h_1h_2\|_{\mathcal B_\varrho}\leq C_\varrho \|h_1\|_{\mathcal B_\varrho}\|h_2\|_{\mathcal B_\varrho}.
\end{equation}

For a known coefficient $h$, we write
\[\mathfrak M_\varrho(h):= \sum_{N=0}^\infty(1+N)^{s_*}\frac{\varrho^N}{N!}\|\nabla^Nh\|_{L^\infty(\Omega)}.\]
For $0<r<r'$ and $h\in\mathcal B_{r'}$, Lemma \ref{Lem.analytic_to_multiplier} gives $\mathfrak M_r(\nabla^qh)\leq C_{q,r,r'}\|h\|_{\mathcal B_{r'}}$ for all $q\in\Z_{\geq0}$. In particular, \eqref{Eq.common_coefficient_radius} implies finite multiplier norms for the background coefficients and every fixed spatial derivative at each analytic radius below $R_v$.

An estimate $h=\Order_{\rm an}(a(\sigma))$ at an analytic radius $\varrho$ means $\mathfrak M_\varrho(h)\leq Ca(\sigma)$ for sufficiently large $\sigma$. Whenever additional derivatives of a known coefficient are needed, we first work at a strictly larger analytic radius.

If $0<\varrho<\varrho'$, every fixed number of derivatives is bounded from $\mathcal B_{\varrho'}$ to $\mathcal B_\varrho$, and the inclusion $\mathcal B_{\varrho'}\hookrightarrow\mathcal B_\varrho$ is compact. In particular,
\begin{equation}\label{Eq.analytic_radius_recovery}
	\|\nabla^qh\|_{\mathcal B_\varrho} \leq \frac{q!}{(\varrho'-\varrho)^q}\|h\|_{\mathcal B_{\varrho'}},\quad\forall\ q\in\Z_{\geq0}.
\end{equation}
The same inclusion recovers any fixed polynomial weight in the derivative order. The inclusion \eqref{Eq.analytic_radius_recovery}, the compact embedding and \eqref{Eq.fixed_algebra} are proved in Appendix \ref{App.analytic_spaces}. We shall use an analytic radius gap to recover the derivatives in the averaging argument.

For later use, define
\begin{equation}\label{Eq.geometry_budget}
	\gamma_v:=\min_{0\leq X\leq X_1}b'(X)>0,\qquad \mathcal L_v(r):=\sum_{q\geq2}\frac{r^{q-1}}{q!}\|\nabla^q\bvec\|_{L^\infty} \quad(0<r<R_v).
\end{equation}
The series converges locally uniformly, and $\mathcal L_v(r)\to0$ as $r\downarrow0$.

For the approximate solution, we choose an outer analytic radius $\overline\varrho>0$ such that
\begin{equation}\label{Eq.phase_radius_condition}
	4\overline\varrho<R_v,\qquad \mathcal L_v(4\overline\varrho)<\gamma_v,
\end{equation}
and define
\begin{equation}\label{Eq.approximation_radii}
	\varrho_1:=2\overline\varrho,\qquad \varrho_2:=4\overline\varrho.
\end{equation}
Since $\mathcal L_v(r)\to0$ as $r\downarrow0$, an analytic radius $\overline\varrho$ satisfying \eqref{Eq.phase_radius_condition} exists. We construct the phase in $\mathcal B_{\varrho_2}$ and estimate the approximate solution in $\mathcal B_{\varrho_1}$.

For an integer $j\geq0$, define the antiperiodic spaces with even and odd subspaces by
\begin{align*}
	C^j_{\operatorname{ap},L}&:=\{h\in C^j(\R;\R):h(\theta+L)=-h(\theta)\ \text{for every }\theta\in\R\},\\
	C^j_{\operatorname{ap},L,\rm e}&:=\{h\in C^j_{\operatorname{ap},L}:h(-\theta)=h(\theta)\ \text{for every }\theta\in\R\},\\
	C^j_{\operatorname{ap},L,\rm o}&:=\{h\in C^j_{\operatorname{ap},L}:h(-\theta)=-h(\theta)\ \text{for every }\theta\in\R\}.
\end{align*}
We equip $C^j_{\operatorname{ap},L}$ with the norm
\[\norm h_{C^j_{\operatorname{ap},L}}:=\max_{0\leq i\leq j}\sup_{0\leq\theta\leq 2L}\left|h^{(i)}(\theta)\right|.\]
Write $C^\infty_{\operatorname{ap},L}:=\bigcap_{j\geq0}C^j_{\operatorname{ap},L}$, with the corresponding definitions for the parity subspaces. We also set
\begin{align*}
	&C^\omega_{\operatorname{ap},L}(\R):=C^\infty_{\operatorname{ap},L}\cap C^\omega(\R),\\
	C^\omega_{\operatorname{ap},L,\rm e}(\R):= C^\omega_{\operatorname{ap},L}&(\R)\cap C^\omega_{\rm e}(\R),\qquad C^\omega_{\operatorname{ap},L,\rm o}(\R):=C^\omega_{\operatorname{ap},L}(\R)\cap C^\omega_{\rm o}(\R).
\end{align*}
Define
\[L_f:=\partial_\theta^2+pf(\theta)^{p-1}:C^{j+2}_{\operatorname{ap},L}\longrightarrow C^j_{\operatorname{ap},L},\qquad \av h:=\frac1{2L}\int_0^{2L}h(\theta)\,\mathrm d\theta.\]
The kernel and the solvability condition are
\begin{equation}\label{Eq.cell_range}
	\ker L_f=\operatorname{span}\{f'\},\qquad\operatorname{Ran}L_f=\left\{h\in C^j_{\operatorname{ap},L}:\av{hf'}=0\right\}.
\end{equation}
Indeed, besides $f'$, differentiating the amplitude-scaled oscillator gives the solution $f+(p-1)\theta f'/2$ of $L_f h=0$, which is clearly not antiperiodic. The two solutions are linearly independent at $\theta=0$. On $[0,L]$, the boundary conditions are $h(L)=-h(0)$ and $h'(L)=-h'(0)$. The operator $L_f$ on $L^2(0,L)$, with domain $\{h\in H^2(0,L):h(L)=-h(0),\ h'(L)=-h'(0)\}$, is self-adjoint. The Fredholm alternative for the antiperiodic boundary-value problem for $L_f$ gives the range in \eqref{Eq.cell_range}. Reflection in $\theta=0$ commutes with $L_f$. Since $f'$ is odd, restriction to the even subspace is invertible, and the range on the odd subspace is precisely the orthogonal complement of $f'$.

By Lemma \ref{Lem.normalized_equation}, the linearization of $\cN$ at $f(\sigma)$ is
\begin{align*}
	\left.\frac{\mathrm d}{\mathrm da}\cN(f(\sigma)+ah)\right|_{a=0} &=D^2h+pf(\sigma)^{p-1}h+\sigma^{-2}(A_0-V)h\\
	&=(\partial_\sigma^2+pf(\sigma)^{p-1})h+2\sigma^{-1}B\partial_\sigma h +\sigma^{-2}(B^2-B+A_0-V)h.
\end{align*}
The terms without a negative power of $\sigma$ therefore give $L_f$ in the fast variable. Since \eqref{Eq.Nf} has leading error $-\sigma^{-2}Vf(\sigma)$, we choose the unique function
\begin{equation}\label{Eq.g}
	g\in C^\infty_{\operatorname{ap},L,\rm e},\qquad L_fg=f.
\end{equation}
By \eqref{Eq.cell_range}, $L_f$ is invertible on the even antiperiodic subspace. Thus \eqref{Eq.g} has a unique solution, and the analytic ODE gives $g\in C^\omega_{\operatorname{ap},L,\rm e}(\R)$. For $a=a(x)$ and $h=h(\sigma)$, we have
\[D^2(\sigma^{-2}ah) =\sigma^{-2}ah''+2\sigma^{-3}(B-2)ah' +\sigma^{-4}(B-3)(B-2)ah.\]
Consequently,
\begin{align*}
	&\cN\bigl(f(\sigma)+\sigma^{-2}Vg(\sigma)\bigr)=2\sigma^{-3}(B-2)Vg'(\sigma)\\
	&\qquad\qquad\qquad+\sigma^{-4}\bigl((B-3)(B-2)+A_0-V\bigr)Vg(\sigma)+\sum_{j=2}^p\binom pj\sigma^{-2j}V^jf(\sigma)^{p-j}g(\sigma)^j.
\end{align*}
Here the order-$\sigma^{-2}$ term vanishes by $L_fg=f$. In particular,
\begin{equation}\label{Eq.first_correction}
	\cN\bigl(f(\sigma)+\sigma^{-2}Vg(\sigma)\bigr) =2\sigma^{-3}(B-2)Vg'(\sigma)+\Order(\sigma^{-4}).
\end{equation}
We measure the remainders in $\mathcal B_{\varrho_1}$, with constants depending on the fixed background. The bounds for the completed approximation are stated in Proposition \ref{Prop.approximation}.

As we can see from the above paragraph, the leading error $-\sigma^{-2}Vf(\sigma)$ in \eqref{Eq.Nf} can be canceled by introducing the corrector $g$ solving \eqref{Eq.g}. Naturally, one may use the same idea to eliminate the leading error $2\sigma^{-3}(B-2)Vg'(\sigma)$ in \eqref{Eq.first_correction}. However, this procedure fails immediately. The obstruction is explicit. The oscillator energy and integration by parts give
\begin{equation}\label{Eq.oscillator_averages}
	\av{(f')^2}=\frac2{p+3},\qquad \av{f'g'}=\frac{\av{f^2}}{p-1}.
\end{equation}
For the second identity, we use $L_ff=(p-1)f^p$ and $L_fg=f$ to obtain
\[\av{f'g'}=\av{f^pg}=\frac1{p-1}\av{(L_ff)g}=\frac1{p-1}\av{fL_fg}=\frac{\av{f^2}}{p-1}.\]
In particular, $g'\notin\operatorname{Ran}L_f$. We fix
\begin{equation}\label{Eq.kappa_H}
	c_a:=\ell^{-1/2},\qquad \kappa_f:=\frac{\av{f'g'}}{\av{(f')^2}}=\frac\ell2\av{f^2},\qquad L_fH=g'-\kappa_ff',\quad \av{Hf'}=0.
\end{equation}
Equation \eqref{Eq.cell_range} gives a unique $H\in C^\omega_{\operatorname{ap},L,\rm o}(\R)$. For the subsequent linearization, we also define
\begin{equation}\label{Eq.K_cancellation}
	K(\theta):=p(p-1)f(\theta)^{p-2}g(\theta),\qquad L_fg'=(1-K)f',\qquad \av{(K-1)(f')^2}=0.
\end{equation}
The second identity follows by differentiating \eqref{Eq.g}, and the third follows by taking its inner product with $f'$.


The nonzero component in the direction $f'$ is removed by a phase correction. Let $\tau\in\R$ be an independent ``slow'' variable and set
\[T_\nu:=\partial_\tau+B-\nu,\qquad \nu\in\R.\]
For every smooth $h=h(\tau,x)$, direct differentiation gives
\begin{equation}\label{Eq.slow_calculus}
	\begin{aligned}
		D\bigl(\sigma^{-\nu}h(\log\sigma,x)\bigr)&=\sigma^{-\nu-1}(T_\nu h)(\log\sigma,x),\\
		\cP\bigl(\sigma^{-\nu}h(\log\sigma,x)\bigr)&=\sigma^{-\nu-2}(T_{\nu+1}T_\nu+A_0)h(\log\sigma,x).
	\end{aligned}
\end{equation}
To derive the equation for the phase correction, let $\Psi(\tau,x)$ initially be a real-valued polynomial in $\tau$ with radial coefficients in $\mathcal B_{\varrho_2}$. Define
\begin{equation}\label{Eq.approximate_phase}
	\theta_\sigma(\sigma,x):=\sigma+\sigma^{-1}\Psi(\log\sigma,x),\qquad U_\Psi:=T_1\Psi,
\end{equation}
and
\begin{equation}\label{Eq.phase_velocity}
	\lambda(\sigma,x):=1+\sigma^{-2}U_\Psi(\log\sigma,x)=D\theta_\sigma.
\end{equation}
The symbol $\theta$ below denotes an independent fast variable, whereas $\theta_\sigma$ denotes the function in \eqref{Eq.approximate_phase}; the subscript in $\theta_\sigma$ is not a derivative. For functions of $(\theta;\tau,x)$, the operators $T_\nu$, $A_0$ and $\nabla_x$ act with $\theta$ fixed. A slow function appearing without arguments in an identity in $(\sigma,x)$ is evaluated at $(\log\sigma,x)$. The phase $\theta_\sigma$ is used only in compositions and does not replace the time coordinate $\sigma$.

The chain rule and \eqref{Eq.slow_calculus} give the exact identity
\begin{equation}\label{Eq.shifted_leading_error}
	\begin{aligned}
		\cN\bigl(f(\theta_\sigma)\bigr)&=-\sigma^{-2}\bigl(2U_\Psi f(\theta_\sigma)^p+Vf(\theta_\sigma)\bigr)+\sigma^{-3}\bigl((T_2T_1+A_0)\Psi\bigr)f'(\theta_\sigma)\\
		&\quad-\sigma^{-4}\bigl(U_\Psi^2-|\nabla\Psi|^2\bigr)f(\theta_\sigma)^p.
	\end{aligned}
\end{equation}
Thus a second-order correction must solve $L_fA_2=2U_\Psi f^p+Vf$. Since $L_ff=(p-1)f^p$ and $L_fg=f$, the unique even antiperiodic solution gives
\begin{equation}\label{Eq.second_approximation_definition}
	\begin{aligned}
		A_2(\theta;\tau,x)&:=\frac2{p-1}U_\Psi(\tau,x)f(\theta)+V(x)g(\theta),\\
		z_2(\sigma,x)&:=f(\theta_\sigma)+\sigma^{-2}A_2(\theta_\sigma;\log\sigma,x).
	\end{aligned}
\end{equation}
The coefficient of $\sigma^{-3}$ in $\cN(z_2)$ is
\begin{equation}\label{Eq.third_residual}
	\begin{aligned}
		R_3(\theta;\tau,x)&:=\bigl((T_2T_1+A_0)\Psi\bigr)(\tau,x)f'(\theta)+2(T_2\partial_\theta A_2)(\theta;\tau,x)\\
		&=\bigl(\ell T_2T_1\Psi+A_0\Psi\bigr)(\tau,x)f'(\theta)+2((B-2)V)(x)g'(\theta).
	\end{aligned}
\end{equation}
The expansion and its analytic remainder are proved in Subsection \ref{Subsec.approximation_proof}. By \eqref{Eq.cell_range}, the equation $L_fA_3=-R_3$ is solvable precisely when $\av{R_3f'}=0$, where the average is taken in $\theta$ with $(\tau,x)$ fixed. Indeed,
\[\av{R_3f'}=(\ell T_2T_1\Psi+A_0\Psi)\av{(f')^2}+2(B-2)V\av{g'f'}.\]
Using \eqref{Eq.oscillator_averages}, its vanishing is equivalent to the phase equation \eqref{Eq.phase_equation} below.

\begin{proposition}[The phase correction]\label{Prop.phase}
	Under Assumption \ref{Ass.background} and \eqref{Eq.phase_radius_condition}, there are $q_*\in\Z_{\geq0}$ and real-valued radial $\Psi_j\in\mathcal B_{\varrho_2}$ such that
	\[\Psi(\tau,x):=\sum_{j=0}^{q_*}\tau^j\Psi_j(x)\]
	satisfies
	\begin{equation}\label{Eq.phase_equation}
		(T_2T_1+c_a^2A_0)\Psi=-\av{f^2}(B-2)V.
	\end{equation}
	The equation holds in spatial distributions for every $\tau\in\R$, and hence pointwise by analyticity.
\end{proposition}

We prove Proposition \ref{Prop.phase} in Subsection \ref{Subsec.phase_proof} by reducing \eqref{Eq.phase_equation} to the averaged acoustic system and using the compact resolvent of $\cA_{\varrho_2}$ to construct a phase correction that is polynomial in $\tau$. From now on, fix a polynomial $\Psi$ given by this proposition and use \eqref{Eq.approximate_phase}--\eqref{Eq.second_approximation_definition} with this choice. Then
\[R_3(\theta;\tau,x)=2((B-2)V)(x)\bigl(g'(\theta)-\kappa_ff'(\theta)\bigr),\qquad \av{R_3f'}=0.\]

\begin{lemma}\label{Lem.second_approximation}
	For the fixed phase $\Psi$, there are $\sigma_2\geq2$, $C>0$ and $M_2\in\Z_{\geq0}$ such that, for $\sigma\geq\sigma_2$ and $x\in\overline\Omega$,
	\begin{equation}\label{Eq.second_approximation_domain}
		\tfrac12\leq\lambda(\sigma,x)\leq\tfrac32,
	\end{equation}
	and
	\begin{equation}\label{Eq.second_approximation_bounds}
		\begin{aligned}
			\|z_2-f(\theta_\sigma)\|_{\mathcal B_{\varrho_1}}&\leq C\sigma^{-2}(1+\log\sigma)^{M_2},\\
			\|\cN(z_2)-\sigma^{-3}R_3(\theta_\sigma;\log\sigma,x)\|_{\mathcal B_{\varrho_1}}&\leq C\sigma^{-4}(1+\log\sigma)^{M_2}.
		\end{aligned}
	\end{equation}
\end{lemma}

The remaining third-order source is removed by an odd antiperiodic correction. The fourth-order source is even and can therefore be removed without another phase correction. The resulting finite approximation has the bounds needed for the nonlinear construction.

\begin{proposition}[The fourth-order approximation]\label{Prop.approximation}
	There are $\sigma_{\app}\geq\sigma_2$, $C>0$, an integer $M\geq0$, and a smooth real-valued radial function $z_{\app}$ on $[\sigma_{\app},\infty)\times\overline\Omega$ such that
	\begin{equation}\label{Eq.approximation_bounds}
		\begin{aligned}
			\|z_{\app}-f(\theta_\sigma)\|_{\mathcal B_{\varrho_1}}&\leq C\sigma^{-2}(1+\log\sigma)^M,\\
			\|\cN(z_{\app})\|_{\mathcal B_{\varrho_1}}&\leq C\sigma^{-5}(1+\log\sigma)^M.
		\end{aligned}
	\end{equation}
	The same estimates hold with $\|\cdot\|_{\mathcal B_{\varrho_1}}$ replaced by $\mathfrak M_{3\overline\varrho}(\cdot)$, after increasing $C$. Inequality \eqref{Eq.second_approximation_domain} holds for $\sigma\geq\sigma_{\app}$. For every $q\in\Z_{\geq0}$ and $0<r<\varrho_1$, there is $C_{q,r}>0$ such that
	\begin{align*}
		\|\nabla^q(z_{\app}-f(\theta_\sigma))\|_{\mathcal B_r}&\leq C_{q,r}\sigma^{-2}(1+\log\sigma)^M,\\
		\|\nabla^q\cN(z_{\app})\|_{\mathcal B_r}&\leq C_{q,r}\sigma^{-5}(1+\log\sigma)^M.
	\end{align*}
\end{proposition}

The third- and fourth-order corrections and the proof of Proposition \ref{Prop.approximation} are given in Subsection \ref{Subsec.approximation_proof}.

\subsection{The perturbation equation and the new unknown}\label{Subsec.formulation}

We derive the equation for the perturbation of $z_{\app}$. We first estimate the phase coefficients and the scalar remainder in Lemmas \ref{Lem.phase_coefficients} and \ref{Lem.remainder_equation}. We then introduce a two-component unknown which separates the periodic frame from the unknown modulation; Lemmas \ref{Lem.state_equation} and \ref{Lem.state_remainder} give the equation for $W$ and the decomposition into averaged, periodic and small-coefficient operators.

Fix $\Psi$ and $z_{\app}$ from Propositions \ref{Prop.phase} and \ref{Prop.approximation}, and use $\lambda$ from \eqref{Eq.phase_velocity}. The identities needed later are
\[D\theta_\sigma=\lambda,\qquad D\lambda=\sigma^{-3}T_2T_1\Psi,\qquad \nabla\theta_\sigma=\sigma^{-1}\nabla\Psi.\]
Since $\theta_\sigma-\sigma=\sigma^{-1}\Psi$, the phase correction tends to zero in the multiplier norm at analytic radius $\varrho_1$. We record the corresponding coefficient bounds, including the positivity of $\lambda$ needed for the change of unknown.

\begin{lemma}\label{Lem.phase_coefficients}
	There are $\sigma_c\geq\sigma_{\app}$ and an integer $M\geq0$ such that, with
	\begin{equation}\label{Eq.nu_sigma}
		\nu(\sigma):=\sigma^{-1}(1+\log\sigma)^M,
	\end{equation}
	one has, for $\sigma\geq\sigma_c$,
	\[\tfrac12\leq\lambda\leq\tfrac32,\quad \mathfrak M_{\varrho_1}(\lambda^{-1})+\mathfrak M_{\varrho_1}(z_{\app})\leq C,\]
	\[\mathfrak M_{\varrho_1}(\lambda-1)+\sigma\mathfrak M_{\varrho_1}(D\lambda) +\mathfrak M_{\varrho_1}(z_{\app}-f(\theta_\sigma))\leq C\sigma^{-1}\nu(\sigma).\]
	For $h\in\{f,f',f'',K\}$,
	\begin{equation}\label{Eq.periodic_composition_bounds}
		\mathfrak M_{\varrho_1}(h(\theta_\sigma))\leq C,\qquad \mathfrak M_{\varrho_1}(h(\theta_\sigma)-h(\sigma))\leq C\nu(\sigma).
	\end{equation}
\end{lemma}

With the phase coefficients controlled, we write the equation for the scalar perturbation. We separate the residual of the approximate solution from the terms which are at least quadratic in $\zeta$. For $z=z_{\app}+\zeta$, define $\cR_{\app}:=\cN(z_{\app})$ and
\begin{equation}\label{Eq.remainder_equation}
	\cL_\sigma:=D^2+\sigma^{-2}(A_0-V)+pz_{\app}^{p-1},\quad \cQ_\sigma(\zeta):=(z_{\app}+\zeta)^p-z_{\app}^p-pz_{\app}^{p-1}\zeta.
\end{equation}
Then $\cN(z)=0$ is equivalent to $\cL_\sigma\zeta=-\cR_{\app}-\cQ_\sigma(\zeta)$. We also define
\begin{equation}\label{Eq.linearized_original}
	R_\sigma:=\sigma^2\bigl(pz_{\app}^{p-1}-\lambda^2pf(\theta_\sigma)^{p-1}\bigr)-K(\theta_\sigma)V,\quad
	\Lambda_\sigma:=A_0+(K(\theta_\sigma)-1)V+R_\sigma.
\end{equation}
Thus
\begin{equation}\label{Eq.linearized_form}
	\cL_\sigma=D^2+\lambda^2pf(\theta_\sigma)^{p-1}+\sigma^{-2}\Lambda_\sigma.
\end{equation}
The term $\lambda^2pf(\theta_\sigma)^{p-1}$ retains the principal periodic potential. The following estimates control the coefficient $R_\sigma$, the prescribed source $\cR_{\app}$, and the nonlinear remainder $\cQ_\sigma$ separately.

\begin{lemma}\label{Lem.remainder_equation}
	After increasing $M$ and $\sigma_c$,
	\[\mathfrak M_{\varrho_1}(R_\sigma)\leq C\nu(\sigma),\qquad \|\cR_{\app}(\sigma)\|_{\mathcal B_{\varrho_1}} \leq C\sigma^{-5}(1+\log\sigma)^M.\]
	For $0<r<\overline\varrho$ and $\|\zeta_i\|_{\mathcal B_r}\leq1$,
	\begin{equation}\label{Eq.quadratic_remainder_bounds}
		\begin{aligned}
			\|\cQ_\sigma(\zeta_1)\|_{\mathcal B_r}&\leq C\|\zeta_1\|_{\mathcal B_r}^2,\\
			\|\cQ_\sigma(\zeta_1)-\cQ_\sigma(\zeta_2)\|_{\mathcal B_r} &\leq C(\|\zeta_1\|_{\mathcal B_r}+\|\zeta_2\|_{\mathcal B_r}) \|\zeta_1-\zeta_2\|_{\mathcal B_r}.
		\end{aligned}
	\end{equation}
	The constants are independent of $\sigma\geq\sigma_c$ and locally uniform for $r\in(0,\overline\varrho)$.
\end{lemma}

Lemma \ref{Lem.remainder_equation} reduces the construction to solving the linearized scalar equation with a small source. The potential $\lambda^2pf(\theta_\sigma)^{p-1}$ depends on the fast phase. To separate the periodic coefficient from the unknown modulation, we use a frame built from the zero-frequency oscillator equation and define
\begin{equation}\label{Eq.Pzero}
	P_0(\theta):=
	\begin{pmatrix}
		f'(\theta)&\frac2{p-1}f(\theta)\\
		f''(\theta)&\frac{p+1}{p-1}f'(\theta)
	\end{pmatrix}
	,\qquad \det P_0=\frac2{p-1}.
\end{equation}
Define
\[\mathsf P_0(\theta):=\begin{pmatrix}0&1\\-pf(\theta)^{p-1}&0\end{pmatrix},\qquad N_0:=\begin{pmatrix}0&1\\0&0\end{pmatrix}.\]
The oscillator equation gives
\begin{equation}\label{Eq.frame_identity}
	P_0'=\mathsf P_0P_0-P_0N_0,\qquad P_0^{-1}\binom01=\binom{-f}{(p-1)f'/2}.
\end{equation}
The first identity changes $\partial_\theta U=\mathsf P_0(\theta)U$ into $\partial_\theta(P_0^{-1}U)=N_0P_0^{-1}U$, removing the periodic coefficients from the homogeneous oscillator. In the proof of Lemma \ref{Lem.state_equation}, the equality $DP_0(\theta_\sigma)=\lambda P_0'(\theta_\sigma)$ gives the same cancellation for $(\zeta,\lambda^{-1}D\zeta)^T$, and the second identity determines the source term in the equation for $W$.
For an arbitrary pair $W=(\psi,\pi)^T$, set
\begin{equation}\label{Eq.zeta-w-W-def}
	\binom{\zeta_W}{w_W}:=P_0(\theta_\sigma)\operatorname{diag}(1,\sigma^{-1})W.
\end{equation}

For a solution of $\cL_\sigma\zeta=\mathcal F$, the equivalent definition is
\begin{equation}\label{Eq.state}
	W:=\operatorname{diag}(1,\sigma)P_0(\theta_\sigma)^{-1} \binom\zeta{\lambda^{-1}D\zeta}.
\end{equation}
Since $\det P_0=2/(p-1)>0$, the change of variables \eqref{Eq.state} is invertible for every phase $\theta_\sigma\in\R$. We define $h:=\sigma^2\mathcal F$ and
\[\mathfrak R_\sigma(W,h):=\lambda^{-1}(h-\Lambda_\sigma\zeta_W) -\sigma^2\frac{D\lambda}{\lambda}w_W.\]
The next lemma gives the equation for $W$ before any averaging is performed. In particular, the identities $D\zeta_W=\lambda w_W$ and $\cL_\sigma\zeta_W=\mathcal F$ will follow from the two-component equation.

\begin{lemma}[The exact system for the new unknowns]\label{Lem.state_equation}
	The scalar equation $\cL_\sigma\zeta=\mathcal F$ is equivalent to
	\begin{equation}\label{Eq.exact_state}
		\begin{cases}
			(\sigma\partial_\sigma+B)\psi=\lambda\pi-\sigma^{-1}f(\theta_\sigma)\mathfrak R_\sigma(W,h),\\
			(\sigma\partial_\sigma+B)\pi=\pi+\dfrac{p-1}{2}f'(\theta_\sigma)\mathfrak R_\sigma(W,h).
		\end{cases}
	\end{equation}
	Conversely, reconstruction by \eqref{Eq.zeta-w-W-def} from a smooth solution of \eqref{Eq.exact_state} gives $D\zeta_W=\lambda w_W$ and $\cL_\sigma\zeta_W=\mathcal F$.
\end{lemma}

We now separate the coefficients in \eqref{Eq.exact_state} into the periodic averages and the mean-zero parts. The averaged operator will be independent of $\sigma$, while the remaining small coefficients will be estimated using Lemma \ref{Lem.phase_coefficients}. Define
\[a_0(\theta):=\tfrac{p-1}{2}(f'(\theta))^2,\qquad a_1(\theta):=a_0(\theta)(K(\theta)-1), \qquad \av{a_0}=c_a^2,\quad\av{a_1}=0,\]
\[\mathscr A:=\begin{pmatrix}-B&I\\-c_a^2A_0&I-B\end{pmatrix},\qquad \cB(\sigma):=\begin{pmatrix}0&0\\-(a_0(\sigma)-c_a^2)A_0-a_1(\sigma)V&0\end{pmatrix},\]
\[\mathbf g_\sigma:=\lambda^{-1} \binom{-\sigma^{-1}f(\theta_\sigma)}{(p-1)f'(\theta_\sigma)/2}.\]

\begin{lemma}\label{Lem.state_remainder}
	There is a linear differential operator $\mathscr E_\sigma$, independent of the source $h$ and of spatial order at most two, such that \eqref{Eq.exact_state} takes the form
	\begin{equation}\label{Eq.state_operator_form}
		\sigma\partial_\sigma W=(\mathscr A+\cB(\sigma)+\mathscr E_\sigma)W+\mathbf g_\sigma h.
	\end{equation}
	After increasing $\sigma_c$, every coefficient of $\mathscr E_\sigma$ has multiplier norm at $\overline\varrho$ bounded by $C\nu(\sigma)$ for $\sigma\geq\sigma_c$. For $0<r<r'<\overline\varrho$,
	\[\|\mathscr E_\sigma W\|_{\cX_r}\leq C_{r,r'}\nu(\sigma)\|W\|_{\cX_{r'}},\qquad \|\mathbf g_\sigma h\|_{\cX_r}\leq C\|h\|_{\mathcal B_r}.\]
	The constants are independent of $\sigma\geq\sigma_c$ and locally uniform for $0<r<r'<\overline\varrho$.
\end{lemma}

Thus the homogeneous system for $W$ consists of the fixed operator $\mathscr A$, the mean-zero periodic operator $\cB(\sigma)$, and the small differential remainder $\mathscr E_\sigma$. The periodic primitive of $\cB(\sigma)$ defines the averaging corrector in Subsection \ref{Subsec.fixedpoint}. We first study the evolution generated by $\mathscr A$. The proofs of Lemmas \ref{Lem.phase_coefficients}--\ref{Lem.state_remainder} are given in Subsection \ref{Subsec.perturbation_proofs}.

\subsection{The averaged equation}\label{Subsec.averaged}

The averages in \eqref{Eq.oscillator_averages} and \eqref{Eq.K_cancellation} give the equation for $W=(\psi,\pi)$. With a source $H=(H_1,H_2)$, we write
\begin{equation}\label{Eq.macro_pair}
	\partial_\tau W=\mathscr AW+H,\qquad (\partial_\tau+B)\psi=\pi+H_1,\qquad (\partial_\tau+B)\pi=\pi-c_a^2A_0\psi+H_2.
\end{equation}
The term $-c_a^2A_0\psi$ contains two spatial derivatives. Its pairing with $\pi$ cancels the highest-order term obtained by differentiating $c_a^2\|\nabla\psi\|_{L^2}^2$. We therefore use the space of pairs
\begin{equation}\label{Eq.acoustic_pair_norm}
	\cH_\varrho:=\{W=(\psi,\pi)\in\cX_\varrho:\nabla\psi\in\mathcal B_\varrho\},\quad \|W\|_{\cH_\varrho}^2:=\|\psi\|_{\mathcal B_\varrho}^2+\|\pi\|_{\mathcal B_\varrho}^2+c_a^2\|\nabla\psi\|_{\mathcal B_\varrho}^2.
\end{equation}
As in $\cX_\varrho$, both components are radial and real-valued unless the space is complexified. For $0<\varrho<\varrho'$, Lemma \ref{Lem.analytic_derivatives} gives
\begin{equation}\label{Eq.acoustic_radius_gap}
	\|W\|_{\cX_\varrho}\leq\|W\|_{\cH_\varrho}\leq\left(1+\frac{c_a^2}{(\varrho'-\varrho)^2}\right)^{1/2}\|W\|_{\cX_{\varrho'}}.
\end{equation}
In particular, the source norm in the energy estimate is
\[\|H\|_{\cH_\varrho}^2=\|H_1\|_{\mathcal B_\varrho}^2+\|H_2\|_{\mathcal B_\varrho}^2+c_a^2\|\nabla H_1\|_{\mathcal B_\varrho}^2.\]

The normalized transport field has the quantitative outgoing property
\begin{equation}\label{Eq.acoustic_outgoing_property}
	\operatorname{sym}\mathrm D_x\bvec\geq\gamma_vI_d\quad\text{on }\overline\Omega,\qquad \bvec\cdot n=1>c_a\quad\text{on }\partial\Omega,
\end{equation}
where $n$ is the outward unit normal. The interior inequality gives a favorable term proportional to the derivative order in the acoustic energy estimate. This term permits the analytic radius to increase. The boundary inequality makes the energy flux nonnegative, so no boundary condition is imposed. The two contributions of the outgoing property are computed directly in the proof of Lemma \ref{Lem.macro_finite_energy}.

\begin{proposition}[The averaged evolution]\label{Prop.macro}
	Let $T_*>0$ and let $\varrho\in C^1([0,T_*];(0,R_v))$ be nondecreasing, with
	\[\delta_{\rm ac}:=\inf_{0\leq\tau\leq T_*}\left(\gamma_v-\mathcal L_v(\varrho(\tau))-\frac{\varrho'(\tau)}{\varrho(\tau)}\right)>0.\]
	Let $W_0\in\cH_{\varrho(0)}$, and let $H:[0,T_*]\to\cH_{\varrho(0)}$ be strongly measurable and satisfy
	\[H(\tau)\in\cH_{\varrho(\tau)}\quad\text{for almost every }\tau,\qquad \int_0^{T_*}\|H(u)\|_{\cH_{\varrho(u)}}\,\mathrm du<\infty.\]
	Then there exists a distributional solution $W=(\psi,\pi)$ of \eqref{Eq.macro_pair} with $W(0)=W_0$ such that
	\begin{equation}\label{Eq.macro_solution_class}
		\begin{gathered}
			W\in\bigcap_{m\in\Z_{\geq0}}C([0,T_*];H^{m+1}(\Omega)\times H^m(\Omega)),\\
			W(t)\in\cH_{\varrho(t)}\quad(0\leq t\leq T_*),\qquad \sup_{0\leq t\leq T_*}\|W(t)\|_{\cH_{\varrho(t)}}<\infty.
		\end{gathered}
	\end{equation}
	The solution satisfies
	\begin{equation}\label{Eq.macro_energy}
		\|W(t)\|_{\cH_{\varrho(t)}}\leq C\left(\|W(s)\|_{\cH_{\varrho(s)}}+\int_s^t\|H(u)\|_{\cH_{\varrho(u)}}\,\mathrm du\right),\qquad 0\leq s\leq t\leq T_*,
	\end{equation}
	where $C$ may depend on the fixed background, $T_*$ and $\varrho$, but is independent of $W_0$ and $H$.
	
	More precisely, uniqueness holds in the class \eqref{Eq.macro_solution_class}: if $W^{(1)}$ and $W^{(2)}$ satisfy \eqref{Eq.macro_solution_class}, solve \eqref{Eq.macro_pair} in distributions with the same source $H$, and satisfy $W^{(1)}(0)=W^{(2)}(0)=W_0$, then
	\[W^{(1)}(t)=W^{(2)}(t)\qquad(0\leq t\leq T_*).\]
\end{proposition}

We next define the solution operators for the homogeneous averaged equation. Fix $r\in(0,R_v)$ with $\mathcal L_v(r)<\gamma_v$. For $W_0\in\cH_r$, apply Proposition \ref{Prop.macro} with the constant radius $\varrho(\tau)=r$ and $H=0$, and we define
\[\cS_r(t):\cH_r\longrightarrow\cH_r,\qquad \cS_r(t)W_0:=W(t),\qquad t\geq0.\]
The proposition applies on every finite time interval, and uniqueness identifies the solutions on overlapping time intervals. Thus $\cS_r(t)$ is defined for every $t\geq0$. Linearity of the equation and \eqref{Eq.macro_energy} show that each $\cS_r(t)$ is a bounded linear operator on $\cH_r$. Since $\mathscr A$ is independent of time, uniqueness also gives
\[\cS_r(0)=I,\qquad \cS_r(t+s)=\cS_r(t)\cS_r(s),\quad \forall\ s,t\geq0.\]

The solution operators are consistent on different analytic spaces. If $0<r<r'<R_v$, $\mathcal L_v(r')<\gamma_v$ and $W_0\in\cH_{r'}$, then the solution constructed in $\cH_{r'}$ also satisfies the hypotheses for uniqueness in $\cH_r$. Consequently,
\[\cS_r(t)W_0=\cS_{r'}(t)W_0,\quad\forall\ t\geq0.\]
We write $\cS(t)W_0$ for this common value when the analytic radius is not displayed. Thus $\cS(t)$ denotes the same homogeneous solution operator on all the spaces considered above, whereas $\cS_r(t)$ specifies the operator acting on the fixed space $\cH_r$. In particular, for $H=0$ and a radius path satisfying Proposition \ref{Prop.macro}, uniqueness and \eqref{Eq.macro_energy} give
\[W(t)=\cS(t-s)W(s),\qquad \|\cS(t-s)\|_{\cH_{\varrho(s)}\to\cH_{\varrho(t)}}\leq C,\qquad 0\leq s\leq t\leq T_*.\]

For $r\in(0,R_v)$ with $\mathcal L_v(r)<\gamma_v$, we define the operator $\cA_r$ on $\cH_r$ by
\begin{equation}\label{Eq.macro_generator_domain}
	\mathcal D(\cA_r):=\{W\in\cH_r:\mathscr AW\in\cH_r\},\qquad \cA_rW:=\mathscr AW,
\end{equation}
where the derivatives in $\mathscr AW$ are understood in distributions. The differential expression $\mathscr A$ is independent of $r$, whereas $\cA_r$ denotes its realization on $\cH_r$ with the domain in \eqref{Eq.macro_generator_domain}.

For spectra and resolvents, $\cH_r$ and $\cA_r$ are complexified without changing notation. We denote by $\mathcal K(\cH_r)$ the compact linear operators on $\cH_r$. For $\lambda\in\C$, we define
\[\mathcal G_\lambda(\cA_r):=\bigcup_{n\geq1}\ker(\cA_r-\lambda I)^n.\]
The next proposition proves strong continuity of $(\cS_r(t))_{t\geq0}$ on $\cH_r$, identifies its infinitesimal generator, and gives the compactness and spectral consequences of the analytic radius gain.

\begin{proposition}\label{Prop.macro_spectrum}
	Fix $r\in(0,R_v)$ with $\mathcal L_v(r)<\gamma_v$. The family $(\cS_r(t))_{t\geq0}$ is a strongly continuous semigroup on $\cH_r$, and
	\[\lim_{t\downarrow0}\frac{\cS_r(t)W-W}{t}=\cA_rW\quad\text{in }\cH_r,\quad \forall\  W\in\mathcal D(\cA_r).\]
	There are $M_r\geq1$ and $\omega_r\in\R$ such that
	\begin{equation}\label{Eq.macro_compact_resolvent}
		\begin{gathered}
			\|\cS_r(t)\|_{\cH_r\to\cH_r}\leq M_r\ee^{\omega_rt},\qquad \cS_r(t)\in\mathcal K(\cH_r)\quad(t>0),\\
			(zI-\cA_r)^{-1}\in\mathcal K(\cH_r)\quad(z\in\C\setminus\spec(\cA_r)),\qquad \{z:\operatorname{Re}z>\omega_r\}\subset\C\setminus\spec(\cA_r).
		\end{gathered}
	\end{equation}
	Each $\lambda\in\spec(\cA_r)$ is isolated. For some $r_+\in(r,R_v)$ with $\mathcal L_v(r_+)<\gamma_v$,
	\begin{equation}\label{Eq.macro_generalized_eigenspaces}
		0<\dim\mathcal G_\lambda(\cA_r)<\infty,\qquad \mathcal G_\lambda(\cA_r)\subset\cH_{r_+},\qquad \forall\ \lambda\in\spec(\cA_r).
	\end{equation}
\end{proposition}

The proofs of both propositions are given in Subsection \ref{Subsec.macro_proof}. All spectra and Riesz projections of real operators are understood after complexification.

\subsection{The periodic matrix family}\label{Subsec.periodic_matrix}

We construct the periodic matrix for the fixed exponent $p$ and record the estimates between sampled frequencies. Subsection \ref{Subsec.energy} uses $S(\theta,\xi,\delta)$ to define the analytic energy.

For the independent oscillator variable $\theta$, set
\[\mathsf P_\xi(\theta):=\begin{pmatrix}0&1\\ -pf(\theta)^{p-1}-\xi^2&0\end{pmatrix}, \qquad\forall\ \xi\geq0.\]
At $\xi=0$, the coefficient matrix $\mathsf P_\xi$ agrees with $\mathsf P_0$ from Subsection \ref{Subsec.formulation}. The family $P_\eta$ in Proposition \ref{Prop.matrix} extends the frame $P_0$ and the first identity in \eqref{Eq.frame_identity} to nonzero frequencies. 

To specify the growth constant, let $\theta_z:=3L/2$, so that $f(\theta_z)=0$ and $f'(\theta_z)>0$, and define the fundamental matrix $Y_p$ by
\[\partial_\theta Y_p(\theta,\xi)=\mathsf P_\xi(\theta+\theta_z)Y_p(\theta,\xi), \qquad Y_p(0,\xi)=I_2.\]
Define
\begin{equation}\label{Eq.cell_growth_constant}
	\begin{aligned}
		\mathsf M_p(\xi)&:=Y_p(L,\xi),\\
		b_{{\rm sc},p}(\xi)&:=\frac12\left(\xi(\mathsf M_p(\xi))_{12} +\xi^{-1}(\mathsf M_p(\xi))_{21}\right)\quad(\xi>0),\\
		\chi_p&:=\sup_{\xi>0}\frac{\arsinh|b_{{\rm sc},p}(\xi)|}{L\xi}.
	\end{aligned}
\end{equation}
Proposition \ref{Prop.matrix} proves $0<\chi_p<\infty$ and constructs the periodic energy matrix.

\begin{proposition}[The periodic matrix family]\label{Prop.matrix}
	Let $p\geq3$ be odd. Then $0<\chi_p<\infty$. For every $\chi>\chi_p$, there exist $\delta_0\in(0,1]$, $0<\xi_{\rm lo}<\xi_{\rm hi}$, $c,C>0$ and a family of symmetric matrices with real-valued entries
	\[S=\begin{pmatrix} S_{11} & S_{12}\\ S_{12} & S_{22} \end{pmatrix}\in C^\infty\bigl(\R\times[0,\infty)\times[0,\delta_0];\R^{2\times2}\bigr)\]
	such that $S(\theta+L,\xi,\delta)=S(\theta,\xi,\delta)$ and $S>0$ for $(\xi,\delta)\ne(0,0)$. For $\theta\in\R$, $\xi\geq0$ and $0<\delta\leq\delta_0$,
	\begin{equation}\label{Eq.S_coercivity}
		\begin{aligned}
			&\qquad\qquad\qquad\qquad \partial_\theta S+\mathsf P_\xi^TS+S\mathsf P_\xi\leq(2\chi\xi+C\delta)S,\\
			&\qquad\qquad\qquad\left\|S^{-1/2}\xi (\partial_\xi S) S^{-1/2}\right\|+\left\|S^{-1/2}\delta (\partial_\delta S) S^{-1/2}\right\|\leq C,\\
			&c\diag(\xi^2+\delta^2,\delta^2)\leq S\leq C\diag(1+\xi^2,1),\qquad S_{22}\leq C,\quad \xi S_{22}\leq\sqrt{\det S}.
		\end{aligned}
	\end{equation}
	More precisely, for $0\leq\xi\leq\xi_{\rm lo}$, we have\footnote{Please notice the notation: $P_{\eta}\neq \mathsf P_{\eta}$. }
	\[S=P_{\xi^2}^{-T}\diag(\omega(\xi^2)+\delta^2,1)P_{\xi^2}^{-1},\qquad \omega(\eta)=c_a^2\eta+\Order(\eta^2),\]
	where $P\in C^\infty(\R\times[0,\xi_{\rm lo}^2];\R^{2\times2})$, $\omega\in C^\omega([0,\xi_{\rm lo}^2])$ and $\tfrac12c_a^2\eta\leq\omega(\eta)\leq2c_a^2\eta$ for $0\leq\eta\leq\xi_{\rm lo}^2$. Each $P_\eta(\theta)$ is invertible, $P$ is analytic in $\eta$, and $P_\eta|_{\eta=0}=P_0$, where $P_0$ is defined in \eqref{Eq.Pzero}. For $\theta\in\R$ and $0\leq\eta\leq\xi_{\rm lo}^2$,
	\begin{align*}
		&P_\eta(\theta+L)=-P_\eta(\theta),\qquad \partial_\theta P_\eta=\mathsf P_{\sqrt\eta}P_\eta-P_\eta
		\begin{pmatrix}
			0&1\\
			-\omega(\eta)&0
		\end{pmatrix}
		,\\
		&\qquad\sup_{\theta\in\R}\bigl(\|P_\eta(\theta)-P_0(\theta)\|+\|\partial_\theta P_\eta(\theta)-P_0'(\theta)\|\bigr)\leq C\eta.
	\end{align*}
	For $\xi\geq\xi_{\rm hi}$, one simply has $S=\diag(1+\xi^2,1)$.
\end{proposition}

For the energy calculation, we compare $S$ at neighboring derivative orders. Fix a compact interval $I_{\rm r}=[r_-,r_+]\Subset(0,\overline\varrho)$. The bounds in the next lemma are uniform for the analytic radius $\varrho\in I_{\rm r}$; the constants may depend on $r_-$ and $r_+$.

\begin{lemma}\label{Lem.sampled_matrices}
	Let $p\geq3$ be odd and fix $\chi>\chi_p$. Fix a family $S$ from Proposition \ref{Prop.matrix} and let $\varrho\in I_{\rm r}$. There exist $\Sigma_{\rm mat}\geq\max\{2,\delta_0^{-1}\}$, $C_\varrho>0$ and $J_{\rm cmp}\in\Z_{\geq1}$ such that, for
	\[\sigma\geq\Sigma_{\rm mat},\qquad c_*\in[1/2,2],\qquad \theta_*\in\R,\qquad \xi_N^\star:=\frac{c_*N}{\varrho\sigma},\qquad S_N^\star:=S(\theta_*,\xi_N^\star,\sigma^{-1}),\]
	we have
	\begin{align*}
		&\left\|(S_N^\star)^{1/2}(S_J^\star)^{-1/2}\right\|\leq1+\frac{C_\varrho(1+N-J)^{J_{\rm cmp}}}{J}\qquad(N\geq J\geq1),\\
		&\left\|(S_N^\star)^{1/2}(S_0^\star)^{-1/2}\right\|\leq C_\varrho(1+N)^{J_{\rm cmp}}\qquad(N\geq0),\\
		&\left\|(S_N^\star)^{1/2}(S_{N+1}^\star)^{-1/2}\right\|+\left\|(S_{N+1}^\star)^{1/2}(S_N^\star)^{-1/2}\right\|\leq C_\varrho\qquad(N\geq0).
	\end{align*}
	For the fixed exponent $s_*$ there is $C>0$ such that
	\begin{equation}\label{Eq.metric_difference}
		\left\|(S_N^\star)^{-1/2}(S_N^\star-t_NS_{N-1}^\star)(S_N^\star)^{-1/2}\right\|\leq\frac{C}N,\qquad t_N:=\left(\frac N{N+1}\right)^{2s_*}\quad(N\geq1).
	\end{equation}
	Moreover,
	\begin{equation}\label{Eq.sampled_boundary}
		\sup_{N\geq0}\xi_{N+1}^\star\sqrt{\frac{(S_N^\star)_{22}(S_{N+1}^\star)_{22}}{\det S_{N+1}^\star}}\leq1+\frac{C_\varrho}\sigma.
	\end{equation}
	The constants are independent of $\sigma,c_*,\theta_*$ and the Euler profile.
\end{lemma}

Proposition \ref{Prop.matrix} and Lemma \ref{Lem.sampled_matrices} are proved in Section \ref{Sec.matrix_proof}. In the energy estimate, \eqref{Eq.metric_difference} controls the change of derivative order, and \eqref{Eq.sampled_boundary} controls the characteristic boundary integral.

\subsection{An increasing analytic radius}\label{Subsec.energy}

We use the fixed exponent $s_*$ in \eqref{Eq.dimension_indices}. Let $\tau=\log(\sigma/\Sigma)$ and let $\varrho(\tau)$ take values in $I_{\rm r}$. Set
\begin{equation}\label{Eq.energy_phase_samples}
	\theta_\sigma:=\sigma+\sigma^{-1}\Psi(\log\sigma,x),\quad \xi_N:=\frac{N}{\varrho(\tau)\sigma\lambda},\quad \delta:=\sigma^{-1},\quad S_N:=S(\theta_\sigma,\xi_N,\delta).
\end{equation}
Choose $\xi_*:=\xi_{\rm lo}/8$ and $\chi_0\in C_c^\infty([0,\infty);[0,1])$ such that $\chi_0=1$ on $[0,\xi_*]$ and $\operatorname{supp}\chi_0\subset[0,2\xi_*)$. Define $U=(\zeta, w)^T$ and
\begin{equation}\label{Eq.adapted_layers}
	\mathsf T_N:=I+\chi_0\!\left(\frac{N}{\varrho(\tau)\sigma}\right) \bigl(P_0(\theta_\sigma)P_0(\sigma)^{-1}-I\bigr),\quad \mathcal U_N:=\bigl(\mathsf T_N\partial_I(\mathsf T_N^{-1}U)\bigr)_{|I|=N}.
\end{equation}
The matrices $\mathsf T_N$ are invertible for sufficiently large $\sigma$, uniformly in $N$ and $\varrho\in I_{\rm r}$; see Lemma \ref{Lem.gauge_bounds}. For an arbitrary pair $W$, reconstruct $U$ by \eqref{Eq.zeta-w-W-def} and set
\begin{equation}\label{Eq.energy_weights}
	\cE^\sharp_\varrho(\sigma;W) :=\sigma^2\sum_{N\geq0}(1+N)^{2s_*}\frac{\varrho^{2N}}{(N!)^2} \|S_N^{1/2}\mathcal U_N\|_{L^2(\Omega)}^2.
\end{equation}
At each fixed $\sigma$, the quantity $\sqrt{\cE^\sharp_\varrho(\sigma;W)}$ is a norm on the pairs $W$ for which \eqref{Eq.energy_weights} is finite. Lemma \ref{Lem.energy_comparison} proves the uniform bounds
\begin{equation}\label{Eq.energy_state_comparison}
	c\|W\|_{\cX_\varrho}^2\leq\cE^\sharp_\varrho(\sigma;W) \leq C\bigl(\|W\|_{\cX_\varrho}^2+\|\nabla W\|_{\cX_\varrho}^2\bigr).
\end{equation}
The right side is allowed to be infinite. In particular, an analytic radius gap gives
\begin{equation}\label{Eq.energy_gap_comparison}
	\sqrt{\cE^\sharp_\varrho(\sigma;W)} \leq C_{\varrho,\varrho'}\|W\|_{\cX_{\varrho'}} \qquad(\varrho<\varrho').
\end{equation}
For the differentiated equation, assume the following local regularity: for every $\tau_0\in[0,T_*]$, there are a relative neighborhood $J$ of $\tau_0$ and $r_J\in(0,2\overline\varrho)$ such that
\begin{equation}\label{Eq.energy_smooth_core}
	\sup_{\tau\in J}\varrho(\tau)<r_J, \qquad (\zeta,D\zeta)\in W^{1,1}(\Sigma\ee^J;\mathcal B_{r_J}^2), \qquad \mathcal F\in L^1(\Sigma\ee^J;\mathcal B_{r_J}).
\end{equation}
Here $\Sigma\ee^J:=\{\Sigma\ee^\tau:\tau\in J\}$. The analytic radius gap in \eqref{Eq.energy_smooth_core} ensures absolute convergence of the differentiated energy series.

\begin{proposition}[The full linearized energy estimate]\label{Prop.energy}
	Fix the background, $\overline\varrho$, the phase $\Psi$ and the matrix $S$ in Proposition \ref{Prop.matrix}. Let $I_{\rm r}\Subset(0,\overline\varrho)$ and $A_*>0$. There are $\Sigma_E,C>0$ such that, for $\Sigma\geq\Sigma_E$, $0<T_*\leq2$ and
	\[\varrho\in C^1([0,T_*];I_{\rm r}),\qquad 0\leq\alpha(\tau):=\varrho'(\tau)/\varrho(\tau)\leq A_*,\]
	every real-valued radial solution of $\cL_\sigma\zeta=\mathcal F$ satisfying \eqref{Eq.energy_smooth_core} obeys, for almost every $\tau\in(0,T_*)$,
	\begin{equation}\label{Eq.full_energy_differential}
		\begin{aligned}
			&\frac{\mathrm d}{\mathrm d\tau}\cE^\sharp_{\varrho(\tau)}(\sigma;W)+2\sigma^2\left[\gamma_v-\mathcal L_v(\varrho)-\frac\chi\varrho-\alpha -C\bigl(\sigma^{-1}+\nu(\sigma)\bigr)\right]\\
			&\qquad\qquad\qquad\qquad\qquad\qquad\qquad\qquad\times\sum_{N\geq0}N(1+N)^{2s_*}\frac{\varrho^{2N}}{(N!)^2} \|S_N^{1/2}\mathcal U_N\|_{L^2}^2\\
			&\leq C\cE^\sharp_{\varrho(\tau)}(\sigma;W) +C\sqrt{\cE^\sharp_{\varrho(\tau)}(\sigma;W)} \|\sigma^2\mathcal F\|_{\mathcal B_{\varrho(\tau)}},\qquad \sigma=\Sigma\ee^\tau.
		\end{aligned}
	\end{equation}
	The constants depend on the fixed coefficients, $I_{\rm r}$ and $A_*$, but not on $\Sigma$ or $T_*$.
\end{proposition}

Section \ref{Sec.energy_proof} proves \eqref{Eq.full_energy_differential}. To obtain a positive coefficient of the derivative sum, the outgoing contribution $\gamma_v$ must dominate the transport commutator loss $\mathcal L_v(\varrho)$ and the periodic growth term $\chi/\varrho$, leaving a positive margin for $\varrho'/\varrho$.

\begin{assumption}\label{Ass.analytic_margin}
	There are a background $(\beta,v)$ as in Assumption \ref{Ass.background}, a common coefficient analytic radius $R_v$ as in \eqref{Eq.qualitative_coefficient_extension}--\eqref{Eq.common_coefficient_radius}, and $\chi>\chi_p$, $\varrho_0>0$, such that
	\begin{equation}\label{Eq.analytic_margin}
		4\varrho_0<R_v,\qquad \mathcal L_v(4\varrho_0)<\gamma_v,\qquad \mathcal L_v(\varrho_0)+\frac\chi{\varrho_0}<\gamma_v.
	\end{equation}
	Here $\chi_p$ is the oscillator constant in \eqref{Eq.cell_growth_constant}.
\end{assumption}

Fix one choice $(\beta,v,R_v,\chi,\varrho_0)$ in Assumption \ref{Ass.analytic_margin}. By the strict inequalities in \eqref{Eq.analytic_margin} and the continuity of $\mathcal L_v$, choose $\overline\varrho>\varrho_0$ still satisfying \eqref{Eq.phase_radius_condition}, and then choose $\eta,\epsilon_{\rm r}>0$ such that
\begin{equation}\label{Eq.working_analytic_margin}
	\begin{gathered}
		I_{\rm r}:=[\varrho_0\ee^{-2\epsilon_{\rm r}},\varrho_0\ee^{2\eta}] \Subset(0,\overline\varrho),\qquad 0<\epsilon_{\rm r}<\eta/3,\\
		g(r):=\gamma_v-\mathcal L_v(r)-\chi/r,\qquad \kappa_*:=\inf_{r\in I_{\rm r}}g(r)>0,\qquad 0<\eta<\kappa_*/4.
	\end{gathered}
\end{equation}
The phase and the approximate solution are now fixed once and for all using $\overline\varrho$ from \eqref{Eq.working_analytic_margin} and \eqref{Eq.approximation_radii}. Define
\begin{equation}\label{Eq.radii}
	\varrho_{\rm F}(\tau):=\varrho_0\ee^{\eta\tau-\epsilon_{\rm r}},\qquad \varrho_{\rm A}(\tau):=\varrho_0\ee^{\eta\tau-2\epsilon_{\rm r}},\qquad 0\leq\tau\leq2.
\end{equation}
Then
\begin{align*}
	\varrho_{\rm A}(0)<\varrho_{\rm F}(0)<\varrho_0<\varrho_{\rm A}(1)<\varrho_{\rm F}(1),\quad \max_{0\leq\tau\leq2}\varrho_{\rm F}(\tau)<\overline\varrho, \qquad 4\overline\varrho<R_v.
\end{align*}
The paths $\varrho_{\rm F}$ and $\varrho_{\rm A}$ will be used for the full equation and the averaging comparison, respectively. By \eqref{Eq.radii}, $\inf_{0\leq\tau\leq2}(\varrho_{\rm F}(\tau)-\varrho_{\rm A}(\tau))>0$.

Choose $\Sigma_E$ larger, if necessary, so that
\[C\sup_{\sigma\geq\Sigma_E}\bigl(\sigma^{-1}+\nu(\sigma)\bigr)\leq\kappa_*/4.\]
For $0\leq\varrho'/\varrho\leq\eta$, the coefficient in square brackets in \eqref{Eq.full_energy_differential} is at least $\kappa_*/2$. Dropping the nonnegative derivative sum gives the estimate used below.

\begin{corollary}[The forward solution in an increasing analytic radius]\label{Cor.forward_energy}
	Under Assumption \ref{Ass.analytic_margin}, let $\Sigma\geq\Sigma_E$, $0<T_*\leq2$ and
	\[\varrho\in C^1([0,T_*];I_{\rm r}),\qquad 0\leq\varrho'/\varrho\leq\eta.\]
	Given a real-valued radial pair $W_0$ and a measurable real-valued radial scalar source $\mathcal F$ with
	\[\cE^\sharp_{\varrho(0)}(\Sigma;W_0)<\infty,\qquad \int_0^{T_*}\|(\Sigma\ee^u)^2\mathcal F(\Sigma\ee^u)\|_{\mathcal B_{\varrho(u)}}\,\mathrm du<\infty,\]
	there is a unique finite-energy solution of $\cL_\sigma\zeta=\mathcal F$ with $W(\Sigma)=W_0$. For $0\leq s\leq t\leq T_*$,
	\begin{equation}\label{Eq.energy_estimate}
		\sqrt{\cE^\sharp_{\varrho(t)}(\Sigma\ee^t;W(\Sigma\ee^t))} \leq C\bigg[\sqrt{\cE^\sharp_{\varrho(s)}(\Sigma\ee^s;W(\Sigma\ee^s))}+\int_s^t\|(\Sigma\ee^u)^2\mathcal F(\Sigma\ee^u)\|_{\mathcal B_{\varrho(u)}}\,\mathrm du\bigg].
	\end{equation}
	The constant is independent of $\Sigma$ and $T_*$. For every integer $q\geq0$,
	\[W\in C([\Sigma,\Sigma\ee^{T_*}];H^q(\Omega;\R^2)).\]
	If $J\subset[0,T_*]$ is a closed interval and $r<\inf_{\tau\in J}\varrho(\tau)$, then
	\[\tau\longmapsto W(\Sigma\ee^\tau)\in C(J;\cX_r).\]
\end{corollary}

We prove Corollary \ref{Cor.forward_energy} by approximating the data and source at larger analytic radii; see Subsection \ref{Subsec.forward_energy}.

At low frequencies, $\mathsf T_N$ places spatial derivatives on the modulation rather than on the fast frame:
\[\mathcal U_N=P_0(\theta_\sigma)\operatorname{diag}(1,\sigma^{-1})\nabla^NW \quad\text{when }\chi_0=1.\]
At high frequencies $\mathsf T_N=I$, and the energy is the usual wave energy. At zero frequency the normalized modulation metric equals $I_2$, whereas the metric of the original wave unknowns degenerates as $\sigma\to\infty$.

For homogeneous solutions, we apply Corollary \ref{Cor.forward_energy} along $\varrho_{\rm F}$ from \eqref{Eq.radii}. The initial analytic radius gap absorbs the derivative in \eqref{Eq.energy_state_comparison}, and the final analytic radius gap gives the following endpoint bound.

\begin{corollary}[The homogeneous endpoint]\label{Cor.closed_energy}
	Under the hypotheses of Corollary \ref{Cor.forward_energy}, we define $X:=\cX_{\varrho_0}$. Let $\mathcal P_\Sigma(\tau)$ denote the homogeneous evolution from $\Sigma$ to $\Sigma\ee^\tau$. For $W_0\in X$,
	\begin{equation}\label{Eq.endpoint_smoothing}
		\sup_{0\leq\tau\leq2}\left( \|\mathcal P_\Sigma(\tau)W_0\|_{\cX_{\varrho_{\rm F}(\tau)}} +\|\cS(\tau)W_0\|_{\cX_{\varrho_{\rm F}(\tau)}}\right) \leq C\|W_0\|_X.
	\end{equation}
	Consequently, defining $\mathsf A W_0:=\cS(1)W_0$ for $W_0\in X$, there holds $\mathcal P_\Sigma(1),\mathsf A\in\mathcal K(X)$.
\end{corollary}

By \eqref{Eq.acoustic_radius_gap} and $\varrho_{\rm F}(0)<\varrho_0$, every $W_0\in X$ belongs to $\cH_{\varrho_{\rm F}(0)}$, with norm at most $C\|W_0\|_X$. Thus $\cS(\tau)W_0$ is defined by the averaged evolution with initial space $\cH_{\varrho_{\rm F}(0)}$ and belongs to $\cH_{\varrho_{\rm F}(\tau)}$. Since $\varrho_{\rm F}(1)>\varrho_0$, the inclusion $\cX_{\varrho_{\rm F}(1)}\hookrightarrow X$ is compact. Together with \eqref{Eq.endpoint_smoothing}, we obtain the compactness of $\mathcal P_\Sigma(1)$ and $\mathsf A$. The proof is given in Subsection \ref{Subsec.averaging_proof}.

\subsection{Averaging and decaying trajectories}\label{Subsec.fixedpoint}

We first construct a decaying solution of $x_{n+1}=Kx_n+R_n(x_n)$ by choosing the initial datum in the finite-dimensional spectral subspace of $K$ outside the prescribed decay circle.

\begin{proposition}[Decaying trajectories]\label{Prop.compact_decay}
	Let $X$ be a real Banach space, let $K:X\to X$ be a compact linear operator, and choose $b\in(0,1)$ such that the spectrum of the complexification of $K$ does not meet $|z|=b$. Let $P_+$ be the Riesz projection of $K$ for $|z|>b$, and define $P_-:=I-P_+$. There are constants $C_G,\epsilon_0>0$, depending only on $K$ and $b$, for which the following recurrence has a unique decaying solution.

	Let $R_n:B_X(R)\to X$ satisfy
	\[\operatorname{Lip}(R_n)\leq\epsilon_0,\qquad \|R_n(0)\|\leq h b^n,\qquad 2C_Gh<R.\]
	Then
	\begin{equation}\label{Eq.abstract_decay_recurrence}
		x_{n+1}=Kx_n+R_n(x_n),\qquad P_-x_0=0
	\end{equation}
	has a unique solution in the ball
	$\sup_{n\geq0}b^{-n}\|x_n\|\leq2C_Gh$.
	The constants depend only on $K$ and $b$, not on $R$ or $h$.
\end{proposition}

The range of $P_+$ is finite dimensional. In Subsection \ref{Subsec.abstract_decay}, we solve the recurrence forward on $\operatorname{Ran}P_-$ and backward on $\operatorname{Ran}P_+$, and determine $P_+x_0$ from the condition $\sup_{n\geq0}b^{-n}\|x_n\|<\infty$.

We apply Proposition \ref{Prop.compact_decay} to the endpoint maps of the wave equation. The first step is to compare the true homogeneous endpoint with the compact averaged endpoint $\mathsf A$. We work on $X:=\cX_{\varrho_0}$ and define
\[\mu(\Sigma):=\sup_{\sigma\geq\Sigma}\bigl(\sigma^{-1}+\nu(\sigma)\bigr), \qquad \epsilon_\Sigma:=\Sigma^{-1}(1+\log\Sigma)^M.\]
After increasing $M$ and the initial threshold once, both functions tend to zero as $\Sigma\to\infty$.

\begin{proposition}[Quantitative averaging]\label{Prop.averaging}
	For all sufficiently large $\Sigma$,
	\begin{equation}\label{Eq.averaging}
		\|\mathcal P_\Sigma(1)-\mathsf A\|_{X\to X}\leq C\mu(\Sigma), \qquad \mathsf A=\cS(1).
	\end{equation}
\end{proposition}

Choose periodic zero-average functions $\alpha_0,\alpha_1$ with
$\alpha_0'=a_0-c_a^2$ and $\alpha_1'=a_1$. The corrector used in the proof is
\begin{equation}\label{Eq.averaging_corrector}
	\mathscr N_\sigma:=\sigma^{-1}
	\begin{pmatrix}
		0&0\\
		-\alpha_0(\sigma)A_0-\alpha_1(\sigma)V&0
	\end{pmatrix}
	.
\end{equation}
In logarithmic time,
$\partial_\tau\mathscr N_\sigma=\cB(\sigma)-\mathscr N_\sigma$ and
$\mathscr N_\sigma^2=\mathscr N_\sigma\cB(\sigma)=0$.
The corrected equation is estimated at $\varrho_{\rm A}(\tau)$, while the true solution is controlled at $\varrho_{\rm F}(\tau)$. The positive difference $\varrho_{\rm F}-\varrho_{\rm A}$ controls the derivatives in $\mathscr N_\sigma$ and $F_\sigma^c$.

Proposition \ref{Prop.averaging} estimates the linear part of the discrete map. We next solve the nonlinear equation on one logarithmic time interval. The normalization of the initial data makes the quadratic term small on a fixed ball in $X$. For $x\in X$, prescribe $W(\Sigma)=\Sigma^{-2}x$ and denote the normalized endpoint by
\begin{equation}\label{Eq.nonlinear_endpoint_definition}
	F_\Sigma(x):=(\ee\Sigma)^2W(\ee\Sigma).
\end{equation}
\begin{proposition}[The nonlinear endpoint]\label{Prop.nonlinear_endpoint}
	There are $R_*>0$, $\Sigma_*>0$ and $C>0$ such that $F_\Sigma$ is defined on $B_X(R_*)$ for every $\Sigma\geq\Sigma_*$. The endpoint map satisfies
	\begin{equation}\label{Eq.nonlinear_endpoint_bounds}
		F_\Sigma(x)=\ee^2\mathcal P_\Sigma(1)x+N_\Sigma(x),\qquad \|N_\Sigma(0)\|_X\leq C\epsilon_\Sigma,
	\end{equation}
	and, for $0<R\leq R_*$,
	\[\operatorname{Lip}(N_\Sigma|_{B_X(R)})\leq C(R+\epsilon_\Sigma).\]
	The solution on $[\Sigma,\ee^2\Sigma]$ obeys
	\begin{equation}\label{Eq.local_nonlinear_energy}
		\sup_{0\leq\tau\leq2}\Sigma^2 \sqrt{\cE^\sharp_{\varrho_{\rm F}(\tau)}(\Sigma \ee^\tau;W)} \leq C(\|x\|_X+\epsilon_\Sigma).
	\end{equation}
\end{proposition}

Propositions \ref{Prop.averaging} and \ref{Prop.nonlinear_endpoint} give the two estimates needed to apply the abstract trajectory result: the linear endpoint approaches a compact operator, and the nonlinear remainder is Lipschitz with a small constant on a fixed ball. We obtain the following decaying perturbation by choosing the finitely many initial components which would otherwise violate the prescribed decay.

\begin{theorem}[The decaying perturbation]\label{Thm.state}
	Under Assumption \ref{Ass.analytic_margin}, there are $2<a<3$, $\Sigma$ sufficiently large, and a smooth real-valued radial solution $z=z_{\app}+\zeta$ of $\cN(z)=0$ on $[\Sigma,\infty)\times\overline\Omega$. The pair $W$ defined by \eqref{Eq.state} satisfies
	\begin{equation}\label{Eq.state_decay}
		\sup_{\sigma\geq\Sigma}(\sigma/\Sigma)^a \|W(\sigma)\|_{\cX_{\varrho_-}} \leq C\Sigma^{-3}(1+\log\Sigma)^M, \qquad \varrho_-:=\varrho_{\rm F}(0)>0.
	\end{equation}
\end{theorem}

The proof in Subsection \ref{Subsec.contraction_proof} applies Proposition \ref{Prop.compact_decay} with
\[K=\ee^2\mathsf A,\qquad b=\ee^{2-a},\qquad x_n=(\Sigma \ee^n)^2W(\Sigma \ee^n).\]
Compactness of $K$ permits a choice $a\in(2,3)$ such that $\operatorname{spec}(K)\cap\{|z|=\ee^{2-a}\}=\varnothing$. Propositions \ref{Prop.averaging} and \ref{Prop.nonlinear_endpoint} give the linear and nonlinear bounds in Proposition \ref{Prop.compact_decay}. Moreover, for $b=\ee^{2-a}>\ee^{-1}$, $\epsilon_{\Sigma\ee^n}\leq C_{a,M}\epsilon_\Sigma b^n$.

\begin{proof}[Proof of Theorem \ref{Thm.normalized}]
	Use Theorem \ref{Thm.state} and set $z=z_{\app}+\zeta$. The reconstruction of $\zeta$ from $W$ is a bounded analytic multiplier at $\varrho_-$. Thus \eqref{Eq.state_decay} and Sobolev embedding give $\|\zeta(\sigma)\|_{L^\infty}\leq C\sigma^{-a}$. Also, \eqref{Eq.approximation_bounds} gives
	\[\|z_{\app}-f(\theta_\sigma)\|_{L^\infty} \leq C\sigma^{-2}(1+\log\sigma)^M.\]
	Since $f'$ is bounded on the real line and $\Psi$ is a polynomial in $\log\sigma$ with bounded spatial coefficients,
	$\|f(\theta_\sigma)-f(\sigma)\|_{L^\infty}\leq C\sigma^{-1}(1+\log\sigma)^M$.
	Adding the bounds for $z_{\rm app}-f(\theta_\sigma)$ and $\zeta$ gives
	\[\sup_{x\in\overline\Omega}|z(\sigma,x)-f(\sigma)| \leq C\sigma^{-1}(1+\log\sigma)^M,\qquad \sigma\geq\Sigma,\]
	which proves Theorem \ref{Thm.normalized}.
\end{proof}

\subsection{Proof of the unconditional theorem}\label{Subsec.unconditional}

We finally state a quantitative construction verifying the two assumptions. The construction yields Theorem \ref{Thm.main_thm} from the conditional Theorem \ref{Thm.conditional}.

For the background family, write
\[R_p:=\rho_*\left(\frac4{p-1}\right)^{3/2}.\]
\begin{proposition}[Quantitative Euler profiles]\label{Prop.quantitative_background}
	There are an odd integer $p_{\rm bg}\geq15$ and constants $\gamma_*,M_*,\rho_*>0$, common to $d=4,5$, such that for each odd $p\geq p_{\rm bg}$ there is a background $(\beta_{d,p},v_{d,p})$ satisfying Assumption \ref{Ass.background}. The common coefficient analytic radius of $v_{d,p}$ may be chosen to be $R_p$, and
	\begin{equation}\label{Eq.background_bounds}
		\sym D_x\bvec\geq\gamma_*I_d,\qquad \|\nabla^j\bvec\|_{L^\infty(\Omega)}\leq M_*(j-1)!R_p^{1-j} \quad(j\geq1).
	\end{equation}
	Moreover, $\beta_{d,p}\to(d-1)/2-1/3$ as odd $p\to\infty$.
\end{proposition}

To verify Assumption \ref{Ass.analytic_margin}, we combine the bounds for $\gamma_v$ and $R_p$ in Proposition \ref{Prop.quantitative_background} with the following estimate for the oscillator constant $\chi_p$.

\begin{proposition}[The oscillator growth constant]\label{Prop.scattering_growth}
	There is $C_{\rm sc}>0$ such that
	\begin{equation}\label{Eq.chi_p_bound}
		0<\chi_p\leq C_{\rm sc}(p+1)^{-2} \qquad(p\geq15,\ p\text{ odd}).
	\end{equation}
\end{proposition}

The two propositions are proved independently in Section \ref{Sec.background_proof}. The bounds $\chi_p\lesssim p^{-2}$ and $R_p\asymp p^{-3/2}$ are compatible because the oscillator growth is of order $p^{-2}$, whereas the available analytic radius is of order $p^{-3/2}$. We now choose a working analytic radius for which the two losses in \eqref{Eq.analytic_margin} are smaller than the positive constant in the outgoing estimate. Define
\begin{equation}\label{Eq.quantitative_radii}
	\vartheta:=\min\left\{\frac1{16},\frac{\gamma_*}{32M_*}\right\},\qquad \varrho_p:=\vartheta R_p,\qquad \chi:=2C_{\rm sc}(p+1)^{-2},
\end{equation}
and let $p_0$ be the least odd integer $P$ with
\begin{equation}\label{Eq.p0_choice}
	P\geq p_{\rm bg},\qquad P-1\geq \left(\frac{2C_{\rm sc}}{\vartheta\rho_*\gamma_*}\right)^2.
\end{equation}
\begin{corollary}[An admissible analytic radius]\label{Cor.quantitative_margin}
	For every odd $p\geq p_0$ and $d\in\{4,5\}$, let $v=v_{d,p}$ be the profile in Proposition \ref{Prop.quantitative_background}. With $\varrho_p,\chi$ from \eqref{Eq.quantitative_radii},
	\begin{align*}
		4\varrho_p&<R_p,\\
		\gamma_v-\mathcal L_v(4\varrho_p)&\geq7\gamma_*/8,\\
		\gamma_v-\mathcal L_v(\varrho_p)-\chi/\varrho_p&\geq3\gamma_*/4.
	\end{align*}
	Hence Assumption \ref{Ass.analytic_margin} holds with $\varrho_0=\varrho_p$; the choices in \eqref{Eq.working_analytic_margin} and \eqref{Eq.radii} are available for $v_{d,p}$.
\end{corollary}

\begin{proof}
	Set
	\[G(t):=\sum_{j\geq2}\frac{t^{j-1}}j,\qquad\forall\ 0<t<1.\]
	Since $G(t)\leq t/(2(1-t))\leq t$ for $0<t\leq1/4$, one has
	\[\mathcal L_v(4\varrho_p)\leq M_*G(4\vartheta) \leq4M_*\vartheta\leq\frac{\gamma_*}{8}.\]
	Moreover, $4\varrho_p\leq R_p/4<R_p$ and $\chi>\chi_p$ by \eqref{Eq.chi_p_bound}. For $p\geq p_0$, \eqref{Eq.p0_choice} gives
	\[\frac\chi{\varrho_p} =\frac{C_{\rm sc}}{4\vartheta\rho_*} \frac{(p-1)^{3/2}}{(p+1)^2} \leq\frac{C_{\rm sc}}{4\vartheta\rho_*}(p-1)^{-1/2} \leq\frac{\gamma_*}{8}.\]
	The radial eigenvalue bound in \eqref{Eq.background_bounds} implies $\gamma_v\geq\gamma_*$. Consequently,
	\[\gamma_v-\mathcal L_v(\varrho_p)-\frac\chi{\varrho_p} \geq\frac{3\gamma_*}{4},\qquad \gamma_v-\mathcal L_v(4\varrho_p)\geq\frac{7\gamma_*}{8}.\]
	Thus Assumptions \ref{Ass.background} and \ref{Ass.analytic_margin} hold.
	The strict inequalities and the continuity of $\mathcal L_v$ give the interval and growth parameters in \eqref{Eq.working_analytic_margin}.
\end{proof}

The two assumptions are now available for all sufficiently large odd exponents in dimensions four and five. We complete the proof in higher dimensions by finite speed of propagation.

\begin{proof}[Proof of Theorem \ref{Thm.main_thm}]
	Let $p_0$ be defined by \eqref{Eq.p0_choice}. If $d\in\{4,5\}$ and $p\geq p_0$ is odd, Corollary \ref{Cor.quantitative_margin} verifies Assumption \ref{Ass.analytic_margin}. Theorem \ref{Thm.conditional} gives real-valued radial compactly supported data for $d=4,5$. For $d>5$, let $U_4$ and $h_0$ be the four-dimensional cone solution and time in the proof of Theorem \ref{Thm.conditional}, and we use the compactly supported extensions $u_{j,4}$.

	Write $y=(y',y'')\in\R^4\times\R^{d-4}$ and define
	\[U_d(t,y',y''):=U_4(t,y'),\qquad 0\leq t<h_0,\quad |y'|^2+|y''|^2\leq(h_0-t)^2.\]
	Since $|y'|\leq h_0-t$, the function $U_d$ is well defined and smooth on each compact truncated cone and is independent of $y''$. In particular,
	\[-\partial_t^2U_d+\Delta_yU_d =-\partial_t^2U_4+\Delta_{y'}U_4=U_4^p=U_d^p.\]
	Take a cut-off function
	\[\chi_{\rm ext}\in C_c^\infty(\R^{d-4};[0,1])\quad\hbox{such that}\quad \supp\chi_{\rm ext}\subset B_{2h_0}^{d-4},\qquad \chi_{\rm ext}|_{\overline B_{h_0}^{d-4}}\equiv1,\]
	and set
	\[u_{j,d}(y',y''):=u_{j,4}(y')\chi_{\rm ext}(y''),\qquad j=0,1.\]
	Then $u_{j,d}\in C_c^\infty(\R^d;\R)$. For $|y|\leq h_0$, both $|y'|\leq h_0$ and $|y''|\leq h_0$, and hence
	\[u_{j,d}(y)=u_{j,4}(y')=\partial_t^jU_4(0,y')=\partial_t^jU_d(0,y),\qquad j=0,1.\]
	Let $u_d$ be the maximal smooth solution with the resulting compactly supported data, and let $T_*$ be the maximal existence time of $u_d$. Finite propagation gives $u_d=U_d$ for $0\leq t<\min\{T_*,h_0\}$ and $|y|\leq h_0-t$. The sequence in \eqref{Eq.origin_blowup_sequence}, for the four-dimensional cone, satisfies $|U_d(t_n,0)|=|U_4(t_n,0)|\to\infty$. If $T_*>h_0$, continuity on $[0,h_0]\times\overline B_{h_0}^d$ contradicts $|U_d(t_n,0)|\to\infty$. Thus $0<T_*\leq h_0<\infty$. The threshold $p_0$ was fixed using only the four- and five-dimensional backgrounds. Hence the same $p_0$ works for every $d>5$, completing the proof of Theorem \ref{Thm.main_thm}.
\end{proof}

\section{The averaged equation and the approximate solution}\label{Sec.approximation_proofs}

We prove Propositions \ref{Prop.macro}, \ref{Prop.macro_spectrum}, \ref{Prop.phase} and \ref{Prop.approximation}. An integer $d\geq4$, an odd exponent $p\geq3$, and a pair $(\beta,v)$ satisfying Assumption \ref{Ass.background} are fixed; $s_d$ and $s_*$ are given by \eqref{Eq.dimension_indices}. The construction of the approximate solution requires solving the slow phase equation and the antiperiodic equations for the corrections. Using the analytic estimates collected in Appendix \ref{App.analytic_spaces}, we first establish the forward evolution for the averaged acoustic system. The compact resolvent of its generator then allows us to solve the phase equation, with polynomial dependence on the slow variable accounting for the possible zero eigenvalue of the shifted generator. Once the phase is fixed, we determine the antiperiodic corrections through order four and prove the residual estimate in Proposition \ref{Prop.approximation}. We conclude by deriving the perturbation equations and estimating their coefficients.

\subsection{The averaged evolution}\label{Subsec.macro_proof}

The properties of $\mathcal B_\varrho$ used below are proved in Appendix \ref{App.analytic_spaces}. We begin with the acoustic energy and compute its time derivative directly. The wave terms cancel in the interior, while the outgoing property \eqref{Eq.acoustic_outgoing_property} gives a nonnegative boundary flux and a favorable high-order transport term. The high-order term allows the analytic radius to increase. This use of outgoing transport is related to the transport estimates in \cite{SWZZ-Euler2026}; here the outgoing property is used to control the analytic derivative weights. We then construct the forward evolution and prove the semigroup and spectral properties needed for the phase equation.

\begin{lemma}[The outgoing property and the analytic energy]\label{Lem.macro_finite_energy}
	Let $W=(\psi,\pi)$ be a smooth solution of \eqref{Eq.macro_pair} on $[0,T_*]\times\overline\Omega$, and let $\varrho$ be as in Proposition \ref{Prop.macro}. For $N\in\Z_{\geq0}$, define
	\begin{align*}
		&a_N(\tau):=(1+N)^{s_*}\frac{\varrho(\tau)^N}{N!},\quad
		e_N(\tau)^2:=a_N^2\left(\|\nabla^N\psi\|_{L^2}^2+\|\nabla^N\pi\|_{L^2}^2+c_a^2\|\nabla^{N+1}\psi\|_{L^2}^2\right),\\
		&\qquad\qquad\qquad E_K(\tau):=\sum_{N=0}^Ke_N(\tau)^2,\qquad d_K(\tau):=\sum_{N=0}^KN e_N(\tau)^2.
	\end{align*}
	There is a constant $C$, independent of $K,W,H$, such that
	\begin{equation}\label{Eq.macro_truncated_energy}
		E_K'+2\left(\gamma_v-\mathcal L_v(\varrho)-\frac{\varrho'}{\varrho}\right)d_K\leq CE_K+2\sqrt{E_K}\|H\|_{\cH_{\varrho(\tau)}}.
	\end{equation}
	The constant may depend on the fixed background, $s_*$ and a compact subinterval of $(0,R_v)$ containing the range of $\varrho$.
\end{lemma}

\begin{proof}
	We use the full ordered derivatives defined in Section \ref{Sec.Introduction}. For $I=(i_1,\ldots,i_N)$ and $S\subset\{1,\ldots,N\}$, let $I_S$ denote the ordered subtuple in the positions belonging to $S$. The empty subtuple corresponds to the identity operator. All spatial integrals below are taken at a fixed time $\tau$.
	
	For each ordered index $I$, define
	\[\mathcal E_I(\tau):=\|\partial_I\psi(\tau)\|_{L^2}^2+\|\partial_I\pi(\tau)\|_{L^2}^2+c_a^2\|\nabla\partial_I\psi(\tau)\|_{L^2}^2.\]
	Then $e_N^2=a_N^2\sum_{|I|=N}\mathcal E_I$. Differentiating $\mathcal E_I$ gives
	\begin{align*}
		\mathcal E_I'=2\operatorname{Re}\int_\Omega\bigl(\overline{\partial_I\psi}\,\partial_\tau\partial_I\psi+\overline{\partial_I\pi}\,\partial_\tau\partial_I\pi+c_a^2\overline{\nabla\partial_I\psi}\cdot\partial_\tau\nabla\partial_I\psi\bigr)\,\mathrm dx.
	\end{align*}
	We substitute the two equations in \eqref{Eq.macro_pair} into this identity. Since $A_0=-\Delta-\nabla\log J_0\cdot\nabla$, the differentiated equations are
	\begin{equation}\label{Eq.acoustic_differentiated_identity}
		\begin{aligned}
			(\partial_\tau+B)\partial_I\psi&=\partial_I\pi-[\partial_I,B]\psi+\partial_IH_1,\\
			(\partial_\tau+B)\partial_I\pi&=c_a^2\Delta\partial_I\psi+\partial_I\pi+c_a^2\partial_I(\nabla\log J_0\cdot\nabla\psi)-[\partial_I,B]\pi+\partial_IH_2,\\
			(\partial_\tau+B)\nabla\partial_I\psi&=\nabla\partial_I\pi-[\partial_I,B]\nabla\psi-\partial_I((\mathrm D_x\bvec)^T\nabla\psi)+\nabla\partial_IH_1.
		\end{aligned}
	\end{equation}
	Here $[\partial_I,B]h:=\partial_I(Bh)-B\partial_Ih$, with $B$ acting componentwise on vector-valued functions.
	
	The Laplacian term in the equation for $\partial_I\pi$ and the term $\nabla\partial_I\pi$ in the equation for $\nabla\partial_I\psi$ contribute
	\begin{align*}
		2c_a^2\operatorname{Re}\int_\Omega\left(\overline{\partial_I\pi}\,\Delta\partial_I\psi+\overline{\nabla\partial_I\psi}\cdot\nabla\partial_I\pi\right)\,\mathrm dx=2c_a^2\operatorname{Re}\int_{\partial\Omega}\overline{\partial_I\pi}\,\partial_n\partial_I\psi\,\mathrm dS.
	\end{align*}
	Thus the highest-order interior terms cancel. For a smooth scalar or Euclidean-valued function $h$, integration by parts and $\bvec\cdot n=1$ give
	\[-2\operatorname{Re}\int_\Omega\overline h\cdot Bh\,\mathrm dx=-\int_{\partial\Omega}|h|^2\,\mathrm dS+\int_\Omega(\operatorname{div}\bvec)|h|^2\,\mathrm dx.\]
	Applying this identity to the three transport terms in $\mathcal E_I'$, the total boundary contribution is
	\[-\int_{\partial\Omega}\left(|\partial_I\psi|^2+|\partial_I\pi|^2+c_a^2|\nabla\partial_I\psi|^2-2c_a^2\operatorname{Re}(\overline{\partial_I\pi}\,\partial_n\partial_I\psi)\right)\,\mathrm dS.\]
	For any smooth complex-valued scalar functions $u,v$, the inequalities $|\partial_nu|\leq|\nabla u|$ and $0<c_a<1$ imply
	\[2c_a^2\operatorname{Re}(\overline v\,\partial_nu)\leq c_a\left(|v|^2+c_a^2|\nabla u|^2\right).\]
	Consequently,
	\begin{equation}\label{Eq.macro_boundary_form}
		|u|^2+|v|^2+c_a^2|\nabla u|^2-2c_a^2\operatorname{Re}(\overline v\,\partial_nu)\geq |u|^2+(1-c_a)(|v|^2+c_a^2|\nabla u|^2)\geq0.
	\end{equation}
	This is the boundary contribution of the outgoing property. In particular, the boundary term in $\mathcal E_I'$ is nonpositive and can be discarded.
	
	The terms containing $\operatorname{div}\bvec$, $\overline{\partial_I\psi}\,\partial_I\pi$ and $|\partial_I\pi|^2$ are bounded by $C\mathcal E_I$. Therefore,
	\begin{align*}
		\mathcal E_I'
		&\leq C\mathcal E_I-2\operatorname{Re}\int_\Omega\bigl(\overline{\partial_I\psi}\,[\partial_I,B]\psi+\overline{\partial_I\pi}\,[\partial_I,B]\pi+c_a^2\overline{\nabla\partial_I\psi}\cdot[\partial_I,B]\nabla\psi\bigr)\,\mathrm dx\\
		&\quad+2c_a^2\operatorname{Re}\int_\Omega\bigl(\overline{\partial_I\pi}\,\partial_I(\nabla\log J_0\cdot\nabla\psi)-\overline{\nabla\partial_I\psi}\cdot\partial_I((\mathrm D_x\bvec)^T\nabla\psi)\bigr)\,\mathrm dx\\
		&\quad+2\operatorname{Re}\int_\Omega\bigl(\overline{\partial_I\psi}\,\partial_IH_1+\overline{\partial_I\pi}\,\partial_IH_2+c_a^2\overline{\nabla\partial_I\psi}\cdot\nabla\partial_IH_1\bigr)\,\mathrm dx.
	\end{align*}
	
	We next use the interior outgoing property to estimate the transport commutators. The ordered Leibniz formula gives
	\[[\partial_I,B]h=\sum_{\emptyset\ne S\subset\{1,\ldots,N\}}\sum_{l=1}^d(\partial_{I_S}b_l)\partial_l\partial_{I_{S^c}}h.\]
	We first consider the terms with $|S|=1$. Fix a position $a\in\{1,\ldots,N\}$ and all indices except $i_a$, and define
	\[V_j:=\partial_{i_1}\cdots\partial_{i_{a-1}}\partial_j\partial_{i_{a+1}}\cdots\partial_{i_N}h,\qquad 1\leq j\leq d.\]
	The inequality $\operatorname{sym}\mathrm D_x\bvec\geq\gamma_vI_d$ implies
	\begin{align*}
		\operatorname{Re}\sum_{j,l=1}^d(\partial_jb_l)\overline V_j\cdot V_l=\sum_{j,l=1}^d(\operatorname{sym}\mathrm D_x\bvec)_{jl}\operatorname{Re}(\overline V_j\cdot V_l)\geq\gamma_v\sum_{j=1}^d|V_j|^2.
	\end{align*}
	Integrating in $x$, summing over the fixed indices and then over the $N$ positions $a$, we obtain
	\begin{equation}\label{Eq.acoustic_transport_coercivity}
		-2\operatorname{Re}\sum_{|I|=N}\sum_{a=1}^N\sum_{l=1}^d\int_\Omega\overline{\partial_Ih}\cdot(\partial_{i_a}b_l)\partial_{i_1}\cdots\partial_{i_{a-1}}\partial_l\partial_{i_{a+1}}\cdots\partial_{i_N}h\,\mathrm dx\leq-2\gamma_vN\|\nabla^Nh\|_{L^2}^2.
	\end{equation}
	The estimate also holds for Euclidean-valued $h$, since the same calculation is summed over the components of $h$. For $N=0$, both sides vanish. Applying \eqref{Eq.acoustic_transport_coercivity} to $\psi$, $\pi$ and $c_a\nabla\psi$ shows that the first-order coefficient derivatives in the commutators contribute at most
	\[-2\gamma_vN\left(\|\nabla^N\psi\|_{L^2}^2+\|\nabla^N\pi\|_{L^2}^2+c_a^2\|\nabla^{N+1}\psi\|_{L^2}^2\right).\]
	The factor $N$ is essential: the outgoing property gives the same favorable contribution at each of the $N$ derivative positions.
	
	For the remaining commutator terms, fix $S\subset\{1,\ldots,N\}$ with $|S|=k\geq2$. Cauchy--Schwarz in the index $l$ and summation over the ordered indices give the pointwise estimate
	\[\sum_{|I|=N}\left|\sum_{l=1}^d(\partial_{I_S}b_l)\partial_l\partial_{I_{S^c}}h\right|^2\leq|\nabla^k\bvec|^2|\nabla^{N-k+1}h|^2.\]
	There is no additional constant in this estimate: the indices in $I_S$ and $I_{S^c}$ vary independently, and the full ordered tensor norms include all their values. There are exactly $\binom Nk$ subsets of size $k$. The triangle inequality over these subsets, followed by Cauchy--Schwarz in the energy pairing, therefore bounds their contribution at level $N$ by
	\[2\sum_{k=2}^N\binom Nk\|\nabla^k\bvec\|_{L^\infty}\frac{e_Ne_{N-k+1}}{a_Na_{N-k+1}}.\]
	Here Cauchy--Schwarz has also been applied to the three summands defining $e_N^2$, again with constant one.
	
	We now multiply the derivative energy inequalities by $a_N^2$ and sum over $0\leq N\leq K$. Since $(a_N^2)'=2N\frac{\varrho'}{\varrho}a_N^2$,	the preceding estimates give
	\[E_K'+2\left(\gamma_v-\frac{\varrho'}{\varrho}\right)d_K\leq CE_K+\mathcal T_K+\mathcal R_K+2\sqrt{E_K}\|H\|_{\cH_{\varrho(\tau)}},\]
	where the nonnegative bounds for the higher transport commutators and the fixed coefficient terms are
	\begin{align*}
		\mathcal T_K&:=2\sum_{N=2}^K\sum_{k=2}^N\binom Nk\frac{a_N}{a_{N-k+1}}\|\nabla^k\bvec\|_{L^\infty}e_Ne_{N-k+1},\\
		\mathcal R_K&:=2c_a^2\sum_{N=0}^Ka_N^2\bigl(\|\nabla^N\pi\|_{L^2}\|\nabla^N(\nabla\log J_0\cdot\nabla\psi)\|_{L^2}+\|\nabla^{N+1}\psi\|_{L^2}\|\nabla^N((\mathrm D_x\bvec)^T\nabla\psi)\|_{L^2}\bigr).
	\end{align*}
	The source estimate follows from Cauchy--Schwarz and
	\[\sum_{N=0}^Ka_N^2\left(\|\nabla^NH_1\|_{L^2}^2+\|\nabla^NH_2\|_{L^2}^2+c_a^2\|\nabla^{N+1}H_1\|_{L^2}^2\right)\leq\|H\|_{\cH_\varrho}^2.\]
	
	Choose $r'\in(\max_{0\leq\tau\leq T_*}\varrho(\tau),R_v)$. By \eqref{Eq.common_coefficient_radius}, the coefficients $\bvec,J_0,J_0^{-1}$ belong to $\mathcal B_{r'}$. Lemma \ref{Lem.analytic_to_multiplier} gives uniformly bounded multiplier norms for $J_0^{-1}$, $\nabla J_0$ and $\mathrm D_x\bvec$ at the radii $\varrho(\tau)$. Since $\nabla\log J_0=J_0^{-1}\nabla J_0$, \eqref{Eq.multiplier_algebra} yields
	\[\sup_{0\leq\tau\leq T_*}\left(\mathfrak M_{\varrho(\tau)}(\nabla\log J_0)+\mathfrak M_{\varrho(\tau)}(\mathrm D_x\bvec)\right)\leq C.\]
	The multiplier estimate in Lemma \ref{Lem.analytic_algebra} applies to the finite derivative sums: a derivative of order at most $K$ of either product in $\mathcal R_K$ uses only derivatives of $\nabla\psi$ of order at most $K$. Thus
	\[\mathcal R_K\leq C\left(\mathfrak M_\varrho(\nabla\log J_0)+\mathfrak M_\varrho(\mathrm D_x\bvec)\right)E_K\leq CE_K.\]
	
	It remains to estimate $\mathcal T_K$ without changing the coefficient of $\mathcal L_v(\varrho)$. The analytic weights satisfy the exact identity
	\[\binom Nk\frac{a_N}{a_{N-k+1}}=\frac{\varrho^{k-1}}{k!}(N-k+1)\left(\frac{1+N}{N-k+2}\right)^{s_*}.\]
	Setting $J=N-k+1$ and
	\[c_k:=\frac{\varrho^{k-1}}{k!}\|\nabla^k\bvec\|_{L^\infty},\]
	we obtain
	\[\mathcal T_K=2\sum_{k=2}^Kc_k\sum_{J=1}^{K-k+1}J\left(\frac{J+k}{J+1}\right)^{s_*}e_{J+k-1}e_J.\]
	For $J\geq1$, the binomial theorem gives
	\begin{align*}
		J\left(\frac{J+k}{J+1}\right)^{s_*}
		&=J+\sum_{a=1}^{s_*}\binom{s_*}{a}\frac{J(k-1)^a}{(J+1)^a}\leq J+C(1+k)^{s_*}.
	\end{align*}
	Moreover,
	\begin{align*}
		\sum_{J=1}^{K-k+1}Je_{J+k-1}e_J
		&\leq\left(\sum_{J=1}^{K-k+1}(J+k-1)e_{J+k-1}^2\right)^{1/2}\left(\sum_{J=1}^{K-k+1}Je_J^2\right)^{1/2}\leq d_K,\\
		\sum_{J=1}^{K-k+1}e_{J+k-1}e_J
		&\leq\left(\sum_{J=1}^{K-k+1}e_{J+k-1}^2\right)^{1/2}\left(\sum_{J=1}^{K-k+1}e_J^2\right)^{1/2}\leq E_K.
	\end{align*}
	By \eqref{Eq.geometry_budget} and \eqref{Eq.vecb-analytic-bounds},
	\begin{align*}
		\sum_{k=2}^\infty c_k=\mathcal L_v(\varrho),\quad
		\sum_{k=2}^\infty(1+k)^{s_*}c_k\leq M_v\sum_{k=2}^\infty\frac{(1+k)^{s_*}}{k}\left(\frac{\varrho}{R_v}\right)^{k-1}\leq C.
	\end{align*}
	The last bound is uniform on the prescribed compact interval of analytic radii. Consequently,
	\[\mathcal T_K\leq2\mathcal L_v(\varrho)d_K+CE_K.\]
	Combining the bounds for $\mathcal T_K$ and $\mathcal R_K$ proves \eqref{Eq.macro_truncated_energy}. The outgoing property contributes $2\gamma_vd_K$, while the higher transport commutators and the increasing analytic weights cost $2\mathcal L_v(\varrho)d_K$ and $2(\varrho'/\varrho)d_K$, respectively. The condition in Proposition \ref{Prop.macro} leaves a positive coefficient of $d_K$.
\end{proof}

The outgoing property has given an energy estimate uniform in the number of derivatives. We now prove Proposition \ref{Prop.macro} by combining finite-order existence with the analytic energy estimate \eqref{Eq.macro_truncated_energy}.

\begin{proof}[Proof of Proposition \ref{Prop.macro}]
	We first construct solutions at finite Sobolev orders. For each integer $m\geq2$, the hypotheses imply
	\[W_0\in H^{m+1}(\Omega)\times H^m(\Omega),\qquad H\in L^1([0,T_*];H^{m+1}(\Omega)\times H^m(\Omega)).\]
	Extend the data and source by bounded Sobolev extension operators, and extend $\bvec$ and $\nabla\log J_0$ smoothly to $\R^d$ with bounded derivatives, keeping $c_a^2\Delta$ unchanged. Standard whole-space linear hyperbolic theory and smooth approximation give a solution whose restriction to $\Omega$ belongs to $C([0,T_*];H^{m+1}(\Omega)\times H^m(\Omega))$; see \cite{Sogge}. 
	
	Now we prove the uniqueness and independence of the extensions. Let $W^{(1)},W^{(2)}\in C([0,T_*];H^3(\Omega)\times H^2(\Omega))$ solve \eqref{Eq.macro_pair} with initial value $W_0$ and source $H$, and write $(u,v):=W^{(1)}-W^{(2)}$. The order-zero energy calculation in Lemma \ref{Lem.macro_finite_energy} and the nonnegative boundary flux \eqref{Eq.macro_boundary_form} give
	\begin{align*}
		\frac{\mathrm d}{\mathrm d\tau}\left(\|u\|_{L^2}^2+\|v\|_{L^2}^2+c_a^2\|\nabla u\|_{L^2}^2\right)\leq C\left(\|u\|_{L^2}^2+\|v\|_{L^2}^2+c_a^2\|\nabla u\|_{L^2}^2\right),
	\end{align*}
	with $(u,v)|_{\tau=0}=(0,0)$. Gronwall's inequality yields $u=v=0$. Thus the restrictions are independent of the Sobolev extensions and agree at all orders $m\geq2$. The solution $W$ satisfies
	\[W\in\bigcap_{m\in\Z_{\geq0}}C([0,T_*];H^{m+1}(\Omega)\times H^m(\Omega)).\]
	
	We now prove the analytic estimate. For $K\in\Z_{\geq0}$, define
	\[E_K(\tau):=\sum_{N=0}^K(1+N)^{2s_*}\frac{\varrho(\tau)^{2N}}{(N!)^2}\left(\|\nabla^NW(\tau)\|_{L^2}^2+c_a^2\|\nabla^{N+1}\psi(\tau)\|_{L^2}^2\right).\]
	For each fixed $K$, the equation and the finite Sobolev regularity imply $E_K\in W^{1,1}(0,T_*)$. The calculation in Lemma \ref{Lem.macro_finite_energy} therefore applies for almost every $\tau$. Since the coefficient of $d_K$ in \eqref{Eq.macro_truncated_energy} is positive, for every $\epsilon>0$,
	\[\frac{\mathrm d}{\mathrm d\tau}\sqrt{E_K+\epsilon}\leq C\sqrt{E_K+\epsilon}+\|H\|_{\cH_{\varrho(\tau)}}.\]
	Integrating and then letting $\epsilon\downarrow0$, we obtain
	\[\sqrt{E_K(t)}\leq\ee^{C(t-s)}\sqrt{E_K(s)}+\int_s^t\ee^{C(t-u)}\|H(u)\|_{\cH_{\varrho(u)}}\,\mathrm du,\qquad 0\leq s\leq t\leq T_*.\]
	Taking $s=0$ and letting $K\to\infty$ shows that
	\[W(t)\in\cH_{\varrho(t)}\quad(0\leq t\leq T_*),\qquad \sup_{0\leq t\leq T_*}\|W(t)\|_{\cH_{\varrho(t)}}<\infty.\]
	For arbitrary $0\leq s\leq t\leq T_*$, monotone convergence as $K\to\infty$ gives
	\begin{equation}\label{Eq.macro_Cauchy_explicit}
		\|W(t)\|_{\cH_{\varrho(t)}}\leq\ee^{C(t-s)}\|W(s)\|_{\cH_{\varrho(s)}}+\int_s^t\ee^{C(t-u)}\|H(u)\|_{\cH_{\varrho(u)}}\,\mathrm du.
	\end{equation}
	Since $t-s\leq T_*$, \eqref{Eq.macro_Cauchy_explicit} implies \eqref{Eq.macro_energy}. Applying \eqref{Eq.macro_Cauchy_explicit} to the difference of two solutions proves continuous dependence on the initial data and source in the norms $\|W_0\|_{\cH_{\varrho(0)}}$ and $\int_0^{T_*}\|H(u)\|_{\cH_{\varrho(u)}}\,\mathrm du$.
	
	Finally, fix
	\[J=[t_0,t_1],\qquad 0\leq t_0\leq t_1\leq T_*,\qquad 0<r<r'<\min_{\tau\in J}\varrho(\tau),\]
	and define
	\[M_J:=\sup_{\tau\in J}\|W(\tau)\|_{\cH_{r'}}<\infty.\]
	For $s,t\in J$ and $K\geq0$, splitting the $\cH_r$ norm at derivative order $K$ gives
	\begin{align*}
		\|W(t)-W(s)\|_{\cH_r}^2		&\leq\sum_{N=0}^K(1+N)^{2s_*}\frac{r^{2N}}{(N!)^2}\|\nabla^N(W(t)-W(s))\|_{L^2}^2\\
		&\qquad+\sum_{N=0}^K(1+N)^{2s_*}\frac{r^{2N}}{(N!)^2}c_a^2\|\nabla^{N+1}(\psi(t)-\psi(s))\|_{L^2}^2+4M_J^2\left(\frac r{r'}\right)^{2(K+1)}.
	\end{align*}
	For fixed $K$, the finite sum tends to zero as $t\to s$ by the Sobolev continuity of $W$. Letting $K\to\infty$ proves $W\in C(J;\cH_r)$.
	
	For radial initial data and source, rotation invariance of \eqref{Eq.macro_pair} and uniqueness imply that $W$ is radial. For real-valued initial data and source, complex conjugation and uniqueness give $\overline W=W$. This proves Proposition \ref{Prop.macro}.
\end{proof}

We next prove strong continuity on each fixed analytic space $\cH_r$ and identify the infinitesimal generator with the differential operator defined in \eqref{Eq.macro_generator_domain}.

\begin{lemma}\label{Lem.macro_generator}
	Fix $r\in(0,R_v)$ with $\mathcal L_v(r)<\gamma_v$. The family $(\cS_r(t))_{t\geq0}$ is a strongly continuous semigroup on $\cH_r$ with infinitesimal generator $\cA_r$ and domain \eqref{Eq.macro_generator_domain}. There are $M_r\geq1$ and $\omega_r\in\R$ such that
	\begin{equation}\label{Eq.macro_semigroup_growth}
		\|\cS_r(t)\|_{\cH_r\to\cH_r}\leq M_r\ee^{\omega_rt},\qquad t\geq0.
	\end{equation}
	For $T_*>0$, $W_0\in\cH_r$ and $H\in L^1([0,T_*];\cH_r)$, the solution of \eqref{Eq.macro_pair} with $W(0)=W_0$ satisfies
	\begin{equation}\label{Eq.macro_duhamel}
		W(t)=\cS_r(t)W_0+\int_0^t\cS_r(t-u)H(u)\,\mathrm du,\qquad 0\leq t\leq T_*.
	\end{equation}
\end{lemma}

\begin{proof}
	Uniqueness in Proposition \ref{Prop.macro} gives
	\[\cS_r(0)=I,\qquad \cS_r(t+s)=\cS_r(t)\cS_r(s),\qquad s,t\geq0.\]
	Applying Proposition \ref{Prop.macro} with $\varrho(t)=r$ and $H=0$ on $[0,1]$, choose $C_r\geq1$ such that
	\[\sup_{0\leq t\leq1}\|\cS_r(t)\|_{\cH_r\to\cH_r}\leq C_r.\]
	For $t\geq0$, write $t=n+s$ with $n\in\Z_{\geq0}$ and $0\leq s<1$. The semigroup property gives
	\begin{align*}
		\|\cS_r(t)\|_{\cH_r\to\cH_r}\leq\|\cS_r(s)\|_{\cH_r\to\cH_r}\|\cS_r(1)\|_{\cH_r\to\cH_r}^{\,n}\leq C_r^{n+1}\leq C_r\ee^{t\log C_r}.
	\end{align*}
	Thus \eqref{Eq.macro_semigroup_growth} holds with $M_r:=C_r$ and $\omega_r:=\log C_r$. The constants $M_r$ and $\omega_r$ depend only on $r$ and the fixed background.
	
	We prove strong continuity at $t=0$. Let $r<r'<R_v$ with $\mathcal L_v(r')<\gamma_v$. The derivative estimate \eqref{Eq.analytic_radius_recovery} and the analytic bounds for the coefficients of $\mathscr A$ give
	\[\|\mathscr AU\|_{\cH_r}\leq C_{r,r'}\|U\|_{\cH_{r'}},\qquad U\in\cH_{r'}.\]
	For $W_0\in\cH_{r'}$, Proposition \ref{Prop.macro} with the constant radius $r'$ and the analytic radius gap allow us to integrate the homogeneous equation in $\cH_r$. Since $\cS_r(t)W_0=\cS_{r'}(t)W_0$, we obtain
	\begin{align*}
		\|\cS_r(t)W_0-W_0\|_{\cH_r}\leq\int_0^t\|\mathscr A\cS_{r'}(u)W_0\|_{\cH_r}\,\mathrm du\leq C_{r,r'}t\|W_0\|_{\cH_{r'}},\quad\forall\  0\leq t\leq1.
	\end{align*}
	The dilation argument in Lemma \ref{Lem.analytic_density} also applies to the norm \eqref{Eq.acoustic_pair_norm}, which implies that
	\[\bigcup_{r<r'<R_v,\ \mathcal L_v(r')<\gamma_v}\cH_{r'}\]
	is dense in $\cH_r$. For $W_0\in\cH_r$, choose $W_{0,n}$ in this union with $\|W_{0,n}-W_0\|_{\cH_r}\to0$. For $0\leq t\leq1$,
	\[\|\cS_r(t)W_0-W_0\|_{\cH_r}\leq(C_r+1)\|W_0-W_{0,n}\|_{\cH_r}+\|\cS_r(t)W_{0,n}-W_{0,n}\|_{\cH_r}.\]
	Taking first $t\downarrow0$ and then $n\to\infty$ proves strong continuity at zero. Strong continuity at every $t\geq0$ follows from the semigroup property and the local operator bound.
	
	Let $A$ denote the infinitesimal generator of $(\cS_r(t))_{t\geq0}$ on $\cH_r$. We prove that $\mathcal D(A)=\mathcal D(\cA_r)$ and $AW=\mathscr AW$ on this domain. The embedding $\cH_r\hookrightarrow H^2(\Omega)\times H^1(\Omega)$ implies that $\mathscr A:\cH_r\longrightarrow L^2(\Omega)^2$	is continuous. Consequently, for every $W\in\cH_r$, the homogeneous equation gives
	\[\cS_r(t)W-W=\int_0^t\mathscr A\cS_r(u)W\,\mathrm du\quad\text{in }L^2(\Omega)^2.\]
	Strong continuity of $\cS_r$ on $\cH_r$ implies
	\[\frac1t\int_0^t\mathscr A\cS_r(u)W\,\mathrm du\longrightarrow\mathscr AW\quad\text{in }L^2(\Omega)^2\quad\text{as }t\downarrow0.\]
	If $W\in\mathcal D(A)$, the definition of the generator gives
	\[\frac{\cS_r(t)W-W}{t}\longrightarrow AW\quad\text{in }\cH_r.\]
	The continuous embedding $\cH_r\hookrightarrow L^2(\Omega)^2$ identifies the two limits, so $AW=\mathscr AW$. In particular, $\mathscr AW\in\cH_r$, and therefore $W\in\mathcal D(\cA_r)$.
	
	Conversely, let $W\in\mathcal D(\cA_r)$, let $F:=\mathscr AW\in\cH_r$, and define
	\[V(t):=W+\int_0^t\cS_r(u)F\,\mathrm du.\]
	Strong continuity implies $V\in C^1([0,T_*];\cH_r)$ for every $T_*>0$, with $V(0)=W$ and $\partial_tV(t)=\cS_r(t)F$. The continuity of $\mathscr A:\cH_r\to L^2(\Omega)^2$ and the integrated homogeneous equation with initial value $F$ give
	\begin{align*}
		\mathscr AV(t)=F+\int_0^t\mathscr A\cS_r(u)F\,\mathrm du=F+\cS_r(t)F-F=\cS_r(t)F=\partial_tV(t)
	\end{align*}
	in $L^2(\Omega)^2$. Thus $V$ solves the homogeneous equation with initial value $W$. On every finite time interval, $V$ satisfies \eqref{Eq.macro_solution_class} with the constant radius $r$. Uniqueness in Proposition \ref{Prop.macro} therefore gives $V(t)=\cS_r(t)W$. It follows that
	\[\left\|\frac{\cS_r(t)W-W}{t}-F\right\|_{\cH_r}\leq\frac1t\int_0^t\|\cS_r(u)F-F\|_{\cH_r}\,\mathrm du\longrightarrow0.\]
	Hence $W\in\mathcal D(A)$ and $AW=F=\cA_rW$. We have proved $\mathcal D(A)=\mathcal D(\cA_r)$ and $A=\cA_r$.
	
	Finally, let $W_0\in\cH_r$ and $H\in L^1([0,T_*];\cH_r)$. The right-hand side of \eqref{Eq.macro_duhamel} defines a function $W\in C([0,T_*];\cH_r)$. The integrated homogeneous equation and Fubini's theorem give
	\[W(t)=W_0+\int_0^t\bigl(\mathscr AW(u)+H(u)\bigr)\,\mathrm du\quad\text{in }L^2(\Omega)^2.\]
	Thus $W$ solves \eqref{Eq.macro_pair} in distributions and has initial value $W_0$. Since $\cH_r$ embeds continuously into $H^{m+1}(\Omega)\times H^m(\Omega)$ for every $m\geq0$, the function $W$ satisfies \eqref{Eq.macro_solution_class} with $\varrho(t)=r$. Uniqueness in Proposition \ref{Prop.macro} proves \eqref{Eq.macro_duhamel}.
\end{proof}

The outgoing property permits a strict increase of the analytic radius at positive times. Combining this radius gain with compact restriction, we prove Proposition \ref{Prop.macro_spectrum} now.

\begin{proof}[Proof of Proposition \ref{Prop.macro_spectrum}]
	Fix $r\in(0,R_v)$ with $\mathcal L_v(r)<\gamma_v$. Lemma \ref{Lem.macro_generator} proves that $(\cS_r(t))_{t\geq0}$ is a strongly continuous semigroup on $\cH_r$, identifies its infinitesimal generator with $\cA_r$ on the domain \eqref{Eq.macro_generator_domain}, and gives the semigroup growth bound \eqref{Eq.macro_semigroup_growth}.	It remains to prove compactness of $\cS_r(t)$ for $t>0$, compactness of the resolvent, the inclusion $\{z\in\C:\operatorname{Re}z>\omega_r\}\subset\C\setminus\spec(\cA_r)$, and the isolation, finite algebraic multiplicity and increased analytic radius of the generalized eigenspaces.
	
	Choose $r_+\in(r,R_v)$ with $\mathcal L_v(r_+)<\gamma_v$, fix $0<\eta_r<\gamma_v-\mathcal L_v(r_+)$, and define
	\[t_r:=\frac1{\eta_r}\log\frac{r_+}{r}>0.\]
	Proposition \ref{Prop.macro}, applied with $\varrho(t)=r\ee^{\eta_rt}$ on $[0,t_r]$ and $H=0$, gives
	\[\sup_{0\leq t\leq t_r}\|\cS_r(t)\|_{\cH_r\to\cH_{r\ee^{\eta_rt}}}<\infty.\]
	For $0<t\leq t_r$, Lemma \ref{Lem.analytic_compactness}, applied to the two components and to the gradient of the first component, gives the compact embedding $\cH_{r\ee^{\eta_rt}}\hookrightarrow\cH_r$.	Hence $\cS_r(t)\in\mathcal K(\cH_r)$ for $0<t\leq t_r$. For $t>t_r$, the identity $\cS_r(t)=\cS_r(t-t_r)\cS_r(t_r)$ implies compactness on $\cH_r$ as well.
	
	We next prove continuity in the operator norm at every positive time. For $t>0$ and $|h|<t/2$,
	\[\cS_r(t+h)-\cS_r(t)=\bigl(\cS_r(t/2+h)-\cS_r(t/2)\bigr)\cS_r(t/2).\]
	The operators $\cS_r(t/2+h)-\cS_r(t/2)$ converge strongly to zero as $h\to0$ and are uniformly bounded in $\|\cdot\|_{\cH_r\to\cH_r}$ for $|h|<t/2$. Since $\cS_r(t/2)$ is compact, the strong convergence is uniform on the compact closure of the image of the unit ball under $\cS_r(t/2)$. Consequently,
	\[\lim_{h\to0}\|\cS_r(t+h)-\cS_r(t)\|_{\cH_r\to\cH_r}=0,\qquad t>0.\]
	
	For spectra and resolvents, we work in the complexification of $\cH_r$. Since $\cA_r$ is the infinitesimal generator of $(\cS_r(t))_{t\geq0}$, the Laplace transform formula and \eqref{Eq.macro_semigroup_growth} give
	\[(zI-\cA_r)^{-1}W=\int_0^\infty\ee^{-zt}\cS_r(t)W\,\mathrm dt,\qquad W\in\cH_r,\quad \operatorname{Re}z>\omega_r,\]
	and
	\[\|(zI-\cA_r)^{-1}\|_{\cH_r\to\cH_r}\leq\frac{M_r}{\operatorname{Re}z-\omega_r}.\]
	In particular, $\{z\in\C:\operatorname{Re}z>\omega_r\}\subset\C\setminus\spec(\cA_r)$.
	
	Fix a real number $\omega>\omega_r$. For $0<\epsilon<T$, operator norm continuity on $[\epsilon,T]$ implies that $\int_\epsilon^T\ee^{-\omega t}\cS_r(t)\,\mathrm dt$	is an operator norm limit of finite sums of compact operators on $\cH_r$. Moreover,
	\begin{align*}
		\left\|\int_0^\epsilon\ee^{-\omega t}\cS_r(t)\,\mathrm dt\right\|_{\cH_r\to\cH_r}&\leq M_r\int_0^\epsilon\ee^{-(\omega-\omega_r)t}\,\mathrm dt\leq M_r\epsilon,\\
		\left\|\int_T^\infty\ee^{-\omega t}\cS_r(t)\,\mathrm dt\right\|_{\cH_r\to\cH_r}&\leq\frac{M_r}{\omega-\omega_r}\ee^{-(\omega-\omega_r)T}.
	\end{align*}
	Letting $\epsilon\downarrow0$ and $T\to\infty$ proves $(\omega I-\cA_r)^{-1}\in\mathcal K(\cH_r)$. For every $z\in\C\setminus\spec(\cA_r)$, the resolvent identity gives
	\begin{align*}
		(zI-\cA_r)^{-1}=(\omega I-\cA_r)^{-1}+(\omega-z)(zI-\cA_r)^{-1}(\omega I-\cA_r)^{-1},
	\end{align*}
	which implies that $(z I-\cA_r)^{-1}\in\mathcal K(\cH_r)$. This proves \eqref{Eq.macro_compact_resolvent}.
	
	On $\mathcal D(\cA_r)$, we have
	\[\lambda I-\cA_r=\left[I-(\omega-\lambda)(\omega I-\cA_r)^{-1}\right](\omega I-\cA_r).\]
	The Riesz--Schauder theorem for the compact operator $(\omega I-\cA_r)^{-1}$ therefore shows that every $\lambda\in\spec(\cA_r)$ is an isolated eigenvalue of finite algebraic multiplicity. In particular,
	\[0<\dim\mathcal G_\lambda(\cA_r)<\infty.\]
	
	We finally prove that every generalized eigenspace belongs to $\cH_{r_+}$. Fix $\lambda\in\spec(\cA_r)$, an integer $m\geq1$, and $U\in\ker(\cA_r-\lambda I)^m$. Each vector $(\cA_r-\lambda I)^jU$ with $0\leq j\leq m-1$ belongs to $\mathcal D(\cA_r)$. We have
	\[\cS_r(t)U=\ee^{\lambda t}\sum_{j=0}^{m-1}\frac{t^j}{j!}(\cA_r-\lambda I)^jU,\qquad 0\leq t\leq t_r.\]
	At $t=t_r$, we obtain
	\begin{equation}\label{Eq.U_eigenspace_identity}
		U=\ee^{-\lambda t_r}\cS_r(t_r)U-\sum_{j=1}^{m-1}\frac{t_r^j}{j!}(\cA_r-\lambda I)^jU.
	\end{equation}
	
	We prove $\ker(\cA_r-\lambda I)^m\subset\cH_{r_+}$ by induction on $m$. For $m=1$, \eqref{Eq.U_eigenspace_identity} reduces to
	\[U=\ee^{-\lambda t_r}\cS_r(t_r)U\in\cH_{r_+},\]
	because $\cS_r(t_r):\cH_r\to\cH_{r_+}$ is bounded. Let $m\geq2$ and suppose that the inclusion has been proved for every order below $m$. For $1\leq j\leq m-1$, we have $(\cA_r-\lambda I)^jU\in\ker(\cA_r-\lambda I)^{m-j}\subset\cH_{r_+}$. Since $\cS_r(t_r)U\in\cH_{r_+}$, every term in \eqref{Eq.U_eigenspace_identity} belongs to $\cH_{r_+}$. This proves the induction step. Taking the union over $m\geq1$ gives $\mathcal G_\lambda(\cA_r)\subset\cH_{r_+}$. The radius $r_+$ and the time $t_r$ were chosen independently of $\lambda$ and $m$. Together with the finite algebraic multiplicity proved above, we obtain \eqref{Eq.macro_generalized_eigenspaces}, which completes the proof.
\end{proof}

\subsection{The phase correction}\label{Subsec.phase_proof}

We prove Proposition \ref{Prop.phase} at the fixed outer analytic radius $4\bar\varrho$. We rewrite the phase equation as a first-order system for $(\Psi,T_1\Psi)$ and verify the analytic regularity of its source. The spectral properties of the averaged generator then allow us to construct a particular solution that is polynomial in $\tau$.

\begin{proof}[Proof of Proposition \ref{Prop.phase}]
	We seek a polynomial $\Psi$ satisfying \eqref{Eq.phase_equation}, namely
	\[(T_2T_1+c_a^2A_0)\Psi=-\av{f^2}(B-2)V.\]
	Introduce the pair $P=(P_1,P_2)=(\Psi,T_1\Psi)$. Since $T_1=\partial_\tau+B-1$ and $T_2=\partial_\tau+B-2$, the phase equation becomes
	\begin{align*}
		(\partial_\tau+B)P_1=P_1+P_2,\quad 
		(\partial_\tau+B)P_2=2P_2-c_a^2A_0P_1-\av{f^2}(B-2)V.
	\end{align*}
	Comparing with \eqref{Eq.macro_pair}, we obtain
	\begin{equation}\label{Eq.phase_pair_system}
		\partial_\tau P=(\mathscr A+I)P+G,\qquad G:=\bigl(0,-\av{f^2}(B-2)V\bigr).
	\end{equation}
	Equivalently, $\partial_\tau(\ee^{-\tau}P)=\mathscr A(\ee^{-\tau}P)+\ee^{-\tau}G$. Thus $\ee^{-\tau}P$ satisfies the averaged equation \eqref{Eq.macro_pair} with source $\ee^{-\tau}G$. We will construct a polynomial solution of \eqref{Eq.phase_pair_system}.
	
	We first verify that $G\in\cH_{\varrho_2}$. Choose $r'\in(\varrho_2,R_v)$. The coefficients $\bvec$ and $V$ occur in \eqref{Eq.background_coefficient_map}, so \eqref{Eq.common_coefficient_radius} and the definition \eqref{Eq.analytic_norm} give $\bvec,V\in\mathcal B_{r'}$. Since $B=\bvec\cdot\nabla$, the derivative estimate \eqref{Eq.analytic_radius_recovery} and Lemma \ref{Lem.analytic_algebra} imply
	\begin{align*}
		\|(B-2)V\|_{\mathcal B_{\varrho_2}}\leq C_{\varrho_2}\|\bvec\|_{\mathcal B_{\varrho_2}}\|\nabla V\|_{\mathcal B_{\varrho_2}}+2\|V\|_{\mathcal B_{\varrho_2}}\leq\left(\frac{C_{\varrho_2}\|\bvec\|_{\mathcal B_{\varrho_2}}}{r'-\varrho_2}+2\right)\|V\|_{\mathcal B_{r'}}<\infty.
	\end{align*}
	The first component of $G$ vanishes. Hence \eqref{Eq.acoustic_pair_norm} gives $\|G\|_{\cH_{\varrho_2}}=\av{f^2}\|(B-2)V\|_{\mathcal B_{\varrho_2}}<\infty$.
	
	We now work in the complexification of the radial space $\cH_{\varrho_2}$ and define
	\[\mathscr L:=\cA_{\varrho_2}+I,\qquad \mathcal D(\mathscr L):=\mathcal D(\cA_{\varrho_2}).\]
	The conditions \eqref{Eq.phase_radius_condition} and \eqref{Eq.approximation_radii} allow us to apply Proposition \ref{Prop.macro_spectrum} at $\varrho_2$. Therefore $\mathscr L$ has compact resolvent. If $0\in\spec(\mathscr L)$, let $E_0$ be the Riesz projection of $\mathscr L$ at zero. The space $E_0\cH_{\varrho_2}$ is finite-dimensional, and there is an integer $\nu_0\geq1$ such that $\mathscr L^{\nu_0}E_0=0$. If $0\notin\spec(\mathscr L)$, define $E_0=0$ and $\nu_0=0$.
	
	The restriction of $\mathscr L$ to $(I-E_0)\cH_{\varrho_2}$, with domain $(I-E_0)\mathcal D(\mathscr L)$, is invertible. Consequently, the source $(I-E_0)G$ in the equation $P'-\mathscr LP=(I-E_0)G$ has the stationary solution
	\[-\bigl(\mathscr L|_{(I-E_0)\cH_{\varrho_2}}\bigr)^{-1}(I-E_0)G.\]
	On $E_0\cH_{\varrho_2}$, the operator $\mathscr L$ is nilpotent. Every vector in $E_0\cH_{\varrho_2}$ belongs to the domain of every power of $\mathscr L$. For $\nu_0\geq1$, we therefore have
	\begin{align*}
		\left(\partial_\tau-\mathscr L\right)\sum_{j=0}^{\nu_0-1}\frac{\tau^{j+1}}{(j+1)!}\mathscr L^jE_0G
		&=\sum_{j=0}^{\nu_0-1}\frac{\tau^j}{j!}\mathscr L^jE_0G-\sum_{j=1}^{\nu_0}\frac{\tau^j}{j!}\mathscr L^jE_0G\\
		&=E_0G-\frac{\tau^{\nu_0}}{\nu_0!}\mathscr L^{\nu_0}E_0G=E_0G.
	\end{align*}
	Thus the source $E_0G$ is handled by a polynomial of degree at most $\nu_0$.
	
	Adding the stationary solution and the polynomial solution, define
	\begin{equation}\label{Eq.polynomial_phase}
		P(\tau):=-\bigl(\mathscr L|_{(I-E_0)\cH_{\varrho_2}}\bigr)^{-1}(I-E_0)G+\sum_{j=0}^{\nu_0-1}\frac{\tau^{j+1}}{(j+1)!}\mathscr L^jE_0G.
	\end{equation}
	The sum is empty when $\nu_0=0$. The constant coefficient belongs to $(I-E_0)\mathcal D(\mathscr L)$, and every coefficient in the finite sum belongs to $E_0\cH_{\varrho_2}\subset\mathcal D(\mathscr L)$. Hence
	\[\partial_\tau P-\mathscr LP=G\quad\text{in }\cH_{\varrho_2},\qquad \tau\in\R.\]
	By \eqref{Eq.macro_generator_domain}, this is precisely \eqref{Eq.phase_pair_system}.
	
	The coefficients of $\mathscr A$ and the source $G$ are real, so replacing $P$ by $\operatorname{Re}P$ preserves \eqref{Eq.phase_pair_system}. The coefficients of $P$ are radial because the construction takes place in the radial space $\cH_{\varrho_2}$. Let $\Psi=\operatorname{Re} P_1$, take $q_*=\nu_0$, and write
	\[\Psi(\tau,x)=\sum_{j=0}^{q_*}\tau^j\Psi_j(x).\]
	Every coefficient of $P$ belongs to $\cH_{\varrho_2}$, so each $\Psi_j$ belongs to $\mathcal B_{\varrho_2}$ and is real-valued and radial. The first component of \eqref{Eq.phase_pair_system} gives $P_2=T_1\Psi$, and the second component gives
	\[T_2P_2+c_a^2A_0\Psi=-\av{f^2}(B-2)V.\]
	Therefore $\Psi$ satisfies \eqref{Eq.phase_equation} in spatial distributions for every $\tau\in\R$. Since the coefficients $\Psi_j$ are analytic in $x$, the phase equation also holds pointwise. This completes the proof.
\end{proof}

\subsection{The fourth-order approximation}\label{Subsec.approximation_proof}

We prove Lemma \ref{Lem.second_approximation} and Proposition \ref{Prop.approximation}. The phase $\Psi$ is fixed by Proposition \ref{Prop.phase}. We first justify composition with $\theta_\sigma$ and record the exact differentiation formula. We then determine the third- and fourth-order corrections from the coefficients of the residual. Throughout this subsection, derivatives of a function of $(\theta;\tau,x)$ are taken before evaluation at $(\theta_\sigma;\log\sigma,x)$.

\begin{lemma}\label{Lem.phase_functional_calculus}
	Let $r_h>3\overline\varrho$ and let
	\begin{equation}\label{Eq.finite_periodic_profile}
		h(\theta;\tau,x)=\sum_{m=1}^{m_h}\sum_{j=0}^{J_h}\tau^j h_{m,j}(x)\phi_m(\theta),\qquad h_{m,j}\in\mathcal B_{r_h},\quad \phi_m\in C^\omega_{\operatorname{ap},L}(\R),
	\end{equation}
	where $m_h\geq1$ and $J_h\geq0$ are integers. There is $\sigma_h\geq2$ such that, for every $a,b,c\in\Z_{\geq0}$, there is $C_{a,b,c}>0$ satisfying
	\begin{equation}\label{Eq.finite_phase_composition_bound}
		\mathfrak M_{3\overline\varrho}\bigl((\partial_\theta^a\partial_\tau^b\nabla_x^ch)(\theta_\sigma;\log\sigma,x)\bigr)\leq C_{a,b,c}(1+\log\sigma)^{J_h},\qquad \sigma\geq\sigma_h.
	\end{equation}
	The same conclusion holds for analytic $2L$-periodic functions $\phi_m$.
\end{lemma}

\begin{proof}
	By Lemma \ref{Lem.analytic_to_multiplier} with $r=3\overline\varrho$ and $r'=r_h$, we have
	\[\mathfrak M_{3\overline\varrho}(\nabla^ch_{m,j})\leq C_{c,r_h}\|h_{m,j}\|_{\mathcal B_{r_h}},\qquad c\in\Z_{\geq0},\quad 1\leq m\leq m_h,\quad 0\leq j\leq J_h.\]
	
	By Proposition \ref{Prop.phase}, the phase correction has the form $\Psi(\tau,x)=\sum_{j=0}^{q_*}\tau^j\Psi_j(x)$ with $\Psi_j\in\mathcal B_{\varrho_2}$ and $\varrho_2=4\overline\varrho$. Applying Lemma \ref{Lem.analytic_to_multiplier} with $r=3\overline\varrho$, $r'=\varrho_2$ and $c=0$, we obtain, for $\sigma\geq2$,
	\begin{align*}
		\mathfrak M_{3\overline\varrho}\bigl(\sigma^{-1}\Psi(\log\sigma,\cdot)\bigr)\leq\sigma^{-1}\sum_{j=0}^{q_*}(\log\sigma)^j\mathfrak M_{3\overline\varrho}(\Psi_j)\leq\frac{C}{\sigma}(1+\log\sigma)^{q_*}\sum_{j=0}^{q_*}\|\Psi_j\|_{\mathcal B_{\varrho_2}}\to0
	\end{align*}
	as $\sigma\to\infty$.
	
	The finitely many functions $\phi_m$ are real analytic and $2L$-periodic. Analyticity on a compact period interval and periodic continuation give $r_c>0$ and $C_\phi>0$ such that every $\phi_m$ extends holomorphically to $\mathbb S_{4r_c}$ and
	\[\max_{1\leq m\leq m_h}\sup_{z\in\mathbb S_{4r_c}}|\phi_m(z)|\leq C_\phi.\]
	See also Lemma \ref{Lem.analytic_extension} for construction. For real $\sigma$ and integers $a,k\geq0$, Cauchy's estimate on the disk of radius $2r_c$ centered at $\sigma$ gives
	\begin{align*}
		\frac{|\phi_m^{(a+k)}(\sigma)|}{k!}\leq\frac{C_\phi(a+k)!}{k!(2r_c)^{a+k}}=\frac{C_\phi a!}{(2r_c)^{a+k}}\binom{a+k}{a}\leq C_\phi a!r_c^{-a-k},
	\end{align*}
	where we used $\binom{a+k}{a}\leq2^{a+k}$. Choose $\sigma_h\geq2$ such that
	\[\sup_{\sigma\geq\sigma_h}\mathfrak M_{3\overline\varrho}\bigl(\sigma^{-1}\Psi(\log\sigma,\cdot)\bigr)\leq\frac{r_c}{2}.\]
	The threshold $\sigma_h$ is independent of $a,b,c$.
	
	For $\sigma\geq\sigma_h$, Taylor's formula gives
	\[\phi_m^{(a)}(\theta_\sigma)=\sum_{k=0}^\infty\frac{\phi_m^{(a+k)}(\sigma)}{k!}\bigl(\sigma^{-1}\Psi(\log\sigma,x)\bigr)^k.\]
	The product estimate \eqref{Eq.multiplier_algebra} implies
	\[\mathfrak M_{3\overline\varrho}\bigl((\sigma^{-1}\Psi(\log\sigma,\cdot))^k\bigr)\leq\mathfrak M_{3\overline\varrho}\bigl(\sigma^{-1}\Psi(\log\sigma,\cdot)\bigr)^k,\qquad k\geq0.\]
	Therefore,
	\begin{align*}
		\sum_{k=0}^\infty\frac{|\phi_m^{(a+k)}(\sigma)|}{k!}\mathfrak M_{3\overline\varrho}\bigl((\sigma^{-1}\Psi(\log\sigma,\cdot))^k\bigr)&\leq C_\phi a!r_c^{-a}\sum_{k=0}^\infty\left(\frac{\mathfrak M_{3\overline\varrho}(\sigma^{-1}\Psi(\log\sigma,\cdot))}{r_c}\right)^k\\
		&\leq C_\phi a!r_c^{-a}\sum_{k=0}^\infty2^{-k}=2C_\phi a!r_c^{-a}.
	\end{align*}
	The Taylor partial sums are Cauchy in every $C^N(\overline\Omega)$ norm, since the multiplier norm controls each spatial derivative in $L^\infty$. Their pointwise limit is $\phi_m^{(a)}(\theta_\sigma)$ because $|\sigma^{-1}\Psi(\log\sigma,x)|\leq r_c/2$. Passing to the limit in each derivative and then in the nonnegative series defining the multiplier norm gives
	\[\mathfrak M_{3\overline\varrho}\bigl(\phi_m^{(a)}(\theta_\sigma)\bigr)\leq2C_\phi a!r_c^{-a},\qquad \sigma\geq\sigma_h,\quad 1\leq m\leq m_h.\]
	
	Finally, the derivatives in \eqref{Eq.finite_phase_composition_bound} are taken with $(\theta,\tau,x)$ independent, before evaluation at $(\theta_\sigma;\log\sigma,x)$. If $b>J_h$, the derivative vanishes. For $0\leq b\leq J_h$, differentiating \eqref{Eq.finite_periodic_profile} gives
	\[(\partial_\theta^a\partial_\tau^b\nabla_x^ch)(\theta_\sigma;\log\sigma,x)=\sum_{m=1}^{m_h}\sum_{j=b}^{J_h}\frac{j!}{(j-b)!}(\log\sigma)^{j-b}(\nabla^ch_{m,j})(x)\phi_m^{(a)}(\theta_\sigma).\]
	Using \eqref{Eq.multiplier_algebra}, \eqref{Eq.analytic_to_multiplier} and the estimate for $\phi_m^{(a)}(\theta_\sigma)$, we obtain
	\begin{align*}
		\mathfrak M_{3\overline\varrho}\bigl((\partial_\theta^a\partial_\tau^b\nabla_x^ch)(\theta_\sigma;\log\sigma,\cdot)\bigr)&\leq\sum_{m=1}^{m_h}\sum_{j=b}^{J_h}\frac{j!}{(j-b)!}(\log\sigma)^{j-b}\mathfrak M_{3\overline\varrho}(\nabla^ch_{m,j})\mathfrak M_{3\overline\varrho}\bigl(\phi_m^{(a)}(\theta_\sigma)\bigr)\\
		&\leq2C_\phi a!r_c^{-a}C_{c,r_h}\sum_{m=1}^{m_h}\sum_{j=b}^{J_h}\frac{j!}{(j-b)!}(\log\sigma)^{j-b}\|h_{m,j}\|_{\mathcal B_{r_h}}\\
		&\leq C_{a,b,c}(1+\log\sigma)^{J_h},\qquad \sigma\geq\sigma_h.
	\end{align*}
	This proves \eqref{Eq.finite_phase_composition_bound}. Only analyticity and $2L$-periodicity of the functions $\phi_m$ were used, so the conclusion also holds for analytic $2L$-periodic functions.
\end{proof}

\begin{lemma}\label{Lem.cell_differentiation}
	Let $h$ be as in \eqref{Eq.finite_periodic_profile}, let $n\in\Z_{\geq0}$, and define
	\[h_n(\sigma,x):=\sigma^{-n}h(\theta_\sigma;\log\sigma,x).\]
	Then
	\begin{equation}\label{Eq.correction_differentiation}
		\begin{aligned}
			\cP h_n={}&\sigma^{-n}\bigl(\lambda^2-\sigma^{-4}|\nabla\Psi|^2\bigr)\partial_\theta^2h+2\sigma^{-n-1}\lambda T_n\partial_\theta h\\
			&+\sigma^{-n-2}(T_{n+1}T_n+A_0)h+\sigma^{-n-3}\bigl[\bigl((T_2T_1+A_0)\Psi\bigr)\partial_\theta h-2\nabla\Psi\cdot\nabla_x\partial_\theta h\bigr].
		\end{aligned}
	\end{equation}
	On the right-hand side, $h$ and its derivatives are evaluated at $(\theta_\sigma;\log\sigma,x)$, and $\Psi$ and its derivatives at $(\log\sigma,x)$.
\end{lemma}

\begin{proof}
	By \eqref{Eq.slow_calculus} and $D\theta_\sigma=\lambda$, we have
	\begin{align*}
		Dh_n&=\sigma^{-n}\lambda\partial_\theta h+\sigma^{-n-1}T_nh,\\
		D^2h_n&=\sigma^{-n}\lambda^2\partial_\theta^2h+2\sigma^{-n-1}\lambda T_n\partial_\theta h+\sigma^{-n}(D\lambda)\partial_\theta h+\sigma^{-n-2}T_{n+1}T_nh.
	\end{align*}
	The spatial chain rule gives
	\begin{align*}
		A_0\bigl(h(\theta_\sigma;\log\sigma,x)\bigr)=(A_0h)(\theta_\sigma;\log\sigma,x)+\sigma^{-1}(A_0\Psi)\partial_\theta h-2\sigma^{-1}\nabla\Psi\cdot\nabla_x\partial_\theta h-\sigma^{-2}|\nabla\Psi|^2\partial_\theta^2h.
	\end{align*}
	Here $A_0h$ on the right is computed with $\theta$ fixed. Adding $\sigma^{-n-2}$ times this identity to the formula for $D^2h_n$, and using $D\lambda=\sigma^{-3}T_2T_1\Psi$, proves \eqref{Eq.correction_differentiation}.
\end{proof}

We use $\Order_{\rm an}$ below at the fixed multiplier radius $3\overline\varrho$. The spatial coefficients of $\Psi$, the background coefficients, and every fixed finite product of their spatial derivatives belong to $\mathcal B_{r_h}$ for each $r_h\in(3\overline\varrho,\varrho_2)$, by \eqref{Eq.common_coefficient_radius}, \eqref{Eq.analytic_radius_recovery} and Lemma \ref{Lem.analytic_algebra}. Lemma \ref{Lem.analytic_to_multiplier} therefore gives finite multiplier norms at $3\overline\varrho$ for the spatial coefficients in the finite expansion. These estimates concern spatial norms at fixed $\sigma$; a time derivative falling on a fast periodic factor need not improve the decay rate.

\begin{proof}[Proof of Lemma \ref{Lem.second_approximation}]
	The coefficients of $A_2$ in \eqref{Eq.second_approximation_definition} are polynomials in $\tau$ with spatial coefficients in $\mathcal B_{r_h}$ for every fixed $r_h\in(0,\varrho_2)$. Choose $r_h\in(3\overline\varrho,\varrho_2)$ and apply Lemma \ref{Lem.phase_functional_calculus} to $f$, $A_2$ and the finitely many derivatives used below. In particular,
	\[z_2-f(\theta_\sigma)=\Order_{\rm an}\bigl(\sigma^{-2}(1+\log\sigma)^M\bigr).\]
	By \eqref{Eq.correction_differentiation}, with $n=2$ and $h=A_2$, and the finite binomial expansion of $z_2^p$, we have
	\begin{align*}
		\cN(z_2)={}&\sigma^{-2}\bigl(L_fA_2-2U_\Psi f^p-Vf\bigr)(\theta_\sigma;\log\sigma,x)\\
		&+\sigma^{-3}\bigl[\bigl((T_2T_1+A_0)\Psi\bigr)f'+2T_2\partial_\theta A_2\bigr](\theta_\sigma;\log\sigma,x)+\Order_{\rm an}\bigl(\sigma^{-4}(1+\log\sigma)^M\bigr).
	\end{align*}
	Every omitted term is a finite product of composed periodic profiles and slow coefficients, with an explicit factor $\sigma^{-j}$ for some $j\geq4$. Lemma \ref{Lem.phase_functional_calculus} bounds these products. The first coefficient vanishes by \eqref{Eq.second_approximation_definition}, and the second is $R_3$ by \eqref{Eq.third_residual}. Since
	\[\|a\|_{\mathcal B_{\varrho_1}}\leq |\Omega|^{1/2}\mathfrak M_{\varrho_1}(a)\leq |\Omega|^{1/2}\mathfrak M_{3\overline\varrho}(a),\]
	we obtain \eqref{Eq.second_approximation_bounds}. Finally, $\lambda-1=\sigma^{-2}U_\Psi$ tends to zero in the multiplier norm. Increasing $\sigma_2$ gives \eqref{Eq.second_approximation_domain} and completes the proof.
\end{proof}

We next define the remaining corrections for independent $(\theta,\tau,x)\in\R\times\R\times\overline\Omega$. By \eqref{Eq.kappa_H} and the phase equation, the third-order equation has the solution
\begin{equation}\label{Eq.third_cell}
	A_3(\theta;\tau,x):=-2((B-2)V)(x)H(\theta),\qquad L_fA_3=-R_3,\qquad \av{A_3f'}=0.
\end{equation}
Define the fourth-order source by
\begin{equation}\label{Eq.R4}
	\begin{aligned}
		R_4(\theta;\tau,x):={}&-\bigl(U_\Psi^2-|\nabla\Psi|^2\bigr)f(\theta)^p+2U_\Psi\partial_\theta^2A_2\\
		&+(T_3T_2+A_0-V)A_2+2T_3\partial_\theta A_3+\frac{p(p-1)}2f(\theta)^{p-2}A_2^2.
	\end{aligned}
\end{equation}
On the right-hand side, $A_2,A_3$ and their derivatives are evaluated at $(\theta;\tau,x)$, and the slow coefficients at $(\tau,x)$. Since $f,A_2$ are even and antiperiodic in $\theta$, and $A_3$ is odd and antiperiodic, each term of $R_4$ is even and antiperiodic. Thus \eqref{Eq.cell_range} gives a unique function $A_4(\cdot;\tau,x)\in C^\omega_{\operatorname{ap},L,\rm e}(\R)$ satisfying
\begin{equation}\label{Eq.fourth_cell}
	L_fA_4(\theta;\tau,x)=-R_4(\theta;\tau,x).
\end{equation}
Finally, define
\begin{align}
	&\qquad\qquad\qquad\quad z_3(\sigma,x):=z_2(\sigma,x)+\sigma^{-3}A_3(\theta_\sigma;\log\sigma,x),\label{Eq.third_approximation}\\
	&z_{\app}(\sigma,x):=z_3(\sigma,x)+\sigma^{-4}A_4(\theta_\sigma;\log\sigma,x)=f(\theta_\sigma)+\sum_{j=2}^4\sigma^{-j}A_j(\theta_\sigma;\log\sigma,x).\label{Eq.zapp}
\end{align}

\begin{proof}[Proof of Proposition \ref{Prop.approximation}]
	We first verify the regularity of $A_4$. The right-hand side of \eqref{Eq.R4} is a finite sum of products of slow polynomials and analytic antiperiodic functions in $\theta$. After collecting terms, write
	\[R_4(\theta;\tau,x)=\sum_{m=1}^{m_4}\sum_{j=0}^{J_4}\tau^j c_{m,j}(x)\phi_m(\theta),\qquad \phi_m\in C^\omega_{\operatorname{ap},L,\rm e}(\R).\]
	For every fixed $r_h\in(0,\varrho_2)$, the coefficients $c_{m,j}$ belong to $\mathcal B_{r_h}$, by \eqref{Eq.fixed_algebra} and \eqref{Eq.analytic_radius_recovery}. For each $m$, let $v_m$ be the unique even antiperiodic solution of $L_fv_m=-\phi_m$. Analytic ODE regularity gives $v_m\in C^\omega_{\operatorname{ap},L,\rm e}(\R)$, and uniqueness in \eqref{Eq.fourth_cell} yields
	\[A_4(\theta;\tau,x)=\sum_{m=1}^{m_4}\sum_{j=0}^{J_4}\tau^j c_{m,j}(x)v_m(\theta).\]
	Thus $A_2,A_3,A_4$ all have the form \eqref{Eq.finite_periodic_profile} with some $r_h>3\overline\varrho$. Their slow coefficients are real-valued and radial. Choose a common threshold $\sigma_h$ in Lemma \ref{Lem.phase_functional_calculus} for the finitely many profiles $f,g,H,K,v_m$ and all terms used in the following calculation.

	Here $A_j$ and their derivatives are evaluated at $(\theta_\sigma;\log\sigma,x)$, while $f$ is evaluated at $\theta_\sigma$ and slow coefficients at $(\log\sigma,x)$. Using \eqref{Eq.correction_differentiation} gives
	\begin{align*}
		\cP(\sigma^{-2}A_2)&=\sigma^{-2}\partial_\theta^2A_2+2\sigma^{-3}T_2\partial_\theta A_2\\
		&\quad+\sigma^{-4}\bigl[(T_3T_2+A_0)A_2+2U_\Psi\partial_\theta^2A_2\bigr]+\Order_{\rm an}\bigl(\sigma^{-5}(1+\log\sigma)^M\bigr),\\
		\cP(\sigma^{-3}A_3)&=\sigma^{-3}\partial_\theta^2A_3+2\sigma^{-4}T_3\partial_\theta A_3+\Order_{\rm an}\bigl(\sigma^{-5}(1+\log\sigma)^M\bigr),\\
		\cP(\sigma^{-4}A_4)&=\sigma^{-4}\partial_\theta^2A_4+\Order_{\rm an}\bigl(\sigma^{-5}(1+\log\sigma)^M\bigr).
	\end{align*}
	The polynomial nonlinearity has the expansion
	\begin{align*}
		z_{\app}^p=f(\theta_\sigma)^p+pf(\theta_\sigma)^{p-1}\sum_{j=2}^4\sigma^{-j}A_j+\frac{p(p-1)}2\sigma^{-4}f(\theta_\sigma)^{p-2}A_2^2+\Order_{\rm an}\bigl(\sigma^{-5}(1+\log\sigma)^M\bigr).
	\end{align*}
	Indeed, among the quadratic products of the corrections, only $\sigma^{-4}A_2^2$ has order below five, and every product of at least three corrections has order at least six. The product estimate \eqref{Eq.multiplier_algebra} controls all the finite products in the expansion. 

	Using \eqref{Eq.shifted_leading_error} and adding $-\sigma^{-2}V(z_{\app}-f(\theta_\sigma))$, we obtain
	\begin{equation}\label{Eq.finite_residual_expansion}
		\begin{aligned}
			\cN(z_{\app})={}&\sigma^{-2}\bigl(L_fA_2-2U_\Psi f^p-Vf\bigr)+\sigma^{-3}(L_fA_3+R_3)+\sigma^{-4}(L_fA_4+R_4)\\
			&+\Order_{\rm an}\bigl(\sigma^{-5}(1+\log\sigma)^M\bigr).
		\end{aligned}
	\end{equation}
	The three displayed coefficients vanish by \eqref{Eq.second_approximation_definition}, \eqref{Eq.third_cell} and \eqref{Eq.fourth_cell}. This proves the required residual bound in $\mathfrak M_{3\overline\varrho}$. The definition \eqref{Eq.zapp} and Lemma \ref{Lem.phase_functional_calculus} also give
	\begin{equation}\label{Eq.zapp-f-multiplier-est}
		\mathfrak M_{3\overline\varrho}(z_{\app}-f(\theta_\sigma))\leq C\sigma^{-2}(1+\log\sigma)^M.
	\end{equation}
	There are only finitely many slow polynomials in the residual, so $M$ can be fixed independently of $\sigma$. The inequality $\|a\|_{\mathcal B_{\varrho_1}}\leq |\Omega|^{1/2}\mathfrak M_{3\overline\varrho}(a)$ proves \eqref{Eq.approximation_bounds}; \eqref{Eq.analytic_radius_recovery} gives all the spatial derivative estimates in Proposition \ref{Prop.approximation} with the same $M$. Increasing the common threshold to $\sigma_{\app}\geq\sigma_2$ preserves \eqref{Eq.second_approximation_domain}. Finally, by \eqref{Eq.zapp}, $z_{\rm app}$ is real-valued, radial, and smooth on $[\sigma_{\app},\infty)\times\overline\Omega$. This completes the proof.
\end{proof}

\subsection{The perturbation equations}\label{Subsec.perturbation_proofs}

We prove the perturbation lemmas in Subsection \ref{Subsec.formulation}. We first estimate the phase coefficients and the scalar remainder, and then apply the invertible frame $P_0$ to derive the exact two-component equation. Averaging the periodic coefficients of the two-component equation gives the acoustic generator and the small differential remainder in Lemma \ref{Lem.state_remainder}. The phase and approximation remain fixed throughout the subsection.

\begin{proof}[Proof of Lemma \ref{Lem.phase_coefficients}]
	By \eqref{Eq.phase_velocity} and \eqref{Eq.slow_calculus},
	\[\lambda-1=\sigma^{-2}U_\Psi,\qquad D\lambda=\sigma^{-3}T_2T_1\Psi.\]
	The spatial coefficients of the polynomials $T_1\Psi$ and $T_2T_1\Psi$ belong to $\mathcal B_{r'}$ for every $0<r'<\varrho_2$, by \eqref{Eq.analytic_radius_recovery} and Lemma \ref{Lem.analytic_algebra}. Their degrees in $\tau$ are at most $q_*$. Applying Lemma \ref{Lem.analytic_to_multiplier} with $3\overline\varrho<r'<\varrho_2$ gives
	\[\mathfrak M_{3\overline\varrho}(\lambda-1)+\sigma\mathfrak M_{3\overline\varrho}(D\lambda)\leq C\sigma^{-2}(1+\log\sigma)^{q_*},\qquad \sigma\geq2.\]
	Choose $\sigma_c\geq\sigma_{\app}$ so large that
	\[\sup_{\sigma\geq\sigma_c}\mathfrak M_{\varrho_1}(\lambda-1)\leq\tfrac12.\]
	The geometric series then gives $\mathfrak M_{\varrho_1}(\lambda^{-1})\leq2$, and the pointwise bounds on $\lambda$ follow as well.

	Apply Lemma \ref{Lem.phase_functional_calculus} to the fixed periodic functions $f,f',f'',K$ and their first derivatives. Their compositions with $\sigma+s\sigma^{-1}\Psi$, $0\leq s\leq1$, have uniform multiplier bounds by the same Taylor-series argument. Hence, for each $h\in\{f,f',f'',K\}$,
	\[h(\theta_\sigma)-h(\sigma)=\sigma^{-1}\Psi\int_0^1h'(\sigma+s\sigma^{-1}\Psi)\,\mathrm ds\]
	implies \eqref{Eq.periodic_composition_bounds}. By \eqref{Eq.zapp-f-multiplier-est}, we have $\mathfrak M_{\varrho_1}(z_{\app}-f(\theta_\sigma))\leq C\sigma^{-2}(1+\log\sigma)^M$. Together with the uniform bound for $f(\theta_\sigma)$, this also bounds $z_{\app}$. Choose $M$ in \eqref{Eq.nu_sigma} to dominate the finitely many logarithmic degrees. The required estimates follow.
\end{proof}

The phase bounds allow us to expand the exact polynomial nonlinearity around the approximate solution. We keep the scalar source and the quadratic remainder separate.

\begin{proof}[Proof of Lemma \ref{Lem.remainder_equation}]
	The equivalence in \eqref{Eq.remainder_equation} is the polynomial identity
	\[\cN(z_\app+\zeta)=\cN(z_\app)+\bigl(D^2+\sigma^{-2}(A_0-V)+pz_\app^{p-1}\bigr)\zeta+\cQ_\sigma(\zeta).\]
	The definitions of $R_\sigma$ and $\Lambda_\sigma$ give \eqref{Eq.linearized_form} exactly.

	By \eqref{Eq.zapp} and Lemma \ref{Lem.phase_functional_calculus}, the finite approximation satisfies
	\begin{align*}
		z_{\app}={}&f(\theta_\sigma)+\sigma^{-2}\left[\frac2{p-1}U_\Psi f(\theta_\sigma)+Vg(\theta_\sigma)\right]+\Order_{\rm an}\bigl(\sigma^{-3}(1+\log\sigma)^M\bigr)
	\end{align*}
	at the multiplier radius $3\overline\varrho$. Expanding the polynomial $pz_{\app}^{p-1}$ and using the definition of $K$ gives
	\begin{align*}
		pz_{\app}^{p-1}=pf(\theta_\sigma)^{p-1}+2p\sigma^{-2}U_\Psi f(\theta_\sigma)^{p-1}+\sigma^{-2}K(\theta_\sigma)V+\Order_{\rm an}\bigl(\sigma^{-3}(1+\log\sigma)^M\bigr).
	\end{align*}
	Since $\lambda^2=1+2\sigma^{-2}U_\Psi+\sigma^{-4}U_\Psi^2$, we obtain
	\begin{equation}\label{Eq.finite_linearized_potential}
		pz_{\app}^{p-1}=\lambda^2pf(\theta_\sigma)^{p-1}+\sigma^{-2}K(\theta_\sigma)V+\Order_{\rm an}\bigl(\sigma^{-3}(1+\log\sigma)^M\bigr).
	\end{equation}
	Thus \eqref{Eq.linearized_original} implies $\mathfrak M_{\varrho_1}(R_\sigma)\leq C\sigma^{-1}(1+\log\sigma)^M$. The bound for $\cR_{\app}$ follows from Proposition \ref{Prop.approximation}.

	For $0<r<\overline\varrho$, the algebra estimate \eqref{Eq.fixed_algebra} has a constant locally uniform in $r$. Finally,
	\[\cQ_\sigma(\zeta)=\sum_{a=2}^p\binom pa z_\app^{p-a}\zeta^a,\qquad \zeta_1^a-\zeta_2^a=(\zeta_1-\zeta_2)\sum_{b=0}^{a-1}\zeta_1^{a-1-b}\zeta_2^b,\quad\forall\ a\in\Z\cap[2,p].\]
	For $\|\zeta_i\|_{\mathcal B_r}\leq1$, the algebra estimate \eqref{Eq.fixed_algebra} gives
	\begin{align*}
		\|\cQ_\sigma(\zeta_1)\|_{\mathcal B_r} &\leq C\sum_{j=2}^p\|\zeta_1\|_{\mathcal B_r}^j \leq C\|\zeta_1\|_{\mathcal B_r}^2,\\
		\|\cQ_\sigma(\zeta_1)-\cQ_\sigma(\zeta_2)\|_{\mathcal B_r} &\leq C(\|\zeta_1\|_{\mathcal B_r}+\|\zeta_2\|_{\mathcal B_r}) \|\zeta_1-\zeta_2\|_{\mathcal B_r}.
	\end{align*}
	Here we used the uniform multiplier bound for $z_\app$.
\end{proof}

We differentiate the transformation \eqref{Eq.zeta-w-W-def} to prove \eqref{Eq.exact_state} and the reconstruction identities $D\zeta_W=\lambda w_W$, $\cL_\sigma\zeta_W=\mathcal F$.

\begin{proof}[Proof of Lemma \ref{Lem.state_equation}]
	The oscillator energy \eqref{Eq.oscillator_energy} gives
	\[\det P_0=\frac{p+1}{p-1}(f')^2-\frac2{p-1}ff''=\frac2{p-1}.\]
	Using $f''=-f^p$, we compute
	\[P_0^{-1}=\begin{pmatrix}(p+1)f'/2&-f\\-(p-1)f''/2&(p-1)f'/2\end{pmatrix},\qquad P_0'=\mathsf P_0P_0-P_0N_0.\]
	Thus $P_0^{-1}(0,1)^T=(-f,(p-1)f'/2)^T$, proving \eqref{Eq.frame_identity}.

	Let $w:=\lambda^{-1}D\zeta$ and $U:=(\zeta,w)^T$. By \eqref{Eq.linearized_form}, we have
	\[DU=\lambda\mathsf P_0(\sigma+\sigma^{-1}\Psi)U+\sigma^{-2}\pair0{\mathfrak R_\sigma(W,h)}.\]
	Indeed, the second component follows from
	\[Dw=\lambda^{-1}D^2\zeta-\frac{D\lambda}{\lambda}w=-\lambda pf(\sigma+\sigma^{-1}\Psi)^{p-1}\zeta+\sigma^{-2}\lambda^{-1}(h-\Lambda_\sigma\zeta)-\frac{D\lambda}{\lambda}w.\]
	Write $D_\sigma:=\diag(1,\sigma^{-1})$ for the reconstruction of $(\zeta,w)$ from $W$. Substitution of $U=P_0(\sigma+\sigma^{-1}\Psi)D_\sigma W$ and \eqref{Eq.frame_identity} gives
	\[D_\sigma DW-\lambda N_0D_\sigma W+(DD_\sigma)W=\sigma^{-2}P_0(\sigma+\sigma^{-1}\Psi)^{-1}\pair0{\mathfrak R_\sigma(W,h)}.\]
	Multiplying by $\sigma D_\sigma^{-1}$ yields exactly \eqref{Eq.exact_state}; here $DD_\sigma=\diag(0,-\sigma^{-2})$ produces the term $\pi$ in the equation for $\pi$. Each operation is reversible. Therefore a solution of the pair equation satisfies both $D\zeta_W=\lambda w_W$ and  $\cL_\sigma\zeta_W=\mathcal F$.
\end{proof}

We now isolate $\mathscr A+\cB(\sigma)$ in \eqref{Eq.exact_state} and estimate the coefficients of $\mathscr E_\sigma$.

\begin{proof}[Proof of Lemma \ref{Lem.state_remainder}]
	Starting from \eqref{Eq.exact_state}, we separate the source by using
	\[\mathfrak R_\sigma(W,h)=\mathfrak R_\sigma(W,0)+\lambda^{-1}h.\]
	Moving $B\psi$ and $B\pi$ to the right gives
	\begin{align*}
		\sigma\partial_\sigma\psi&=-B\psi+\lambda\pi-\sigma^{-1}f(\theta_\sigma)\mathfrak R_\sigma(W,0)-\sigma^{-1}\lambda^{-1}f(\theta_\sigma)h,\\
		\sigma\partial_\sigma\pi&=(1-B)\pi+\frac{p-1}{2}f'(\theta_\sigma)\mathfrak R_\sigma(W,0)+\frac{p-1}{2}\lambda^{-1}f'(\theta_\sigma)h.
	\end{align*}
	By the definitions of $\mathscr A$ and $\cB(\sigma)$,
	\[(\mathscr A+\cB(\sigma))W=\binom{-B\psi+\pi}{(1-B)\pi-a_0(\sigma)A_0\psi-a_1(\sigma)V\psi}.\]
	We therefore define
	\begin{equation}\label{Eq.state_remainder_definition}
		\begin{aligned}
			(\mathscr E_\sigma W)_1&:=(\lambda-1)\pi-\sigma^{-1}f(\theta_\sigma)\mathfrak R_\sigma(W,0),\\
			(\mathscr E_\sigma W)_2&:=\frac{p-1}{2}f'(\theta_\sigma)\mathfrak R_\sigma(W,0)+a_0(\sigma)A_0\psi+a_1(\sigma)V\psi.
		\end{aligned}
	\end{equation}
	The two equations for $\psi,\pi$ now give \eqref{Eq.state_operator_form}, with
	\[\mathbf g_\sigma=\lambda^{-1}\binom{-\sigma^{-1}f(\theta_\sigma)}{(p-1)f'(\theta_\sigma)/2}.\]
	Since $\mathfrak R_\sigma(W,0)=-\lambda^{-1}\Lambda_\sigma\zeta_W-\sigma^2\frac{D\lambda}{\lambda}w_W$, the operator $\mathscr E_\sigma$ is linear and independent of $h$. The reconstruction \eqref{Eq.zeta-w-W-def} contains no spatial derivative of $W$, while $\Lambda_\sigma$ has order two. Hence $\mathscr E_\sigma$ has spatial order at most two. The identities \eqref{Eq.oscillator_averages} and \eqref{Eq.K_cancellation} also give
	\[\av{a_0}=\frac{p-1}{2}\av{(f')^2}=\frac{p-1}{p+3}=c_a^2,\qquad \av{a_1}=\frac{p-1}{2}\av{(f')^2(K-1)}=0.\]
	
	We now estimate $\mathscr E_\sigma$. Fix $0<r<r'<\overline\varrho$ and $W\in\cX_{r'}$. Increase $\sigma_c$ so that $\nu(\sigma)\leq1$ for $\sigma\geq\sigma_c$. The first component in \eqref{Eq.state_remainder_definition} requires a uniform coefficient bound for $\mathfrak R_\sigma(W,0)$. The second component requires cancellation between its principal part and $a_0(\sigma)A_0\psi+a_1(\sigma)V\psi$.
	
	The reconstruction formula \eqref{Eq.zeta-w-W-def} gives
	\begin{align*}
		\zeta_W-f'(\sigma)\psi&=\bigl(f'(\theta_\sigma)-f'(\sigma)\bigr)\psi+\frac2{p-1}\sigma^{-1}f(\theta_\sigma)\pi,\\
		w_W&=f''(\theta_\sigma)\psi+\frac{p+1}{p-1}\sigma^{-1}f'(\theta_\sigma)\pi.
	\end{align*}
	We use $f'(\sigma)\psi$ as the principal part of $\zeta_W$ because $f'(\sigma)$ is independent of $x$. Since $\varrho_1=2\overline\varrho$, Lemma \ref{Lem.analytic_to_multiplier} and the definitions of the norms imply
	\[\sum_{j=0}^2\mathfrak M_{\overline\varrho}(\nabla^ja)\leq C\|a\|_{\mathcal B_{\varrho_1}}\leq C|\Omega|^{1/2}\mathfrak M_{\varrho_1}(a)\]
	for every scalar coefficient $a$ with finite $\mathfrak M_{\varrho_1}(a)$. Applying this estimate to $f'(\theta_\sigma)-f'(\sigma)$ and $f(\theta_\sigma)$, and using \eqref{Eq.periodic_composition_bounds}, we obtain
	\[\sum_{j=0}^2\left(\mathfrak M_{\overline\varrho}\bigl(\nabla^j(f'(\theta_\sigma)-f'(\sigma))\bigr)+\sigma^{-1}\mathfrak M_{\overline\varrho}\bigl(\nabla^jf(\theta_\sigma)\bigr)\right)\leq C\nu(\sigma).\]
	
	By \eqref{Eq.common_coefficient_radius}, Lemma \ref{Lem.analytic_to_multiplier} and \eqref{Eq.multiplier_algebra}, we have $\mathfrak M_{\overline\varrho}(J_0^{-1}\nabla J_0)+\mathfrak M_{\overline\varrho}(V)\leq C$. Thus the coefficients of $A_0=-\Delta-J_0^{-1}\nabla J_0\cdot\nabla$ have bounded multiplier norms at $\overline\varrho$. Lemmas \ref{Lem.phase_coefficients} and \ref{Lem.remainder_equation}, together with \eqref{Eq.multiplier_algebra}, give
	\begin{align*}
		&\mathfrak M_{\overline\varrho}(\lambda^{-1})+\sum_{j=0}^2\mathfrak M_{\overline\varrho}\bigl(f^{(j)}(\theta_\sigma)\bigr)+\mathfrak M_{\overline\varrho}\bigl(K(\theta_\sigma)\bigr)\leq C,\quad 
		\mathfrak M_{\overline\varrho}(\lambda-1)+\mathfrak M_{\overline\varrho}(\lambda^{-1}-1)\leq C\sigma^{-1}\nu(\sigma),\\
		&\qquad\qquad\qquad\qquad\mathfrak M_{\overline\varrho}\bigl(K(\theta_\sigma)-K(\sigma)\bigr)+\mathfrak M_{\overline\varrho}(R_\sigma)+\sigma^2\mathfrak M_{\overline\varrho}(D\lambda/\lambda)\leq C\nu(\sigma).
	\end{align*}
	In particular, the coefficients of $\Lambda_\sigma=A_0+(K(\theta_\sigma)-1)V+R_\sigma$ have uniformly bounded multiplier norms at $\overline\varrho$.
	
	Since $f'(\sigma)$ is independent of $x$, the definition of $\mathfrak R_\sigma(W,0)$ gives
	\begin{equation}\label{Eq.state_R_principal_remainder}
		\begin{aligned}
			&\mathfrak R_\sigma(W,0)+f'(\sigma)\bigl(A_0+(K(\sigma)-1)V\bigr)\psi\\
			=&\,-\lambda^{-1}\Lambda_\sigma\bigl(\zeta_W-f'(\sigma)\psi\bigr)-(\lambda^{-1}-1)f'(\sigma)\bigl(A_0+(K(\sigma)-1)V\bigr)\psi\\
			&\qquad-\lambda^{-1}f'(\sigma)\bigl((K(\theta_\sigma)-K(\sigma))V+R_\sigma\bigr)\psi-\sigma^2\frac{D\lambda}{\lambda}w_W.
		\end{aligned}
	\end{equation}
	We estimate the first term by the product formula $\Lambda_\sigma(au)=a\Lambda_\sigma u-2\nabla a\cdot\nabla u+(A_0a)u$. The bounds for the coefficients of $\zeta_W-f'(\sigma)\psi$ and their first two derivatives imply that $W\mapsto\Lambda_\sigma(\zeta_W-f'(\sigma)\psi)$ has order at most two, with every coefficient having multiplier norm at $\overline\varrho$ bounded by $C\nu(\sigma)$. In particular,
	\[\|\Lambda_\sigma(\zeta_W-f'(\sigma)\psi)\|_{\mathcal B_r}\leq C\nu(\sigma)\sum_{j=0}^2\|\nabla^jW\|_{\cX_r}.\]
	The term containing $\lambda^{-1}-1$ has coefficient bounds $C\sigma^{-1}\nu(\sigma)$. The term containing $(K(\theta_\sigma)-K(\sigma))V+R_\sigma$ has coefficient bounds $C\nu(\sigma)$. Finally, the coefficients of $w_W$ are uniformly bounded, so the term containing $\sigma^2D\lambda/\lambda$ has coefficient bounds $C\nu(\sigma)$ as well. Thus the scalar operator on the left of \eqref{Eq.state_R_principal_remainder} has order at most two and coefficient multiplier norms bounded by $C\nu(\sigma)$. We obtain
	\begin{align*}
		&\left\|\mathfrak R_\sigma(W,0)+f'(\sigma)\bigl(A_0+(K(\sigma)-1)V\bigr)\psi\right\|_{\mathcal B_r}
		\leq C\nu(\sigma)\sum_{j=0}^2\|\nabla^jW\|_{\cX_r},\\
		&\qquad\qquad\qquad\|\mathfrak R_\sigma(W,0)\|_{\mathcal B_r}\leq C\sum_{j=0}^2\|\nabla^jW\|_{\cX_r}.
	\end{align*}
	The second estimate follows by adding back the principal part, whose coefficient multiplier norms are uniformly bounded.
	
	We apply these bounds to \eqref{Eq.state_remainder_definition}. The first component has coefficient multiplier norms bounded by $C\nu(\sigma)$ because $\mathfrak M_{\overline\varrho}(\lambda-1)\leq C\sigma^{-1}\nu(\sigma)$ and the term containing $\mathfrak R_\sigma(W,0)$ has the additional factor $\sigma^{-1}$. For $(\mathscr E_\sigma W)_2$, the definitions of $a_0,a_1$ give the exact rearrangement
	\begin{align*}
		(\mathscr E_\sigma W)_2&=\frac{p-1}{2}f'(\theta_\sigma)\left[\mathfrak R_\sigma(W,0)+f'(\sigma)\bigl(A_0+(K(\sigma)-1)V\bigr)\psi\right]\\
		&\qquad\qquad+\frac{p-1}{2}f'(\sigma)\bigl(f'(\sigma)-f'(\theta_\sigma)\bigr)\bigl(A_0+(K(\sigma)-1)V\bigr)\psi.
	\end{align*}
	The operator inside the brackets has coefficient multiplier norms bounded by $C\nu(\sigma)$, by \eqref{Eq.state_R_principal_remainder}. The second term has the same bound because $\mathfrak M_{\overline\varrho}\bigl(f'(\sigma)-f'(\theta_\sigma)\bigr)\leq C\nu(\sigma)$. Hence every coefficient of $\mathscr E_\sigma$ has multiplier norm at $\overline\varrho$ bounded by $C\nu(\sigma)$. Lemma \ref{Lem.analytic_algebra} and \eqref{Eq.analytic_radius_recovery} yield
	\begin{align*}
		\|\mathscr E_\sigma W\|_{\cX_r}\leq C\nu(\sigma)\sum_{j=0}^2\|\nabla^jW\|_{\cX_r}\leq C\nu(\sigma)\left(1+\frac1{r'-r}+\frac2{(r'-r)^2}\right)\|W\|_{\cX_{r'}}.
	\end{align*}
	
	Finally, Lemma \ref{Lem.phase_coefficients} and the multiplier estimate in Lemma \ref{Lem.analytic_algebra} give
	\begin{align*}
		\|\mathbf g_\sigma h\|_{\cX_r}\leq\mathfrak M_r(\lambda^{-1})\left(\sigma^{-1}\mathfrak M_r(f(\theta_\sigma))+\frac{p-1}{2}\mathfrak M_r(f'(\theta_\sigma))\right)\|h\|_{\mathcal B_r}\leq C\|h\|_{\mathcal B_r}.
	\end{align*}
	The constants are independent of $\sigma\geq\sigma_c$ and locally uniform for $0<r<r'<\overline\varrho$. This completes the proof.
\end{proof}

\section{Construction of the periodic matrices}\label{Sec.matrix_proof}

Fix $p\in2\Z_{\geq1}+1$. We prove Proposition \ref{Prop.matrix} and Lemma \ref{Lem.sampled_matrices} by treating low, intermediate and high frequencies separately. Let $\theta_z:=3L/2$, and define
\[f_c(\theta):=f(\theta_z+\theta),\qquad c_f:=\sqrt{\frac2{p+1}},\qquad \mathsf P^c_\xi(\theta):=\begin{pmatrix}0&1\\-pf_c(\theta)^{p-1}-\xi^2&0\end{pmatrix}.\]
The oscillator equation \eqref{Eq.oscillator_intro}, the energy identity \eqref{Eq.oscillator_energy} and the antiperiodicity of $f$ give
\begin{equation}\label{Eq.centered_oscillator}
	\begin{gathered}
		f_c''+f_c^p=0,\qquad f_c(\theta+L)=-f_c(\theta),\\
		f_c(0)=f_c(L)=0,\qquad f_c'(0)=c_f,\qquad f_c'(L)=-c_f.
	\end{gathered}
\end{equation}
Moreover, $f_c(\theta)>0$ for $0<\theta<L$. We use $Y(\theta,\xi):=Y_p(\theta,\xi)$ and $\mathsf M(\xi):=\mathsf M_p(\xi)$ from \eqref{Eq.cell_growth_constant}. Thus
\begin{equation}\label{Eq.centered_fundamental_matrix}
	\partial_\theta Y(\theta,\xi)=\mathsf P^c_\xi(\theta)Y(\theta,\xi),\qquad Y(0,\xi)=I_2,\qquad \mathsf M(\xi)=Y(L,\xi).
\end{equation}
The frames and energy matrices constructed using $f_c$ will be translated back to the original phase.

\subsection{The low-frequency equation}\label{Subsec.low_matrix}

For small $\xi\geq0$, we seek an invertible matrix $P_{\xi^2}(\theta)$ such that, for every solution $h$ of
\begin{equation}\label{Eq.periodic_oscillator_model}
	h''+\bigl(\xi^2+pf(\theta)^{p-1}\bigr)h=0,\qquad \xi\geq0,
\end{equation}
the new unknowns $Z=(Z_1,Z_2)^T$ defined by $Z:=P_{\xi^2}(\theta)^{-1}(h, h')^T$ solves
\begin{equation}\label{Eq.low_frame_unknowns}
	Z_1'=Z_2,\qquad Z_2'=-\omega(\xi^2)Z_1,
\end{equation}
a system whose coefficients are independent of $\theta$. We will obtain $\omega(\xi^2)=c_a^2\xi^2+\Order(\xi^4)$ and then use the matrix
\[S_{\rm lo}(\theta,\xi,\delta):=P_{\xi^2}(\theta)^{-T}\diag(\omega(\xi^2)+\delta^2,1)P_{\xi^2}(\theta)^{-1}\]
to define the energy of \eqref{Eq.periodic_oscillator_model}. The auxiliary parameter $\delta>0$ is introduced to ensure the positive-definiteness of the energy at frequency $\xi=0$. Lemma \ref{Lem.zero_monodromy} computes $\mathsf M(0)$ and $\left.\partial_\eta \mathsf M_{21}(\sqrt\eta)\right|_{\eta=0}$. Lemma \ref{Lem.low_frame} uses \eqref{Eq.zero_cell_matrix} to construct $P_{\xi^2}$ and determine the coefficient $c_a^2$. Lemma \ref{Lem.low_metric} then proves the required bounds for $S_{\rm lo}$.

\begin{lemma}\label{Lem.zero_monodromy}
	For $\xi\geq0$, one has $\det \mathsf M(\xi)=1$ and $\mathsf M_{11}(\xi)=\mathsf M_{22}(\xi)$. The matrix depends analytically on $\eta:=\xi^2$. At $\xi=0$,
	\begin{equation}\label{Eq.zero_cell_matrix}
		\mathsf M(0)=-
		\begin{pmatrix}
			1&L(p-1)/(p+1)\\
			0&1
		\end{pmatrix}
		,\qquad \left.\partial_\eta \mathsf M_{21}(\sqrt\eta)\right|_{\eta=0}=L\frac{p+1}{p+3}.
	\end{equation}
	Moreover,
	\[b_{{\rm sc},p}(\xi)=\frac{2L}{(p+1)(p+3)}\xi+\Order(\xi^3),\qquad 0<\chi_p<\infty.\]
\end{lemma}

\begin{proof}
	By \eqref{Eq.centered_fundamental_matrix}, the coefficient matrix $\mathsf P^c_\xi$ has trace zero and $Y(0,\xi)=I_2$. Liouville's formula gives $\partial_\theta\det Y(\theta,\xi)=\operatorname{tr}(\mathsf P^c_\xi(\theta))\det Y(\theta,\xi)=0$. Hence $\det Y(\theta,\xi)=1$ and $\det \mathsf M(\xi)=1$.
	
	The functions $\theta\mapsto f_c(L-\theta)$ and $\theta\mapsto f_c(\theta)$ solve the oscillator equation in \eqref{Eq.centered_oscillator} and have the same initial values $(0,c_f)$ at $\theta=0$. Uniqueness therefore gives $f_c(L-\theta)=f_c(\theta)$. Define $J_c:=\diag(1,-1)$. Then $\mathsf P^c_\xi(L-\theta)=\mathsf P^c_\xi(\theta)$ and $J_c\mathsf P^c_\xi(\theta)J_c=-\mathsf P^c_\xi(\theta)$. Using \eqref{Eq.centered_fundamental_matrix}, we obtain
	\begin{align*}
		\partial_\theta Y(L-\theta,\xi)&=-\mathsf P^c_\xi(\theta)Y(L-\theta,\xi),\\
		\partial_\theta\bigl(J_cY(\theta,\xi)J_c\mathsf M(\xi)\bigr)&=-\mathsf P^c_\xi(\theta)J_cY(\theta,\xi)J_c\mathsf M(\xi).
	\end{align*}
	Both matrix functions have value $\mathsf M(\xi)$ at $\theta=0$. Uniqueness yields
	\begin{equation}\label{Eq.monodromy_reflection}
		Y(L-\theta,\xi)=J_cY(\theta,\xi)J_c\mathsf M(\xi),\qquad 0\leq\theta\leq L.
	\end{equation}
	Taking $\theta=L$ in \eqref{Eq.monodromy_reflection} gives $\mathsf M(\xi)^{-1}=J_c\mathsf M(\xi)J_c$. Since $\det \mathsf M(\xi)=1$,
	\[\begin{pmatrix}\mathsf M_{22}(\xi)&-\mathsf M_{12}(\xi)\\-\mathsf M_{21}(\xi)&\mathsf M_{11}(\xi)\end{pmatrix}=\mathsf M(\xi)^{-1}=J_c\mathsf M(\xi)J_c=\begin{pmatrix}\mathsf M_{11}(\xi)&-\mathsf M_{12}(\xi)\\-\mathsf M_{21}(\xi)&\mathsf M_{22}(\xi)\end{pmatrix}.\]
	Comparing the diagonal entries proves $\mathsf M_{11}(\xi)=\mathsf M_{22}(\xi)$.
	
	We next consider dependence on $\eta=\xi^2$. For complex $\eta$, the notation $Y(\theta,\sqrt\eta)$ denotes the solution of \eqref{Eq.centered_fundamental_matrix} with $\xi^2$ replaced by $\eta$. Its integral equation is
	\[Y(\theta,\sqrt\eta)=I_2+\int_0^\theta\begin{pmatrix}0&1\\-pf_c(s)^{p-1}-\eta&0\end{pmatrix}Y(s,\sqrt\eta)\,\mathrm ds.\]
	For every $R>0$, the coefficient matrix is bounded by a constant $C_R$ on $[0,L]\times\{|\eta|\leq R\}$. The uniformly convergent Picard series on $[0,L]$ implies that $Y(\theta,\sqrt\eta)$ and $\mathsf M(\sqrt\eta)=Y(L,\sqrt\eta)$ depend analytically on $\eta$. For real $\eta$, the coefficient matrix and the initial value are real, so the matrices have real-valued entries. The determinant and reflection calculations also apply with $\xi^2$ replaced by real $\eta$.
	
	At $\xi=0$, define
	\begin{equation}\label{Eq.zero_basis}
		u(\theta):=c_f^{-1}f_c'(\theta),\qquad v_c(\theta):=c_ff_c(\theta)+\frac{p-1}{p+1}\theta u(\theta).
	\end{equation}
	The functions $u,v_c$ solve $h''+pf_c^{p-1}h=0$, as follows by differentiating the oscillator equation \eqref{Eq.centered_oscillator}. The endpoint values in \eqref{Eq.centered_oscillator} now give
	\begin{align*}
		u(0)&=1,\qquad u'(0)=0,\qquad v_c(0)=0,\qquad v_c'(0)=c_f^2+\frac{p-1}{p+1}=1,\\
		u(L)&=-1,\qquad u'(L)=0,\qquad v_c(L)=-\frac{L(p-1)}{p+1},\qquad v_c'(L)=-c_f^2-\frac{p-1}{p+1}=-1.
	\end{align*}
	Consequently, uniqueness in \eqref{Eq.centered_fundamental_matrix} identifies
	\[Y(\theta,0)=\begin{pmatrix}u(\theta)&v_c(\theta)\\u'(\theta)&v_c'(\theta)\end{pmatrix}.\]
	Evaluating at $\theta=L$ proves the formula for $\mathsf M(0)$ in \eqref{Eq.zero_cell_matrix}.
	
	To compute the derivative of $\mathsf M_{21}(\sqrt\eta)$, let $u_\eta$ solve
	\begin{equation}\label{Eq.monodromy_parameter_solution}
		u_\eta''+(pf_c^{p-1}+\eta)u_\eta=0,\qquad u_\eta(0)=1,\qquad u_\eta'(0)=0.
	\end{equation}
	The relation with the matrix in \eqref{Eq.centered_fundamental_matrix} is
	\[\binom{u_\eta(\theta)}{u_\eta'(\theta)}=Y(\theta,\sqrt\eta)\binom10,\qquad \mathsf M_{21}(\sqrt\eta)=u_\eta'(L).\]
	At $\eta=0$, uniqueness gives $u_0=u$, where $u$ is defined in \eqref{Eq.zero_basis}. Analytic dependence on $\eta$ allows us to define $w(\theta):=\left.\partial_\eta u_\eta(\theta)\right|_{\eta=0}$ and differentiate \eqref{Eq.monodromy_parameter_solution}. We obtain
	\begin{equation}\label{Eq.monodromy_parameter_derivative}
		w''+pf_c^{p-1}w=-u,\qquad w(0)=w'(0)=0.
	\end{equation}
	The equations for $u$ and $w$ give $(uw'-u'w)'=uw''-u''w=-u^2$.
	Integrating over $[0,L]$, and using $u(L)=-1$, $u'(L)=0$ and the initial values in \eqref{Eq.monodromy_parameter_derivative}, yields
	\[-w'(L)=[uw'-u'w]_0^L=-\int_0^Lu(\theta)^2\,\mathrm d\theta.\]
	Since $(f_c')^2$ is $L$-periodic, \eqref{Eq.zero_basis} and \eqref{Eq.oscillator_averages} imply
	\begin{align*}
		\left.\partial_\eta \mathsf M_{21}(\sqrt\eta)\right|_{\eta=0}
		&=w'(L)=\int_0^Lu(\theta)^2\,\mathrm d\theta=\frac{p+1}{2}\int_0^L(f_c'(\theta))^2\,\mathrm d\theta\\
		&=\frac{L(p+1)}2\av{(f')^2}=\frac{L(p+1)}{p+3}.
	\end{align*}
	This proves the derivative formula in \eqref{Eq.zero_cell_matrix}.
	
	We now expand $b_{{\rm sc},p}$. Analyticity in $\eta=\xi^2$ and \eqref{Eq.zero_cell_matrix} give
	\begin{align*}
		\mathsf M_{12}(\xi)=-\frac{L(p-1)}{p+1}+\Order(\xi^2),\quad 	\mathsf M_{21}(\xi)=\frac{L(p+1)}{p+3}\xi^2+\Order(\xi^4)
	\end{align*}
	as $\xi\downarrow0$. Substituting both expansions into the definition \eqref{Eq.cell_growth_constant}, we obtain
	\begin{equation}\label{Eq.zero_sc_expansion}
		\begin{aligned}
			b_{{\rm sc},p}(\xi)&=\frac12\left(\xi \mathsf M_{12}(\xi)+\xi^{-1}\mathsf M_{21}(\xi)\right)=\frac L2\left(-\frac{p-1}{p+1}+\frac{p+1}{p+3}\right)\xi+\Order(\xi^3)\\
			&=\frac{2L}{(p+1)(p+3)}\xi+\Order(\xi^3).
		\end{aligned}
	\end{equation}
	Since the coefficient of $\xi$ is positive and $\arsinh t=t+\Order(t^3)$ as $t\to0$,
	\[\lim_{\xi\downarrow0}\frac{\arsinh|b_{{\rm sc},p}(\xi)|}{L\xi}=\frac2{(p+1)(p+3)}>0.\]
	
	It remains to bound the quotient in \eqref{Eq.cell_growth_constant} at positive frequencies. For $\xi>0$, define
	\[D_\xi:=\diag(\xi^2,1),\qquad \widehat{\mathsf M}(\xi):=D_\xi^{1/2}\mathsf M(\xi)D_\xi^{-1/2}=\begin{pmatrix}\mathsf M_{11}(\xi)&\xi \mathsf M_{12}(\xi)\\\xi^{-1}\mathsf M_{21}(\xi)&\mathsf M_{22}(\xi)\end{pmatrix}.\]
	Using $\det \mathsf M=1$, we compute
	\begin{align*}
		\operatorname{tr}(\widehat{\mathsf M}^T\widehat{\mathsf M})-2
		&=\mathsf M_{11}^2+\mathsf M_{22}^2+\xi^2\mathsf M_{12}^2+\xi^{-2}\mathsf M_{21}^2-2(\mathsf M_{11}\mathsf M_{22}-\mathsf M_{12}\mathsf M_{21})\\
		&=(\mathsf M_{11}-\mathsf M_{22})^2+(\xi \mathsf M_{12}+\xi^{-1}\mathsf M_{21})^2=4|b_{{\rm sc},p}(\xi)|^2,
	\end{align*}
	where all entries of $\mathsf M$ are evaluated at $\xi$, and the last equality uses $\mathsf M_{11}=\mathsf M_{22}$ and \eqref{Eq.cell_growth_constant}. Since $\det\widehat{\mathsf M}=1$, its singular values have the form $\ee^{a_\xi},\ee^{-a_\xi}$ with $a_\xi\geq0$. Hence
	\[4\sinh^2a_\xi=\ee^{2a_\xi}+\ee^{-2a_\xi}-2=4|b_{{\rm sc},p}(\xi)|^2.\]
	We have therefore proved
	\begin{equation}\label{Eq.sc_matrix_norm}
		a_\xi=\arsinh|b_{{\rm sc},p}(\xi)|,\qquad \|\widehat{\mathsf M}(\xi)\|=\ee^{a_\xi}.
	\end{equation}
	
	Let $h$ solve $h''+(\xi^2+pf_c^{p-1})h=0$	with arbitrary real initial values. Direct differentiation gives $\frac{\mathrm d}{\mathrm d\theta}\bigl(\xi^2h^2+(h')^2\bigr)=-2pf_c^{p-1}hh'$. Since $2\xi|hh'|\leq\xi^2h^2+(h')^2$, it follows that
	\[\left|\frac{\mathrm d}{\mathrm d\theta}\bigl(\xi^2h^2+(h')^2\bigr)\right|\leq\frac{pf_c^{p-1}}\xi\bigl(\xi^2h^2+(h')^2\bigr).\]
	Gronwall's inequality on $[0,L]$ yields
	\begin{equation}\label{Eq.incoming-outgoing_energies}
		\xi^2h(L)^2+h'(L)^2\leq\exp\left(\frac1\xi\int_0^Lpf_c(\theta)^{p-1}\,\mathrm d\theta\right)\bigl(\xi^2h(0)^2+h'(0)^2\bigr).
	\end{equation}
	For $v=(h(0),h'(0))^T$, \eqref{Eq.centered_fundamental_matrix} gives $\binom{h(L)}{h'(L)}=\mathsf M(\xi)v$. Since $D_\xi^{1/2}$ is invertible for $\xi>0$, the change of variables $z=D_\xi^{1/2}v$ gives
	\begin{align*}
		\sup_{(h(0),h'(0))\ne(0,0)}\frac{\xi^2h(L)^2+h'(L)^2}{\xi^2h(0)^2+h'(0)^2}=\sup_{v\in\R^2\setminus\{0\}}\frac{|D_\xi^{1/2}\mathsf M(\xi)v|^2}{|D_\xi^{1/2}v|^2}=\sup_{z\in\R^2\setminus\{0\}}\frac{|\widehat{\mathsf M}(\xi)z|^2}{|z|^2}=\|\widehat{\mathsf M}(\xi)\|^2.
	\end{align*}
	The energy estimate on $[0,L]$ and \eqref{Eq.sc_matrix_norm} therefore give
	\[0\leq a_\xi\leq\frac1{2\xi}\int_0^Lpf_c(\theta)^{p-1}\,\mathrm d\theta\leq\frac{pL}{2\xi},\]
	where the last inequality follows from $|f_c|\leq1$, by \eqref{Eq.oscillator_energy}. In particular,
	\[\frac{\arsinh|b_{{\rm sc},p}(\xi)|}{L\xi}=\frac{a_\xi}{L\xi}\leq\frac{p}{2\xi^2}\longrightarrow0\qquad\text{as }\xi\to\infty.\]
	The quotient is continuous on $(0,\infty)$, has a finite positive limit at zero by \eqref{Eq.zero_sc_expansion}, and tends to zero at infinity. Its supremum is therefore finite and positive. By \eqref{Eq.cell_growth_constant}, this proves $0<\chi_p<\infty$.
\end{proof}

To construct an antiperiodic frame, we seek a matrix $C(\eta)$ independent of $\theta$ such that $\ee^{LC(\eta)}=-\mathsf M(\sqrt\eta)$, so that $Y(\theta,\sqrt\eta)\ee^{-\theta C(\eta)}$ is antiperiodic in $\theta$. Lemma \ref{Lem.low_frame} constructs the normalized frame and uses the derivative in \eqref{Eq.zero_cell_matrix} to obtain $\omega(\eta)=c_a^2\eta+\Order(\eta^2)$.

\begin{lemma}[The analytic low-frequency frame]\label{Lem.low_frame}
	There is $\eta_0>0$ and there are real-valued analytic functions $a,b_0,\omega$ on $[0,\eta_0]$ such that
	\[a(0)=\frac{p-1}{p+1},\qquad b_0(0)=0,\qquad b_0'(0)=\frac{p+1}{p+3},\qquad \omega:=ab_0.\]
	For $0\leq\eta\leq\eta_0$, define
	\[C(\eta):=\begin{pmatrix}0&a(\eta)\\ -b_0(\eta)&0\end{pmatrix},\qquad P^c_\eta(\theta):=\frac{\sqrt{2(p+1)}}{p-1}Y(\theta,\sqrt\eta)\ee^{-\theta C(\eta)}\diag(a(\eta),1),\]
	one has
	\[(P^c_\eta)'=\mathsf P^c_{\sqrt\eta}P^c_\eta-P^c_\eta\begin{pmatrix}0&1\\ -\omega(\eta)&0\end{pmatrix},\qquad P^c_\eta(\theta+L)=-P^c_\eta(\theta),\qquad P^c_0(\theta)=P_0(\theta+\theta_z).\]
	The frame is invertible. On $0\leq\eta\leq\eta_0$, one has $\tfrac12c_a^2\eta\leq\omega(\eta)\leq2c_a^2\eta$ and
	\[\sup_{\theta\in\R}\bigl(\|P^c_\eta-P^c_0\|+\|(P^c_\eta)'-(P^c_0)'\|\bigr)\leq C\eta.\]
\end{lemma}

\begin{proof}
	By \eqref{Eq.zero_cell_matrix}, the spectrum of $-\mathsf M(0)$ is $\{1\}$. Lemma \ref{Lem.zero_monodromy} proves analytic dependence of $\mathsf M(\sqrt\eta)$ on the parameter $\eta$ at zero. Choose $\eta_1>0$ such that
	\[\spec(-\mathsf M(\sqrt\eta))\subset\{z\in\C:|z-1|<1/2\},\qquad 0\leq\eta\leq\eta_1.\]
	For $0\leq\eta\leq\eta_1$, use the branch of $\log z$ on $|z-1|<1$ with $\log1=0$, and define
	\begin{equation}\label{Eq.low_frame_logarithm}
		C(\eta):=\frac1L\log[-\mathsf M(\sqrt\eta)]=\frac1{2\pi\mathrm iL}\int_{|z-1|=3/4}\log z\,(zI_2+\mathsf M(\sqrt\eta))^{-1}\,\mathrm dz.
	\end{equation}
	The contour is positively oriented. The analytic dependence proved in Lemma \ref{Lem.zero_monodromy} and the fixed contour in \eqref{Eq.low_frame_logarithm} imply that $C(\eta)$ is analytic, including at $\eta=0$, and $\ee^{LC(\eta)}=-\mathsf M(\sqrt\eta)$. For real $\eta$, the matrix $\mathsf M(\sqrt\eta)$ is real and the logarithm satisfies $\overline{\log z}=\log\overline z$ on its domain. Therefore $C(\eta)$ has real-valued entries.
	
	We determine the form of $C(\eta)$. Since $\mathsf M_{11}(\sqrt\eta)=\mathsf M_{22}(\sqrt\eta)$, the inverse
	\[(zI_2+\mathsf M(\sqrt\eta))^{-1}=\frac1{\det(zI_2+\mathsf M(\sqrt\eta))}\begin{pmatrix}z+\mathsf M_{22}(\sqrt\eta)&-\mathsf M_{12}(\sqrt\eta)\\-\mathsf M_{21}(\sqrt\eta)&z+\mathsf M_{11}(\sqrt\eta)\end{pmatrix}\]
	has equal diagonal entries. Integrating in \eqref{Eq.low_frame_logarithm} gives $C_{11}(\eta)=C_{22}(\eta)$. On the other hand,
	\[\exp\bigl(L\operatorname{tr}C(\eta)\bigr)=\det\bigl(\ee^{LC(\eta)}\bigr)=\det(-\mathsf M(\sqrt\eta))=1.\]
	Since $\operatorname{tr}C(\eta)$ is real for real $\eta$, we obtain $\operatorname{tr}C(\eta)=0$. Thus
	\[C(\eta)=\begin{pmatrix}0&a(\eta)\\-b_0(\eta)&0\end{pmatrix},\qquad a(\eta):=C_{12}(\eta),\qquad b_0(\eta):=-C_{21}(\eta).\]
	
	At $\eta=0$, \eqref{Eq.zero_cell_matrix} gives
	\[-\mathsf M(0)=I_2+\begin{pmatrix}0&L(p-1)/(p+1)\\0&0\end{pmatrix}.\]
	The matrix added to $I_2$ has square zero, so its logarithm equals that matrix. Consequently,
	\[C(0)=\begin{pmatrix}0&(p-1)/(p+1)\\0&0\end{pmatrix},\qquad a(0)=\frac{p-1}{p+1},\qquad b_0(0)=0.\]
	To compute $b_0'(0)$, we differentiate $\mathsf M(\sqrt\eta)=-\ee^{LC(\eta)}$ to obtain
	\[\left.\partial_\eta \mathsf M(\sqrt\eta)\right|_{\eta=0}=-\int_0^L\ee^{(L-s)C(0)}C'(0)\ee^{sC(0)}\,\mathrm ds.\]
	Since $\ee^{sC(0)}=\begin{pmatrix}1&s(p-1)/(p+1)\\0&1\end{pmatrix}$, the $(2,1)$ entry of the integrand is
	\[\left(\ee^{(L-s)C(0)}C'(0)\ee^{sC(0)}\right)_{21}=C_{21}'(0)=-b_0'(0).\]
	It follows from \eqref{Eq.zero_cell_matrix} that
	\[Lb_0'(0)=\left.\partial_\eta \mathsf M_{21}(\sqrt\eta)\right|_{\eta=0}=\frac{L(p+1)}{p+3}.\]
	With $\omega(\eta):=a(\eta)b_0(\eta)$, we therefore have
	\begin{equation}\label{Eq.low_dispersion_expansion}
		\omega(0)=0,\qquad \omega'(0)=a(0)b_0'(0)=\frac{p-1}{p+3}=c_a^2,\qquad \omega(\eta)=c_a^2\eta+\Order(\eta^2).
	\end{equation}
	Here $c_a$ is defined in \eqref{Eq.kappa_H}.
	
	Choose $0<\eta_0<\eta_1$ such that
	\[\sup_{0\leq\eta\leq\eta_0}|a(\eta)-a(0)|\leq\frac{a(0)}2,\qquad \sup_{0\leq\eta\leq\eta_0}|\omega'(\eta)-c_a^2|\leq\frac{c_a^2}2.\]
	Then $a(\eta)\geq a(0)/2>0$ on $[0,\eta_0]$. Since $\omega(0)=0$, integration of the bound for $\omega'$ gives
	\[\frac12c_a^2\eta\leq\omega(\eta)\leq\frac32c_a^2\eta\leq2c_a^2\eta,\qquad 0\leq\eta\leq\eta_0.\]
	
	Define the frame by
	\begin{equation}\label{Eq.low_frame_construction}
		P_\eta^c(\theta):=\frac{\sqrt{2(p+1)}}{p-1}Y(\theta,\sqrt\eta)\ee^{-\theta C(\eta)}\diag(a(\eta),1).
	\end{equation}
	Every factor in \eqref{Eq.low_frame_construction} is invertible: $\det Y=1$, a matrix exponential is invertible, and $a(\eta)>0$. The frame has real-valued entries for real $\eta$ and is analytic in $\eta$.
	
	The identity $\omega=ab_0$ gives
	\[C(\eta)\diag(a(\eta),1)=\diag(a(\eta),1)\begin{pmatrix}0&1\\-\omega(\eta)&0\end{pmatrix}.\]
	Differentiating \eqref{Eq.low_frame_construction} and using \eqref{Eq.centered_fundamental_matrix}, we obtain
	\begin{equation}\label{Eq.low_frame_equation}
		\partial_\theta P_\eta^c=\mathsf P^c_{\sqrt\eta}P_\eta^c-P_\eta^c\begin{pmatrix}0&1\\-\omega(\eta)&0\end{pmatrix}.
	\end{equation}
	
	The coefficient matrix $\mathsf P^c_{\sqrt\eta}$ is $L$-periodic by \eqref{Eq.centered_oscillator}. Uniqueness in \eqref{Eq.centered_fundamental_matrix} therefore gives
	\[Y(\theta+L,\sqrt\eta)=Y(\theta,\sqrt\eta)\mathsf M(\sqrt\eta).\]
	Combining this identity with $\mathsf M(\sqrt\eta)=-\ee^{LC(\eta)}$ and \eqref{Eq.low_frame_construction}, we find
	\begin{align*}
		P_\eta^c(\theta+L)&=\frac{\sqrt{2(p+1)}}{p-1}Y(\theta,\sqrt\eta)\mathsf M(\sqrt\eta)\ee^{-(\theta+L)C(\eta)}\diag(a(\eta),1)\\
		&=-\frac{\sqrt{2(p+1)}}{p-1}Y(\theta,\sqrt\eta)\ee^{-\theta C(\eta)}\diag(a(\eta),1)=-P_\eta^c(\theta).
	\end{align*}
	
	At $\eta=0$, by \eqref{Eq.Pzero}, \eqref{Eq.centered_oscillator} and \eqref{Eq.low_frame_construction}, we have
	\[P_0^c(0)=\diag\left(c_f,\frac{p+1}{p-1}c_f\right)=P_0(\theta_z).\]
	By \eqref{Eq.low_frame_equation} and $\omega(0)=0$, the matrix $P_0^c$ satisfies $\partial_\theta P_0^c=\mathsf P_0^cP_0^c-P_0^cN_0$. The translated frame $\theta\mapsto P_0(\theta+\theta_z)$ satisfies the same equation by \eqref{Eq.frame_identity}, with the same value at $\theta=0$. Uniqueness gives $P_0^c(\theta)=P_0(\theta+\theta_z)$ for all $\theta\in\R$.
	
	Finally, analyticity in $\eta$ and smoothness in $\theta$ in \eqref{Eq.low_frame_construction} imply
	\begin{align*}
		&\sup_{\theta\in[0,L]}\left(\|P_\eta^c(\theta)-P_0^c(\theta)\|+\|\partial_\theta P_\eta^c(\theta)-\partial_\theta P_0^c(\theta)\|\right)\\
		&\qquad\qquad\leq\eta\sup_{\substack{\theta\in[0,L]\\0\leq s\leq\eta_0}}\left(\|\partial_sP_s^c(\theta)\|+\|\partial_\theta\partial_sP_s^c(\theta)\|\right)\leq C\eta.
	\end{align*}
	The antiperiodicity of $P_\eta^c$ and $\partial_\theta P_\eta^c$ extends this bound to all $\theta\in\R$. 
\end{proof}

We return to the original phase by defining $P_\eta(\theta):=P^c_\eta(\theta-\theta_z)$. For a real solution $U=(h,h')^T$ of \eqref{Eq.periodic_oscillator_model}, let $Z=P_{\xi^2}^{-1}U$, then $Z$ solves \eqref{Eq.low_frame_unknowns}. For fixed $\xi$, there holds
\begin{equation}\label{Eq.low_oscillator_energy}
	\frac{\mathrm d}{\mathrm d\theta}\left(\omega(\xi^2)Z_1^2+Z_2^2\right)=2\omega(\xi^2)Z_1Z_2-2\omega(\xi^2)Z_2Z_1=0.
\end{equation}
By \eqref{Eq.low_dispersion_expansion}, the coefficient $\omega(\xi^2)$ is positive for small $\xi>0$, but $\omega(0)=0$. At $\xi=0$, the expression in \eqref{Eq.low_oscillator_energy} is $Z_2^2$ and does not control $Z_1$.

For $\delta>0$, we instead use $(\omega(\xi^2)+\delta^2)Z_1^2+Z_2^2$. Substituting $Z=P_{\xi^2}^{-1}U$ gives
\begin{equation}\label{Eq.low_regularized_energy}
	(\omega(\xi^2)+\delta^2)Z_1^2+Z_2^2
	=U^TP_{\xi^2}^{-T}\diag(\omega(\xi^2)+\delta^2,1)P_{\xi^2}^{-1}U=U^TS_{\rm lo}U.
\end{equation}
For fixed $\xi,\delta$, it follows from \eqref{Eq.low_frame_unknowns} that
\begin{align*}
	\frac{\mathrm d}{\mathrm d\theta}(U^TS_{\rm lo}U)=2\delta^2Z_1Z_2\leq\delta\left((\omega(\xi^2)+\delta^2)Z_1^2+Z_2^2\right)=\delta U^TS_{\rm lo}U.
\end{align*}
Thus the additional term $\delta^2Z_1^2$ makes the energy positive definite at $\xi=0$, at the cost of an energy growth rate bounded by $\delta$. Lemma \ref{Lem.low_metric} proves the parameter derivative bounds and the strict boundary inequality for the matrix $S_{\rm lo}$ in \eqref{Eq.low_regularized_energy}.

\begin{lemma}\label{Lem.low_metric}
	There exist $a_*,\delta_*,\epsilon_b,c,C>0$, with $4a_*^2\leq\eta_0$, such that the symmetric matrix
	\[S_{\rm lo}:=P_{\xi^2}(\theta)^{-T}\diag(\omega(\xi^2)+\delta^2,1)P_{\xi^2}(\theta)^{-1}\]
	is smooth for $0\leq\xi\leq2a_*$, $0\leq\delta\leq\delta_*$, is positive for $(\xi,\delta)\ne(0,0)$, and satisfies
	\begin{align*}
		&\qquad\qquad\qquad\qquad\qquad\quad\partial_\theta S_{\rm lo}+\mathsf P_\xi^TS_{\rm lo}+S_{\rm lo}\mathsf P_\xi\leq\delta S_{\rm lo},\\
		&\qquad\qquad\qquad\quad\|S_{\rm lo}^{-1/2}(\xi\partial_\xi S_{\rm lo})S_{\rm lo}^{-1/2}\| +\|S_{\rm lo}^{-1/2}(\delta\partial_\delta S_{\rm lo})S_{\rm lo}^{-1/2}\|\leq C,\\
		&c\diag(\xi^2+\delta^2,\delta^2)\leq S_{\rm lo} \leq C\diag(1+\xi^2,1),\qquad (S_{\rm lo})_{22}\leq C,\quad \frac{\xi(S_{\rm lo})_{22}}{\sqrt{\det S_{\rm lo}}}\leq1-\epsilon_b
	\end{align*}
	for $0\leq\xi\leq2a_*$ and $0<\delta\leq\delta_*$.
\end{lemma}

\begin{proof}
	Choose $a_*,\delta_*>0$ with $4a_*^2\leq\eta_0$, to be reduced below. For $0\leq\xi\leq2a_*$ and $0<\delta\leq\delta_*$, define
	\begin{equation}\label{Eq.low_metric_coordinates}
		\kappa:=\bigl(\omega(\xi^2)+\delta^2\bigr)^{1/2},\qquad D_{\xi,\delta}:=\diag(\kappa^2,1),\qquad S_{\rm lo}=P_{\xi^2}^{-T}D_{\xi,\delta}P_{\xi^2}^{-1}.
	\end{equation}
	Lemma \ref{Lem.low_frame} gives $\kappa^2\asymp\xi^2+\delta^2$. We first choose $a_*,\delta_*$ small enough that $\kappa\leq1$. The matrices $P_{\xi^2}(\theta)$ and $P_{\xi^2}(\theta)^{-1}$ are uniformly bounded on this parameter range, by Lemma \ref{Lem.low_frame} and antiperiodicity. Thus, for $z\in\R^2$,
	\begin{align*}
		z^TS_{\rm lo}z=\kappa^2|(P_{\xi^2}^{-1}z)_1|^2+|(P_{\xi^2}^{-1}z)_2|^2\geq\kappa^2|P_{\xi^2}^{-1}z|^2\geq c(\xi^2+\delta^2)|z|^2,
	\end{align*}
	while $z^TS_{\rm lo}z\leq C|z|^2$. Consequently,
	\[c\diag(\xi^2+\delta^2,\delta^2)\leq S_{\rm lo}\leq CI_2\leq C\diag(1+\xi^2,1),\qquad (S_{\rm lo})_{22}\leq C.\]
	The definition \eqref{Eq.low_metric_coordinates} also gives positivity whenever $(\xi,\delta)\ne(0,0)$, including $\delta=0$ and $\xi>0$.
	
	We next compute the energy derivative. Since $P_\eta(\theta)=P_\eta^c(\theta-\theta_z)$, \eqref{Eq.low_frame_equation} gives
	\[\partial_\theta P_{\xi^2}=\mathsf P_\xi P_{\xi^2}-P_{\xi^2}\begin{pmatrix}0&1\\-\omega(\xi^2)&0\end{pmatrix}.\]
	The matrix $D_{\xi,\delta}$ is independent of $\theta$. Differentiating \eqref{Eq.low_metric_coordinates} therefore yields
	\begin{equation}\label{Eq.low_metric_defect}
		\begin{aligned}
			&P_{\xi^2}^T\left(\partial_\theta S_{\rm lo}+\mathsf P_\xi^TS_{\rm lo}+S_{\rm lo}\mathsf P_\xi\right)P_{\xi^2}\\
			&\qquad\qquad\qquad=\begin{pmatrix}0&-\omega(\xi^2)\\1&0\end{pmatrix}D_{\xi,\delta}+D_{\xi,\delta}\begin{pmatrix}0&1\\-\omega(\xi^2)&0\end{pmatrix}=\begin{pmatrix}0&\delta^2\\\delta^2&0\end{pmatrix}.
		\end{aligned}
	\end{equation}
	After multiplication on both sides by $D_{\xi,\delta}^{-1/2}$, the last matrix becomes
	\[D_{\xi,\delta}^{-1/2}\begin{pmatrix}0&\delta^2\\\delta^2&0\end{pmatrix}D_{\xi,\delta}^{-1/2}=\begin{pmatrix}0&\delta^2/\kappa\\\delta^2/\kappa&0\end{pmatrix}.\]
	Its Euclidean operator norm is $\delta^2/\kappa\leq\delta$, because $\kappa\geq\delta$. Hence \eqref{Eq.low_metric_defect} implies
	\[\partial_\theta S_{\rm lo}+\mathsf P_\xi^TS_{\rm lo}+S_{\rm lo}\mathsf P_\xi\leq\delta S_{\rm lo}.\]
	
	We now prove the strict boundary inequality. It follows from \eqref{Eq.low_metric_coordinates} that
	\begin{equation}\label{Eq.low_metric_entries}
		(S_{\rm lo})_{22}=\frac{(P_{\xi^2})_{11}^2+\kappa^2(P_{\xi^2})_{12}^2}{(\det P_{\xi^2})^2},\qquad \sqrt{\det S_{\rm lo}}=\frac{\kappa}{|\det P_{\xi^2}|}.
	\end{equation}
	By Lemma \ref{Lem.low_frame} and \eqref{Eq.Pzero}, uniformly in $\theta$,
	\begin{align*}
		(P_{\xi^2})_{11}=f'(\theta)+\Order(\xi^2),\quad		(P_{\xi^2})_{12}=\frac2{p-1}f(\theta)+\Order(\xi^2),\quad\det P_{\xi^2}=\frac2{p-1}+\Order(\xi^2).
	\end{align*}
	Thus
	\[\frac{(P_{\xi^2})_{11}^2}{|\det P_{\xi^2}|}=\frac{p-1}{2}(f'(\theta))^2+\Order(\xi^2),\qquad \frac{(P_{\xi^2})_{12}^2}{|\det P_{\xi^2}|}\leq C.\]
	For $\xi>0$, the expansion \eqref{Eq.low_dispersion_expansion} implies
	\[\frac{\xi}{\kappa}\leq\frac{\xi}{\sqrt{\omega(\xi^2)}}\leq\frac1{c_a}(1+C\xi^2),\qquad \xi\kappa\leq C(\xi^2+\delta^2).\]
	Substitution into \eqref{Eq.low_metric_entries} gives
	\begin{equation}\label{Eq.low_metric_boundary_ratio}
		\begin{aligned}
			\frac{\xi(S_{\rm lo})_{22}}{\sqrt{\det S_{\rm lo}}}=\frac{\xi}{\kappa}\frac{(P_{\xi^2})_{11}^2}{|\det P_{\xi^2}|}+\xi\kappa\frac{(P_{\xi^2})_{12}^2}{|\det P_{\xi^2}|}\leq\frac{p-1}{2c_a}(f'(\theta))^2+C_b(\xi^2+\delta^2).
		\end{aligned}
	\end{equation}
	For $\xi=0$ and $\delta>0$, the left-hand side vanishes, so the inequality also holds at $\xi=0$.
	
	The energy identity \eqref{Eq.oscillator_energy} gives $\max_\theta(f'(\theta))^2=2/(p+1)$. Since $c_a^2=(p-1)/(p+3)$, define
	\[r_p:=\frac{p-1}{(p+1)c_a}=\frac{\sqrt{(p-1)(p+3)}}{p+1}=\sqrt{1-\frac4{(p+1)^2}}<1.\]
	By adjusting $a_*,\delta_*$ to be smaller if necessary, we may assume without loss of generality that
	\[C_b(4a_*^2+\delta_*^2)\leq\frac{1-r_p}{2},\qquad \epsilon_b:=\frac{1-r_p}{2}>0.\]
	Then \eqref{Eq.low_metric_boundary_ratio} yields
	\[\frac{\xi(S_{\rm lo})_{22}}{\sqrt{\det S_{\rm lo}}}\leq r_p+\frac{1-r_p}{2}=1-\epsilon_b.\]
	
	It remains to estimate the relative parameter derivatives. Define $F_\xi:=P_{\xi^2}^{-1}\xi\partial_\xi P_{\xi^2}$.	The frame $P$ is analytic in $\eta=\xi^2$ by Lemma \ref{Lem.low_frame}. Since $\xi\partial_\xi P_{\xi^2}=2\xi^2\left.(\partial_\eta P_\eta)\right|_{\eta=\xi^2}$, uniform bounds for $\partial_\eta P_\eta$ and $P_\eta^{-1}$ give $\|F_\xi\|\leq C\xi^2$. Differentiating \eqref{Eq.low_metric_coordinates} at fixed $\theta,\delta$, we obtain
	\begin{equation}\label{Eq.partial_xiS}
		P_{\xi^2}^T(\xi\partial_\xi S_{\rm lo})P_{\xi^2}=\xi\partial_\xi D_{\xi,\delta}-F_\xi^TD_{\xi,\delta}-D_{\xi,\delta}F_\xi.
	\end{equation}
	
	For a positive definite matrix $S$ and a real symmetric matrix $H$, we have
	\[\|S^{-1/2}HS^{-1/2}\|=\sup_{v\in\R^2\setminus\{0\}}\frac{|v^TS^{-1/2}HS^{-1/2}v|}{|v|^2}=\sup_{z\in\R^2\setminus\{0\}}\frac{|z^THz|}{z^TSz}.\]
	By \eqref{Eq.low_metric_coordinates}, $P_{\xi^2}^TS_{\rm lo}P_{\xi^2}=D_{\xi,\delta}$. Since $P_{\xi^2}$ is invertible, changing variables by $z=P_{\xi^2}y$ yields
	\begin{equation}\label{Eq.relative_quardratic-form}
		\begin{aligned}
			\|S_{\rm lo}^{-1/2}HS_{\rm lo}^{-1/2}\|
			&=\sup_{z\in\R^2\setminus\{0\}}\frac{|z^THz|}{z^TS_{\rm lo}z}=\sup_{y\in\R^2\setminus\{0\}}\frac{|y^TP_{\xi^2}^THP_{\xi^2}y|}{y^TD_{\xi,\delta}y}\\
			&=\|D_{\xi,\delta}^{-1/2}P_{\xi^2}^THP_{\xi^2}D_{\xi,\delta}^{-1/2}\|.
		\end{aligned}
	\end{equation}
	Then it follows from \eqref{Eq.partial_xiS} and \eqref{Eq.relative_quardratic-form} that 
	\begin{align*}
		\|S_{\rm lo}^{-1/2}(\xi\partial_\xi S_{\rm lo})S_{\rm lo}^{-1/2}\|\leq\|D_{\xi,\delta}^{-1/2}(\xi\partial_\xi D_{\xi,\delta})D_{\xi,\delta}^{-1/2}\|+2\|D_{\xi,\delta}^{1/2}F_\xi D_{\xi,\delta}^{-1/2}\|.
	\end{align*}
	The first term is bounded because
	\[D_{\xi,\delta}^{-1/2}(\xi\partial_\xi D_{\xi,\delta})D_{\xi,\delta}^{-1/2}=\diag\left(\frac{2\xi^2\omega'(\xi^2)}{\kappa^2},0\right),\qquad \frac{2\xi^2|\omega'(\xi^2)|}{\kappa^2}\leq C.\]
	For the second term, $\kappa\leq1$, $\kappa\geq c\xi$ and $\|F_\xi\|\leq C\xi^2$ give
	\[\|D_{\xi,\delta}^{1/2}F_\xi D_{\xi,\delta}^{-1/2}\|\leq\frac{C\xi^2}{\kappa}\leq C\xi\leq C.\]
	This proves the relative $\xi$-derivative bound.
	
	The frame $P_{\xi^2}$ is independent of $\delta$, so $P_{\xi^2}^T(\delta\partial_\delta S_{\rm lo})P_{\xi^2}=\delta\partial_\delta D_{\xi,\delta}=\diag(2\delta^2,0)$.
	Using \eqref{Eq.relative_quardratic-form} implies that
	\[\|S_{\rm lo}^{-1/2}(\delta\partial_\delta S_{\rm lo})S_{\rm lo}^{-1/2}\|=\frac{2\delta^2}{\kappa^2}\leq2.\]
	This completes the proof.
\end{proof}

\subsection{Frequencies bounded away from zero}\label{Subsec.middle_matrix}

Lemma \ref{Lem.low_metric} constructs $S_{\rm lo}$ near $\xi=0$. For frequencies in a compact interval $K\Subset(0,\infty)$, we use the solution matrix $Y$ in \eqref{Eq.centered_fundamental_matrix} to define
\begin{equation}\label{Eq.middle_endpoint_matrices}
	S_{\rm in}(\theta,\xi):=Y(\theta,\xi)^{-T}D_\xi Y(\theta,\xi)^{-1},\quad S_{\rm out}(\theta,\xi):=Y(\theta,\xi)^{-T}\mathsf M(\xi)^TD_\xi \mathsf M(\xi)Y(\theta,\xi)^{-1},
\end{equation}
where $D_\xi:=\diag(\xi^2,1)$. For $U(\theta)=Y(\theta,\xi)U(0)$, the definitions give
\[U(\theta)^TS_{\rm in}(\theta,\xi)U(\theta)=U(0)^TD_\xi U(0),\qquad U(\theta)^TS_{\rm out}(\theta,\xi)U(\theta)=U(L)^TD_\xi U(L).\]
We join $S_{\rm in}$ near $\theta=0$ to $S_{\rm out}$ near $\theta=L$ to obtain an $L$-periodic matrix. The boundary inequality $\xi S_{22}\leq\sqrt{\det S}$ requires control of the $(2,2)$ entries during the interpolation from $S_{\rm in}$ to $S_{\rm out}$. Lemma \ref{Lem.half_cell} proves the necessary bounds, and Lemma \ref{Lem.middle_metric} gives the interpolation. At sufficiently large frequencies, Lemma \ref{Lem.high_metric} uses $S_{\rm hi}=\diag(1+\xi^2,1)$ instead.

\begin{lemma}\label{Lem.half_cell}
	Fix a compact interval $K\Subset(0,\infty)$. The matrices in \eqref{Eq.middle_endpoint_matrices} are smooth and positive definite on $[0,L]\times K$, and, for $S_\bullet\in\{S_{\rm in},S_{\rm out}\}$,
	\begin{equation}\label{Eq.endpoint_matrix_identity}
		\partial_\theta S_\bullet+(\mathsf P_\xi^c)^TS_\bullet+S_\bullet\mathsf P_\xi^c=0,\qquad \det S_\bullet=\xi^2.
	\end{equation}
	For every $\xi\in K$,
	\begin{equation}\label{Eq.half_period_bounds}
		(S_{\rm in})_{22}(\theta,\xi)\leq1\quad(0\leq\theta\leq L/2),\qquad (S_{\rm out})_{22}(\theta,\xi)\leq1\quad(L/2\leq\theta\leq L).
	\end{equation}
	There are $\eta_K\in(0,L/4)$ and $\epsilon_K>0$ such that
	\begin{equation}\label{Eq.half_period_overlap}
		\max\{(S_{\rm in})_{22}(\theta,\xi),(S_{\rm out})_{22}(\theta,\xi)\}\leq1-\epsilon_K\qquad(|\theta-L/2|\leq\eta_K,\ \xi\in K).
	\end{equation}
\end{lemma}

\begin{proof}
	Smoothness and positivity follow from \eqref{Eq.middle_endpoint_matrices} and invertibility of $Y$ and $\mathsf M$. By \eqref{Eq.centered_fundamental_matrix}, $\partial_\theta(Y^{-1})=-Y^{-1}\mathsf P_\xi^c$. Thus, for either $C=D_\xi$ or $C=\mathsf M(\xi)^TD_\xi \mathsf M(\xi)$,
	\[\partial_\theta(Y^{-T}CY^{-1})=-(\mathsf P_\xi^c)^TY^{-T}CY^{-1}-Y^{-T}CY^{-1}\mathsf P_\xi^c.\]
	Since $\det Y=\det \mathsf M=1$ by Lemma \ref{Lem.zero_monodromy}, we compute that $\det S_{\rm in}=\xi^2=\det S_{\rm out}$. This proves \eqref{Eq.endpoint_matrix_identity}.
	
	To estimate $(S_{\rm in})_{22}$, fix $\xi\in K$ and solve
	\begin{equation}\label{Eq.half_period_solution}
		\psi''+(\xi^2+pf_c^{p-1})\psi=0,\qquad \psi(0)=1,\qquad \psi'(0)=\mathrm i\xi.
	\end{equation}
	Writing $Y=\left(\begin{smallmatrix}u&v\\u'&v'\end{smallmatrix}\right)$, equations \eqref{Eq.centered_fundamental_matrix} and \eqref{Eq.half_period_solution} give $\psi=u+\mathrm i\xi v$. The inverse formula for $Y$, together with $\det Y=1$, yields
	\[(S_{\rm in})_{22}=\left|D_\xi^{1/2}Y^{-1}\binom01\right|^2=u^2+\xi^2v^2=|\psi|^2.\]
	Moreover, \eqref{Eq.half_period_solution} gives $\partial_\theta\operatorname{Im}(\overline\psi\psi')=0$, so $\operatorname{Im}(\overline\psi\psi')=\xi$ and $|\psi'|^2\geq\frac{\xi^2}{|\psi|^2}$. In particular, $\psi$ never vanishes.
	
	Define
	\[\Omega_\xi^2:=\xi^2+pf_c^{p-1},\qquad E:=|\psi'|^2+\Omega_\xi^2|\psi|^2.\]
	The bound for $|\psi'|$ implies
	\[\Omega_\xi^2|\psi|^4-E|\psi|^2+\xi^2\leq0.\]
	Hence $E^2-4\xi^2\Omega_\xi^2\geq0$, and $|\psi|^2$ does not exceed the larger root
	\begin{equation}\label{Eq.half_period_envelope}
		\mathcal A:=\frac{E+\sqrt{E^2-4\xi^2\Omega_\xi^2}}{2\Omega_\xi^2},\qquad |\psi|^2\leq\mathcal A,\qquad E=\Omega_\xi^2\mathcal A+\frac{\xi^2}{\mathcal A}.
	\end{equation}
	The initial values in \eqref{Eq.half_period_solution} give $\mathcal A(0)=1$. By \eqref{Eq.centered_oscillator} and \eqref{Eq.oscillator_energy}, $f_c$ increases from zero to one on $[0,L/2]$, so $(\Omega_\xi^2)'=(pf_c^{p-1})'>0$ on $(0,L/2)$. Differentiating $E$ by \eqref{Eq.half_period_solution} gives $E'=(\Omega_\xi^2)'|\psi|^2$. On the open set $\{\mathcal A>1\}$, the larger root in \eqref{Eq.half_period_envelope} is simple and
	\[\left(\Omega_\xi^2-\frac{\xi^2}{\mathcal A^2}\right)\mathcal A'=(\Omega_\xi^2)'(|\psi|^2-\mathcal A)\leq0.\]
	Since $\Omega_\xi^2\geq\xi^2>0$, the coefficient of $\mathcal A'$ is positive when $\mathcal A>1$. Thus $\mathcal A'\leq0$ on $\{\mathcal A>1\}$. Every connected component of $\{\mathcal A>1\}$ in $[0,L/2]$ would start with value $\mathcal A=1$ by continuity, contradicting $\mathcal A'\leq0$. Consequently $|\psi|^2\leq\mathcal A\leq1$ on $[0,L/2]$, which proves the bound for $S_{\rm in}$ in \eqref{Eq.half_period_bounds}.
	
	We need strict inequality at $\theta=L/2$ to allow interpolation on an interval. Since $E(0)=2\xi^2$ and $\Omega_\xi(0)^2=\xi^2$,
	\[E(L/2)-\Omega_\xi(L/2)^2-\xi^2=\int_0^{L/2}(pf_c^{p-1})'(|\psi|^2-1)\,\mathrm d\theta\leq0.\]
	Equality would force $|\psi|=1$ on $(0,L/2)$. Then $\operatorname{Re}(\overline\psi\psi')=0$ and $\operatorname{Im}(\overline\psi\psi')=\xi$ imply $\psi'=\mathrm i\xi\psi$ and $\psi''=-\xi^2\psi$, contradicting \eqref{Eq.half_period_solution} because $f_c>0$ there. Thus $E(L/2)<\Omega_\xi(L/2)^2+\xi^2$. If $|\psi(L/2)|=1$, the inequality $|\psi'|^2\geq\xi^2/|\psi|^2$ would give the opposite bound. Therefore $(S_{\rm in})_{22}(L/2,\xi)=|\psi(L/2)|^2<1$.
	
	For $J_c=\diag(1,-1)$, the reflection identity \eqref{Eq.monodromy_reflection} and $J_cD_\xi J_c=D_\xi$ give
	\[J_cS_{\rm in}(L-\theta,\xi)J_c=Y(\theta,\xi)^{-T}J_c\mathsf M(\xi)^{-T}D_\xi \mathsf M(\xi)^{-1}J_cY(\theta,\xi)^{-1}=S_{\rm out}(\theta,\xi),\]
	where we used $\mathsf M^{-1}=J_c\mathsf MJ_c$, proved in Lemma \ref{Lem.zero_monodromy}. This proves the bound for $S_{\rm out}$ in \eqref{Eq.half_period_bounds} and gives $(S_{\rm out})_{22}(L/2,\xi)=(S_{\rm in})_{22}(L/2,\xi)<1$. By continuity and compactness,
	\[\epsilon_K^0:=\min_{\xi\in K}\bigl(1-(S_{\rm in})_{22}(L/2,\xi)\bigr)>0.\]
	Uniform continuity of both $(2,2)$ entries on $[0,L]\times K$ gives $\eta_K\in(0,L/4)$ such that their values at $\theta$ and $L/2$ differ by at most $\epsilon_K^0/2$ when $|\theta-L/2|\leq\eta_K$. Taking $\epsilon_K=\epsilon_K^0/2$ proves \eqref{Eq.half_period_overlap}.
\end{proof}

The endpoint values in \eqref{Eq.middle_endpoint_matrices} satisfy $S_{\rm in}(0,\xi)=S_{\rm out}(L,\xi)=D_\xi$. We interpolate between $S_{\rm in}$ and $S_{\rm out}$ only where \eqref{Eq.half_period_overlap} holds. The interpolation must preserve the determinant $\xi^2$ and the bound $S_{22}\leq1$, while the change from $D_\xi$ to $\mathsf M^TD_\xi \mathsf M$ determines the energy growth that must be allowed.

\begin{lemma}\label{Lem.middle_metric}
	For a compact interval $K\Subset(0,\infty)$ and $\chi>\chi_p$, there is a smooth, symmetric, positive definite, $L$-periodic matrix $S_{\rm mid}(\theta,\xi)$ with real-valued entries on a neighborhood of $\R\times K$ such that
	\begin{equation}\label{Eq.middle_matrix_bounds}
		\partial_\theta S_{\rm mid}+(\mathsf P^c_\xi)^TS_{\rm mid}+S_{\rm mid}\mathsf P^c_\xi\leq2\widehat\chi\xi S_{\rm mid},\qquad \xi(S_{\rm mid})_{22}\leq\sqrt{\det S_{\rm mid}},
	\end{equation}
	where $\widehat\chi:=(\chi+\chi_p)/2$. The inverse and every fixed-order derivative of $S_{\rm mid}$ are bounded on $[0,L]\times K$.
\end{lemma}

\begin{proof}
	Choose a compact interval $K_1\Subset(0,\infty)$ with $K\subset\operatorname{int}K_1$. Apply Lemma \ref{Lem.half_cell} on $K_1$, and choose $h\in C^\infty([0,L];[0,1])$ with
	\[h'\geq0,\qquad h=0\text{ on }[0,L/2-\eta_{K_1}],\qquad h=1\text{ on }[L/2+\eta_{K_1},L].\]
	Let $\widehat{\mathsf M}=D_\xi^{1/2}\mathsf M D_\xi^{-1/2}$ and define $C_\xi:=\widehat{\mathsf M}^T\widehat{\mathsf M}$. By \eqref{Eq.sc_matrix_norm}, $C_\xi$ has eigenvalues $\ee^{2a_\xi}$ and $\ee^{-2a_\xi}$, where $a_\xi=\arsinh|b_{{\rm sc},p}(\xi)|$. In particular, $\det C_\xi=1$ and $\log C_\xi\leq2a_\xi I_2$. Define
	\begin{equation}\label{Eq.middle_interpolation}
		S_0(\theta,\xi):=Y^{-T}D_\xi^{1/2}C_\xi^{h(\theta)}D_\xi^{1/2}Y^{-1},\qquad C_\xi^t:=\exp(t\log C_\xi).
	\end{equation}
	The positive definite matrix $C_\xi$ has a smooth logarithm, including when its eigenvalues coincide. Thus $S_0$ is smooth and positive definite. At $h=0$, formula \eqref{Eq.middle_interpolation} gives $S_0=S_{\rm in}$; at $h=1$, the identity $D_\xi^{1/2}C_\xi D_\xi^{1/2}=\mathsf M^TD_\xi \mathsf M$ gives $S_0=S_{\rm out}$.
	
	For $c>0$ and $0\leq t\leq1$, one has $c^t\leq1-t+tc$. Applying $c^t\leq1-t+tc$ to the eigenvalues of $C_\xi$ and then using \eqref{Eq.middle_interpolation}, we obtain
	\[S_0\leq(1-h)S_{\rm in}+hS_{\rm out},\qquad \det S_0=(\det Y)^{-2}\det D_\xi\,(\det C_\xi)^h=\xi^2.\]
	On $|\theta-L/2|\leq\eta_{K_1}$, inequality \eqref{Eq.half_period_overlap} bounds both $(S_{\rm in})_{22}$ and $(S_{\rm out})_{22}$ by $1-\epsilon_{K_1}$. On $[0,L/2-\eta_{K_1}]$, we have $S_0=S_{\rm in}$; on $[L/2+\eta_{K_1},L]$, we have $S_0=S_{\rm out}$. Applying \eqref{Eq.half_period_bounds} on these two intervals proves $(S_0)_{22}\leq1$ on $[0,L]\times K_1$, and hence $\xi(S_0)_{22}\leq\sqrt{\det S_0}$.
	
	Using \eqref{Eq.centered_fundamental_matrix} to differentiate \eqref{Eq.middle_interpolation}, we find
	\begin{align*}
		\partial_\theta S_0+(\mathsf P^c_\xi)^TS_0+S_0\mathsf P^c_\xi=h'Y^{-T}D_\xi^{1/2}C_\xi^{h/2}(\log C_\xi)C_\xi^{h/2}D_\xi^{1/2}Y^{-1}\leq2a_\xi h'S_0.
	\end{align*}
	To replace the coefficient $2a_\xi h'$ by a bound independent of the width of the interpolation interval, define
	\begin{equation}\label{Eq.middle_matrix_definition}
		S_{\rm mid}(\theta,\xi):=\ee^{2\widehat\chi\xi(\theta-Lh(\theta))}S_0(\theta,\xi).
	\end{equation}
	By \eqref{Eq.cell_growth_constant}, $a_\xi\leq L\xi\chi_p<L\xi\widehat\chi$. Differentiating the scalar factor in \eqref{Eq.middle_matrix_definition} gives
	\begin{align*}
		\partial_\theta S_{\rm mid}+(\mathsf P^c_\xi)^TS_{\rm mid}+S_{\rm mid}\mathsf P^c_\xi\leq\bigl[2\widehat\chi\xi+2(a_\xi-L\widehat\chi\xi)h'\bigr]S_{\rm mid}\leq2\widehat\chi\xi S_{\rm mid}.
	\end{align*}
	Multiplication by a positive scalar leaves $\xi S_{22}/\sqrt{\det S}$ unchanged in dimension two. Therefore
	\[\frac{\xi(S_{\rm mid})_{22}}{\sqrt{\det S_{\rm mid}}}=\frac{\xi(S_0)_{22}}{\sqrt{\det S_0}}\leq1,\]
	which proves \eqref{Eq.middle_matrix_bounds}.
	
	We verify smooth periodicity before extending $S_{\rm mid}$ beyond $[0,L]$. Periodicity of $\mathsf P^c_\xi$ and \eqref{Eq.centered_fundamental_matrix} give $Y(\theta+L,\xi)=Y(\theta,\xi)\mathsf M(\xi)$, and hence $S_{\rm out}(\theta+L,\xi)=S_{\rm in}(\theta,\xi)$. The definition \eqref{Eq.middle_matrix_definition} equals $\ee^{2\widehat\chi\xi\theta}S_{\rm in}(\theta,\xi)$ near zero and $\ee^{2\widehat\chi\xi(\theta-L)}S_{\rm out}(\theta,\xi)$ near $L$. The identity $S_{\rm out}(\theta+L,\xi)=S_{\rm in}(\theta,\xi)$ shows that the endpoint formulas agree after translation by $L$, including all derivatives in $\theta$ and $\xi$. Thus the $L$-periodic extension is smooth. Positivity and compactness of $[0,L]\times K$ give the uniform inverse and derivative bounds.
\end{proof}

For large $\xi$, the fixed quadratic form $(1+\xi^2)h^2+(h')^2$ suffices. Along \eqref{Eq.periodic_oscillator_model}, its derivative is $2(1-pf^{p-1})hh'$, whose absolute value relative to $(1+\xi^2)h^2+(h')^2$ is at most $(p-1)/\sqrt{1+\xi^2}$. We now record the corresponding matrix bounds.

\begin{lemma}\label{Lem.high_metric}
	Let $\chi>0$, and define $R_{\rm hi}:=\max\{1,\sqrt{(p-1)/\chi}\}$ and $S_{\rm hi}(\xi):=\diag(1+\xi^2,1)$. For $\xi\geq R_{\rm hi}$ and $0<\delta\leq1$,
	\begin{align*}
		&\mathsf P_\xi^TS_{\rm hi}+S_{\rm hi}\mathsf P_\xi\leq\chi\xi S_{\rm hi},\quad\|S_{\rm hi}^{-1/2}(\xi\partial_\xi S_{\rm hi})S_{\rm hi}^{-1/2}\|\leq2,\quad \partial_\delta S_{\rm hi}=0,\\
		&\diag(\xi^2+\delta^2,\delta^2)\leq S_{\rm hi},\qquad (S_{\rm hi})_{22}=1,\qquad \frac{\xi(S_{\rm hi})_{22}}{\sqrt{\det S_{\rm hi}}}=\frac\xi{\sqrt{1+\xi^2}}<1.
	\end{align*}
\end{lemma}

\begin{proof}
	The energy identity \eqref{Eq.oscillator_energy} gives $|f|\leq1$, and $p$ is odd, so $0\leq pf^{p-1}\leq p$. Direct multiplication yields
	\[S_{\rm hi}^{-1/2}(\mathsf P_\xi^TS_{\rm hi}+S_{\rm hi}\mathsf P_\xi)S_{\rm hi}^{-1/2}=\begin{pmatrix}0&(1-pf^{p-1})/\sqrt{1+\xi^2}\\(1-pf^{p-1})/\sqrt{1+\xi^2}&0\end{pmatrix}.\]
	Its operator norm satisfies
	\[\frac{|1-pf^{p-1}|}{\sqrt{1+\xi^2}}\leq\frac{p-1}{\xi}\leq\chi\xi,\qquad \xi\geq R_{\rm hi}.\]
	Multiplying the resulting quadratic-form inequality by $S_{\rm hi}^{1/2}$ proves the first bound. The parameter derivatives are
	\[S_{\rm hi}^{-1/2}(\xi\partial_\xi S_{\rm hi})S_{\rm hi}^{-1/2}=\diag\left(\frac{2\xi^2}{1+\xi^2},0\right),\qquad \partial_\delta S_{\rm hi}=0.\]
	Finally, $\delta\leq1$ implies $\diag(\xi^2+\delta^2,\delta^2)\leq\diag(1+\xi^2,1)$, and $\det S_{\rm hi}=1+\xi^2$ gives the boundary ratio in the statement.
\end{proof}

\subsection{Gluing and comparison of the derivative levels}
\label{Subsec.matrix_gluing}

We combine $S_{\rm lo}$, $S_{\rm mid}$ and $S_{\rm hi}$ with nonnegative weights depending only on $\xi$. Lemma \ref{Lem.global_metric} proves that the resulting matrix satisfies both the energy inequality and $\xi S_{22}\leq\sqrt{\det S}$ in \eqref{Eq.S_coercivity}. In \eqref{Eq.energy_phase_samples}, the energy \eqref{Eq.energy_weights} uses $S$ at $\xi_N=c_*N/(\varrho\sigma)$ with $c_*=\lambda^{-1}$ and $\delta=\sigma^{-1}$. We then compare $S(\theta_*,\xi_N,\delta)$ at different values of $N$, keeping the coefficient $1+C/\sigma$ in the boundary estimate \eqref{Eq.sampled_boundary}.

\begin{lemma}\label{Lem.global_metric}
	Let $p\geq3$ be odd and $\chi>\chi_p$. There are $\delta_0\in(0,1]$, $0<\xi_{\rm lo}<\xi_{\rm hi}$ and a smooth symmetric matrix $S(\theta,\xi,\delta)$ with real-valued entries on $\R\times[0,\infty)\times[0,\delta_0]$ such that $S(\theta+L,\xi,\delta)=S(\theta,\xi,\delta)$ and $S>0$ for $(\xi,\delta)\ne(0,0)$. For $\xi\geq0$ and $0<\delta\leq\delta_0$,
	\begin{equation}\label{Eq.global_matrix_bounds}
		\begin{aligned}
			&\qquad\qquad\qquad\qquad\qquad\partial_\theta S+\mathsf P_\xi^TS+S\mathsf P_\xi\leq(2\chi\xi+C\delta)S,\\
			&\qquad\qquad\qquad\|S^{-1/2}\xi(\partial_\xi S)S^{-1/2}\|+\|S^{-1/2}\delta(\partial_\delta S)S^{-1/2}\|\leq C,\\
			&c\diag(\xi^2+\delta^2,\delta^2)\leq S\leq C\diag(1+\xi^2,1),\quad S_{22}\leq C,\quad \xi S_{22}\leq\sqrt{\det S}.
		\end{aligned}
	\end{equation}
	Moreover, $S=S_{\rm lo}$ for $0\leq\xi\leq\xi_{\rm lo}$, and $S=S_{\rm hi}$ for $\xi\geq\xi_{\rm hi}$.
\end{lemma}

\begin{proof}
	Fix the constants $a_*,\delta_*$ in Lemma \ref{Lem.low_metric} and define $R:=\max\{2a_*,R_{\rm hi}\}$, where $R_{\rm hi}$ is defined in Lemma \ref{Lem.high_metric}. Apply Lemma \ref{Lem.middle_metric} with $K=[a_*,2R]$, and choose smooth functions $\alpha_{\rm lo},\alpha_{\rm hi}:[0,\infty)\to[0,1]$ such that
	\begin{align*}
		\alpha_{\rm lo}=1\text{ on }[0,a_*],\quad \supp\alpha_{\rm lo}\subset[0,3a_*/2),\quad		\alpha_{\rm hi}=0\text{ on }[0,R],\quad \alpha_{\rm hi}=1\text{ on }[2R,\infty).
	\end{align*}
	Since $R\geq2a_*$, the supports of $\alpha_{\rm lo}$ and $\alpha_{\rm hi}$ are disjoint. Thus $\alpha_{\rm mid}:=1-\alpha_{\rm lo}-\alpha_{\rm hi}\geq0$. Set
	\begin{equation}\label{Eq.global_matrix_definition}
		\begin{aligned}
			S(\theta,\xi,\delta)&:=\alpha_{\rm lo}(\xi)S_{\rm lo}(\theta,\xi,\delta)+\alpha_{\rm mid}(\xi)S_{\rm mid}(\theta-\theta_z,\xi)+\alpha_{\rm hi}(\xi)S_{\rm hi}(\xi),\\
			\delta_0&:=\min\{1,\delta_*\},\qquad \xi_{\rm lo}:=a_*,\qquad \xi_{\rm hi}:=2R.
		\end{aligned}
	\end{equation}
	Each matrix is used only where it is defined; products with zero weights are extended by zero. The shift $\theta-\theta_z$ changes the coefficient $\mathsf P_\xi^c$ in \eqref{Eq.middle_matrix_bounds} to $\mathsf P_\xi$. Formula \eqref{Eq.global_matrix_definition} is therefore smooth and $L$-periodic, is positive definite for $(\xi,\delta)\ne(0,0)$, and equals the required matrices for $\xi\leq a_*$ and $\xi\geq2R$.
	
	For fixed $(\theta,\xi,\delta)$, write $S_j$ for the three matrices evaluated at the arguments in \eqref{Eq.global_matrix_definition}, with $j\in\{\mathrm{lo},\mathrm{mid},\mathrm{hi}\}$ and $\alpha_j>0$. Since $\partial_\theta\alpha_j=0$, Lemmas \ref{Lem.low_metric}, \ref{Lem.middle_metric} and \ref{Lem.high_metric} give
	\[\partial_\theta S+\mathsf P_\xi^TS+S\mathsf P_\xi=\sum_j\alpha_j(\partial_\theta S_j+\mathsf P_\xi^TS_j+S_j\mathsf P_\xi)\leq(2\chi\xi+C\delta)\sum_j\alpha_jS_j.\]
	Here $\widehat\chi<\chi$ and all $\alpha_j$ are nonnegative. This proves the energy inequality in \eqref{Eq.global_matrix_bounds}.
	
	The boundary inequality is also preserved, but requires a separate calculation. For $\delta>0$, write $S_j=\left(\begin{smallmatrix}a_j&h_j\\h_j&d_j\end{smallmatrix}\right)$, so $d_j>0$. The bound $\xi(S_j)_{22}\leq\sqrt{\det S_j}$ is equivalent to $a_j\geq h_j^2/d_j+\xi^2d_j$. Cauchy--Schwarz gives
	\begin{align*}
		S_{11}=\sum_j\alpha_ja_j&\geq\sum_j\alpha_j\frac{h_j^2}{d_j}+\xi^2\sum_j\alpha_jd_j\\
		&\geq\frac{(\sum_j\alpha_jh_j)^2}{\sum_j\alpha_jd_j}+\xi^2\sum_j\alpha_jd_j=\frac{S_{12}^2}{S_{22}}+\xi^2S_{22}.
	\end{align*}
	Multiplication by $S_{22}>0$ yields $\det S\geq\xi^2S_{22}^2$, as required.
	
	On the compact interval $a_*\leq\xi\leq2R$, each matrix with nonzero weight in \eqref{Eq.global_matrix_definition} and its first parameter derivatives are bounded uniformly for $\theta\in[0,L]$ and $0\leq\delta\leq\delta_0$. Positivity, including at $\delta=0$, gives $cI_2\leq S_j\leq CI_2$ on the relevant supports, and hence $S\geq cI_2$. Differentiating \eqref{Eq.global_matrix_definition}, we obtain
	\begin{align*}
		&\|S^{-1/2}\xi(\partial_\xi S)S^{-1/2}\|\leq c^{-1}\sum_j\bigl(|\xi\alpha_j'|\|S_j\|+\alpha_j\|\xi\partial_\xi S_j\|\bigr)\leq C,\\
		&\qquad\qquad\|S^{-1/2}\delta(\partial_\delta S)S^{-1/2}\|\leq c^{-1}\sum_j\alpha_j\|\delta\partial_\delta S_j\|\leq C.
	\end{align*}
	The uniform positive bounds on $[a_*,2R]$ also imply the matrix upper and lower bounds and $S_{22}\leq C$ in \eqref{Eq.global_matrix_bounds}. For $\xi\leq a_*$ and $\xi\geq2R$, the bounds in \eqref{Eq.global_matrix_bounds} follow directly from Lemmas \ref{Lem.low_metric} and \ref{Lem.high_metric}. This proves \eqref{Eq.global_matrix_bounds} on the full frequency range.
\end{proof}

\begin{proof}[Proof of Proposition \ref{Prop.matrix}]
	Lemma \ref{Lem.zero_monodromy} proves $0<\chi_p<\infty$. The matrix defined in \eqref{Eq.global_matrix_definition} satisfies \eqref{Eq.S_coercivity} by \eqref{Eq.global_matrix_bounds}. For $\xi\leq\xi_{\rm lo}$, its formula is \eqref{Eq.low_regularized_energy}, and Lemma \ref{Lem.low_frame} gives $P_\eta|_{\eta=0}=P_0$, the frame equation \eqref{Eq.low_frame_equation} translated to the original phase, and the uniform $\Order(\eta)$ bounds for $P_\eta-P_0$ and its phase derivative. The expansion \eqref{Eq.low_dispersion_expansion} gives $\omega(\eta)=c_a^2\eta+\Order(\eta^2)$. For $\xi\geq\xi_{\rm hi}$, formula \eqref{Eq.global_matrix_definition} gives $S=\diag(1+\xi^2,1)$. This proves Proposition \ref{Prop.matrix}.
\end{proof}

The spatial derivative estimates compare the quadratic forms defined by $S(\theta_*,\xi_N,\delta)$ and $S(\theta_*,\xi_J,\delta)$. The relative $\xi$-derivative bound in \eqref{Eq.global_matrix_bounds} controls the ratios of the quadratic forms when $N,J\geq1$. Comparisons with $N=0$ require the explicit formula \eqref{Eq.low_metric_coordinates}, since a bound for $\xi\partial_\xi S$ cannot be integrated through $\xi=0$.

\begin{lemma}\label{Lem.matrix_comparisons}
	Let $I_{\rm r}=[r_-,r_+]\Subset(0,\overline\varrho)$, $\varrho\in I_{\rm r}$, $\sigma\geq\max\{2,\delta_0^{-1}\}$, $c_*\in[1/2,2]$ and $\theta_*\in\R$. Define
	\begin{equation}\label{Eq.sampled_matrix_definition}
		\delta:=\sigma^{-1},\qquad \xi_N:=\frac{c_*N}{\varrho\sigma},\qquad S_N^\star:=S(\theta_*,\xi_N,\delta),\qquad \forall\ N\in\Z_{\geq0}.
	\end{equation}
	There are $C>0$ and $H\in\Z_{\geq1}$, depending only on $p,\chi,I_{\rm r}$ and $s_*$, such that
	\begin{equation}\label{Eq.matrix_level_comparisons}
		\begin{aligned}
			&\|(S_N^\star)^{1/2}(S_J^\star)^{-1/2}\|\leq1+\frac{C(1+N-J)^H}{J}\qquad(N\geq J\geq1),\\
			&\|(S_N^\star)^{1/2}(S_0^\star)^{-1/2}\|\leq C(1+N)^H\qquad(N\geq0),\\
			&\|(S_N^\star)^{1/2}(S_{N+1}^\star)^{-1/2}\|+\|(S_{N+1}^\star)^{1/2}(S_N^\star)^{-1/2}\|\leq C\qquad(N\geq0).
		\end{aligned}
	\end{equation}
	For the fixed exponent $s_*$,
	\begin{equation}\label{Eq.matrix_weighted_difference}
		\|(S_N^\star)^{-1/2}(S_N^\star-t_NS_{N-1}^\star)(S_N^\star)^{-1/2}\|\leq\frac C N,\qquad t_N:=\left(\frac N{N+1}\right)^{2s_*},\qquad N\geq1.
	\end{equation}
\end{lemma}

\begin{proof}
	Fix $\theta_*,\delta$ and write $S(\xi)=S(\theta_*,\xi,\delta)$ in this proof. By \eqref{Eq.global_matrix_bounds}, there is $C_0>0$, independent of $\theta_*,\delta$, such that
	\[-C_0S(\xi)\leq\xi\partial_\xi S(\xi)\leq C_0S(\xi),\qquad \xi>0.\]
	For each fixed $v\in\R^2\setminus\{0\}$, it follows that $|\partial_\xi\log(v^TS(\xi)v)|\leq C_0/\xi$. Integrating between $0<\xi_1\leq\xi_2$ and then allowing $v$ to vary gives
	\begin{equation}\label{Eq.matrix_frequency_comparison}
		\left(\frac{\xi_2}{\xi_1}\right)^{-C_0}S(\xi_1)\leq S(\xi_2)\leq\left(\frac{\xi_2}{\xi_1}\right)^{C_0}S(\xi_1).
	\end{equation}
	For $N\geq J\geq1$, the definition of the operator norm and \eqref{Eq.sampled_matrix_definition} therefore imply
	\begin{align*}
		\|(S_N^\star)^{1/2}(S_J^\star)^{-1/2}\|^2=\sup_{v\ne0}\frac{v^TS_N^\star v}{v^TS_J^\star v}\leq\left(\frac{\xi_N}{\xi_J}\right)^{C_0}=\left(\frac N J\right)^{C_0}.
	\end{align*}
	The mean value theorem gives
	\[\left(\frac N J\right)^{C_0/2}\leq1+C\frac{N-J}{J}\left(1+\frac{N-J}{J}\right)^{\max\{C_0/2-1,0\}}\leq1+\frac{C(1+N-J)^H}{J}\]
	for an integer $H>1+C_0/2$. This proves the first bound in \eqref{Eq.matrix_level_comparisons}. Taking consecutive positive indices in \eqref{Eq.matrix_frequency_comparison} also bounds both adjacent norm ratios, because $(N+1)/N\leq2$ for $N\geq1$.
	
	For $N\geq2$, formula \eqref{Eq.matrix_frequency_comparison} implies
	\[\left(\frac N{N-1}\right)^{-C_0}I_2\leq(S_N^\star)^{-1/2}S_{N-1}^\star(S_N^\star)^{-1/2}\leq\left(\frac N{N-1}\right)^{C_0}I_2.\]
	Both scalar bounds differ from one by at most $C/N$. Since $|1-t_N|\leq C/N$, we obtain
	\begin{align*}
		\|(S_N^\star)^{-1/2}(S_N^\star-t_NS_{N-1}^\star)(S_N^\star)^{-1/2}\|\leq|1-t_N|+t_N\|I_2-(S_N^\star)^{-1/2}S_{N-1}^\star(S_N^\star)^{-1/2}\|\leq C/N.
	\end{align*}
	This proves \eqref{Eq.matrix_weighted_difference} except at $N=1$.
	
	To compare with level zero, suppose first that $\xi_N\leq\xi_{\rm lo}$ and define $\kappa_N:=\sqrt{\omega(\xi_N^2)+\delta^2}$. The formula \eqref{Eq.low_metric_coordinates} gives
	\[v^TS_N^\star v=|\diag(\kappa_N,1)P_{\xi_N^2}(\theta_*)^{-1}v|^2,\qquad v^TS_0^\star v=|\diag(\delta,1)P_0(\theta_*)^{-1}v|^2.\]
	Taking the quotient $v^TS_N^\star v/(v^TS_0^\star v)$ and substituting $z=\diag(\delta,1)P_0(\theta_*)^{-1}v$ yields
	\[\|(S_N^\star)^{1/2}(S_0^\star)^{-1/2}\|=\|\diag(\kappa_N,1)P_{\xi_N^2}(\theta_*)^{-1}P_0(\theta_*)\diag(\delta^{-1},1)\|.\]
	Lemma \ref{Lem.low_frame} gives $P_{\xi_N^2}^{-1}P_0=I_2+\Order(\xi_N^2)$ uniformly in $\theta_*$. Since $\kappa_N\asymp\xi_N+\delta$, $\delta\leq1$ and $\xi_N/\delta=c_*N/\varrho$, we obtain
	\[\|(S_N^\star)^{1/2}(S_0^\star)^{-1/2}\|\leq C\left(1+\frac{\xi_N}{\delta}+\frac{\xi_N^2}{\delta}\right)\leq C(1+N)\qquad(\xi_N\leq\xi_{\rm lo}).\]
	If $\xi_N>\xi_{\rm lo}$, the bounds $S_N^\star\leq C(1+\xi_N^2)I_2$ and $S_0^\star\geq c\delta^2I_2$ in \eqref{Eq.global_matrix_bounds} instead give
	\[\|(S_N^\star)^{1/2}(S_0^\star)^{-1/2}\|\leq C\delta^{-1}(1+\xi_N)\leq C\left(\frac1{\xi_{\rm lo}}+1\right)\frac{c_*N}{\varrho}\leq CN.\]
	The second bound in \eqref{Eq.matrix_level_comparisons} follows in both cases.
	
	It remains to bound $\|(S_0^\star)^{1/2}(S_1^\star)^{-1/2}\|$. If $\xi_1\leq\xi_{\rm lo}$, use \eqref{Eq.low_metric_coordinates} and substitute 
\[z=\diag(\kappa_1,1)P_{\xi_1^2}(\theta_*)^{-1}v\]
 in $v^TS_0^\star v/(v^TS_1^\star v)$ to obtain
	\begin{align*}
		\|(S_0^\star)^{1/2}(S_1^\star)^{-1/2}\|=\|\diag(\delta,1)P_0(\theta_*)^{-1}P_{\xi_1^2}(\theta_*)\diag(\kappa_1^{-1},1)\|\leq C\left(1+\frac{\xi_1^2}{\kappa_1}\right)\leq C.
	\end{align*}
	If $\xi_1>\xi_{\rm lo}$, then $\delta=\varrho\xi_1/c_*>r_-\xi_{\rm lo}/2$. The bounds $S_0^\star\leq CI_2$ and $S_1^\star\geq c\delta^2I_2$ in \eqref{Eq.global_matrix_bounds} therefore imply $\|(S_0^\star)^{1/2}(S_1^\star)^{-1/2}\|\leq C/\delta\leq C$. This completes the adjacent comparison at $N=0$. Finally,
	\[\|(S_1^\star)^{-1/2}(S_1^\star-t_1S_0^\star)(S_1^\star)^{-1/2}\|\leq1+t_1\|(S_0^\star)^{1/2}(S_1^\star)^{-1/2}\|^2\leq C,\]
	which proves \eqref{Eq.matrix_weighted_difference} at $N=1$. All constants are uniform for $\varrho\in I_{\rm r}$.
\end{proof}

The boundary calculation in Lemma \ref{Lem.boundary_energy} requires
\[\xi_{N+1}\sqrt{\frac{(S_N^\star)_{22}(S_{N+1}^\star)_{22}}{\det S_{N+1}^\star}}\leq1+\frac C\sigma.\]
A uniform constant larger than one would not give the error of order $\sigma^{-2}$ in \eqref{Eq.total_boundary_estimate}. We prove the displayed bound using Lemma \ref{Lem.low_metric} near zero, smoothness on compact positive frequency intervals, and the formula $S_{\rm hi}=\diag(1+\xi^2,1)$ at large frequencies.

\begin{lemma}\label{Lem.cross_boundary}
	Fix $I_{\rm r}=[r_-,r_+]\Subset(0,\overline\varrho)$. There are $\Sigma_{\rm mat}\geq\max\{2,\delta_0^{-1}\}$ and $C>0$ such that, for $\varrho\in I_{\rm r}$, $\sigma\geq\Sigma_{\rm mat}$, $c_*\in[1/2,2]$ and $\theta_*\in\R$, the matrices and frequencies defined in \eqref{Eq.sampled_matrix_definition} satisfy
	\begin{equation}\label{Eq.cross_level_boundary_bound}
		\sup_{N\geq0}\xi_{N+1}\sqrt{\frac{(S_N^\star)_{22}(S_{N+1}^\star)_{22}}{\det S_{N+1}^\star}}\leq1+\frac C\sigma.
	\end{equation}
\end{lemma}

\begin{proof}
	For the parameters in \eqref{Eq.sampled_matrix_definition}, define
	\[\mathfrak R_N:=\xi_{N+1}\sqrt{\frac{(S_N^\star)_{22}(S_{N+1}^\star)_{22}}{\det S_{N+1}^\star}}.\]
	We first consider small $\xi_{N+1}$, so that both matrices equal $S_{\rm lo}$. For $j=N,N+1$, let $\kappa_j:=\sqrt{\omega(\xi_j^2)+\delta^2}$. By \eqref{Eq.low_metric_entries} and Lemma \ref{Lem.low_frame}, uniformly in $\theta_*$,
	\begin{align*}
		&\qquad\qquad\quad(S_j^\star)_{22}\leq\frac{(f'(\theta_*))^2+C(\xi_{N+1}^2+\delta^2)}{|\det P_{\xi_j^2}(\theta_*)|^2},\\
		&\sqrt{\det S_j^\star}=\frac{\kappa_j}{|\det P_{\xi_j^2}(\theta_*)|},\qquad |\det P_{\xi_j^2}(\theta_*)|=\frac2{p-1}+\Order(\xi_j^2).
	\end{align*}
	The expansion \eqref{Eq.low_dispersion_expansion} gives $\xi_{N+1}/\kappa_{N+1}\leq c_a^{-1}(1+C\xi_{N+1}^2)$. Substituting the two bounds for $(S_j^\star)_{22}$ into $\mathfrak R_N$, the determinant factor at level $N+1$ cancels, and we obtain
	\begin{align*}
		\mathfrak R_N&\leq\frac{\xi_{N+1}}{\kappa_{N+1}}\frac{(f'(\theta_*))^2+C(\xi_{N+1}^2+\delta^2)}{|\det P_{\xi_N^2}(\theta_*)|}\\
		&\leq\frac{p-1}{2c_a}(f'(\theta_*))^2+C(\xi_{N+1}^2+\delta^2)\leq r_p+C(\xi_{N+1}^2+\delta^2),
	\end{align*}
	where \eqref{Eq.oscillator_energy} gives
	\[r_p:=\frac{p-1}{(p+1)c_a}=\sqrt{1-\frac4{(p+1)^2}}<1.\]
	Choose $a_b\in(0,\xi_{\rm lo}/2]$ small and $\Sigma_{\rm mat}$ large enough that
	\[C(a_b^2+\Sigma_{\rm mat}^{-2})\leq\frac{1-r_p}{2}.\]
	Then $\mathfrak R_N\leq(1+r_p)/2<1$ whenever $\xi_{N+1}\leq a_b$.
	
	Increase $\Sigma_{\rm mat}$ so that the frequency increment satisfies
	\[\xi_{N+1}-\xi_N=\frac{c_*}{\varrho\sigma}\leq\frac2{r_-\Sigma_{\rm mat}}\leq\min\{a_b/2,1\}.\]
	If $\xi_{N+1}>a_b$ and $\xi_N<\xi_{\rm hi}$, both frequencies lie in $[a_b/2,\xi_{\rm hi}+1]$. Smoothness and positivity of the matrix in \eqref{Eq.global_matrix_definition} on this compact interval give
	\[\sup_{\substack{\theta\in[0,L],\ a_b/2\leq\xi\leq\xi_{\rm hi}+1\\0\leq\delta\leq\delta_0}}|\partial_\xi\log S_{22}(\theta,\xi,\delta)|\leq C.\]
	Consequently,
	\[\left|\log\frac{(S_N^\star)_{22}}{(S_{N+1}^\star)_{22}}\right|\leq\int_{\xi_N}^{\xi_{N+1}}|\partial_\xi\log S_{22}(\theta_*,\xi,\delta)|\,\mathrm d\xi\leq\frac C\sigma.\]
	Using the boundary inequality in \eqref{Eq.global_matrix_bounds} at $\xi_{N+1}$, we obtain
	\[\mathfrak R_N=\frac{\xi_{N+1}(S_{N+1}^\star)_{22}}{\sqrt{\det S_{N+1}^\star}}\sqrt{\frac{(S_N^\star)_{22}}{(S_{N+1}^\star)_{22}}}\leq\ee^{C/(2\sigma)}\leq1+\frac C\sigma.\]
	Finally, if $\xi_N\geq\xi_{\rm hi}$, both matrices equal $S_{\rm hi}$ and
	\[\mathfrak R_N=\frac{\xi_{N+1}}{\sqrt{1+\xi_{N+1}^2}}<1.\]
	The three frequency ranges cover every $N\geq0$, including pairs crossing $\xi_{\rm hi}$. This proves \eqref{Eq.cross_level_boundary_bound}.
\end{proof}

\begin{proof}[Proof of Lemma \ref{Lem.sampled_matrices}]
	For $\varrho\in I_{\rm r}$, $c_*\in[1/2,2]$ and $\theta_*\in\R$, the definition \eqref{Eq.sampled_matrix_definition} agrees with the matrices in Lemma \ref{Lem.sampled_matrices}, with $\xi_N=\xi_N^\star=c_*N/(\varrho\sigma)$. Lemma \ref{Lem.matrix_comparisons} gives the three comparisons in \eqref{Eq.matrix_level_comparisons} and the bound \eqref{Eq.metric_difference} through \eqref{Eq.matrix_weighted_difference}. Taking $J_{\rm cmp}=H$ and using \eqref{Eq.cross_level_boundary_bound} gives \eqref{Eq.sampled_boundary}. The constants and the starting time are uniform for $\varrho\in I_{\rm r}$ and depend only on $p,\chi,I_{\rm r}$ and $s_*$, not on the Euler profile.
\end{proof}

\section{The analytic energy estimate}
\label{Sec.energy_proof}

We prove the estimate \eqref{Eq.full_energy_differential} for the energy \eqref{Eq.energy_weights}. The matrix inequality \eqref{Eq.S_coercivity} controls the periodic potential; the outgoing condition \eqref{Eq.acoustic_outgoing_property} must dominate the resulting loss $\chi/\varrho$ and the transport loss $\mathcal L_v(\varrho)$. Fix $I_{\rm r}\Subset(0,\overline\varrho)$, $A_*>0$, $\tau=\log(\sigma/\Sigma)$ and $\varrho\in C^1([0,T_*];I_{\rm r})$ with $0\leq\alpha:=\varrho'/\varrho\leq A_*$. We first calculate on the class \eqref{Eq.energy_smooth_core}, retain the boundary flux in \eqref{Eq.total_boundary_estimate}, and then extend the estimate to finite-energy data.

Set
\[\mathbf g_0:=\na\log J_0,\qquad Q_\sigma:=(K(\theta_\sigma)-1)V+R_\sigma.\]
Then $A_0=-\Delta-\mathbf g_0\cdot\na$, and the linearized equation is exactly
\begin{equation}\label{Eq.energy_scalar_equation}
	D^2\zeta-\sigma^{-2}(\Delta+\mathbf g_0\cdot\na)\zeta+\lambda^2pf(\theta_\sigma)^{p-1}\zeta+\sigma^{-2}Q_\sigma\zeta=\mathcal F.
\end{equation}
The coefficient bounds \eqref{Eq.common_coefficient_radius} and Lemma \ref{Lem.analytic_to_multiplier} give finite multiplier norms for $V,J_0^{-1}$ and $\nabla J_0$ at $2\overline\varrho$. Since $\mathbf g_0=J_0^{-1}\nabla J_0$ and $Q_\sigma=(K(\theta_\sigma)-1)V+R_\sigma$, \eqref{Eq.multiplier_algebra} and Lemmas \ref{Lem.phase_coefficients} and \ref{Lem.remainder_equation} imply
\[\mathfrak M_{2\overline\varrho}(\mathbf g_0)+\mathfrak M_{2\overline\varrho}(Q_\sigma)\leq C\]
for all sufficiently large $\sigma$. With $w:=\lambda^{-1}D\zeta$ and $U:=(\zeta,w)^T$, the corresponding first-order equation is
\begin{equation}\label{Eq.energy_first_order}
	DU=\lambda\mathsf P_0(\theta_\sigma)U+\sigma^{-2}\lambda^{-1}E_{21}(\Delta+\mathbf g_0\cdot\na-Q_\sigma)U-\frac{D\lambda}{\lambda}E_{22}U+\lambda^{-1}e_2\mathcal F,
\end{equation}
where $e_1=(1,0)^T$, $e_2=(0,1)^T$ and $E_{ab}:=e_ae_b^T$. The definitions in \eqref{Eq.energy_phase_samples} and $D\theta_\sigma=\lambda$ give
\begin{equation}\label{Eq.sampled_metric_chain_rule}
	\begin{aligned}
		D\xi_N/\xi_N&=-(1+\alpha)\sigma^{-1}-D\lambda/\lambda\quad(N\geq1),\qquad D\delta/\delta=-\sigma^{-1},\\
		DS_N&=\lambda\partial_\theta S(\theta_\sigma,\xi_N,\delta)+(D\xi_N)\partial_\xi S(\theta_\sigma,\xi_N,\delta)+(D\delta)\partial_\delta S(\theta_\sigma,\xi_N,\delta).
	\end{aligned}
\end{equation}
The terms $(D\xi_N)\partial_\xi S$ and $(D\delta)\partial_\delta S$ in \eqref{Eq.sampled_metric_chain_rule} are bounded in relative quadratic form by $CS_N/\sigma$, using \eqref{Eq.S_coercivity} and $|D\lambda|\leq C\sigma^{-2}\nu(\sigma)$ from Lemma \ref{Lem.phase_coefficients}. Since $\lambda\xi_N=N/(\varrho\sigma)$, we obtain
\begin{equation}\label{Eq.frozen_original_defect}
	DS_N+\lambda(\mathsf P_{\xi_N}^TS_N+S_N\mathsf P_{\xi_N})\leq\left(\frac{2\chi N}{\varrho\sigma}+\frac C\sigma\right)S_N.
\end{equation}
Here $\mathsf P_{\xi_N}$ is evaluated at $\theta_\sigma$, and $C$ may depend on $A_*$. At $N=0$, $\xi_0=0$ is constant and the frequency derivative in \eqref{Eq.sampled_metric_chain_rule} vanishes.

\subsection{Differentiated variables and matrix norms}
\label{Subsec.energy_calculus}

The differentiated unknowns are $\mathcal U_N=\mathsf T_N\nabla^N(\mathsf T_N^{-1}U)$ from \eqref{Eq.adapted_layers}, rather than $\nabla^NU$. Differentiation produces coefficients involving $\mathsf T_N^{\pm1}$ and $S_N$, while the product rule produces inputs $\mathsf T_N\nabla^J(\mathsf T_N^{-1}U)$ at orders $J\ne N$. We bound the coefficients in Lemmas \ref{Lem.gauge_bounds} and \ref{Lem.covariant_metric}, and compare the inputs with $\mathcal U_J$ in Lemma \ref{Lem.gauge_transfer}. The factorial identity \eqref{Eq.energy_factorial_shift} then makes the sums in $N,J$ explicit.

Use $S_N,\mathsf T_N,\mathcal U_N$ from \eqref{Eq.energy_phase_samples}--\eqref{Eq.adapted_layers}. For the energy calculation, we define
\begin{equation}\label{Eq.energy_level_sums}
	a_N=(1+N)^{s_*}\frac{\varrho(\tau)^N}{N!},\qquad e_N=a_N\|S_N^{1/2}\mathcal U_N\|_{L^2},\qquad E=\sum_Ne_N^2,\qquad \mathfrak d=\sum_NNe_N^2.
\end{equation}
Thus $\cE^\sharp_\varrho=\sigma^2E$. The sum $\mathfrak d$ records the factor $N$ in the outgoing estimate \eqref{Eq.energy_transport_bounds}. Define
\begin{equation}\label{Eq.energy_covariant_derivatives}
	\Gamma_{N,i}:=(\partial_i\mathsf T_N)\mathsf T_N^{-1},\quad \nabla_i^{[N]}:=\partial_i-\Gamma_{N,i},\quad \nabla_I^{[N]}:=\mathsf T_N\partial_I\mathsf T_N^{-1},\quad \Delta^{[N]}:=\mathsf T_N\Delta \mathsf T_N^{-1}.
\end{equation}
The expressions $\mathsf T_N\partial_I\mathsf T_N^{-1}$ and $\mathsf T_N\Delta \mathsf T_N^{-1}$ in \eqref{Eq.energy_covariant_derivatives} are operator compositions. In particular, $\mathcal U_N=(\nabla_I^{[N]}U)_{|I|=N}$. All subsequent tensor norms and contractions include the full ordered Cartesian indices.

\begin{lemma}\label{Lem.gauge_bounds}
	There are $\Sigma_0\geq\sigma_c$, $R_c\in(\overline\varrho,2\overline\varrho)$ and $C>0$ such that $\nu(\sigma)\leq1$, $1/2\leq\lambda\leq3/2$, and $\mathsf T_N,\mathsf T_N^{-1}$ are uniformly bounded for $\sigma\geq\Sigma_0$. For all integers $r,N\geq0$,
	\begin{equation}\label{Eq.energy_gauge_coefficients}
		\begin{aligned}
			\|\nabla^r(\mathsf T_N-I)\|_{L^\infty}+\|\nabla^r(\mathsf T_N^{-1}-I)\|_{L^\infty}&\leq C\nu(\sigma)r!R_c^{-r},\\
			\|\nabla^r(\lambda-1)\|_{L^\infty}+\|\nabla^r(\lambda^{-1}-1)\|_{L^\infty}&\leq C\sigma^{-1}\nu(\sigma)r!R_c^{-r},\\
			\|\nabla^rD\lambda\|_{L^\infty}&\leq C\sigma^{-2}\nu(\sigma)r!R_c^{-r}.
		\end{aligned}
	\end{equation}
	For every $r,N\geq0$,
	\begin{equation}\label{Eq.energy_adjacent_gauges}
		\|\nabla^r(\mathsf T_{N+1}-\mathsf T_N)\|_{L^\infty}+ \|\nabla^r(\mathsf T_{N+1}^{-1}-\mathsf T_N^{-1})\|_{L^\infty} \leq C\sigma^{-1}\nu(\sigma)r!R_c^{-r}.
	\end{equation}
	For fixed $H\geq0$ and $\max I_{\rm r}<R_c'<R_c$, the bounds \eqref{Eq.energy_gauge_coefficients}--\eqref{Eq.energy_adjacent_gauges} also hold with $R_c^{-r}$ replaced by $(1+r)^{-H}(R_c')^{-r}$, after increasing $C$.
\end{lemma}

\begin{proof}
	Fix analytic radii $\sup I_{\rm r}<\overline\varrho<R_c<R_{\rm aux}<2\overline\varrho$. We first work in the multiplier norm $\mathfrak M_{3\overline\varrho}$. Since $\varrho_2=4\overline\varrho$, Lemma \ref{Lem.analytic_to_multiplier} gives $\mathfrak M_{3\overline\varrho}(\nabla^c\Psi_j)\leq C_c\|\Psi_j\|_{\mathcal B_{\varrho_2}}$ for all $c\in\Z_{\geq0}$ and $0\leq j\leq q_*$. For $a,c\in\Z_{\geq0}$ and $\sigma\geq2$, the polynomial formula for $\Psi$ therefore implies
	\begin{align*}
		\mathfrak M_{3\overline\varrho}\bigl((\partial_\tau^a\nabla_x^c\Psi)(\log\sigma,\cdot)\bigr)\leq\sum_{j=a}^{q_*}\frac{j!}{(j-a)!}(\log\sigma)^{j-a}\mathfrak M_{3\overline\varrho}(\nabla^c\Psi_j)\leq C_{a,c}(1+\log\sigma)^{q_*}.
	\end{align*}
	The sum is empty when $a>q_*$.
	
	We use \eqref{Eq.multiplier_algebra} for scalar and matrix products. The periodic matrix $P_0$ is holomorphic and uniformly bounded in a fixed complex strip. Apply the Taylor-series estimate in Lemma \ref{Lem.phase_functional_calculus} to the entries of $P_0$ in \eqref{Eq.Pzero}, subtracting the constant Taylor term. Since $\mathfrak M_{3\overline\varrho}(\sigma^{-1}\Psi)\to0$, after increasing $\Sigma_0$ we obtain
	\[\sup_{N\geq0}\mathfrak M_{3\overline\varrho}(\mathsf T_N-I)\leq C\nu(\sigma)\leq\tfrac12.\]
	The cutoff $\chi_0(N/(\varrho\sigma))$ in \eqref{Eq.adapted_layers} is independent of $x$. The geometric series for $(I+(\mathsf T_N-I))^{-1}$ and the product bound \eqref{Eq.multiplier_algebra} give
	\[\sup_{N\geq0}\mathfrak M_{3\overline\varrho}(\mathsf T_N^{-1}-I)\leq2\sup_{N\geq0}\mathfrak M_{3\overline\varrho}(\mathsf T_N-I)\leq C\nu(\sigma).\]
	Each term of the multiplier series for $\mathsf T_N-I$ and $\mathsf T_N^{-1}-I$ is bounded by $C\nu(\sigma)$, which proves the first derivative estimate, initially with analytic radius $3\overline\varrho$ and therefore with analytic radius $R_{\rm aux}$. The identities $\lambda-1=\sigma^{-2}T_1\Psi$ and $D\lambda=\sigma^{-3}T_2T_1\Psi$, from \eqref{Eq.phase_velocity} and \eqref{Eq.slow_calculus}, give the last two bounds in \eqref{Eq.energy_gauge_coefficients}; the inverse of $\lambda$ is obtained from a geometric series. The formulas for $\mathsf T_N^{\pm1}$, $\lambda^{\pm1}$ and $D\lambda$ involve only a fixed finite set of slow derivatives of $\Psi$. Thus one integer $M$ in the definition of $\nu$ bounds all the logarithmic degrees used here.
	
	By the mean value theorem,
	\[\left|\chi_0\left(\frac{N+1}{\varrho\sigma}\right)-\chi_0\left(\frac{N}{\varrho\sigma}\right)\right|\leq\frac{\|\chi_0'\|_{L^\infty}}{\varrho\sigma}.\]
	The formula for $\mathsf T_N$ gives $\mathfrak M_{R_{\rm aux}}(\mathsf T_{N+1}-\mathsf T_N)\leq C\nu(\sigma)/\sigma$. For the inverse difference, we use
	\[\mathsf T_{N+1}^{-1}-\mathsf T_N^{-1}=\mathsf T_{N+1}^{-1}(\mathsf T_N-\mathsf T_{N+1})\mathsf T_N^{-1}.\]
	For $H\geq0$ and $\max I_{\rm r}<R_c'<R_c$, define
	\[C_{H,R_c,R_c'}:=\sup_{r\geq0}(1+r)^H(R_c'/R_c)^r<\infty.\]
	Then $r!R_c^{-r}\leq C_{H,R_c,R_c'}r!(1+r)^{-H}(R_c')^{-r}$ for every $r\geq0$. Increase $\Sigma_0$ once more to ensure $\nu\leq1$ and $1/2\leq\lambda\leq3/2$.
\end{proof}

Integration by parts differentiates $\mathcal U_N^*S_N\mathcal U_N$. Substituting $\partial_i\mathcal U_N=\nabla_i^{[N]}\mathcal U_N+\Gamma_{N,i}\mathcal U_N$ from \eqref{Eq.energy_covariant_derivatives} combines the coefficient derivatives into $\partial_iS_N+\Gamma_{N,i}^TS_N+S_N\Gamma_{N,i}$. Lemma \ref{Lem.covariant_metric} estimates this matrix relative to $S_N$, including the low-frequency range where $S_N^{-1}$ is large.

\begin{lemma}[The spatial metric derivative]\label{Lem.covariant_metric}
	We define
	\begin{equation}\label{Eq.energy_metric_derivative}
		\mathcal Q_{N,i}:=\partial_iS_N+\Gamma_{N,i}^TS_N+S_N\Gamma_{N,i}.
	\end{equation}
	For $\sigma\geq\Sigma_0$ one has
	\begin{equation}\label{Eq.energy_metric_derivative_bound}
		\|S_N^{-1/2}\mathcal Q_{N,i}S_N^{-1/2}\|\leq C\nu(\sigma)\min\{1,\xi_N+\delta\},
	\end{equation}
	and
	\begin{equation}\label{Eq.rank_one_metric}
		\|S_N^{1/2}E_{21}S_N^{-1/2}\|=\frac{(S_N)_{22}}{\sqrt{\det S_N}}\leq\frac{C}{\xi_N+\delta},\quad\|S_N^{1/2}E_{22}S_N^{-1/2}\|\leq C\left(1+\frac1{\xi_N+\delta}\right).
	\end{equation}
\end{lemma}

\begin{proof}
	When $\chi_0(N/(\varrho\sigma))=1$, define $P:=P_{\xi_N^2}(\theta_\sigma)$, $P_*:=P_0(\theta_\sigma)$ and $R:=P^{-1}P_*$. The estimates in Lemma \ref{Lem.low_frame} give $\|R-I\|+\|\partial_\theta R\|\leq C\xi_N^2$ and $\|\xi_N\partial_{\xi_N}R\|\leq C\xi_N^2$, with phase and frequency treated as independent parameters. The low-frequency formula in Proposition \ref{Prop.matrix} and \eqref{Eq.adapted_layers} imply
	\begin{equation}\label{Eq.energy_metric_conjugation}
			\mathsf T_N^TS_N\mathsf T_N=P_0(\sigma)^{-T}R^T\diag(\omega(\xi_N^2)+\delta^2,1)RP_0(\sigma)^{-1},\quad
			\partial_i(\mathsf T_N^TS_N\mathsf T_N)=\mathsf T_N^T\mathcal Q_{N,i}\mathsf T_N.
	\end{equation}
	The cancellation in \eqref{Eq.energy_metric_conjugation} is essential near $\xi_N=0$: $P_0(\sigma)$ has no spatial dependence, so spatial differentiation of the conjugated matrix differentiates $R$, not $P_0(\theta_\sigma)$ separately. In particular, at $N=0$ we have $R=I_2$, $\xi_0=0$ and
\[\mathsf T_0^TS_0\mathsf T_0=P_0(\sigma)^{-T}\diag(\delta^2,1)P_0(\sigma)^{-1},\qquad \mathcal Q_{0,i}=0.\]
For $N>0$, write $\kappa_N=(\omega(\xi_N^2)+\delta^2)^{1/2}\asymp\xi_N+\delta$. Holding $\xi_N,\delta$ fixed when differentiating in phase, the relative quadratic-form identity \eqref{Eq.relative_quardratic-form} gives
\begin{align*}
&\left\|\diag(\kappa_N^{-1},1)R^{-T}\partial_\theta\!\left(R^T\diag(\kappa_N^2,1)R\right)R^{-1}\diag(\kappa_N^{-1},1)\right\|\\
&\qquad\leq2\left\|\diag(\kappa_N,1)(\partial_\theta R)R^{-1}\diag(\kappa_N^{-1},1)\right\|\leq\frac{C\xi_N^2}{\kappa_N}.
\end{align*}
Indeed, differentiating $R^T\diag(\kappa_N^2,1)R$ produces the displayed matrix and its transpose after normalization. Multiplication by $|\partial_i\theta_\sigma|=\sigma^{-1}|\partial_i\Psi|\leq C\nu(\sigma)$ therefore gives $C\nu(\sigma)\xi_N$. The factor $\xi_N^2$ from $\partial_\theta R$ compensates the possible factor $\kappa_N^{-1}$ from the degenerate energy; estimating $\partial_i S_N$ without the two terms containing $\Gamma_{N,i}$ would not reveal this cancellation.

For the frequency contribution, $\mathsf T_N$ is independent of the frequency parameter where the cutoff equals one. Congruence by $\mathsf T_N$ therefore leaves the relative frequency estimate in \eqref{Eq.S_coercivity} unchanged. Using \eqref{Eq.energy_gauge_coefficients},
\[\left|\frac{\partial_i\xi_N}{\xi_N}\right|=|\partial_i\log\lambda|\leq C\nu(\sigma)/\sigma,\qquad \partial_i\delta=0.\]
The sum of the phase and frequency contributions is bounded by $C\nu(\sigma)(\xi_N+\delta)$. This proves \eqref{Eq.energy_metric_derivative_bound} where the cutoff equals one, including $N=0$.

	Where $\chi_0(N/(\varrho\sigma))\ne1$ and $\xi_N$ remains bounded, $\xi_N$ is bounded below by a positive constant. Smoothness and positivity of $S$ on the resulting compact parameter set, together with \eqref{Eq.energy_gauge_coefficients}, give $\|S_N^{-1/2}\mathcal Q_{N,i}S_N^{-1/2}\|\leq C\nu(\sigma)$. At sufficiently large $\xi_N$, $\mathsf T_N=I$ and $S_N=\diag(1+\xi_N^2,1)$, so
	\[\|S_N^{-1/2}(\partial_iS_N)S_N^{-1/2}\|=\frac{2\xi_N^2}{1+\xi_N^2}|\partial_i\log\lambda|\leq C\nu(\sigma)/\sigma.\]
	These frequency ranges cover every $N\geq0$ and prove \eqref{Eq.energy_metric_derivative_bound}.
	
	For a positive two-by-two matrix $S$, the rank-one norm is
	\[\|S^{1/2}e_2e_1^TS^{-1/2}\|^2=(e_2^TSe_2)(e_1^TS^{-1}e_1)=S_{22}^2/\det S.\]
	At low frequency, \eqref{Eq.low_metric_entries} gives $S_{22}^2/\det S\leq C/(\xi+\delta)^2$. On positive compact intervals use positivity, and at high frequency use $S=\diag(1+\xi^2,1)$. For $E_{22}=e_2e_2^T$, the rank-one formula gives
	\[\|S^{1/2}E_{22}S^{-1/2}\|^2=(e_2^TSe_2)(e_2^TS^{-1}e_2)=\frac{S_{11}S_{22}}{\det S}.\]
	The low-frequency bound is again $C/(\xi+\delta)^2$, whereas the expression equals one at high frequency. Compactness handles the intermediate region and proves \eqref{Eq.rank_one_metric}.
\end{proof}

Before estimating the differentiated equation, we record the factorial identities used in the sums over derivative orders. Lemma \ref{Lem.energy_shifts} identifies the factors of the input derivative order that remain after cancellation of the analytic weights.

\begin{lemma}\label{Lem.energy_shifts}
	Let $q\in\{0,1,2\}$, $N\geq r\geq0$ and $J:=N-r+q$. Then
	\begin{equation}\label{Eq.energy_factorial_shift}
		\frac{a_N}{a_J}\binom Nr=\left(\frac{1+N}{1+J}\right)^{s_*}\frac{\varrho^{r-q}}{r!}\frac{J!}{(J-q)!}.
	\end{equation}
	For a nonnegative sequence $e=(e_N)_{N\geq0}$ and $r\geq0$,
	\begin{equation}\label{Eq.energy_sequence_sums}
		\sum_{N\geq r}e_Ne_{N-r+q}\leq\sum_{N\geq0}e_N^2,\quad\sum_{N\geq r}(1+N-r+q)e_Ne_{N-r+q}\leq\sqrt3\sum_{N\geq0}(1+N)e_N^2.
	\end{equation}
	If $r\geq2$ and $q=1$, then
	\[\sum_{N\geq r}(N-r+1)e_Ne_{N-r+1}\leq\sum_{N\geq0}Ne_N^2.\]
	Each inequality remains valid after summation against a nonnegative $\ell^1$ sequence in $r$, with the $\ell^1$ norm of the coefficient sequence as an additional factor on the right.
\end{lemma}

\begin{proof}
	Since $J-q=N-r$, the definitions \eqref{Eq.energy_level_sums} give
	\[\frac{a_N}{a_J}\binom Nr=\frac{(1+N)^{s_*}}{(1+J)^{s_*}}\varrho^{N-J}\frac{J!}{r!(N-r)!}=\frac{(1+N)^{s_*}}{(1+J)^{s_*}}\frac{\varrho^{r-q}}{r!}\frac{J!}{(J-q)!},\]
	which proves \eqref{Eq.energy_factorial_shift}. Cauchy--Schwarz gives the first bound in \eqref{Eq.energy_sequence_sums}. For the second bound, define $f_N:=\sqrt{1+N}e_N$. Since $J\leq N+2$,
	\[(1+J)e_Ne_J=\sqrt{\frac{1+J}{1+N}}f_Nf_J\leq\sqrt3 f_Nf_J.\]
	Consequently,
	\[\sum_{N\geq r}(1+J)e_Ne_J\leq\sqrt3\left(\sum_{N\geq r}f_N^2\right)^{1/2}\left(\sum_{N\geq r}f_{N-r+q}^2\right)^{1/2}\leq\sqrt3\sum_{N\geq0}(1+N)e_N^2.\]
	When $r\geq2$ and $q=1$, one has $1\leq J<N$, and $Je_Ne_J\leq\sqrt{NJ}e_Ne_J$. The same argument with $f_N=\sqrt N e_N$ gives the final inequality. 
\end{proof}

At output derivative order $N$, the product rule in the differentiated equation produces $\mathsf T_N\nabla^J(\mathsf T_N^{-1}U)$, where $J$ need not equal $N$. Lemma \ref{Lem.gauge_transfer} expresses $\mathsf T_N\nabla^J(\mathsf T_N^{-1}U)$ in terms of $\mathcal U_J$ and lower derivative orders.

\begin{lemma}[Changes of differentiated frame]\label{Lem.gauge_transfer}
	For $N,J\geq0$, define
	\begin{equation}\label{Eq.gauge_transfer_definitions}
		\mathcal V_{N,J}:=\mathsf T_N\na^J(\mathsf T_N^{-1}U),\qquad t_N:=\chi_0\!\left(\frac{N}{\varrho\sigma}\right),\qquad d_{N,J}:=|t_N-t_J|.
	\end{equation}
	There exist $C>0$, $H\in\Z_{\geq1}$ and $t\in(0,1)$ such that, for sufficiently large $\sigma$,
	\begin{equation}\label{Eq.gauge_transfer_bound}
		a_J\|S_J^{1/2}(\mathcal V_{N,J}-\mathcal U_J)\|_{L^2}\leq C\nu(\sigma)d_{N,J}(1+|N-J|)^H\sum_{r=1}^Jt^re_{J-r}.
	\end{equation}
	The empty sum is zero. In particular, the difference vanishes if $t_N=t_J$, and the right side has an additional factor $\sigma^{-1}$ when $|N-J|=1$.
\end{lemma}

\begin{proof}
	We define $K:=P_0(\theta_\sigma)P_0(\sigma)^{-1}-I$. By \eqref{Eq.gauge_transfer_definitions}, $\mathsf T_N=I+t_NK$, and
	\[G_{N,J}:=\mathsf T_N^{-1}\mathsf T_J=I+(t_J-t_N)\mathsf T_N^{-1}K.\]
	Fix $\varrho<R_2<R_c$. Lemma \ref{Lem.gauge_bounds} and the convergent inverse series give, for $r\geq1$,
	\[\|\mathsf T_N\na^rG_{N,J}\mathsf T_J^{-1}\|_{L^\infty}\leq C\nu d_{N,J}r!R_2^{-r}.\]
	We first record the corresponding two-metric estimate:
	\begin{equation}\label{Eq.gauge_two_metric}
		\frac{a_J}{a_{J-r}}\binom Jr\|S_J^{1/2}\mathsf T_N\na^rG_{N,J}\mathsf T_J^{-1}S_{J-r}^{-1/2}\|_{L^\infty}\leq C\nu d_{N,J}(1+|N-J|+r)^H t_0^r,
	\end{equation}
	where $t_0:=\varrho/R_2$. Here and below $1\leq r\leq J$. For the low-frequency comparison, we use
	\[\|S_J^{1/2}\|\leq C(1+J/\sigma),\qquad \|S_{J-r}^{-1/2}\|\leq C\left(1+\frac{\sigma}{1+J-r}\right).\]
	If $d_{N,J}\ne0$, then $\max\{N,J\}\geq c\sigma$ and $\min\{N,J\}\leq C\sigma$, with constants uniform for $\varrho\in I_{\rm r}$. If $|N-J|+r\leq c\sigma/4$, the indices $J$ and $J-r$ are bounded above and below by positive multiples of $\sigma$, so both matrix norms are bounded. Otherwise $\sigma\leq C(|N-J|+r)$ and $J\leq C\sigma+|N-J|$; substitution into the two norm bounds gives $C(1+|N-J|+r)^2$. Thus the product of the matrix norms is bounded by a fixed polynomial in $1+|N-J|+r$. Lemma \ref{Lem.energy_shifts} absorbs the scalar index weight into $C(1+|N-J|+r)^H$ after increasing $H$ and proves \eqref{Eq.gauge_two_metric}.
	
	For an ordered index $I=(i_1,\ldots,i_J)$ and $S\subset\{1,\ldots,J\}$, let $I_S$ be the subindex in the positions belonging to $S$, in their original order. Apply the ordered product rule to $\mathsf T_N\partial_I(G_{N,J}\mathsf T_J^{-1}U)$ and use $\mathsf T_NG_{N,J}=\mathsf T_J$. The term with no derivative on $G_{N,J}$ is $(\mathcal U_J)_I$, so
	\begin{equation}\label{Eq.gauge_transfer_product}
		(\mathcal V_{N,J}-\mathcal U_J)_I=\sum_{\varnothing\ne S\subset\{1,\ldots,J\}}\mathsf T_N(\partial_{I_S}G_{N,J})\mathsf T_J^{-1}(\mathcal V_{J,J-|S|})_{I_{S^c}}.
	\end{equation}
	For $|S|=r$, the rearrangement of derivative indices preserves the Euclidean tensor norm, and there are $\binom Jr$ terms. In \eqref{Eq.gauge_transfer_product}, write
\[\mathcal V_{J,J-r}=\mathcal U_{J-r}+(\mathcal V_{J,J-r}-\mathcal U_{J-r})\]
and apply \eqref{Eq.gauge_transfer_product} again to the difference. After $a$ substitutions, a term ending at derivative order $J-r$ has lost $r_1,\ldots,r_a\geq1$ derivatives with $r_1+\cdots+r_a=r$. Thus $a\leq r\leq J$, so the expansion is finite before any estimate or infinite summation is used.

The analytic weights telescope along the derivative orders $J,J-r_1,J-r_1-r_2,\ldots,J-r$. Formula \eqref{Eq.gauge_two_metric} therefore bounds the coefficient of each such term by a product of scalar coefficient bounds. The factor $d_{N,J}$ occurs in the first factor and is never discarded. At every subsequent substitution we bound the cutoff difference by one. Also,
\[(1+|N-J|+r_1)^H\leq(1+|N-J|)^H(1+r_1)^H.\]
After enlarging $C,H$, all additional polynomial factors in the derivative losses can be included in
\[k_r:=C\nu(\sigma)(1+r)^{2H}t_0^r\quad(r\geq1),\qquad k_0:=0.\]
The coefficient of $\mathcal U_{J-r}$ in the fully expanded difference consequently has the bound
\begin{equation}\label{Eq.gauge_transfer_compositions}
C d_{N,J}(1+|N-J|)^H\sum_{a=1}^r\ \sum_{\substack{r_1+\cdots+r_a=r\\r_1,\ldots,r_a\geq1}}k_{r_1}\cdots k_{r_a}
=C d_{N,J}(1+|N-J|)^H\sum_{a=1}^r(k^{*a})_r.
\end{equation}
Choose $t\in(t_0,1)$ and increase the starting time until $\sum_{r\geq1}t^{-r}k_r\leq1/2$. Since all coefficients are nonnegative, summing the convolutions in \eqref{Eq.gauge_transfer_compositions} gives
\[\sum_{r\geq1}t^{-r}\sum_{a\geq1}(k^{*a})_r=\frac{\sum_{r\geq1}t^{-r}k_r}{1-\sum_{r\geq1}t^{-r}k_r}\leq C\nu(\sigma).\]
In particular, $\sum_{a=1}^r(k^{*a})_r\leq C\nu(\sigma)t^r$. Multiplying \eqref{Eq.gauge_transfer_compositions} by $e_{J-r}$ and summing over $1\leq r\leq J$ proves \eqref{Eq.gauge_transfer_bound}, with constants independent of $J$. For $|N-J|=1$, use $d_{N,J}\leq\|\chi_0'\|_{L^\infty}/(\varrho\sigma)$ in \eqref{Eq.gauge_transfer_bound} to obtain
	\[a_J\|S_J^{1/2}(\mathcal V_{N,J}-\mathcal U_J)\|_{L^2} \leq C\sigma^{-1}\nu(\sigma)\sum_{r=1}^Jt^re_{J-r}.\]
	This completes the proof.
\end{proof}

We next estimate the contribution of $(a\cdot\nabla+c)U$. The factor $(\xi_J+\delta)^{-1}$ in \eqref{Eq.rank_one_metric} combines with $\sigma^{-2}J$ to give a bound of order $\sigma^{-1}$.

\begin{lemma}\label{Lem.energy_lower_order}
	Let $a=a(\sigma,x)\in\R^d$ and $c=c(\sigma,x)\in\R$ have uniformly bounded multiplier norms at $2\overline\varrho$. For $U\in\mathcal B_{r'}(\Omega;\C^2)$ with $\varrho<r'<3\overline\varrho/2$, define
	\[\mathcal R_N[a,c]:=\sigma^{-2}\mathsf T_N\na^N\bigl(\mathsf T_N^{-1}\lambda^{-1}E_{21}(a\cdot\na+c)U\bigr).\]
	For all sufficiently large $\sigma$,
	\begin{equation}\label{Eq.lower_order_energy_bound}
		\left|2\operatorname{Re}\sum_{N\geq0}(a_N)^2(S_N\mathcal U_N,\mathcal R_N[a,c])_{L^2}\right|\leq\frac{C}\sigma E.
	\end{equation}
	The constant depends on the coefficient bounds but not on $\sigma$ or $U$.
\end{lemma}

\begin{proof}
	Write $V_N:=\mathsf T_N^{-1}U$. Then the expression inside the differentiated brackets is
	\[\mathsf T_N^{-1}\lambda^{-1}E_{21}(a\cdot\na+c)U=\sum_{i=1}^d M_N^i\partial_iV_N+M_N^0V_N,\]
	where
	\begin{equation}\label{Eq.energy_lower_coefficients}
		M_N^i:=\lambda^{-1}a_i\mathsf T_N^{-1}E_{21}\mathsf T_N,\qquad M_N^0:=\lambda^{-1}\mathsf T_N^{-1}E_{21}(a\cdot\na \mathsf T_N+c\mathsf T_N).
	\end{equation}
	The matrices $M_N^i,M_N^0$ in \eqref{Eq.energy_lower_coefficients} have derivative bounds $Cr!R_c^{-r}$ for some $R_c>\max I_{\rm r}$, by \eqref{Eq.energy_gauge_coefficients} and \eqref{Eq.multiplier_algebra}. For the first-order terms, we define $J:=N-r+1$, and for the zeroth-order terms we define $J:=N-r$. Lemmas \ref{Lem.energy_shifts} and \ref{Lem.gauge_transfer} express the first-order and zeroth-order contributions in the tensors $\mathcal U_J$.
	
	We next compare the input and output matrix norms. In the low-frequency region the frame formula gives
	\[\|S_N^{1/2}\|\leq C,\qquad \|S_J^{-1/2}\|\leq\frac{C}{\xi_J+\delta}.\]
	Consequently every bounded differentiated coefficient between the balanced input and output metrics costs at most $C/(\xi_J+\delta)$. When both gauges are the identity, the coefficient matrices are multiples of $E_{21}$; the same bound follows from \eqref{Eq.rank_one_metric} and the sampled matrix comparisons, with a fixed polynomial in $1+|N-J|$. In the cutoff transition the frequency is bounded below and the metrics are uniformly nondegenerate. If two levels belong to different regions and are separated by a large gap, the exponential coefficient tail absorbs the polynomial comparison factors, as in \eqref{Eq.gauge_two_metric}. Here $J\leq N+1$, so adjacent levels that cross a threshold are covered by the same compact-frequency bound. The expansions in Lemma \ref{Lem.gauge_transfer} satisfy the same bounds.
	
	The remaining derivative factors are
	\[\frac{\sigma^{-2}J}{\xi_J+\delta}=\frac1\sigma\frac{J}{1+J/(\varrho\lambda)}\leq\frac{C_\varrho}\sigma,\qquad \frac{\sigma^{-2}}{\xi_J+\delta}\leq\frac1\sigma.\]
	Thus, for a summable sequence $b_r\leq C(1+r)^H(\varrho/R_c)^r$, the absolute value in \eqref{Eq.lower_order_energy_bound} is at most
	\[\frac C\sigma\sum_{q=0}^1\sum_{r\geq0}b_r\sum_{N\geq r}e_Ne_{N-r+q}\leq\frac C\sigma\left(\sum_{r\geq0}b_r\right)\sum_{N\geq0}e_N^2.\]
	The last inequality is Cauchy--Schwarz in each shifted sequence. The coefficients $a=\mathbf g_0$ and $c=-Q_\sigma$ satisfy the hypotheses by the bounds following \eqref{Eq.energy_scalar_equation}.
\end{proof}

\subsection{The interior identity}
\label{Subsec.energy_interior}

We apply $\mathsf T_N\partial_I\mathsf T_N^{-1}$ to \eqref{Eq.energy_first_order} for $|I|=N$ and derive \eqref{Eq.differentiated_energy_equation}. The terms with one derivative of $\bvec$ give $-2\gamma_v\mathfrak d/\sigma$; integration by parts in the Laplacian produces the matrix contribution $-\lambda\xi_N^2E_{21}$. The remaining terms are estimated in \eqref{Eq.energy_remainder_bounds}, retaining the dependence of $\mathsf T_N$ on the moving radius $\varrho$.

For an ordered tensor $Z_I$ of length $N$, define
\begin{equation}\label{Eq.energy_transport_tensor}
	(\mathscr B_NZ)_I:=\sum_{a=1}^N(\partial_{i_a}b_k)Z_{i_1\cdots i_{a-1}k i_{a+1}\cdots i_N},\qquad \mathscr B_0:=0.
\end{equation}
The coefficient matrix acts only on tensor indices and commutes with every two-by-two matrix acting on the unknowns. Define also
\begin{equation}\label{Eq.energy_transport_higher}
	\mathcal C_N:=\sum_{r=2}^N\binom Nr(\nabla^r\bvec)\cdot \mathsf T_N\nabla^{N-r+1}(\mathsf T_N^{-1}U),
\end{equation}
where the contraction is the sum of the terms in $\mathsf T_N[\partial_I,B](\mathsf T_N^{-1}U)$ with at least two derivatives on $\bvec$; define $\mathcal C_0=\mathcal C_1:=0$.

For $V_N:=\mathsf T_N^{-1}U$, define the coefficients
\begin{equation}\label{Eq.energy_conjugated_coefficients}
	\begin{aligned}
		A_N&:=\mathsf T_N^{-1}\lambda\mathsf P_0(\sigma+\sigma^{-1}\Psi)\mathsf T_N-\mathsf T_N^{-1}D\mathsf T_N,\\
		K_N&:=\lambda^{-1}\mathsf T_N^{-1}E_{21}\mathsf T_N,\qquad K_N^i:=\lambda^{-1}\mathsf T_N^{-1}E_{21}(2\partial_i\mathsf T_N+(\mathbf g_0)_i\mathsf T_N),\\
		K_N^0&:=\lambda^{-1}\mathsf T_N^{-1}E_{21}(\Delta \mathsf T_N+\mathbf g_0\cdot\na \mathsf T_N-Q_\sigma \mathsf T_N)-\sigma^2\frac{D\lambda}{\lambda}\mathsf T_N^{-1}E_{22}\mathsf T_N.
	\end{aligned}
\end{equation}

\begin{lemma}\label{Lem.differentiated_energy_equation}
	For a smooth solution of \eqref{Eq.energy_scalar_equation}, set
	\[\mathcal F_N:=\mathsf T_N\nabla^N(\mathsf T_N^{-1}\lambda^{-1}e_2\mathcal F).\]
	For $|I|=N$, define
	\begin{equation}\label{Eq.exact_energy_remainder}
		(\mathcal R_N)_I=\mathsf T_N[\partial_I,A_N]V_N+\sigma^{-2}\mathsf T_N[\partial_I,K_N]\Delta V_N+\sigma^{-2}\mathsf T_N\partial_I\left(\sum_iK_N^i\partial_iV_N+K_N^0V_N\right).
	\end{equation}
	Then, for every $N\geq0$,
	\begin{equation}\label{Eq.differentiated_energy_equation}
		D\mathcal U_N=\lambda\mathsf P_0(\theta_\sigma)\mathcal U_N +\sigma^{-2}\lambda^{-1}E_{21}\Delta^{[N]}\mathcal U_N -\sigma^{-1}\mathscr B_N\mathcal U_N-\sigma^{-1}\mathcal C_N+\mathcal R_N+\mathcal F_N.
	\end{equation}
\end{lemma}

\begin{proof}
	For $V_N=\mathsf T_N^{-1}U$, substitute $U=\mathsf T_NV_N$ into \eqref{Eq.energy_first_order} and use the coefficients \eqref{Eq.energy_conjugated_coefficients}. The product rule gives
	\begin{equation}\label{Eq.energy_conjugated_system}
		DV_N=A_NV_N+\sigma^{-2}\left(K_N\Delta V_N+\sum_iK_N^i\partial_iV_N+K_N^0V_N\right)+\mathsf T_N^{-1}\lambda^{-1}e_2\mathcal F.
	\end{equation}
	To differentiate \eqref{Eq.energy_conjugated_system}, use, for $|I|=N$,
	\[D\partial_IV_N=\partial_IDV_N-\sigma^{-1}[\partial_I,B]V_N.\]
	Multiplication by $\mathsf T_N$ and \eqref{Eq.energy_transport_tensor}--\eqref{Eq.energy_transport_higher} give
	\[\mathsf T_N[\partial_I,B]V_N=(\mathscr B_N\mathcal U_N)_I+(\mathcal C_N)_I.\]
	In the terms with no derivative on $A_N,K_N$, the definitions \eqref{Eq.energy_conjugated_coefficients} imply
	\begin{align*}
		(D\mathsf T_N)\partial_IV_N+\mathsf T_NA_N\partial_IV_N=\lambda\mathsf P_0(\theta_\sigma)(\mathcal U_N)_I,\quad 
		\mathsf T_NK_N\Delta\partial_IV_N=\lambda^{-1}E_{21}\Delta^{[N]}(\mathcal U_N)_I.
	\end{align*}
	The coefficient derivatives are exactly \eqref{Eq.exact_energy_remainder}, because
	\[[\partial_I,a]V_N=\sum_{\varnothing\ne J\subset\{1,\ldots,N\}}(\partial_{I_J}a)\partial_{I_{J^c}}V_N.\]
	The source is $\mathcal F_N=\mathsf T_N\nabla^N(\mathsf T_N^{-1}\lambda^{-1}e_2\mathcal F)$. Combining the four contributions proves \eqref{Eq.differentiated_energy_equation}.
\end{proof}

We next estimate the terms in $D\mathsf T_N$ containing $\alpha=\varrho'/\varrho$. The factor $q_N\chi_0'(q_N)$ restricts the output derivative order to $N\asymp\sigma$.

\begin{lemma}\label{Lem.moving_gauge}
	We define $q_N=N/(\varrho\sigma)$ and $H_\sigma=P_0(\theta_\sigma)P_0(\sigma)^{-1}-I$. The part of the differentiated equation caused by the variation of the analytic radius is
	\begin{equation}\label{Eq.moving_gauge_remainder}
		\mathcal R_N^{\rm mov}=\frac{\alpha}{\sigma} \mathsf T_N[\nabla^N,Z_N](\mathsf T_N^{-1}U),\qquad Z_N=q_N\chi_0'(q_N)\mathsf T_N^{-1}H_\sigma, \qquad \alpha=\varrho'/\varrho.
	\end{equation}
	For $0\leq\alpha\leq A_*$, the contribution of $\mathcal R_N^{\rm mov}$ is bounded by
	\[\left|2\operatorname{Re}\sum_Na_N^2 (S_N\mathcal U_N,\mathcal R_N^{\rm mov})_{L^2}\right| \leq \frac{C\alpha\nu(\sigma)}{\sigma}E.\]
\end{lemma}

\begin{proof}
	Since $Dq_N=-(1+\alpha)q_N/\sigma$, the additional term in $D\mathsf T_N$ is $-\alpha q_N\chi_0'(q_N)H_\sigma/\sigma$. Thus $A_N$ in \eqref{Eq.energy_conjugated_coefficients} acquires $\alpha Z_N/\sigma$. After differentiating the conjugated equation and multiplying by $\mathsf T_N$, the term $\alpha \mathsf T_NZ_N\nabla^N(\mathsf T_N^{-1}U)/\sigma$ cancels with the additional term in $D\mathsf T_N$. The remainder is exactly \eqref{Eq.moving_gauge_remainder}.
	
	By \eqref{Eq.energy_gauge_coefficients} and \eqref{Eq.multiplier_algebra}, $\|\nabla^kZ_N\|_{L^\infty}\leq C\nu k!R_c^{-k}$, and $Z_N$ vanishes unless $q_N$ lies in the fixed cutoff transition. On $\operatorname{supp}\chi_0'$ $N\asymp\sigma$ and the output metric is uniformly positive and bounded. Write $J=N-k$ in a commutator term. If $J\geq N/2$, $\xi_J$ is bounded away from zero, so the input metric is also uniformly nondegenerate. If $J<N/2$, the general input estimate costs at most $C\sigma\leq C(1+k)$; the factor $1+k$ is absorbed by the exponentially decaying analytic coefficient sequence. Lemma \ref{Lem.gauge_transfer} permits the same estimates after replacing the $N$-gauge input by $\mathcal U_J$. The factorial coefficient is
	\[\frac{a_N}{a_J}\binom Nk =\left(\frac{1+N}{1+J}\right)^{s_*}\frac{\varrho^k}{k!}.\]
	Consequently the absolute energy contribution is bounded by
	\[\frac{C\alpha\nu}{\sigma} \sum_{k\geq1}(1+k)^H(\varrho/R_c)^k\sum_{N\geq k}e_Ne_{N-k}.\]
	Since $\sum_{k\geq1}(1+k)^H(\varrho/R_c)^k\leq C$ uniformly for $\varrho\in I_{\rm r}$, Cauchy--Schwarz gives
	\[\frac{C\alpha\nu}{\sigma}\sum_{k\geq1}(1+k)^H(\varrho/R_c)^k\sum_{N\geq k}e_Ne_{N-k}\leq\frac{C\alpha\nu}{\sigma}\sum_{N\geq0}e_N^2=\frac{C\alpha\nu}{\sigma}E.\qedhere\]
\end{proof}

The source in \eqref{Eq.energy_first_order} is $\lambda^{-1}e_2\mathcal F$. The bound $e_2^TS_Ne_2=(S_N)_{22}\leq C$ in \eqref{Eq.S_coercivity} permits a source estimate without a factor $(\xi_N+\delta)^{-1}$. Lemma \ref{Lem.scalar_source_energy} retains this bound after spatial differentiation.

\begin{lemma}\label{Lem.scalar_source_energy}
	For $\mathcal F_N=\mathsf T_N\nabla^N(\mathsf T_N^{-1}\lambda^{-1}e_2\mathcal F)$ one has
	\begin{equation}\label{Eq.scalar_source_matrix_bound}
		\left(\sum_{N\geq0}a_N^2\|S_N^{1/2}\mathcal F_N\|_{L^2}^2\right)^{1/2}\leq C\|\mathcal F\|_{\mathcal B_\varrho}.
	\end{equation}
	The constant is uniform for $\varrho\in I_{\rm r}$ and sufficiently large $\sigma$.
\end{lemma}

\begin{proof}
	For $r=0$, set $C_{N,0}:=\lambda^{-1}e_2$, the coefficient after multiplication by $\mathsf T_N$. Then
	\[\|S_N^{1/2}\lambda^{-1}e_2\|^2=\lambda^{-2}(S_N)_{22}\leq C.\]
	For $r\geq1$ let $C_{N,r}=\mathsf T_N\nabla^r(\mathsf T_N^{-1}\lambda^{-1}e_2)$. On the support of the cutoff, $\xi_N$ is bounded and $\|S_N^{1/2}\|\leq C$; Lemma \ref{Lem.gauge_bounds} therefore gives $\|S_N^{1/2}C_{N,r}\|_{L^\infty}\leq Cr!R_c^{-r}$. Outside the cutoff, $\mathsf T_N=I$ and $C_{N,r}=(\nabla^r\lambda^{-1})e_2$, so the same estimate follows from $(S_N)_{22}\leq C$. The ordered Leibniz rule and the factorial weights yield
	\begin{align*}
		a_N\|S_N^{1/2}\mathcal F_N\|_{L^2} \leq\sum_{r=0}^N\frac{a_N}{a_{N-r}}\binom Nr\|S_N^{1/2}C_{N,r}\|_{L^\infty}\,a_{N-r}\|\nabla^{N-r}\mathcal F\|_{L^2}\leq\sum_{r=0}^Nb_r a_{N-r}\|\nabla^{N-r}\mathcal F\|_{L^2},
	\end{align*}
	where $b_r=C(1+r)^{s_*}(\varrho/R_c)^r$ and $\sum_rb_r\leq C$. Young's inequality gives
	\[\left(\sum_Na_N^2\|S_N^{1/2}\mathcal F_N\|_{L^2}^2\right)^{1/2} \leq\|b\|_{\ell^1}\|\mathcal F\|_{\mathcal B_\varrho} \leq C\|\mathcal F\|_{\mathcal B_\varrho}.\qedhere\]
\end{proof}

We next estimate each coefficient commutator in \eqref{Eq.exact_energy_remainder}, keeping the two-derivative terms separate from the lower-order terms. The source estimate is already available from Lemma \ref{Lem.scalar_source_energy}.

\begin{lemma}\label{Lem.energy_remainders}
	The remainder $\mathcal R_N$ in \eqref{Eq.exact_energy_remainder} and the source $\mathcal F_N=\mathsf T_N\nabla^N(\mathsf T_N^{-1}\lambda^{-1}e_2\mathcal F)$ satisfy
	\begin{equation}\label{Eq.energy_remainder_bounds}
		\begin{aligned}
			\left|2\operatorname{Re}\sum_N(a_N)^2(S_N\mathcal U_N,\mathcal R_N)_{L^2}\right|&\leq\frac C\sigma E+\frac{C\nu(\sigma)}\sigma(E+\mathfrak d),\\
			\left|2\operatorname{Re}\sum_N(a_N)^2(S_N\mathcal U_N,\mathcal F_N)_{L^2}\right|&\leq C\sqrt{E}\|\mathcal F\|_{\mathcal B_\varrho}.
		\end{aligned}
	\end{equation}
	The constants are uniform for sufficiently large $\sigma$. The series are absolutely convergent when the solution has a strictly larger analytic radius.
\end{lemma}

\begin{proof}
	Lemma \ref{Lem.moving_gauge} bounds the additional part $\mathcal R_N^{\rm mov}$. Lemma \ref{Lem.energy_lower_order}, with $(a,c)=(\mathbf g_0,-Q_\sigma)$, bounds the contributions of $\mathbf g_0$ and $Q_\sigma$ by $CE/\sigma$. Using \eqref{Eq.energy_conjugated_coefficients} and \eqref{Eq.moving_gauge_remainder}, define
	\begin{equation}\label{Eq.energy_remainder_split}
		\begin{aligned}
			(\mathcal R_N^A)_I&:=\mathsf T_N[\partial_I,A_N-\alpha\sigma^{-1}Z_N]V_N,\qquad (\mathcal R_N^K)_I:=\sigma^{-2}\mathsf T_N[\partial_I,K_N]\Delta V_N,\\
			(\mathcal R_N^1)_I&:=2\sigma^{-2}\mathsf T_N\partial_I\left(\sum_i\lambda^{-1}\mathsf T_N^{-1}E_{21}(\partial_i\mathsf T_N)\partial_iV_N\right),\\
			(\mathcal R_N^0)_I&:=\sigma^{-2}\mathsf T_N\partial_I\left(\lambda^{-1}\mathsf T_N^{-1}E_{21}(\Delta \mathsf T_N)V_N-\sigma^2\frac{D\lambda}{\lambda}\mathsf T_N^{-1}E_{22}\mathsf T_NV_N\right).
		\end{aligned}
	\end{equation}
	The sum of the four terms in \eqref{Eq.energy_remainder_split}, $\mathcal R_N^{\rm mov}$ and $\mathcal R_N[\mathbf g_0,-Q_\sigma]$ equals \eqref{Eq.exact_energy_remainder}.
	
	Where $\chi_0(N/(\varrho\sigma))=1$, the matrix $Z_N$ in \eqref{Eq.moving_gauge_remainder} vanishes. The identity \eqref{Eq.frame_identity} and $D\theta_\sigma=\lambda$ give
\begin{equation}\label{Eq.energy_low_principal_matrix}
A_N=\mathsf P_0(\sigma)+(\lambda-1)P_0(\sigma)N_0P_0(\sigma)^{-1}.
\end{equation}
The matrix $P_0(\sigma)$ is constant in $x$. For every ordered spatial derivative $\partial_I$ of length $r\geq1$, differentiating \eqref{Eq.energy_low_principal_matrix} and using $\mathsf T_N=P_0(\theta_\sigma)P_0(\sigma)^{-1}$ therefore yields
\begin{equation}\label{Eq.energy_low_principal_derivative}
\mathsf T_N(\partial_I A_N)\mathsf T_N^{-1}=(\partial_I\lambda)P_0(\theta_\sigma)N_0P_0(\theta_\sigma)^{-1}.
\end{equation}
Thus no spatial derivative of the rapidly oscillating matrix $\mathsf P_0(\theta_\sigma)$ remains. The coefficient $\partial_I\lambda$ has the extra factor $\sigma^{-1}\nu(\sigma)$ by \eqref{Eq.energy_gauge_coefficients}. We must still measure the matrix on the right of \eqref{Eq.energy_low_principal_derivative} in the possibly degenerate energy.
	For $R=P_{\xi^2}(\theta)^{-1}P_0(\theta)=\left(\begin{smallmatrix}a&b\\c&d\end{smallmatrix}\right)$, Lemma \ref{Lem.low_frame} gives $R=I+\Order(\xi^2)$, and direct multiplication yields
	\[RN_0R^{-1}=\frac1{\det R}\begin{pmatrix}-ac&a^2\\-c^2&ac\end{pmatrix}.\]
	Writing $\kappa=(\omega(\xi^2)+\delta^2)^{1/2}$, the entries of $\diag(\kappa,1)RN_0R^{-1}\diag(\kappa^{-1},1)$ satisfy
	\[\left\|\operatorname{diag}(\kappa,1)RN_0R^{-1}\operatorname{diag}(\kappa^{-1},1)\right\| \leq C\left(\xi^2+\kappa+\frac{\xi^4}{\kappa}\right)\leq C(\xi+\delta).\]
	The lower-left entry contains $c^2=\Order(\xi^4)$, rather than merely $\Order(\xi^2)$; this follows from the explicit conjugation of $N_0$ and is what makes the factor $\kappa^{-1}$ harmless. Combining this bound with \eqref{Eq.energy_low_principal_derivative}, $\nabla^rA_N$, $r\geq1$, contributes $C\nu\sigma^{-1}(\xi_J+\delta)r!R_c^{-r}$ at an input level $J$, up to the polynomial comparison factors of Lemma \ref{Lem.matrix_comparisons}. Where $\chi_0(N/(\varrho\sigma))\ne1$, we have $N\geq c\sigma$; the differentiated entries of $A_N$ are bounded by $C\nu r!R_c^{-r}$, and hence by $C\nu\sigma^{-1}(1+N)r!R_c^{-r}$. The cutoff derivative has the smaller bound $C\nu\sigma^{-1}r!R_c^{-r}$. Using Lemma \ref{Lem.gauge_transfer} and enlarging a summable sequence $b_r\leq C(1+r)^H(\varrho/R_c)^r$, we obtain
	\[\left|2\operatorname{Re}\sum_Na_N^2(S_N\mathcal U_N,\mathcal R_N^A)_{L^2}\right| \leq\frac{C\nu}\sigma\sum_{r\geq1}b_r\sum_{N\geq r}(1+N)e_Ne_{N-r}\leq\frac{C\nu}\sigma(E+\mathfrak d).\]
	In the last estimate, $1+N\leq(1+r)(1+N-r)$ and the additional factor $1+r$ is absorbed into $b_r$.
	
	For $\mathcal R_N^K$, we define $J=N-r+2$. The differentiated coefficient has a balanced bound $C\nu r!R_c^{-r}/(\xi_J+\delta)$, with the same type of polynomial comparison factors. The factorial shift gives the factor $J(J-1)$, and
	\begin{equation}\label{Eq.energy_derivative_costs}
		\frac{\sigma^{-2}J(J-1)}{\xi_J+\delta} =\frac{J(J-1)}{\sigma(1+J/(\varrho\lambda))}\leq\frac{C J}{\sigma},\quad\frac{\sigma^{-2}J}{\xi_J+\delta}\leq\frac C\sigma,\quad \frac{\sigma^{-2}}{\xi_J+\delta}\leq\frac1\sigma.
	\end{equation}
	For example, using \eqref{Eq.energy_factorial_shift} with $q=2$ and $J=N-r+2$, the contribution of $\mathcal R_N^K$ is bounded by
\begin{align*}
&\left|2\operatorname{Re}\sum_Na_N^2(S_N\mathcal U_N,\mathcal R_N^K)_{L^2}\right|\leq C\nu(\sigma)\sum_{r\geq1}b_r\sum_{N\geq r}\frac{\sigma^{-2}J(J-1)}{\xi_J+\delta}e_Ne_J\\
&\qquad\leq\frac{C\nu(\sigma)}\sigma\sum_{r\geq1}b_r\sum_{N\geq r}(1+J)e_Ne_J\leq\frac{C\nu(\sigma)}\sigma(E+\mathfrak d),
\end{align*}
where $b_r\leq C(1+r)^H(\varrho/R_c)^r$ includes the coefficient and frame comparisons. The last inequality is \eqref{Eq.energy_sequence_sums}. Thus the two derivatives in $\Delta V_N$ leave only one derivative-order weight after the factor $(\xi_J+\delta)^{-1}$ is used; the resulting sum is exactly the one controlled by outgoing transport. The three bounds in \eqref{Eq.energy_derivative_costs} estimate $\mathcal R_N^K$, $\mathcal R_N^1$ and the $\Delta \mathsf T_N$ term of $\mathcal R_N^0$ in \eqref{Eq.energy_remainder_split}, respectively. For the remaining part, \eqref{Eq.rank_one_metric} and $\sigma^2D\lambda=\Order_{\rm an}(\nu)$ give
	\[\sigma^{-2}\nu\left(1+\frac1{\xi_J+\delta}\right)\leq\frac{C\nu}\sigma.\]
	Each normalized sum is therefore bounded by 
\[C\nu\sigma^{-1}\sum_{q=0}^2\sum_{r\geq0}b_r\sum_{N\geq r}(1+N-r+q)e_Ne_{N-r+q},\]
 hence by $C\nu(E+\mathfrak d)/\sigma$, using Lemma \ref{Lem.energy_shifts}. We obtain the first inequality in \eqref{Eq.energy_remainder_bounds}.
	
	Finally, \eqref{Eq.scalar_source_matrix_bound} and Cauchy--Schwarz give
	\[2\left|\operatorname{Re}\sum_Na_N^2(S_N\mathcal U_N,\mathcal F_N)_{L^2}\right| \leq2\sqrt E\left(\sum_Na_N^2\|S_N^{1/2}\mathcal F_N\|_{L^2}^2\right)^{1/2} \leq C\sqrt E\|\mathcal F\|_{\mathcal B_\varrho}.\]
	For solutions with a strict additional analytic radius, the derivative sequences have exponentially decreasing tails, which justify all rearrangements and integrations.
\end{proof}

We now estimate the transport commutators. The terms with $\nabla\bvec$ contribute $-2\gamma_v\mathfrak d/\sigma$. The terms containing $\nabla^r\bvec$, $r\geq2$, are bounded by $2\mathcal L_v(\varrho)\mathfrak d/\sigma$ and lower-order energy terms.

\begin{lemma}\label{Lem.transport_energy}
	For $\mathscr B_N$ and $\mathcal C_N$ defined in \eqref{Eq.energy_transport_tensor}--\eqref{Eq.energy_transport_higher},
	\begin{equation}\label{Eq.energy_transport_bounds}
		\begin{aligned}
			-\frac2\sigma\operatorname{Re}\sum_{N\geq0}a_N^2 (S_N\mathcal U_N,\mathscr B_N\mathcal U_N)_{L^2} &\leq-\frac{2\gamma_v}{\sigma}\mathfrak d,\\
			\left|\frac2\sigma\operatorname{Re}\sum_{N\geq0}a_N^2 (S_N\mathcal U_N,\mathcal C_N)_{L^2}\right| &\leq\frac{2\mathcal L_v(\varrho)}\sigma\mathfrak d +\frac C\sigma E+\frac{C\nu(\sigma)}\sigma(E+\mathfrak d).
		\end{aligned}
	\end{equation}
\end{lemma}

\begin{proof}
	We define $V_I=S_N^{1/2}(\mathcal U_N)_I$ for $|I|=N$. Since $S_N$ acts only on the two components of the unknowns,
	\begin{align*}
		-\frac2\sigma\operatorname{Re}\sum_Na_N^2(S_N\mathcal U_N,\mathscr B_N\mathcal U_N)_{L^2} &=-\frac2\sigma\operatorname{Re}\sum_Na_N^2\sum_{a=1}^N\sum_I\int_\Omega\overline V_I\cdot(\partial_{i_a}b_k)V_{i_1\cdots k\cdots i_N}\,\mathrm dx\\
		&\leq-\frac{2\gamma_v}{\sigma}\sum_NNa_N^2\sum_I\|V_I\|_{L^2}^2 =-\frac{2\gamma_v}{\sigma}\mathfrak d.
	\end{align*}
	For each $a$, fix all indices except $i_a$ and apply $\operatorname{sym}\mathrm D_x\bvec\geq\gamma_vI_d$ from \eqref{Eq.acoustic_outgoing_property}; summing over $a=1,\ldots,N$ produces the factor $N$.
	
	For the higher commutator set $J=N-r+1$ and $c_r=\varrho^{r-1}\|\nabla^r\bvec\|_{L^\infty}/r!$. After replacing $\mathsf T_N\nabla^J(\mathsf T_N^{-1}U)$ by $\mathcal U_J$, Lemmas \ref{Lem.energy_shifts} and \ref{Lem.matrix_comparisons} give the bound
	\[\frac2\sigma\sum_{r\geq2}c_r\sum_{J\geq1}J\left(\frac{J+r}{J+1}\right)^{s_*}\left(1+\frac{C(1+r)^H}{J}\right)e_{J+r-1}e_J.\]
	Increasing the fixed integer $H$ if necessary, we have
	\[J\left(\frac{J+r}{J+1}\right)^{s_*}\left(1+\frac{C(1+r)^H}{J}\right)\leq J+C'(1+r)^{H+s_*}.\]
	Consequently,
	\begin{align*}
		&\sum_{J\geq1}Je_{J+r-1}e_J \leq\left(\sum_{J\geq1}(J+r-1)e_{J+r-1}^2\right)^{1/2}\left(\sum_{J\geq1}Je_J^2\right)^{1/2}\leq\mathfrak d,\quad \sum_{J\geq1}e_{J+r-1}e_J\leq E,\\
		&\qquad\qquad\qquad\frac2\sigma\sum_{r\geq2}c_r\sum_{J\geq1}\bigl[J+C'(1+r)^{H+s_*}\bigr]e_{J+r-1}e_J \leq\frac{2\mathcal L_v(\varrho)}\sigma\mathfrak d+\frac C\sigma E.
	\end{align*}
	By \eqref{Eq.geometry_budget}, $\sum_{r\geq2}c_r=\mathcal L_v(\varrho)$; the coefficient bound \eqref{Eq.common_coefficient_radius} also gives $\sum_{r\geq2}(1+r)^{H+s_*}c_r\leq C$. The replacement error in Lemma \ref{Lem.gauge_transfer} has an additional factor $\nu$. The coefficients of the replacement error are bounded by a summable sequence $b_r$ times $(1+J)\nu/\sigma$. Therefore Lemma \ref{Lem.energy_shifts} gives
	\[\frac{C\nu}\sigma\sum_{r\geq0}b_r\sum_{N\geq r}(1+N-r+1)e_Ne_{N-r+1}\leq\frac{C\nu}\sigma(E+\mathfrak d).\]
	Adding the bounds for the unchanged frame and the frame difference gives
	\[\left|\frac2\sigma\operatorname{Re}\sum_Na_N^2(S_N\mathcal U_N,\mathcal C_N)_{L^2}\right| \leq\frac{2\mathcal L_v(\varrho)}\sigma\mathfrak d +\frac C\sigma E+\frac{C\nu(\sigma)}\sigma(E+\mathfrak d).\qedhere\]
\end{proof}

The equation \eqref{Eq.differentiated_energy_equation} contains $\lambda\mathsf P_0$ and a Laplacian, whereas \eqref{Eq.frozen_original_defect} concerns $\lambda\mathsf P_{\xi_N}=\lambda\mathsf P_0-\lambda\xi_N^2E_{21}$. Integration by parts and the factorial relation \eqref{Eq.energy_laplacian_shift} produce precisely this frequency term, up to \eqref{Eq.metric_difference} and the boundary integral defined below.

\begin{lemma}[The Laplacian cancellation]\label{Lem.laplacian_energy}
	Define
	\begin{equation}\label{Eq.energy_laplacian_forms}
		\begin{aligned}
			\mathcal I_\Delta&:=2\sigma^{-2}\operatorname{Re}\sum_{N\geq0}a_N^2 (S_N\mathcal U_N,\lambda^{-1}E_{21}\Delta^{[N]}\mathcal U_N)_{L^2},\\
			\mathcal I_\xi&:=-2\operatorname{Re}\sum_{N\geq0}a_N^2 (S_N\mathcal U_N,\lambda\xi_N^2E_{21}\mathcal U_N)_{L^2},\\
			\mathcal B_\Delta&:=2\sigma^{-2}\operatorname{Re}\sum_{N\geq0}a_N^2 \int_{\partial\Omega}\lambda^{-1}(e_2^TS_N\mathcal U_N)^* e_1^T(n_i\mathcal U_{N+1,i})\,\mathrm dS.
		\end{aligned}
	\end{equation}
	Then
	\begin{equation}\label{Eq.laplacian_form_error}
		|\mathcal I_\Delta-\mathcal I_\xi-\mathcal B_\Delta| \leq\frac C\sigma E+\frac{C\nu(\sigma)}\sigma(E+\mathfrak d).
	\end{equation}
	The additional boundary terms generated by changing $\nabla^{[N]}\mathcal U_N$ to $\mathcal U_{N+1}$ have total absolute value at most $C\nu(\sigma)(E+\mathfrak d)/\sigma$.
\end{lemma}

\begin{proof}
	For $Z=\mathcal U_N$, covariant integration by parts reads
	\begin{align*}
		&\int_\Omega Z^*\lambda^{-1}S_NE_{21}\Delta^{[N]}Z\,\mathrm dx =\int_{\partial\Omega}Z^*\lambda^{-1}S_NE_{21}\nabla_n^{[N]}Z\,\mathrm dS\\
		&\qquad\qquad-\sum_i\int_\Omega(\nabla_i^{[N]}Z)^*\lambda^{-1}S_NE_{21}\nabla_i^{[N]}Z\,\mathrm dx-\sum_i\int_\Omega Z^*L_{N,i}\nabla_i^{[N]}Z\,\mathrm dx,
	\end{align*}
	where
	\begin{equation}\label{Eq.energy_laplacian_coefficient}
		L_{N,i}:=\lambda^{-1}\mathcal Q_{N,i}E_{21}+\lambda^{-1}S_N[E_{21},\Gamma_{N,i}]+(\partial_i\lambda^{-1})S_NE_{21}.
	\end{equation}
	The coefficient \eqref{Eq.energy_laplacian_coefficient} follows by substituting $\partial_iZ=\nabla_i^{[N]}Z+\Gamma_{N,i}Z$ and $\partial_iS_N=\mathcal Q_{N,i}-\Gamma_{N,i}^TS_N-S_N\Gamma_{N,i}$ in the ordinary integration-by-parts formula. We estimate the three terms of $L_{N,i}$ separately. For the term containing $\mathcal Q_{N,i}$, \eqref{Eq.energy_metric_derivative_bound} and \eqref{Eq.rank_one_metric} give
\begin{align*}
\|S_N^{-1/2}\mathcal Q_{N,i}E_{21}S_N^{-1/2}\|
&\leq\|S_N^{-1/2}\mathcal Q_{N,i}S_N^{-1/2}\|\|S_N^{1/2}E_{21}S_N^{-1/2}\|\\
&\leq C\nu(\sigma)\frac{\min\{1,\xi_N+\delta\}}{\xi_N+\delta}\leq C\nu(\sigma).
\end{align*}
The small factor in the spatial metric derivative cancels the singular factor in the rank-one matrix norm. For the commutator term, $\Gamma_{N,i}$ vanishes outside the cutoff, while on the cutoff support
\[\|S_N^{1/2}\|\leq C,\qquad \|S_N^{-1/2}\|\leq C(\xi_N+\delta)^{-1},\qquad \|\Gamma_{N,i}\|\leq C\nu(\sigma).\]
Hence $\|S_N^{1/2}[E_{21},\Gamma_{N,i}]S_N^{-1/2}\|\leq C\nu(\sigma)/(\xi_N+\delta)$. Finally, $|\partial_i\lambda^{-1}|\leq C\nu(\sigma)/\sigma$ by \eqref{Eq.energy_gauge_coefficients}. We conclude that
\begin{equation}\label{Eq.laplacian_coefficient_bound}
\|S_N^{-1/2}L_{N,i}S_N^{-1/2}\|\leq C\nu(\sigma)\left(1+(\xi_N+\delta)^{-1}\right).
\end{equation}
The adjacent matrix comparison in Lemma \ref{Lem.sampled_matrices} and $a_N/a_{N+1}\leq(N+1)/\varrho$ now give, for the term with $\mathcal U_{N+1,i}$,
\begin{align*}
&\sigma^{-2}\sum_{N,i}a_N^2\left|\int_\Omega\mathcal U_N^*L_{N,i}\mathcal U_{N+1,i}\,\mathrm dx\right|\\
&\qquad\leq C\nu(\sigma)\sum_N\frac{N+1}{\sigma^2}\left(1+(\xi_N+\delta)^{-1}\right)e_Ne_{N+1}\leq\frac{C\nu(\sigma)}\sigma(E+\mathfrak d).
\end{align*}
Indeed, $\delta=\sigma^{-1}$ and $\xi_N=N/(\varrho\sigma\lambda)$ imply
\[\frac{N+1}{\sigma^2(\xi_N+\delta)}=\frac{N+1}{\sigma(1+N/(\varrho\lambda))}\leq\frac C\sigma,\qquad \sum_N(N+1)e_Ne_{N+1}\leq C(E+\mathfrak d).\]
The correction $\nabla_i^{[N]}\mathcal U_N-\mathcal U_{N+1,i}$ is treated by \eqref{Eq.gauge_transfer_bound} and \eqref{Eq.energy_sequence_sums}, retaining the bound \eqref{Eq.laplacian_coefficient_bound}. The resulting coefficient sequence is summable and has the extra factor $\nu/\sigma$ from \eqref{Eq.energy_adjacent_gauges}; its contribution satisfies the same bound $C\nu(E+\mathfrak d)/\sigma$.
	
	In the quadratic term $(\nabla_i^{[N]}\mathcal U_N)^*S_NE_{21}\nabla_i^{[N]}\mathcal U_N$, replace $\nabla_i^{[N]}\mathcal U_N$ by $\mathcal U_{N+1,i}$ using Lemma \ref{Lem.gauge_transfer}. The difference $\nabla_i^{[N]}\mathcal U_N-\mathcal U_{N+1,i}$ contains $\nu/\sigma$ before differentiation. Apply \eqref{Eq.energy_derivative_costs} and \eqref{Eq.energy_sequence_sums} to the adjacent-gauge estimate \eqref{Eq.gauge_transfer_bound}; the resulting error is bounded by $C\nu(E+\mathfrak d)/\sigma$. Thus the principal interior term is
	\[-2\sigma^{-2}\operatorname{Re}\sum_Na_N^2\int_\Omega\lambda^{-1}\mathcal U_{N+1}^*S_NE_{21}\mathcal U_{N+1}\,\mathrm dx.\]
	Replace $N+1$ by $N$ in the principal interior integral. The exact coefficient relation is
	\begin{equation}\label{Eq.energy_laplacian_shift}
		\sigma^{-2}\lambda^{-1}(a_{N-1})^2=\lambda\xi_N^2(a_N)^2\left(\frac N{N+1}\right)^{2s_*}\quad(N\geq1).
	\end{equation}
	The change from $t_NS_{N-1}$ to $S_N$ must be made before estimating the frequency matrix $E_{21}$. Indeed,
\begin{equation}\label{Eq.laplacian_relative_factorization}
\begin{aligned}
&S_N^{-1/2}(t_NS_{N-1}-S_N)E_{21}S_N^{-1/2}\\
&\quad=\left[S_N^{-1/2}(t_NS_{N-1}-S_N)S_N^{-1/2}\right]\left[S_N^{1/2}E_{21}S_N^{-1/2}\right].
\end{aligned}
\end{equation}
By \eqref{Eq.metric_difference} and \eqref{Eq.rank_one_metric}, the norm in \eqref{Eq.laplacian_relative_factorization} is at most $C/[N(\xi_N+\delta)]$. Hence replacing $t_NS_{N-1}$ by $S_N$ in the principal interior integral has error bounded by
	\begin{align*}
		C\sum_{N\geq1}a_N^2\int_\Omega\frac{\lambda\xi_N^2}{N(\xi_N+\delta)}|S_N^{1/2}\mathcal U_N|^2\,\mathrm dx \leq\frac C\sigma\sum_{N\geq1}a_N^2\int_\Omega\frac{\xi_N}{\varrho(\xi_N+\delta)}|S_N^{1/2}\mathcal U_N|^2\,\mathrm dx\leq\frac C\sigma E,
	\end{align*}
	by \eqref{Eq.metric_difference} and \eqref{Eq.rank_one_metric}. The frequency contribution at $N=0$ vanishes since $\xi_0=0$.
	
	Replacing $\nabla_i^{[N]}\mathcal U_N$ by $\mathcal U_{N+1,i}$ in the boundary integral gives $\mathcal B_\Delta$. By Lemma \ref{Lem.gauge_transfer} and \eqref{Eq.rank_one_metric}, the remaining boundary integral $\mathcal R_\partial$ satisfies
	\[|\mathcal R_\partial|\leq\frac{C\nu}{\sigma}\sum_Na_N^2\int_{\partial\Omega}\mathcal U_N^*S_N\mathcal U_N\,\mathrm dS\leq\frac{C\nu}{\sigma}(E+\mathfrak d),\]
	where the last inequality is \eqref{Eq.energy_trace_bound}. Combining the interior and boundary bounds gives \eqref{Eq.laplacian_form_error}. This completes the proof.
\end{proof}

\subsection{The characteristic boundary}
\label{Subsec.energy_boundary}

The boundary term in \eqref{Eq.total_boundary_estimate} combines the negative transport integral with the Laplacian integral $\mathcal B_\Delta$ defined in \eqref{Eq.energy_laplacian_forms}. We estimate the combined flux before taking absolute values. The cross-level bound \eqref{Eq.sampled_boundary} preserves the leading constant one in the boundary square; \eqref{Eq.energy_trace_bound} controls the remaining error of order $\sigma^{-2}$.

\begin{lemma}\label{Lem.energy_trace}
	For a solution satisfying \eqref{Eq.energy_smooth_core}, with $a_N,E,\mathfrak d$ from \eqref{Eq.energy_level_sums},
	\begin{equation}\label{Eq.energy_trace_bound}
		\sum_N(a_N)^2\int_{\partial\Omega}\mathcal U_N^*S_N\mathcal U_N\,\mathrm dS\leq C(E+\mathfrak d),
	\end{equation}
	with $C$ independent of $\sigma\geq\Sigma_0$.
\end{lemma}

\begin{proof}
	Since $\Omega=B_{X_1}(0)$, the divergence theorem applied to $x(\mathcal U_N^*S_N\mathcal U_N)$ gives
	\begin{align*}
		X_1\int_{\partial\Omega}\mathcal U_N^*S_N\mathcal U_N\,\mathrm dS ={}&d\int_\Omega\mathcal U_N^*S_N\mathcal U_N\,\mathrm dx\\
		&+2\operatorname{Re}\sum_i\int_\Omega x_i\mathcal U_N^*S_N\nabla_i^{[N]}\mathcal U_N\,\mathrm dx +\sum_i\int_\Omega x_i\mathcal U_N^*\mathcal Q_{N,i}\mathcal U_N\,\mathrm dx.
	\end{align*}
	By \eqref{Eq.energy_metric_derivative}--\eqref{Eq.energy_metric_derivative_bound}, we have $\left|\sum_{N,i}a_N^2\int_\Omega x_i\mathcal U_N^*\mathcal Q_{N,i}\mathcal U_N\,\mathrm dx\right|\leq C\nu E$.
	
	Replace $\nabla_i^{[N]}\mathcal U_N$ by $\mathcal U_{N+1,i}$. The adjacent metric comparison and $a_N/a_{N+1}\leq(N+1)/\varrho$ give
	\begin{align*}
		&\left|2\operatorname{Re}\sum_{N,i}a_N^2\int_\Omega x_i\mathcal U_N^*S_N\mathcal U_{N+1,i}\,\mathrm dx\right|\leq C\sum_{N\geq0}(N+1)e_Ne_{N+1}\\
		&\qquad\leq C\left(\sum_{N\geq0}(N+1)e_N^2\right)^{1/2} \left(\sum_{N\geq0}(N+1)e_{N+1}^2\right)^{1/2} \leq C(E+\mathfrak d).
	\end{align*}
	Set $V_{N,i}=\nabla_i^{[N]}\mathcal U_N-\mathcal U_{N+1,i}$. The product expansion \eqref{Eq.gauge_transfer_product} and the adjacent-cutoff bound in \eqref{Eq.energy_adjacent_gauges} give a factor $C\nu/\sigma$. The one-derivative factorial coefficient is $(N+1)\varrho^{r-1}/r!$. On the support of $\mathsf T_N-\mathsf T_{N+1}$, $N+1\asymp\sigma$. The two-metric estimate \eqref{Eq.gauge_two_metric} therefore gives a sequence $b_r\leq C(1+r)^H(\varrho/R_c)^r$ such that
	\begin{align*}
		\left|2\operatorname{Re}\sum_{N,i}a_N^2\int_\Omega x_i\mathcal U_N^*S_NV_{N,i}\,\mathrm dx\right|\leq C\nu\sum_{r\geq1}b_r\sum_{N\geq r-1}e_Ne_{N+1-r} \leq C\nu\left(\sum_{r\geq1}b_r\right)E\leq C\nu E.
	\end{align*}
	For indices outside a common positive-frequency interval, the polynomial comparison factor is absorbed into $(1+r)^H$ exactly as in \eqref{Eq.gauge_two_metric}. The bounds for $\mathcal Q_{N,i}$, $\mathcal U_{N+1,i}$ and $V_{N,i}$ give
	\[X_1\sum_Na_N^2\int_{\partial\Omega}\mathcal U_N^*S_N\mathcal U_N\,\mathrm dS \leq C E+C(E+\mathfrak d)+C\nu E\leq C(E+\mathfrak d).\]
	Division by $X_1>0$ proves the lemma.
\end{proof}

The trace estimate controls an error of order $\sigma^{-2}$. To obtain that error, we combine the exact boundary fluxes and apply Lemma \ref{Lem.cross_boundary} before taking absolute values.

\begin{lemma}[The sign of the principal boundary form]\label{Lem.boundary_energy}
	For the boundary integral $\mathcal B_\Delta$ in \eqref{Eq.energy_laplacian_forms},
	\begin{equation}\label{Eq.total_boundary_estimate}
		-\sigma^{-1}\sum_{N\geq0}a_N^2 \int_{\partial\Omega}\mathcal U_N^*S_N\mathcal U_N\,\mathrm dS +\mathcal B_\Delta \leq C\sigma^{-2}(E+\mathfrak d).
	\end{equation}
\end{lemma}

\begin{proof}
	For every ordered tensor component write $\mathcal U_N=(u_N^\partial,w_N^\partial)^T$ and set
	\begin{align*}
		A_N^\partial=a_N\sqrt{(S_N)_{22}}\left(w_N^\partial+\frac{(S_N)_{12}}{(S_N)_{22}}u_N^\partial\right),\quad B_N^\partial=a_N\sqrt{\frac{\det S_N}{(S_N)_{22}}}u_N^\partial.
	\end{align*}
	Completing the square gives
	\[a_N^2\mathcal U_N^*S_N\mathcal U_N=|A_N^\partial|^2+|B_N^\partial|^2.\]
	Use $\bvec\cdot n=1$ from \eqref{Eq.acoustic_outgoing_property}, the definition of $\mathcal B_\Delta$ in \eqref{Eq.energy_laplacian_forms}, and $a_N/a_{N+1}=((N+1)/(N+2))^{s_*}(N+1)/\varrho$. The left-hand side of \eqref{Eq.total_boundary_estimate} is
	\[\mathfrak B=\frac1\sigma\int_{\partial\Omega}\left[-\sum_{N\geq0}(|A_N^\partial|^2+|B_N^\partial|^2) +2\operatorname{Re}\sum_{N\geq0}\rho_N(A_N^\partial)^*n_iB_{N+1,i}^\partial\right]\,\mathrm dS,\]
	where
	\[\rho_N=\left(\frac{N+1}{N+2}\right)^{s_*}\xi_{N+1} \sqrt{\frac{(S_N)_{22}(S_{N+1})_{22}}{\det S_{N+1}}}\leq1+\frac C\sigma.\]
	Apply \eqref{Eq.sampled_boundary} pointwise with $\theta_*=\theta_\sigma$ and $c_*=1/\lambda\in[1/2,2]$; the estimate is uniform in both parameters and gives the bound for $\rho_N$. Since $\sum_i n_i^2=1$,
	\[|n_iB_{N+1,i}^\partial|^2\leq\sum_i|B_{N+1,i}^\partial|^2=|B_{N+1}^\partial|^2.\]
	The constant one in \eqref{Eq.sampled_boundary} is needed before applying the trace estimate. To see the sign, replace $\rho_N$ temporarily by one. Pointwise on $\partial\Omega$, completing the square and shifting the derivative index give
\begin{align*}
&-\sum_{N\geq0}(|A_N^\partial|^2+|B_N^\partial|^2)+2\operatorname{Re}\sum_{N\geq0}(A_N^\partial)^*n_iB_{N+1,i}^\partial\\
=&\,-|B_0^\partial|^2-\sum_{N\geq0}|A_N^\partial-n_iB_{N+1,i}^\partial|^2-\sum_{N\geq0}\bigl(|B_{N+1}^\partial|^2-|n_iB_{N+1,i}^\partial|^2\bigr)\leq0.
\end{align*}
The last sum is nonnegative before the minus sign, by $|n_iB_{N+1,i}^\partial|\leq|B_{N+1}^\partial|$. For the actual coefficients, $0\leq\rho_N\leq1+C/\sigma$ is sufficient; no lower bound close to one is needed. Cauchy--Schwarz therefore yields
\begin{align*}
&\qquad2\operatorname{Re}\sum_{N\geq0}\rho_N(A_N^\partial)^*n_iB_{N+1,i}^\partial \leq\left(1+\frac C\sigma\right)\sum_{N\geq0}(|A_N^\partial|^2+|B_{N+1}^\partial|^2),\\
		&\mathfrak B \leq-\frac1\sigma\int_{\partial\Omega}|B_0^\partial|^2\,\mathrm dS +\frac C{\sigma^2}\int_{\partial\Omega}\sum_{N\geq0}(|A_N^\partial|^2+|B_N^\partial|^2)\,\mathrm dS\leq\frac C{\sigma^2}(E+\mathfrak d),
	\end{align*}
	where the last inequality is \eqref{Eq.energy_trace_bound}.
\end{proof}

We now differentiate $E$ and combine \eqref{Eq.energy_remainder_bounds}, \eqref{Eq.laplacian_form_error} and \eqref{Eq.total_boundary_estimate} with the transport estimate in Lemma \ref{Lem.transport_energy}.

\begin{lemma}\label{Lem.moving_energy_identity}
	For solutions satisfying \eqref{Eq.energy_smooth_core}, the quantities $E,\mathfrak d$ in \eqref{Eq.energy_level_sums} satisfy
	\begin{equation}\label{Eq.local_energy_budget}
		\frac{\mathrm dE}{\mathrm d\sigma}+\frac2\sigma \left[\gamma_v-\mathcal L_v(\varrho)-\frac\chi\varrho-\alpha -C\bigl(\sigma^{-1}+\nu(\sigma)\bigr)\right]\mathfrak d \leq\frac C\sigma E+C\sqrt E\|\mathcal F\|_{\mathcal B_\varrho}.
	\end{equation}
\end{lemma}

\begin{proof}
	For $Z=\mathcal U_N$, use $\partial_\sigma a_N=\alpha Na_N/\sigma$ from \eqref{Eq.energy_level_sums} and integrate $B=\bvec\cdot\nabla$ by parts, with $\bvec\cdot n=1$ from \eqref{Eq.acoustic_outgoing_property}. We obtain
	\begin{align*}
		&\frac{\mathrm d}{\mathrm d\sigma}\left(a_N^2\int_\Omega Z^*S_NZ\,\mathrm dx\right) =\frac{2\alpha N}{\sigma}a_N^2\int_\Omega Z^*S_NZ\,\mathrm dx +2a_N^2\operatorname{Re}\int_\Omega Z^*S_NDZ\,\mathrm dx\\
		&\qquad\qquad+a_N^2\int_\Omega Z^*(DS_N)Z\,\mathrm dx +\frac{a_N^2}{\sigma}\int_\Omega(\operatorname{div}\bvec)Z^*S_NZ\,\mathrm dx-\frac{a_N^2}{\sigma}\int_{\partial\Omega}Z^*S_NZ\,\mathrm dS.
	\end{align*}
	Insert \eqref{Eq.differentiated_energy_equation} and sum in $N$. Lemma \ref{Lem.laplacian_energy} changes the interior principal matrix from $\mathsf P_0$ to $\mathsf P_{\xi_N}=\mathsf P_0-\xi_N^2E_{21}$. The associated matrix contribution is bounded by
	\[\sum_Na_N^2\int_\Omega Z^*\left[DS_N+\lambda(\mathsf P_{\xi_N}^TS_N+S_N\mathsf P_{\xi_N})\right]Z\,\mathrm dx \leq\frac{2\chi}{\varrho\sigma}\mathfrak d+\frac C\sigma E\]
	by \eqref{Eq.frozen_original_defect}. Lemmas \ref{Lem.transport_energy}, \ref{Lem.energy_remainders} and \ref{Lem.boundary_energy}, together with the Laplacian error, now give
	\begin{align*}
		E'\leq{}&\frac{2\alpha}{\sigma}\mathfrak d+\frac{2\chi}{\varrho\sigma}\mathfrak d -\frac{2\gamma_v}{\sigma}\mathfrak d+\frac{2\mathcal L_v(\varrho)}\sigma\mathfrak d +\frac C\sigma E\\
		&\qquad+C\left(\frac{\nu(\sigma)}\sigma+\frac1{\sigma^2}\right)(E+\mathfrak d) +C\sqrt E\|\mathcal F\|_{\mathcal B_\varrho}.
	\end{align*}
	Since $\nu\leq1$, the terms multiplying $E$ are bounded by $CE/\sigma$. Move the terms multiplying $\mathfrak d$ to the left and enlarge $C$ to obtain \eqref{Eq.local_energy_budget}. All identities hold first for solutions with a strict additional analytic radius. Absolute convergence then permits termwise differentiation and summation; Lemma \ref{Lem.energy_completion} gives the finite-energy extension.
\end{proof}

\subsection{Construction of the forward solution}
\label{Subsec.forward_energy}

We prove $\|W\|_{\cX_\varrho}^2\lesssim\cE^\sharp_\varrho\lesssim\|W\|_{\cX_\varrho}^2+\|\nabla W\|_{\cX_\varrho}^2$ in Lemma \ref{Lem.energy_comparison}, so that the estimate \eqref{Eq.local_energy_budget} controls the original modulation. At low derivative orders, the reconstruction \eqref{Eq.zeta-w-W-def} gives a direct comparison; at high orders, Lemma \ref{Lem.sequence_tail} controls the lower-order terms produced by $\mathsf T_N$. We then solve \eqref{Eq.natural_forward_system} for $(\zeta,D\zeta,\sigma^{-1}\nabla\zeta)$ and use approximation to extend \eqref{Eq.energy_estimate} to finite-energy data.

\begin{lemma}\label{Lem.sequence_tail}
	Let $t\in(0,1)$, $n_*\in\Z_{\geq1}$, and let $k=(k_r)_{r\geq0}$ be nonnegative with $k_0=0$ and $\sum_{r\geq1}t^{-r}k_r\leq1/2$. If nonnegative sequences $X,e\in\ell^2(\Z_{\geq0})$ satisfy
	\[X_N\leq e_N+\sum_{r=1}^Nk_rX_{N-r}\qquad(N\geq n_*),\]
	then
	\[\left(\sum_{N\geq2n_*}X_N^2\right)^{1/2}\leq2\|e\|_{\ell^2}+2t^{n_*}\left(\sum_{N<n_*}X_N^2\right)^{1/2}.\]
\end{lemma}

\begin{proof}
	Set $Y_N:=\mathbf1_{N\geq n_*}X_N$ and $Z_N:=\mathbf1_{N<n_*}X_N$. Componentwise, $Y\leq e+k*Y+k*Z$. The convolution inverse is the nonnegative sequence
	\[r:=\delta^{(0)}+\sum_{a\geq1}k^{*a},\qquad \sum_{n\geq0}t^{-n}r_n\leq2,\qquad r*k=r-\delta^{(0)},\]
	where $\delta^{(0)}=(1,0,0,\ldots)$. Iteration gives $Y\leq r*e+(r*k)*Z$; the remaining iterate tends to zero in $\ell^2$ because $\|k\|_{\ell^1}\leq1/2$. For $N\geq2n_*$ and $J<n_*$, one has $N-J>n_*$. Thus $\|\mathbf1_{N\geq2n_*}((r*k)*Z)\|_{\ell^2}\leq\|\mathbf1_{n>n_*}(r*k)_n\|_{\ell^1}\|Z\|_{\ell^2}\leq2t^{n_*}\|Z\|_{\ell^2}$. Together with $\|r*e\|_{\ell^2}\leq2\|e\|_{\ell^2}$, we obtain the required bound for $\mathbf1_{N\geq2n_*}X_N$.
\end{proof}

The sequence estimate allows us to remove the small high-order frame error. We apply Lemma \ref{Lem.sequence_tail} to compare $\cE^\sharp_\varrho(\sigma;W)$ with the norm of $W$.

\begin{lemma}[The energy of an arbitrary pair]\label{Lem.energy_comparison}
	After one increase of $\Sigma_0$, uniformly for $\varrho\in I_{\rm r}$ and $\sigma\geq\Sigma_0$,
	\[c\|W\|_{\cX_\varrho}^2\leq\cE^\sharp_\varrho(\sigma;W) \leq C\bigl(\|W\|_{\cX_\varrho}^2+\|\nabla W\|_{\cX_\varrho}^2\bigr).\]
	Here $W\in\cX_\varrho$, and the upper bound is interpreted as $+\infty$ when $\nabla W\notin\cX_\varrho$.
\end{lemma}

\begin{proof}
	Fix $\sigma$ and define $\delta:=\sigma^{-1}$, $\mathsf D:=\diag(1,\delta)$, and $U=P_0(\theta_\sigma)\mathsf DW$. Set
	\[q_N:=a_N\|\na^NW\|_{L^2},\qquad X_N:=\sigma a_N\|\na^NU\|_{L^2},\qquad z_N:=\sigma e_N.\]
	Then $\|q\|_{\ell^2}=\|W\|_{\cX_\varrho}$ and $\|z\|_{\ell^2}^2=\cE^\sharp_\varrho$. The uniform multiplier bound for $P_0(\theta_\sigma)$ gives $\|X\|_{\ell^2}\leq C\sigma\|q\|_{\ell^2}$.
	
	Let $c_*:=\varrho\xi_*/16$ and $L_\sigma:=\lceil c_*\sigma\rceil$. Increase $\Sigma_0$ so that $4L_\sigma\leq\varrho\xi_*\sigma$. At levels $N\leq4L_\sigma$ the cutoff equals one, and
	\[\mathcal U_N=P_0(\theta_\sigma)\mathsf D\na^NW,\qquad R_N:=P_{\xi_N^2}^{-1}P_0=I+\Order(\xi_N^2).\]
	By \eqref{Eq.energy_weights}, the contribution $z_N^2$ equals $a_N^2\|M_N\nabla^NW\|_{L^2}^2$, where
	\[M_N:=\delta^{-1}\diag(\sqrt{\omega(\xi_N^2)+\delta^2},1)R_N\mathsf D.\]
	At $N=0$, $R_0=I_2$ and $M_0=I_2$, so the normalization in \eqref{Eq.energy_weights} gives exactly
	\[\sigma^2U^TS(\theta_\sigma,0,\sigma^{-1})U=|W|^2.\]
	For the other low derivative orders, $R_N=I_2+\Order(\xi_N^2)$ and $\sqrt{\omega(\xi_N^2)+\delta^2}\asymp\xi_N+\delta$ imply $\|M_N\|\leq C(1+N)$. To obtain the lower bound without losing a factor $\sigma$, write the inverse explicitly:
	\begin{equation}\label{Eq.energy_comparison_inverse}
		M_N^{-1}=\frac1{\det R_N}\begin{pmatrix}
			\dfrac{\delta(R_N)_{22}}{\sqrt{\omega(\xi_N^2)+\delta^2}}&-\delta(R_N)_{12}\\
			-\dfrac{(R_N)_{21}}{\sqrt{\omega(\xi_N^2)+\delta^2}}&(R_N)_{11}
		\end{pmatrix}.
	\end{equation}
	The determinant of $R_N$ stays away from zero. The off-diagonal entries of $R_N$ are $\Order(\xi_N^2)$, so the four entries of \eqref{Eq.energy_comparison_inverse} are bounded respectively by
	\[\frac{C\delta}{\xi_N+\delta},\qquad C\delta\xi_N^2,\qquad\frac{C\xi_N^2}{\xi_N+\delta},\qquad C.\]
	All four bounds are uniform in $\sigma$ in the low-frequency region. In particular, the factor $(\xi_N+\delta)^{-1}$ multiplies an off-diagonal entry vanishing to second order. Thus $\|M_N^{-1}\|\leq C$, and
	\begin{equation}\label{Eq.energy_low_comparison}
		cq_N\leq z_N\leq C(1+N)q_N\qquad(N\leq4L_\sigma).
	\end{equation}
	For $N\geq L_\sigma$, the frequency is bounded below by a fixed positive constant, so
	\[S_N\geq cI_2,\qquad \sigma\|S_N^{1/2}\|\leq C(1+N).\]
	The product formula for $\mathcal U_N=\mathsf T_N\na^N(\mathsf T_N^{-1}P_0(\theta_\sigma)\mathsf DW)$ therefore gives
	\[z_N\leq C\sum_{r=0}^Nb_r(1+N-r)q_{N-r}\qquad(N\geq L_\sigma),\]
	where $b_r\leq C(1+r)^H(\varrho/R_c)^r$ for fixed $H$ and $R_c>\varrho$. In passing from the Leibniz formula to this bound, the output factor $1+N$ is replaced by $(1+r)(1+N-r)$, with $1+r$ included in $b_r$. Since $\|b\|_{\ell^1}\leq C$, Young's inequality and \eqref{Eq.energy_low_comparison} yield
\[\|z\|_{\ell^2}\leq C\left(\sum_{N\geq0}(1+N)^2q_N^2\right)^{1/2}.\]
For $N\geq1$, the exact weight ratio is
\[Na_N=\varrho\left(\frac{N+1}{N}\right)^{s_*}a_{N-1}\leq 2^{s_*}\varrho a_{N-1}.\]
Using the ordered tensor identity $\|\nabla^NW\|_{L^2}=\|\nabla^{N-1}(\nabla W)\|_{L^2}$, we obtain
\[\sum_{N\geq0}(1+N)^2q_N^2\leq C\left(\|W\|_{\cX_\varrho}^2+\|\nabla W\|_{\cX_\varrho}^2\right).\]
This proves the upper bound without requiring a loss of analytic radius.
	
	To prove $\|W\|_{\cX_\varrho}^2\leq C\cE^\sharp_\varrho(\sigma;W)$, we may assume $\|z\|_{\ell^2}<\infty$. Expanding $\mathsf T_N\na^N(\mathsf T_N^{-1}U)$ and using Lemma \ref{Lem.gauge_bounds} yields
	\[X_N\leq Cz_N+\sum_{r=1}^Nk_rX_{N-r}\qquad(N\geq L_\sigma),\qquad k_r\leq C\nu(\sigma)(1+r)^H(\varrho/R_c)^r.\]
	Choose $\varrho/R_c<t<1$ and then increase $\Sigma_0$ so that $\sum_{r\geq1}t^{-r}k_r\leq1/2$. Lemma \ref{Lem.sequence_tail} gives
	\[\|\mathbf1_{N\geq2L_\sigma}X_N\|_{\ell^2}\leq C\|z\|_{\ell^2}+C\sigma t^{L_\sigma}\|q\|_{\ell^2}.\]
	The inverse identity $W=\mathsf D^{-1}P_0(\theta_\sigma)^{-1}U$ and the multiplier bound for $P_0^{-1}$ give
\[q_N\leq C\sum_{J=0}^Nb_{N-J}X_J,\qquad b_r\leq C(1+r)^H(\varrho/R_c)^r.\]
For $N\geq4L_\sigma$, separate $J\geq2L_\sigma$ from $J<2L_\sigma$. Increasing $t<1$ if necessary, while keeping $t>\varrho/R_c$, gives $\sum_{r\geq2L_\sigma}b_r\leq Ct^{2L_\sigma}$. Young's inequality on the two parts then yields
\begin{align*}
\|\mathbf1_{N\geq4L_\sigma}q_N\|_{\ell^2}
&\leq C\|\mathbf1_{J\geq2L_\sigma}X_J\|_{\ell^2}+Ct^{2L_\sigma}\|X\|_{\ell^2}\leq C\|z\|_{\ell^2}+C\sigma(t^{L_\sigma}+t^{2L_\sigma})\|q\|_{\ell^2}.
\end{align*}
Consequently,
	\[\|\mathbf1_{N\geq4L_\sigma}q_N\|_{\ell^2}\leq C\|z\|_{\ell^2}+C\sigma t^{L_\sigma}\|q\|_{\ell^2}.\]
	This estimate avoids inverting the reconstruction with a uniform bound on all derivative orders: that inverse alone would cost $\delta^{-1}=\sigma$. Instead, \eqref{Eq.energy_low_comparison} controls $N<4L_\sigma$ without loss, while any contribution from $N<2L_\sigma$ to $N\geq4L_\sigma$ has at least $2L_\sigma$ derivatives on an analytic coefficient. Its exponentially small coefficient tail absorbs the factor $\sigma$, since $L_\sigma\geq c\sigma$ uniformly for $\varrho\in I_{\rm r}$. Use \eqref{Eq.energy_low_comparison} for the remaining indices. Take $\Sigma_0$ so that $C\sigma t^{L_\sigma}\leq1/2$ for $\sigma\geq\Sigma_0$. Then $\|q\|_{\ell^2}\leq C\|z\|_{\ell^2}+\tfrac12\|q\|_{\ell^2}$, hence $\|W\|_{\cX_\varrho}\leq2C\sqrt{\cE^\sharp_\varrho(\sigma;W)}$.
\end{proof}

We now combine Lemmas \ref{Lem.energy_comparison} and \ref{Lem.moving_energy_identity} to obtain the normalized energy estimate.

\begin{proof}[Proof of Proposition \ref{Prop.energy}]
	By Lemma \ref{Lem.moving_energy_identity},
	\[\frac{\mathrm dE}{\mathrm d\sigma}+\frac2\sigma \left[g(\varrho)-\alpha-C(\sigma^{-1}+\nu(\sigma))\right]\mathfrak d \leq\frac C\sigma E+C\sqrt E\|\mathcal F\|_{\mathcal B_\varrho}.\]
	Use $\cE^\sharp_\varrho=\sigma^2E$, $\partial_\tau=\sigma\partial_\sigma$, and
	\[\partial_\tau(\sigma^2E)=2\sigma^2E+\sigma^3\partial_\sigma E,\qquad \sigma^3\sqrt E\|\mathcal F\|_{\mathcal B_\varrho} =\sqrt{\cE^\sharp_\varrho}\|\sigma^2\mathcal F\|_{\mathcal B_\varrho}.\]
	Multiplication by $\sigma^3$ gives \eqref{Eq.full_energy_differential}. Lemma \ref{Lem.energy_comparison} gives \eqref{Eq.energy_state_comparison}, and \eqref{Eq.analytic_radius_recovery} gives \eqref{Eq.energy_gap_comparison}.
\end{proof}

For the Cauchy problem, fix the background and analytic radii in \eqref{Eq.working_analytic_margin}. We first solve for $(\zeta,D\zeta,\sigma^{-1}\nabla\zeta)$ and then pass to the normalized energy by approximation.

\begin{lemma}\label{Lem.physical_forward}
	Let $0<T_*\leq2$ and let $R\in C^1([0,T_*];(0,2\overline\varrho))$ satisfy
	\[0\leq R'/R,\qquad \inf_{0\leq\tau\leq T_*}\bigl(\gamma_v-\mathcal L_v(R)-R'/R\bigr)>0.\]
	Let $\mathcal F$ be measurable and suppose
	\[\zeta_0,\nabla\zeta_0,z_{1,0}\in\mathcal B_{R(0)},\qquad \int_0^{T_*}\|\mathcal F(\Sigma\ee^u)\|_{\mathcal B_{R(u)}}\,\mathrm du<\infty.\]
	The equation \eqref{Eq.energy_scalar_equation}, with $(\zeta,D\zeta)(\Sigma)=(\zeta_0,z_{1,0})$, has a unique solution on $[\Sigma,\Sigma\ee^{T_*}]$. For
	\[V(\sigma):=(\zeta,D\zeta,\sigma^{-1}\nabla\zeta),\qquad V_0:=(\zeta_0,z_{1,0},\Sigma^{-1}\nabla\zeta_0),\]
	we have
	\begin{equation}\label{Eq.natural_forward_bound}
		\sup_{0\leq t\leq T_*}\|V(\Sigma\ee^t)\|_{\mathcal B_{R(t)}} \leq C_{\Sigma,T_*}\left[\|V_0\|_{\mathcal B_{R(0)}} +\int_0^{T_*}\|\mathcal F(\Sigma\ee^u)\|_{\mathcal B_{R(u)}}\,\mathrm du\right].
	\end{equation}
	Uniqueness is taken in the class continuous in all finite spatial Sobolev norms; the constants may depend on the finite time interval.
\end{lemma}

\begin{proof}
	Set $z_1=D\zeta$ and $G=\sigma^{-1}\nabla\zeta$. The natural symmetric system is
	\begin{equation}\label{Eq.natural_forward_system}
		\begin{cases}
			D\zeta=z_1,\\
			Dz_1-\sigma^{-1}\operatorname{div}G =-\lambda^2pf(\theta_\sigma)^{p-1}\zeta +\sigma^{-1}\mathbf g_0\cdot G-\sigma^{-2}Q_\sigma\zeta+\mathcal F,\\
			DG-\sigma^{-1}\nabla z_1=-\sigma^{-1}(I+(\mathrm D_x\mathbf b)^T)G.
		\end{cases}
	\end{equation}
	By $\bvec\cdot n=1$ in \eqref{Eq.acoustic_outgoing_property}, the outward boundary form of \eqref{Eq.natural_forward_system} is
	\[\sigma^{-1}\bigl(|\zeta|^2+|z_1-G\cdot n|^2+|G|^2-|G\cdot n|^2\bigr)\geq0.\]
	Extend the smooth symmetric principal matrices to $\mathbb R^d$ and extend the data and source at each finite Sobolev order. The symmetric hyperbolic Cauchy theorem gives a solution; see \cite{Sogge}. The nonnegative outward boundary form gives uniqueness on $\Omega$. The solutions constructed at different Sobolev orders therefore coincide on $\Omega$.
	
	For $V=(\zeta,z_1,G)$, use the finite sum
	\[H_K=\sum_{N=0}^K(1+N)^{2s_*}\frac{R(\tau)^{2N}}{(N!)^2} \|\nabla^NV\|_{L^2}^2.\]
	Here $H_K$ is differentiated in $\tau=\log(\sigma/\Sigma)$. The $\nabla$ and $\operatorname{div}$ terms in \eqref{Eq.natural_forward_system} cancel in the interior after integration by parts at every order, leaving the nonnegative boundary form. The ordered transport estimate \eqref{Eq.acoustic_transport_coercivity} and the commutator summation in Lemma \ref{Lem.macro_finite_energy}, with $R'/R$ included from the start, give
	\[\frac{\mathrm dH_K}{\mathrm d\tau}+2\bigl(\gamma_v-\mathcal L_v(R)-R'/R\bigr) \sum_{N=0}^KN(1+N)^{2s_*}\frac{R^{2N}}{(N!)^2}\|\nabla^NV\|_{L^2}^2 \leq C_{\Sigma,T_*}H_K+C_{\Sigma,T_*}\sqrt{H_K}\|\mathcal F\|_{\mathcal B_R}.\]
	The coefficients $\lambda^2pf(\theta_\sigma)^{p-1}$, $\mathbf g_0$, $Q_\sigma$ and $\mathrm D_x\mathbf b$ have bounded analytic multiplier norms on $[\Sigma,\Sigma\ee^{T_*}]$, by Lemmas \ref{Lem.phase_coefficients} and \ref{Lem.analytic_to_multiplier} and the bounds following \eqref{Eq.energy_scalar_equation}. The constant is independent of $K$. For admissible $R$ the coefficient on the left is positive. Gronwall and monotone convergence give \eqref{Eq.natural_forward_bound}. For $Q:=G-\sigma^{-1}\nabla\zeta$, we have $DQ=-\sigma^{-1}(I+(\mathrm D_x\mathbf b)^T)Q$ and $Q(\Sigma)=0$. The $L^2$ energy inequality gives $Q=0$, hence \eqref{Eq.natural_forward_system} implies \eqref{Eq.energy_scalar_equation}. Rotation invariance and conjugation show that the solution is radial and real-valued.
\end{proof}

We approximate finite-energy data by functions with a larger analytic radius. Lemma \ref{Lem.physical_forward} gives the approximating solutions, and \eqref{Eq.local_energy_budget} controls the differences in the moving analytic energy.

\begin{lemma}\label{Lem.energy_completion}
	Assume \eqref{Eq.working_analytic_margin}, let $0<T_*\leq2$, and let $\varrho\in C^1([0,T_*];I_{\rm r})$ satisfy $0\leq\varrho'/\varrho\leq\eta$, and let $\Sigma$ be sufficiently large. For a measurable scalar source $\mathcal F$ and initial data satisfying
	\[\cE^\sharp_{\varrho(0)}(\Sigma;W_0)<\infty,\qquad \int_0^{T_*}\|(\Sigma\ee^u)^2\mathcal F(\Sigma\ee^u)\|_{\mathcal B_{\varrho(u)}}\,\mathrm du<\infty,\]
	there is a unique finite-energy pair $W$ with $W(\Sigma)=W_0$ whose reconstruction satisfies $\cL_\sigma\zeta_W=\mathcal F$ and $D\zeta_W=\lambda w_W$. The unknown $W$ satisfies
	\[\sup_{0\leq t\leq T_*}\sqrt{\cE^\sharp_{\varrho(t)}(\Sigma\ee^t;W)} \leq C\left[\sqrt{\cE^\sharp_{\varrho(0)}(\Sigma;W_0)} +\int_0^{T_*}\|(\Sigma\ee^u)^2\mathcal F(\Sigma\ee^u)\|_{\mathcal B_{\varrho(u)}}\,\mathrm du\right].\]
	Moreover, $W\in C([\Sigma,\Sigma\ee^{T_*}];H^q(\Omega;\R^2))$ for every integer $q\geq0$. For every closed interval $J\subset[0,T_*]$ and $r<\inf_J\varrho$,
	\[\tau\longmapsto W(\Sigma\ee^\tau)\in C(J;\cX_r).\]
\end{lemma}

\begin{proof}
	At fixed $\sigma$, the norm $(\cE^\sharp_\varrho)^{1/2}$ is equivalent, with constants allowed to depend on $\sigma$, to the natural wave norm
	\[\|\zeta\|_{\mathcal B_\varrho} +\|\nabla\zeta\|_{\mathcal B_\varrho} +\|D\zeta\|_{\mathcal B_\varrho}.\]
	Indeed, above a finite derivative level the cutoff is zero and $S_N=\operatorname{diag}(1+\xi_N^2,1)$. The finitely many lower levels have positive metrics and triangular changes of derivatives with bounded coefficients. Induction over the finite lower levels and the factorial shift at the higher levels give equivalence with the natural wave norm. Since the natural wave space is complete, the norm $(\cE^\sharp_\varrho)^{1/2}$ is complete.
	
	Reconstruct $(\zeta_0,w_0)$ from $W_0$ by \eqref{Eq.zeta-w-W-def} at $\sigma=\Sigma$, and set $z_{1,0}=\lambda(\Sigma)w_0$, $G_0=\Sigma^{-1}\nabla\zeta_0$. For $0<a<1$, approximate these data and the source by
	\[(\zeta_{0,a},z_{1,0,a},G_{0,a})(x)=(\zeta_0(ax),z_{1,0}(ax),aG_0(ax)),\qquad \mathcal F_a(\sigma,x)=\mathcal F(\sigma,ax).\]
	Then $G_{0,a}=\Sigma^{-1}\nabla\zeta_{0,a}$. The corresponding $W_{0,a}$ is obtained from $(\zeta_{0,a},\lambda(\Sigma)^{-1}z_{1,0,a})$ by \eqref{Eq.state}; the background coefficients are not dilated. The dilated data and sources converge in the initial natural norm and in the integrable moving source norm by the dilation estimates in Lemma \ref{Lem.analytic_density} and dominated convergence. The dilated functions have analytic radius $R(\tau)=\varrho(\tau)/a$. Since $R\to\varrho$ uniformly and $R'/R=\varrho'/\varrho$, the strict inequalities in \eqref{Eq.working_analytic_margin} imply $R<2\overline\varrho$ and $\inf_\tau(\gamma_v-\mathcal L_v(R)-R'/R)>0$ for $a$ sufficiently close to one. Lemma \ref{Lem.physical_forward} therefore gives the approximating solutions.
	
	The positive gap between $R$ and $\varrho$ makes the refined energy sums and differentiated identities absolutely convergent. Let $W_a,W_b$ be two approximating solutions with sources $\mathcal F_a,\mathcal F_b$. Lemma \ref{Lem.moving_energy_identity}, the positive margin in \eqref{Eq.working_analytic_margin}, and the square-root inequality give
	\begin{align*}
		\sup_{0\leq\tau\leq T_*}\sqrt{\cE^\sharp_{\varrho(\tau)}(\Sigma \ee^\tau;W_a-W_b)} \leq{}&C\sqrt{\cE^\sharp_{\varrho(0)}(\Sigma;W_a(\Sigma)-W_b(\Sigma))}\\
		&+C\int_0^{T_*}\|(\Sigma \ee^\tau)^2(\mathcal F_a-\mathcal F_b)(\Sigma \ee^\tau)\|_{\mathcal B_{\varrho(\tau)}}\,\mathrm d\tau \longrightarrow0.
	\end{align*}
	Thus the approximations are Cauchy in the supremum of the moving refined energy norms. The lower comparison in Lemma \ref{Lem.energy_comparison} gives a common distributional limit; closedness of each spatial derivative and Fatou's lemma identify $\cE^\sharp_{\varrho(\tau)}(\Sigma\ee^\tau;W)$. The scalar equation passes to the limit at every finite Sobolev order. The same estimate gives uniqueness. Continuity at a smaller analytic radius follows from the equation and the gap between the working and target analytic radii, and continuity at the initial time follows by approximation. Integrable sources are handled first for simple functions in time and then by approximation in the integral source norm.
\end{proof}

We finish the proof of Corollary \ref{Cor.forward_energy} by choosing a common late-time threshold and integrating the square-root energy inequality.

\begin{proof}[Proof of Corollary \ref{Cor.forward_energy}]
	Choose the threshold so that $C(\sigma^{-1}+\nu(\sigma))\leq\kappa_*/4$. Since $g(\varrho)\geq\kappa_*$ and $\alpha\leq\eta<\kappa_*/4$, the derivative term in Lemma \ref{Lem.moving_energy_identity} is nonnegative. Multiply \eqref{Eq.local_energy_budget} by $\sigma^3$, using
	\[\partial_\tau(\sigma^2E)=2\sigma^2E+\sigma^3\partial_\sigma E, \qquad \sigma^3\sqrt E\|\mathcal F\|_{\mathcal B_\varrho} =\sqrt{\cE^\sharp_\varrho}\|\sigma^2\mathcal F\|_{\mathcal B_\varrho}.\]
	Discard the nonnegative derivative sum. Applying the resulting inequality to $\sqrt{\cE^\sharp_\varrho+\varepsilon}$ and letting $\varepsilon\downarrow0$ gives
	\[\partial_\tau\sqrt{\cE^\sharp_{\varrho(\tau)}} \leq C\sqrt{\cE^\sharp_{\varrho(\tau)}} +C\|\sigma^2\mathcal F\|_{\mathcal B_{\varrho(\tau)}}.\]
	Gronwall's inequality on any subinterval $[s,t]\subset[0,T_*]$ gives \eqref{Eq.energy_estimate}. Lemma \ref{Lem.energy_completion} proves existence and uniqueness for the data and sources in Corollary \ref{Cor.forward_energy}. Lemma \ref{Lem.energy_comparison} and \eqref{Eq.analytic_radius_recovery} give \eqref{Eq.energy_state_comparison} and \eqref{Eq.energy_gap_comparison}.
\end{proof}

\section{Averaging and the nonlinear solution}
\label{Sec.state_proof}

We construct a decaying solution of $\cN(z_{\app}+\zeta)=0$ by choosing the initial components that would otherwise violate the decay rate. The abstract recurrence is \eqref{Eq.abstract_decay_recurrence}. For the wave equation, \eqref{Eq.averaging} compares $\mathcal P_\Sigma(1)$ with the compact map $\cS(1)$, and \eqref{Eq.nonlinear_endpoint_bounds} separates the nonlinear endpoint into its linear part and a small Lipschitz remainder. After applying Proposition \ref{Prop.compact_decay}, we reconstruct $u^{\mathrm c}$ and prove the cone estimates \eqref{Eq.cone_rate}--\eqref{Eq.peaks}.

\subsection{An abstract decaying trajectory}
\label{Subsec.abstract_decay}

Let $X$ be a real Banach space, let $K:X\to X$ be compact, and fix $b\in(0,1)$ with $\spec(K)\cap\{|z|=b\}=\varnothing$, where the spectrum is taken after complexification. Let $P_+$ be the Riesz projection for $|z|>b$ and $P_-=I-P_+$. We seek $x_{n+1}=Kx_n+R_n(x_n)$ with $P_-x_0=0$ and $\sup_nb^{-n}\|x_n\|_X<\infty$. Forward iteration controls $P_-x_n$; the decay requirement determines $P_+x_n$ by backward iteration, as made explicit in \eqref{Eq.discrete_green_proof}.

\begin{lemma}\label{Lem.abstract_spectral}
	The space $P_+X$ is finite dimensional and $K_+:=K|_{P_+X}$ is invertible. There are $C>0$ and $\theta\in(0,1)$ such that, for $B=b^{-1}K$,
	\begin{equation}\label{Eq.abstract_spectral_powers}
		\|B^nP_-\|+\|(B|_{P_+X})^{-n}P_+\|\leq C\theta^n \qquad(n\geq0).
	\end{equation}
	The second term is omitted when $P_+=0$.
\end{lemma}

\begin{proof}
	Compactness implies that the part of the spectrum in $|z|>b$ consists of finitely many nonzero eigenvalues with finite algebraic multiplicities. The Riesz projection $P_+$ therefore has finite rank, and $0\notin\operatorname{spec}(K_+)$. After scaling by $B=b^{-1}K$, the spectrum on $P_-X$ is strictly inside the unit circle and the spectrum of $(B|_{P_+X})^{-1}$ is also strictly inside. Thus, with $r(A):=\sup\{|z|:z\in\spec(A)\}$,
	\[r_-:=r(B|_{P_-X})<1,\qquad r_+:=r((B|_{P_+X})^{-1})<1.\]
	When $P_+=0$ set $r_+=0$. Choose $\max\{r_-,r_+\}<\theta_0<\theta<1$. The spectral-radius formula gives an integer $n_0$ such that
	\[\|B^nP_-\|\leq C_0\theta_0^n,\qquad\|(B|_{P_+X})^{-n}P_+\|\leq C_0\theta_0^n\quad(n\geq n_0),\]
	where the fixed projections can be included in $C_0$. Taking
	\[C=2C_0+\max_{0\leq n<n_0}\theta^{-n}\left(\|B^nP_-\|+\|(B|_{P_+X})^{-n}P_+\|\right)\]
	gives $\|B^nP_-\|+\|(B|_{P_+X})^{-n}P_+\|\leq C\theta^n$ for every $n\geq0$. The inside and outside spectral sets are invariant under complex conjugation. The resolvent identity $\overline{(z-K)^{-1}u}=(\overline z-K)^{-1}\overline u$ in the contour formula therefore shows that $P_\pm$ commute with conjugation. Thus the restrictions to $X$ are real operators.
\end{proof}

The projection $P_+$ is defined by the prescribed decay rate $b$, not by the unit circle. The generalized eigenspace of an eigenvalue with $b<|\lambda|<1$ still lies in $P_+X$: the homogeneous sequence $\lambda^n v$ decays, but $b^{-n}|\lambda|^n\|v\|$ is unbounded. Compactness makes the sum of all generalized eigenspaces with $|\lambda|>b$ finite-dimensional, without requiring $K$ to be diagonalizable.

For sequences define $\|x\|_b:=\sup_{n\geq0}b^{-n}\|x_n\|_X$. Let $\ell_b(X)$ be the resulting Banach space.

The condition $P_-x_0=0$ determines $P_-x_n$ by forward iteration of $x_{n+1}=Kx_n+f_n$. On $P_+X$, imposing $\|x_n\|_X=O(b^n)$ determines the solution by a sum over $j\geq n$. The power bounds \eqref{Eq.abstract_spectral_powers} control both sums in the norm $\|\cdot\|_b$.

\begin{lemma}[The discrete Green operator]\label{Lem.discrete_green}
	For $f\in\ell_b(X)$, the recurrence
	\[x_{n+1}=Kx_n+f_n,\qquad P_-x_0=0\]
	has a unique solution in $\ell_b(X)$. The solution is $x=\mathcal Gf$, where
	\begin{equation}\label{Eq.discrete_green_proof}
		(\mathcal Gf)_n =\sum_{j=0}^{n-1}K^{n-1-j}P_-f_j -\sum_{j=n}^{\infty}K_+^{-(j-n+1)}P_+f_j.
	\end{equation}
	Moreover, $\|\mathcal Gf\|_b\leq C_G\|f\|_b$.
\end{lemma}

\begin{proof}
	Forward iteration gives $P_-x_n=\sum_{j=0}^{n-1}K^{n-1-j}P_-f_j$. On $P_+X$, iterate from $n$ to $m>n$ and solve for $P_+x_n$:
	\[P_+x_n=K_+^{-(m-n)}P_+x_m-\sum_{j=n}^{m-1}K_+^{-(j-n+1)}P_+f_j.\]
	For $x\in\ell_b(X)$, \eqref{Eq.abstract_spectral_powers} bounds the first term by $Cb^n\theta^{m-n}\|x\|_b$, which tends to zero as $m\to\infty$. This yields the formula \eqref{Eq.discrete_green_proof}.
	
	Set $v_n=b^{-n}x_n$ and $g_n=b^{-(n+1)}f_n$. The equation becomes $v_{n+1}=Bv_n+g_n$, and $\|g\|_{\ell^\infty}\leq b^{-1}\|f\|_b$. Lemma \ref{Lem.abstract_spectral} gives
	\begin{align*}
		\left\|\sum_{j=0}^{n-1}B^{n-1-j}P_-g_j\right\| &\leq C\|g\|_{\ell^\infty}\sum_{k=0}^{n-1}\theta^k\leq\frac C{1-\theta}\|g\|_{\ell^\infty},\\
		\left\|\sum_{j=n}^\infty(B|_{P_+X})^{-(j-n+1)}P_+g_j\right\| &\leq C\|g\|_{\ell^\infty}\sum_{k=1}^\infty\theta^k\leq\frac{C\theta}{1-\theta}\|g\|_{\ell^\infty}.
	\end{align*}
	Hence the forward series on $P_-X$ minus the backward series on $P_+X$ defines a bounded sequence $v$, and rescaling gives \eqref{Eq.discrete_green_proof} with $C_G=2C/[b(1-\theta)]$. Subtracting $Bv_n$ from $v_{n+1}$ cancels every summand except $P_-g_n$ from the first series and $P_+g_n$ from the second. Thus $v_{n+1}-Bv_n=g_n$. At $n=0$ the first series is empty and the second lies in $P_+X$, so $P_-v_0=0$.
	
	For uniqueness, let $v$ be a bounded homogeneous solution with $P_-v_0=0$. Then
	\[P_-v_n=B^nP_-v_0=0,\qquad P_+v_n=(B|_{P_+X})^{-j}P_+v_{n+j}.\]
	Therefore $\|P_+v_n\|\leq C\theta^j\sup_k\|v_k\|\to0$ as $j\to\infty$. Thus $v_n=0$ for every $n$. The spectral projections and both series preserve reality by Lemma \ref{Lem.abstract_spectral}.
\end{proof}

We now apply the discrete Green operator to $R_n(x_n)$. The proof of Proposition \ref{Prop.compact_decay} gives an invariant ball in the weighted sequence space and a contraction on that ball.

\begin{proof}[Proof of Proposition \ref{Prop.compact_decay}]
	Choose $\epsilon_0$ with $C_G\epsilon_0\leq1/2$. On the closed ball $\|x\|_b\leq2C_Gh$, we have $\|x_n\|_X\leq b^n2C_Gh<R$, so every $R_n(x_n)$ is defined. Because $\|R_n(0)\|\leq hb^n$ and $\operatorname{Lip}(R_n)\leq\epsilon_0$, the map $\mathcal T(x)=\mathcal G((R_n(x_n))_n)$ satisfies
	\[\|\mathcal T(x)\|_b\leq C_Gh+C_G\epsilon_0\|x\|_b\leq2C_Gh, \qquad \|\mathcal T(x)-\mathcal T(y)\|_b\leq\tfrac12\|x-y\|_b.\]
	Banach's fixed-point theorem gives a unique fixed point. Lemma \ref{Lem.discrete_green} identifies the fixed point $x$ with a solution of \eqref{Eq.abstract_decay_recurrence}. Conversely, let $x,y$ satisfy \eqref{Eq.abstract_decay_recurrence} and $\|x\|_b,\|y\|_b\leq2C_Gh$. The identities $x=\mathcal G(R_n(x_n))_n$ and $y=\mathcal G(R_n(y_n))_n$ give
	\[\|x-y\|_b\leq C_G\epsilon_0\|x-y\|_b\leq\tfrac12\|x-y\|_b,\]
	so $x=y$.
\end{proof}

\subsection{Comparison with the averaged evolution}
\label{Subsec.averaging_proof}

Fix the background and the radii in \eqref{Eq.working_analytic_margin}--\eqref{Eq.radii}, and let $X=\cX_{\varrho_0}$. The homogeneous equation is $\partial_\tau W=(\mathscr A+\cB(\sigma)+\mathscr E_\sigma)W$, with $\sigma=\Sigma\ee^\tau$. We use $W^c=(I-\mathscr N_\sigma)W$ from \eqref{Eq.averaging_corrector} to cancel $\cB(\sigma)$. The full solution is bounded at $\varrho_{\rm F}(\tau)$, while the equation for $W^c-\cS(\tau)W_0$ is estimated at $\varrho_{\rm A}(\tau)<\varrho_{\rm F}(\tau)$. Since $\varrho_{\rm A}(1)>\varrho_0$, the endpoint estimate returns to the initial space $X$.

\begin{proof}[Proof of Corollary \ref{Cor.closed_energy}]
	By \eqref{Eq.energy_gap_comparison},
	\[\sqrt{\cE^\sharp_{\varrho_{\rm F}(0)}(\Sigma;W_0)}\leq C\|W_0\|_X.\]
	For the homogeneous solution $W(\Sigma\ee^\tau)=\mathcal P_\Sigma(\tau)W_0$, Corollary \ref{Cor.forward_energy} and \eqref{Eq.energy_state_comparison} give
\begin{align*}
\|\mathcal P_\Sigma(\tau)W_0\|_{\cX_{\varrho_{\rm F}(\tau)}}
&\leq C\sqrt{\cE^\sharp_{\varrho_{\rm F}(\tau)}(\Sigma\ee^\tau;W(\Sigma\ee^\tau))}\\
&\leq C\sqrt{\cE^\sharp_{\varrho_{\rm F}(0)}(\Sigma;W_0)}\leq C\|W_0\|_X,\qquad 0\leq\tau\leq2.
\end{align*}
The initial inequality $\varrho_{\rm F}(0)<\varrho_0$ pays for the extra spatial derivative in the upper energy bound. The constant is independent of $\Sigma$, because \eqref{Eq.energy_gap_comparison} and \eqref{Eq.energy_estimate} are uniform at these fixed analytic-radius gaps. This proves \eqref{Eq.endpoint_smoothing} for the true evolution. For the averaged evolution, \eqref{Eq.acoustic_radius_gap} and Proposition \ref{Prop.macro} give
	\[\|\cS(\tau)W_0\|_{\cH_{\varrho_{\rm F}(\tau)}}\leq C\|W_0\|_{\cH_{\varrho_{\rm F}(0)}}\leq C\|W_0\|_X.\]
	Since $\|W\|_{\cX_r}\leq\|W\|_{\cH_r}$, the averaged bound in \eqref{Eq.endpoint_smoothing} follows. At $\tau=1$, both maps factor through a strictly larger analytic radius:
\[X\xrightarrow{\ \mathcal P_\Sigma(1)\ \text{or}\ \cS(1)\ }\cX_{\varrho_{\rm F}(1)}\lhook\joinrel\longrightarrow X,\qquad \varrho_{\rm F}(1)>\varrho_0.\]
The first arrow is bounded by \eqref{Eq.endpoint_smoothing}; the last inclusion is compact by Lemma \ref{Lem.analytic_compactness}. This proves compactness of the endpoint maps on $X$. The argument does not require regularization from arbitrary Sobolev data: the improvement used here is the increase from the initial analytic radius $\varrho_0$ to $\varrho_{\rm F}(1)$. Approximation by data of larger analytic radius extends the endpoint maps to $X$, and uniqueness identifies the extensions with the distributional evolutions.
\end{proof}

The operator $\cB(\sigma)$ is not small in its coefficients. Its periodic primitive in \eqref{Eq.averaging_corrector} has the additional factor $\sigma^{-1}$, and satisfies $\partial_\tau\mathscr N_\sigma=\cB(\sigma)-\mathscr N_\sigma$ and $\mathscr N_\sigma\cB(\sigma)=0$. Differentiating $(I-\mathscr N_\sigma)W$ therefore cancels $\cB(\sigma)W$. Since $\mathscr N_\sigma$ contains $A_0$, we do not regard $I-\mathscr N_\sigma$ as a bounded change of unknowns on the single space $X$. The larger analytic radius of $W$ controls the derivatives in $\mathscr N_\sigma W$, and the equation for the corrected unknown is estimated at the smaller radius $\varrho_{\rm A}$.

\begin{lemma}[The corrector identity]\label{Lem.weak_averaging}
	Let $W$ solve \eqref{Eq.state_operator_form} with $h=0$ and $W(\Sigma)=W_0\in X$. We define $W^c=(I-\mathscr N_\sigma)W$. Then
	\begin{equation}\label{Eq.corrected_equation}
		(\partial_\tau-\mathscr A)W^c=F_\sigma^c, \qquad F_\sigma^c=\bigl([\mathscr A,\mathscr N_\sigma] +\mathscr N_\sigma+\mathscr E_\sigma -\mathscr N_\sigma\mathscr E_\sigma\bigr)W.
	\end{equation}
	For $0\leq\tau\leq1$,
	\begin{equation}\label{Eq.moving_corrector_bound}
		\|F_{\Sigma \ee^\tau}^c\|_{\cH_{\varrho_{\rm A}(\tau)}} \leq C\mu(\Sigma)\|W_0\|_X.
	\end{equation}
	The two endpoint corrections satisfy
	\[\|\mathscr N_\Sigma W_0\|_{\cH_{\varrho_{\rm A}(0)}} +\|\mathscr N_{\Sigma\ee}W(\Sigma\ee)\|_{\cX_{\varrho_{\rm A}(1)}} \leq C\mu(\Sigma)\|W_0\|_X.\]
\end{lemma}

\begin{proof}
	Differentiate $W^c=(I-\mathscr N_\sigma)W$ using \eqref{Eq.state_operator_form} with $h=0$ and the identities from \eqref{Eq.averaging_corrector}, namely $\partial_\tau\mathscr N_\sigma=\cB(\sigma)-\mathscr N_\sigma$ and $\mathscr N_\sigma\cB(\sigma)=0$. Then
	\begin{align*}
		\partial_\tau W^c &=(I-\mathscr N_\sigma)(\mathscr A+\cB+\mathscr E_\sigma)W-(\cB-\mathscr N_\sigma)W\\
		&=\mathscr A(I-\mathscr N_\sigma)W +\bigl([\mathscr A,\mathscr N_\sigma]+\mathscr N_\sigma+\mathscr E_\sigma-\mathscr N_\sigma\mathscr E_\sigma\bigr)W.
	\end{align*}
	Thus $W^c$ satisfies \eqref{Eq.corrected_equation}.
	
	For $C_\sigma=-\alpha_0(\sigma)A_0-\alpha_1(\sigma)V$, direct multiplication gives
	\begin{equation}\label{Eq.corrector_operator_commutator}
		[\mathscr A,\mathscr N_\sigma]=\sigma^{-1}\begin{pmatrix}C_\sigma&0\\
			C_\sigma+[-B,C_\sigma]&-C_\sigma\end{pmatrix}.
	\end{equation}
	The commutator $[-B,C_\sigma]$ in \eqref{Eq.corrector_operator_commutator} is determined by
	\begin{equation}\label{Eq.corrector_coefficient_commutators}
		\begin{aligned}
			[-B,\Delta]&=2\sum_{i,k}(\partial_i b_k)\partial_i\partial_k+(\Delta\bvec)\cdot\nabla,\\
			[-B,\mathbf g_0\cdot\nabla]&=(\mathbf g_0\cdot\nabla\bvec-B\mathbf g_0)\cdot\nabla, &[-B,V]&=-BV.
		\end{aligned}
	\end{equation}
	The coefficient bounds \eqref{Eq.common_coefficient_radius}, Lemma \ref{Lem.analytic_to_multiplier} and \eqref{Eq.multiplier_algebra} bound the coefficients of \eqref{Eq.corrector_coefficient_commutators} in $\mathfrak M_{\overline\varrho}$. Consequently, $[\mathscr A,\mathscr N_\sigma]$ and $\mathscr N_\sigma$ have order at most two, with coefficient multiplier norms at $\overline\varrho$ bounded by $C\sigma^{-1}$. Lemma \ref{Lem.state_remainder} gives order at most two and coefficient multiplier bounds $C\nu(\sigma)$ for $\mathscr E_\sigma$.
	
	We now estimate \eqref{Eq.corrected_equation} with the derivative losses included. For fixed $0<r<R<\overline\varrho$, the coefficient bounds from \eqref{Eq.corrector_coefficient_commutators} and Lemma \ref{Lem.state_remainder}, together with \eqref{Eq.analytic_radius_recovery} give
\begin{align*}
\|([\mathscr A,\mathscr N_\sigma]+\mathscr N_\sigma)U\|_{\cH_r}&\leq C\sigma^{-1}\left(1+(R-r)^{-3}\right)\|U\|_{\cX_R},\\
\|\mathscr E_\sigma U\|_{\cH_r}&\leq C\nu(\sigma)\left(1+(R-r)^{-3}\right)\|U\|_{\cX_R}.
\end{align*}
There are at most two derivatives in each operator and one additional derivative on the first component in \eqref{Eq.acoustic_pair_norm}. For the composition, use the triangular form of $\mathscr N_\sigma$:
\[\mathscr N_\sigma\mathscr E_\sigma U=\sigma^{-1}\binom{0}{C_\sigma(\mathscr E_\sigma U)_1}.\]
The first component vanishes, so this term does not require the additional gradient in the $\cH_r$ norm. Choose $r_m=(r+R)/2$. Applying the second-order estimates successively, with coefficients controlled by \eqref{Eq.multiplier_algebra}, gives
\begin{align*}
\|\mathscr N_\sigma\mathscr E_\sigma U\|_{\cH_r}
&=\sigma^{-1}\|C_\sigma(\mathscr E_\sigma U)_1\|_{\mathcal B_r}\leq C\sigma^{-1}\left(1+(r_m-r)^{-2}\right)\|\mathscr E_\sigma U\|_{\cX_{r_m}}\\
&\leq C\sigma^{-1}\nu(\sigma)\left(1+(R-r)^{-4}\right)\|U\|_{\cX_R}.
\end{align*}
Only $\mathscr N_\sigma\mathscr E_\sigma$ can have order four; the explicit commutator \eqref{Eq.corrector_operator_commutator} has order two. Thus no fourth-order principal term is being added to the averaged equation.

By \eqref{Eq.radii}, the difference $\varrho_{\rm F}(\tau)-\varrho_{\rm A}(\tau)$ has a positive minimum independent of $\Sigma$ on $[0,1]$. Taking $r=\varrho_{\rm A}(\tau)$, $R=\varrho_{\rm F}(\tau)$ in these estimates yields
\begin{align*}
\|([\mathscr A,\mathscr N_\sigma]+\mathscr N_\sigma)W\|_{\cH_{\varrho_{\rm A}(\tau)}}&\leq C\sigma^{-1}\|W\|_{\cX_{\varrho_{\rm F}(\tau)}},\\
\|\mathscr E_\sigma W\|_{\cH_{\varrho_{\rm A}(\tau)}}&\leq C\nu(\sigma)\|W\|_{\cX_{\varrho_{\rm F}(\tau)}},\\
\|\mathscr N_\sigma\mathscr E_\sigma W\|_{\cH_{\varrho_{\rm A}(\tau)}}&\leq C\sigma^{-1}\nu(\sigma)\|W\|_{\cX_{\varrho_{\rm F}(\tau)}}.
\end{align*}
	The constants depend on the positive minimum of $\varrho_{\rm F}-\varrho_{\rm A}$, not on $\Sigma$. Since $\sigma\geq\Sigma$, \eqref{Eq.endpoint_smoothing} proves \eqref{Eq.moving_corrector_bound}. The inequalities $\varrho_{\rm A}(0)<\varrho_0$ and $\varrho_{\rm A}(1)<\varrho_{\rm F}(1)$ in \eqref{Eq.radii} likewise give
	\begin{align*}
		\|\mathscr N_\Sigma W_0\|_{\cH_{\varrho_{\rm A}(0)}}&\leq C\Sigma^{-1}\|W_0\|_X,\\
		\|\mathscr N_{\Sigma \ee}W(\Sigma \ee)\|_{\cX_{\varrho_{\rm A}(1)}}&\leq C\Sigma^{-1}\|W(\Sigma \ee)\|_{\cX_{\varrho_{\rm F}(1)}}\leq C\Sigma^{-1}\|W_0\|_X.
	\end{align*}
	For $W_0\in X$, approximate $W_0$ by data with a larger analytic radius and use \eqref{Eq.moving_corrector_bound} to pass to the limit in \eqref{Eq.corrected_equation}.
\end{proof}

For the difference $W^c-\cS(\tau)W_0$, the initial correction contains two derivatives of $W_0$. The inequalities
\[\varrho_{\rm A}(0)<\varrho_{\rm F}(0)<\varrho_0<\varrho_{\rm A}(1)<\varrho_{\rm F}(1)\]
from \eqref{Eq.radii} serve different purposes: $\varrho_{\rm A}(0)<\varrho_0$ controls the initial correction, $\varrho_{\rm A}(\tau)<\varrho_{\rm F}(\tau)$ controls the source \eqref{Eq.corrected_equation}, and $\varrho_{\rm A}(1)>\varrho_0$ returns the endpoint estimate to $X=\cX_{\varrho_0}$. We compare $W^c$ with $\cS(\tau)W_0$ before using $W=W^c+\mathscr N_\sigma W$ at $\tau=1$. Thus the small operator-norm difference comes from the corrector equation, not from compactness alone.

\begin{proof}[Proof of Proposition \ref{Prop.averaging}]
	Let $\overline W(\tau)=\cS(\tau)W_0$ and $Y(\tau)=W^c(\Sigma\ee^\tau)-\overline W(\tau)$. By \eqref{Eq.corrected_equation},
	\[(\partial_\tau-\mathscr A)Y=F_{\Sigma\ee^\tau}^c,\qquad Y(0)=-\mathscr N_\Sigma W_0.\]
	The initial difference belongs to $\cH_{\varrho_{\rm A}(0)}$, and the source is integrable in the moving $\cH_{\varrho_{\rm A}}$ norm by Lemma \ref{Lem.weak_averaging}. For data of a strictly larger analytic radius, the integral equation for $Y$ reads
\begin{equation}\label{Eq.averaging_endpoint_identity}
\mathcal P_\Sigma(1)W_0-\cS(1)W_0=-\cS(1)\mathscr N_\Sigma W_0+\int_0^1\cS(1-s)F_{\Sigma\ee^s}^c\,\mathrm ds+\mathscr N_{\Sigma\ee}W(\Sigma\ee).
\end{equation}
The terms in \eqref{Eq.averaging_endpoint_identity} are evaluated between different analytic radii, not by treating $\mathscr N_\Sigma$ as an operator $X\to X$. The first term starts at $\varrho_{\rm A}(0)$, the integral starts at $\varrho_{\rm A}(s)$, and the last correction is estimated using the gap $\varrho_{\rm F}(1)-\varrho_{\rm A}(1)$. Proposition \ref{Prop.macro} propagates the first two terms to $\varrho_{\rm A}(1)$, with a constant uniform in $s\in[0,1]$ because the shifted radius paths satisfy the same strict margin. Applying Proposition \ref{Prop.macro} directly to $Y$ therefore gives
	\begin{align*}
		\|Y(1)\|_{\cH_{\varrho_{\rm A}(1)}}&\leq C\|\mathscr N_\Sigma W_0\|_{\cH_{\varrho_{\rm A}(0)}}+C\int_0^1\|F_{\Sigma\ee^s}^c\|_{\cH_{\varrho_{\rm A}(s)}}\,\mathrm ds\leq C\mu(\Sigma)\|W_0\|_X.
	\end{align*}
	Since $\varrho_{\rm A}(1)>\varrho_0$, restriction to $X$ and the endpoint correction estimate yield
	\begin{align*}
		\|\mathcal P_\Sigma(1)W_0-\mathsf A W_0\|_X\leq\|Y(1)\|_{\cH_{\varrho_{\rm A}(1)}}+\|\mathscr N_{\Sigma\ee}W(\Sigma\ee)\|_{\cX_{\varrho_{\rm A}(1)}}\leq C\mu(\Sigma)\|W_0\|_X.
	\end{align*}
	The estimate for $\mathcal P_\Sigma(1)W_0-\mathsf A W_0$ is linear in $\|W_0\|_X$, with a constant independent of the extra radius used to justify \eqref{Eq.averaging_endpoint_identity}. Approximate an arbitrary $W_0\in X$ by data of larger analytic radius, using Lemma \ref{Lem.analytic_density}, and pass to the limit with the uniform endpoint bounds \eqref{Eq.endpoint_smoothing}. Taking the supremum over $\|W_0\|_X\leq1$ proves \eqref{Eq.averaging}. The rate $\mu(\Sigma)$ is therefore an operator-norm rate, not merely convergence for each fixed initial datum.
\end{proof}

\subsection{The nonlinear endpoint and the decaying solution}
\label{Subsec.contraction_proof}

We solve $\cL_\sigma\zeta=-\cR_{\app}-\cQ_\sigma(\zeta)$ with $W(\Sigma)=\Sigma^{-2}x$, where $\cR_{\app}=\cN(z_{\app})$ and $\cQ_\sigma$ is defined in \eqref{Eq.remainder_equation}. The scaling $\Sigma^{-2}$ keeps the quadratic source small in the path norm \eqref{Eq.nonlinear_path_norm}. After solving on $[\Sigma,\ee^2\Sigma]$, we evaluate at $\ee\Sigma$ to obtain $F_\Sigma(x)$ from \eqref{Eq.nonlinear_endpoint_definition}. The bounds on $F_\Sigma$ and \eqref{Eq.averaging} then permit the application of Proposition \ref{Prop.compact_decay}.

For $u,v\in\mathcal B_\varrho$ with norms at most one, the exact expansion
\[\cQ_\sigma(u)=\sum_{k=2}^p\binom pk z_\app^{p-k}u^k,\qquad u^k-v^k=(u-v)\sum_{h=0}^{k-1}u^{k-1-h}v^h\]
and Lemma \ref{Lem.analytic_algebra} give
\begin{equation}\label{Eq.nonlinear_algebra_bounds}
	\begin{aligned}
		\|\cQ_\sigma(u)\|_{\mathcal B_\varrho} &\leq C\sum_{k=2}^p\|u\|_{\mathcal B_\varrho}^k\leq C\|u\|_{\mathcal B_\varrho}^2,\\
		\|\cQ_\sigma(u)-\cQ_\sigma(v)\|_{\mathcal B_\varrho} &\leq C(\|u\|_{\mathcal B_\varrho}+\|v\|_{\mathcal B_\varrho})\|u-v\|_{\mathcal B_\varrho}.
	\end{aligned}
\end{equation}
The coefficient bound for $z_\app$ makes $C$ uniform on the chosen interval of analytic radii and for sufficiently large $\sigma$.

For a pair of unknowns depending on time on $[0,2]$, we define
\begin{equation}\label{Eq.nonlinear_path_norm}
	\|W\|_{\mathfrak E_\Sigma} :=\sup_{0\leq\tau\leq2}\Sigma^2 \sqrt{\cE^\sharp_{\varrho_{\rm F}(\tau)}(\Sigma \ee^\tau;W(\Sigma \ee^\tau))}.
\end{equation}
We use paths continuous at every smaller analytic radius. The normed path space $\mathfrak E_\Sigma$ is complete: a Cauchy sequence converges uniformly in the smaller analytic spaces, each finite spatial derivative is closed, and Fatou's lemma identifies the limit in $\mathfrak E_\Sigma$. At a fixed time, completeness follows from the equivalence of $\sqrt{\cE^\sharp_\varrho}$ and the natural wave norm proved in Lemma \ref{Lem.energy_completion}.

\begin{lemma}\label{Lem.local_nonlinear}
	For $\|x\|_X$ sufficiently small and $\Sigma$ sufficiently large, the equation
	\[\cL_\sigma\zeta=-\cR_{\app}-\cQ_\sigma(\zeta),\qquad W(\Sigma)=\Sigma^{-2}x\]
	has a unique small solution on $[\Sigma,\ee^2\Sigma]$. The solution $W$ satisfies
	\[\|W\|_{\mathfrak E_\Sigma}\leq C(\|x\|_X+\epsilon_\Sigma).\]
	The solution is real-valued and radial and is smooth in space and time.
\end{lemma}

\begin{proof}
	For $V\in\mathfrak E_\Sigma$, reconstruct $\zeta_V$ by \eqref{Eq.zeta-w-W-def} and apply Corollary \ref{Cor.forward_energy} with source $-\cR_{\app}-\cQ_\sigma(\zeta_V)$ and initial value $W(\Sigma)=\Sigma^{-2}x$. Denote the resulting pair of unknowns by $\mathcal T_xV$. Since $\varrho_{\rm F}(0)<\varrho_0$, \eqref{Eq.energy_gap_comparison} gives $\Sigma^2\sqrt{\cE^\sharp_{\varrho_{\rm F}(0)}(\Sigma;\Sigma^{-2}x)}\leq C\|x\|_X$. Proposition \ref{Prop.approximation} and Corollary \ref{Cor.forward_energy} imply
	\begin{equation}\label{Eq.nonlinear_iteration_bound}
		\|\mathcal T_xV\|_{\mathfrak E_\Sigma} \leq C(\|x\|_X+\epsilon_\Sigma)+C\|V\|_{\mathfrak E_\Sigma}^2.
	\end{equation}
	Indeed, on $0\leq\tau\leq2$ the scaled residual satisfies
	\[\Sigma^2\|\sigma^2\cR_{\app}(\sigma)\|_{\mathcal B_{\varrho_{\rm F}(\tau)}} \leq C\Sigma^{-1}(1+\log\Sigma)^M=C\epsilon_\Sigma.\]
	For the quadratic term, the fixed algebra estimate and the energy lower bound give
	\[\|\cQ_\sigma(\zeta_V)\|_{\mathcal B_{\varrho_{\rm F}(\tau)}} \leq C\cE^\sharp_{\varrho_{\rm F}(\tau)}(\sigma;V).\]
	More explicitly, \eqref{Eq.zeta-w-W-def}, \eqref{Eq.energy_state_comparison} and \eqref{Eq.nonlinear_path_norm} imply
\[\|\zeta_V(\sigma)\|_{\mathcal B_{\varrho_{\rm F}(\tau)}}\leq C\Sigma^{-2}\|V\|_{\mathfrak E_\Sigma},\qquad \sigma=\Sigma\ee^\tau.\]
For $\|V\|_{\mathfrak E_\Sigma}\leq1$ and $\Sigma$ sufficiently large, \eqref{Eq.nonlinear_algebra_bounds} therefore gives
\begin{align*}
\Sigma^2\int_0^2\|\sigma^2\cQ_\sigma(\zeta_V)\|_{\mathcal B_{\varrho_{\rm F}(\tau)}}\,\mathrm d\tau\leq C\Sigma^2\int_0^2\sigma^2\Sigma^{-4}\|V\|_{\mathfrak E_\Sigma}^2\,\mathrm d\tau\leq C\|V\|_{\mathfrak E_\Sigma}^2.
\end{align*}
The large starting time makes the residual $\epsilon_\Sigma$ small, but does not make this quadratic constant small. The nonlinear contraction will therefore use a small fixed ball, independently of how late the construction starts. Apply the second inequality in \eqref{Eq.nonlinear_algebra_bounds} to the difference of the sources; Corollary \ref{Cor.forward_energy} yields
	\begin{equation}\label{Eq.nonlinear_iteration_difference}
		\|\mathcal T_xV-\mathcal T_x\widetilde V\|_{\mathfrak E_\Sigma} \leq C(\|V\|_{\mathfrak E_\Sigma}+\|\widetilde V\|_{\mathfrak E_\Sigma}) \|V-\widetilde V\|_{\mathfrak E_\Sigma}.
	\end{equation}
	Let $C_0,C_1\geq1$ be constants for the linear and quadratic bounds in \eqref{Eq.nonlinear_iteration_bound}--\eqref{Eq.nonlinear_iteration_difference}, and choose $r_*>0$ with $2C_1r_*\leq1/2$. Require $C_0(\|x\|_X+\epsilon_\Sigma)\leq r_*/2$. Then, on the closed ball $\|V\|_{\mathfrak E_\Sigma}\leq r_*$,
	\begin{align*}
		\|\mathcal T_xV\|_{\mathfrak E_\Sigma}&\leq r_*/2+C_1r_*^2\leq3r_*/4,\\
		\|\mathcal T_xV-\mathcal T_x\widetilde V\|_{\mathfrak E_\Sigma}&\leq2C_1r_*\|V-\widetilde V\|_{\mathfrak E_\Sigma}\leq\tfrac12\|V-\widetilde V\|_{\mathfrak E_\Sigma}.
	\end{align*}
	The contraction theorem gives a unique fixed point $W$. The inequality $\|W\|_{\mathfrak E_\Sigma}\leq C_0(\|x\|_X+\epsilon_\Sigma)+C_1r_*\|W\|_{\mathfrak E_\Sigma}$ implies $\|W\|_{\mathfrak E_\Sigma}\leq2C_0(\|x\|_X+\epsilon_\Sigma)$. The estimates apply first to bounded measurable paths; the linear solver gives continuity at smaller analytic radii, so $W$ is continuous at every smaller analytic radius. Real-valued radial paths are preserved.
	
	The positive lower bound on the analytic radius gives all finite spatial derivatives. The scalar equation and \eqref{Eq.exact_state} then give all time derivatives on compact intervals, successively using a smaller analytic radius when a spatial derivative occurs. Thus the solution is smooth up to the boundary. Uniqueness in the small energy ball follows from the same difference estimate.
\end{proof}

To estimate $N_\Sigma$ in \eqref{Eq.nonlinear_endpoint_bounds}, compare $W_x-W_y$ first. Subtracting the homogeneous difference $Z_x-Z_y$ leaves zero initial data and only the source $-\cQ_\sigma(\zeta_{W_x})+\cQ_\sigma(\zeta_{W_y})$, which gains the small factor $R+\epsilon_\Sigma$ through \eqref{Eq.nonlinear_algebra_bounds}.

\begin{lemma}\label{Lem.local_difference}
	Let $W_x,W_y$ be the solutions in Lemma \ref{Lem.local_nonlinear}, with $\|x\|_X,\|y\|_X\leq R$. For fixed small $R$ and sufficiently large $\Sigma$,
	\[\|W_x-W_y\|_{\mathfrak E_\Sigma}\leq C\|x-y\|_X.\]
	If $Z_x$ denotes the homogeneous linear solution with $Z_x(\Sigma)=\Sigma^{-2}x$, then
	\[\|(W_x-Z_x)-(W_y-Z_y)\|_{\mathfrak E_\Sigma} \leq C(R+\epsilon_\Sigma)\|x-y\|_X.\]
\end{lemma}

\begin{proof}
	Write $\Delta W:=W_x-W_y$. The initial value of $\Delta W$ is $\Sigma^{-2}(x-y)$ and the scalar source for $\zeta_{\Delta W}$ is $-\cQ_\sigma(\zeta_{W_x})+\cQ_\sigma(\zeta_{W_y})$. At the current analytic radius $\varrho=\varrho_{\rm F}(\tau)$, \eqref{Eq.nonlinear_algebra_bounds} gives
	\begin{align*}
		\|\cQ_\sigma(\zeta_{W_x})-\cQ_\sigma(\zeta_{W_y})\|_{\mathcal B_\varrho} &\leq C(\|\zeta_{W_x}\|_{\mathcal B_\varrho}+\|\zeta_{W_y}\|_{\mathcal B_\varrho})\|\zeta_{\Delta W}\|_{\mathcal B_\varrho}\\
		&\leq C\Sigma^{-4}(\|W_x\|_{\mathfrak E_\Sigma}+\|W_y\|_{\mathfrak E_\Sigma})\|\Delta W\|_{\mathfrak E_\Sigma}.
	\end{align*}
	Since $\sigma/\Sigma\leq \ee^2$, Corollary \ref{Cor.forward_energy} and Lemma \ref{Lem.local_nonlinear} imply
	\begin{align*}
		\|\Delta W\|_{\mathfrak E_\Sigma} &\leq C\|x-y\|_X+C\Sigma^2\int_0^2\sigma^2\|\cQ_\sigma(\zeta_{W_x})-\cQ_\sigma(\zeta_{W_y})\|_{\mathcal B_\varrho}\,\mathrm d\tau\\
		&\leq C\|x-y\|_X+C(R+\epsilon_\Sigma)\|\Delta W\|_{\mathfrak E_\Sigma}.
	\end{align*}
	Choose $R$ and the initial threshold so that $C(R+\epsilon_\Sigma)\leq1/2$. Then $\|\Delta W\|_{\mathfrak E_\Sigma}\leq2C\|x-y\|_X$.
	
	For $H=(W_x-Z_x)-(W_y-Z_y)$ we have $H(\Sigma)=0$ and $\cL_\sigma\zeta_H=-\cQ_\sigma(\zeta_{W_x})+\cQ_\sigma(\zeta_{W_y})$. A second application of Corollary \ref{Cor.forward_energy} therefore gives
	\[\|H\|_{\mathfrak E_\Sigma}\leq C(R+\epsilon_\Sigma)\|\Delta W\|_{\mathfrak E_\Sigma}\leq C(R+\epsilon_\Sigma)\|x-y\|_X.\]
	The constants in the bounds for $\Delta W$ and $H$ are independent of $x,y\in B_X(R)$.
\end{proof}

The endpoint lies at an analytic radius strictly larger than $\varrho_0$. Restriction therefore turns the local energy estimates into estimates for a map on the single fixed space $X$.

\begin{proof}[Proof of Proposition \ref{Prop.nonlinear_endpoint}]
	Lemma \ref{Lem.local_nonlinear} gives the solution and \eqref{Eq.local_nonlinear_energy}. At $\sigma=\ee\Sigma$, the energy lower bound and $\varrho_{\rm F}(1)>\varrho_0$ give, for any path $V$,
	\[\|(\ee\Sigma)^2V(\ee\Sigma)\|_X \leq C\ee^2\Sigma^2\sqrt{\cE^\sharp_{\varrho_{\rm F}(1)}(\ee\Sigma;V(\ee\Sigma))} \leq C\|V\|_{\mathfrak E_\Sigma}.\]
	For the homogeneous solution, $(\ee\Sigma)^2Z_x(\ee\Sigma)=\ee^2\mathcal P_\Sigma(1)x$. Define $N_\Sigma(x)=(\ee\Sigma)^2(W_x-Z_x)(\ee\Sigma)$. Lemma \ref{Lem.local_difference} yields
	\begin{align*}
		\|N_\Sigma(x)-N_\Sigma(y)\|_X &\leq C\|(W_x-Z_x)-(W_y-Z_y)\|_{\mathfrak E_\Sigma}\leq C(R+\epsilon_\Sigma)\|x-y\|_X.
	\end{align*}
	For the solution with $x=0$, $Z_0=0$ and Lemma \ref{Lem.local_nonlinear} gives $\|N_\Sigma(0)\|_X\leq C\epsilon_\Sigma$. Together with $F_\Sigma(x)=\ee^2\mathcal P_\Sigma(1)x+N_\Sigma(x)$ and the difference estimate, this proves \eqref{Eq.nonlinear_endpoint_bounds}.
\end{proof}

For $\sigma_n=\Sigma\ee^n$, the normalization $x_n=\sigma_n^2W(\sigma_n)$ makes the linear endpoint approach $K=\ee^2\cS(1)$, rather than $\cS(1)$. If $\|x_n\|_X\lesssim b^n$ with $b=\ee^{2-a}$, then $\|W(\sigma_n)\|_X\lesssim\Sigma^{-2}\ee^{-an}$. We choose $2<a<3$: $a>2$ gives $b<1$, while $a<3$ gives $b>\ee^{-1}$, so the forcing $\epsilon_{\sigma_n}=\sigma_n^{-1}(1+\log\sigma_n)^M$ belongs to the required weighted sequence space. Compactness of $K$ allows the circle $|z|=b$ to avoid the spectrum. No spectral stability of the full solution operator is assumed.

\begin{proof}[Proof of Theorem \ref{Thm.state}]
	The operator $K=\ee^2\mathsf A$ is compact on $X$. The spectrum of $K$ has only finitely many points in the annulus $\ee^{-1}\leq|z|\leq1$. Choose $a\in(2,3)$ so that $\operatorname{spec}(K)\cap\{|z|=b\}=\varnothing$, where $b=\ee^{2-a}$, and fix the constants in Proposition \ref{Prop.compact_decay}. Choose $R\leq R_*$ sufficiently small, and then $\Sigma$ sufficiently large. For $\sigma_n=\Sigma\ee^n$, we define
	\[R_n(x)=\ee^2(\mathcal P_{\sigma_n}(1)-\mathsf A)x+N_{\sigma_n}(x).\]
	Propositions \ref{Prop.averaging} and \ref{Prop.nonlinear_endpoint} give
	\[\operatorname{Lip}(R_n|_{B_X(R)}) \leq C\mu(\Sigma)+C(R+\sup_{n\geq0}\epsilon_{\Sigma\ee^n}).\]
	Choose $R$ and $\Sigma$ so that $C\mu(\Sigma)+C(R+\sup_{n\geq0}\epsilon_{\Sigma\ee^n})\leq\epsilon_0$, where $\epsilon_0$ is the constant in Proposition \ref{Prop.compact_decay}. Moreover,
	\[\|R_n(0)\|_X\leq C \Sigma^{-1}\ee^{-n}(1+\log \Sigma+n)^M \leq C_{b,M}\Sigma^{-1}(1+\log \Sigma)^M b^n.\]
	The last inequality follows because $b>\ee^{-1}$ and $\sup_{n\geq0}(\ee^{-1}/b)^n(1+n)^M<\infty$. Increasing $\Sigma$ once more gives $2C_Gh<R$, where $h=C_{b,M}\Sigma^{-1}(1+\log \Sigma)^M$.
	
	Proposition \ref{Prop.compact_decay} yields a sequence with $\|x_n\|_X\leq2C_Ghb^n$ and $x_{n+1}=F_{\sigma_n}(x_n)$. The Green formula \eqref{Eq.discrete_green_proof} makes the choice of the initial data explicit:
\[P_-x_0=0,\qquad P_+x_0=-\sum_{j=0}^\infty K_+^{-j-1}P_+R_j(x_j).\]
The series converges by \eqref{Eq.abstract_spectral_powers}, since $\|R_j(x_j)\|_X\leq C h b^j$. Thus the components in $P_+X$ are selected by the decay requirement; no arbitrary initial value is prescribed in these directions. In particular, the construction gives one decaying solution and does not claim decay for every small initial datum.

Solve the local nonlinear equation with data $W(\sigma_n)=\sigma_n^{-2}x_n$. The identity $x_{n+1}=F_{\sigma_n}(x_n)$ gives $W(\sigma_{n+1})=\sigma_{n+1}^{-2}x_{n+1}$. Formula \eqref{Eq.zeta-w-W-def} and $D\zeta=\lambda w$ from Lemma \ref{Lem.state_equation} identify $\zeta$ and $D\zeta$ at each common endpoint. Uniqueness for the scalar equation gives a single smooth solution on $[\Sigma,\infty)$. On $[\sigma_n,\sigma_{n+1}]$, the local estimate gives
	\[\|W(\sigma)\|_{\cX_{\varrho_-}} \leq C\sigma_n^{-2}(\|x_n\|_X+\epsilon_{\sigma_n}) \leq C \Sigma^{-3}(1+\log \Sigma)^M \ee^{-an}.\]
	Here $\varrho_{\rm F}(\tau)\geq\varrho_-:=\varrho_{\rm F}(0)$, and $\ee^{-2n}b^n=\ee^{-an}$. Since $\sigma/\Sigma\in[\ee^n,\ee^{n+1}]$, we obtain \eqref{Eq.state_decay}. Equation \eqref{Eq.remainder_equation} yields $\cN(z_{\app}+\zeta)=0$. The solution is real-valued and radial on every interval and hence on $[\Sigma,\infty)$.
\end{proof}

\subsection{The asymptotic behavior in the cone}
\label{Subsec.cone_rates}

The bounds \eqref{Eq.approximation_bounds} and \eqref{Eq.state_decay} give $u^{\mathrm c}/\rho_0=f(\theta_\sigma)+o(1)$. An upper bound for $|u^{\mathrm c}|$ follows from $|f|\leq1$. For the lower bound and the signed limits in \eqref{Eq.peaks}, we show that the physical phase varies by at least one full period on radii $|y|\lesssim(T-t)^{(\beta+1)/2}$. Fix the initial threshold $\Sigma$ and choose $t_2\in[t_0,T)$ with $T-t_2\leq1/2$. We use radial notation for $u^{\mathrm c}$ and $\rho_0$. For a fixed $\omega\in\mathbb S^{d-1}$, define
\begin{equation}\label{Eq.physical_phase_definition}
	\Theta(t,Z):=(T-t)^{1-\beta}\widehat\phi_0(Z)+\frac{(T-t)^{\beta-1}}{\widehat\phi_0(Z)}\Psi\!\left(\log\bigl((T-t)^{1-\beta}\widehat\phi_0(Z)\bigr),X(Z)\omega\right)
\end{equation}
for $0\leq Z\leq1$. By \eqref{Eq.sigma_z} and \eqref{Eq.approximate_phase}, $\Theta(t,Z)$ is the value of $\theta_\sigma$ at $y=(T-t)Z\omega$. Radiality of $\Psi$ makes $\Theta$ independent of $\omega$.

\begin{lemma}\label{Lem.cone_phase_variation}
	There are $C_1,C_2,C_R>0$ such that, for $t_2\leq t<T$ and $0\leq Z\leq1$,
	\begin{equation}\label{Eq.physical_phase_asymptotic}
		\left|\rho_0(t,(T-t)Z)^{-1}u^{\mathrm c}(t,(T-t)Z)-f(\Theta(t,Z))\right|\leq C_2(T-t)^{2\beta-2}(1+|\log(T-t)|)^M.
	\end{equation}
	Moreover, $|\widehat\rho_0(Z)-1|\leq C_RZ^2$, and
	\begin{equation}\label{Eq.physical_phase_variation}
		\begin{aligned}
			\Theta(t,Z)-\Theta(t,0)&=\frac{\beta\ell}{2d}(T-t)^{1-\beta}Z^2(1+E(t,Z)),\\
			|E(t,Z)|&\leq C_1\left(Z^2+(T-t)^{2\beta-2}(1+|\log(T-t)|)^M\right).
		\end{aligned}
	\end{equation}
	The exponent $M$ is the fixed integer in \eqref{Eq.nu_sigma}.
\end{lemma}

\begin{proof}
	On the cone, $\sigma=(T-t)^{1-\beta}\widehat\phi_0(Z)\geq\Sigma$. The positive lower bound for $\widehat\phi_0$ on $[0,1]$ implies
	\[\sigma^{-2}\leq C(T-t)^{2\beta-2},\qquad 1+|\log\sigma|\leq C(1+|\log(T-t)|).\]
	The reconstruction \eqref{Eq.zeta-w-W-def}, \eqref{Eq.state_decay} and Sobolev embedding give $|\zeta|\leq C_\Sigma\sigma^{-a}$ with $a>2$. Adding this bound to \eqref{Eq.approximation_bounds} proves \eqref{Eq.physical_phase_asymptotic}.
	
	The profile formulas \eqref{Eq.profiles}, with $v'(0)=\beta\ell/d$ from \eqref{Eq.velocity}, give
	\[\widehat\phi_0(Z)=\frac1{\beta-1}+\frac{\beta\ell}{2d}Z^2+\Order(Z^4),\qquad \widehat\rho_0(Z)=1+\Order(Z^2).\]
	Taylor's formula on $[0,1]$, evenness of both profiles and $X(0)=0$ give
	\begin{equation}\label{Eq.profile_Taylor_bounds}
		\begin{aligned}
			&\left|\widehat\phi_0(Z)-\widehat\phi_0(0)-\frac{\beta\ell}{2d}Z^2\right|\leq CZ^4,\qquad |\widehat\rho_0(Z)-1|\leq C_RZ^2,\\
			&\qquad|X(Z)|\leq CZ,\qquad |\log\widehat\phi_0(Z)-\log\widehat\phi_0(0)|\leq CZ^2.
		\end{aligned}
	\end{equation}
	The last bound uses $\min_{[0,1]}\widehat\phi_0>0$. Write $\tau_Z=\log((T-t)^{1-\beta}\widehat\phi_0(Z))$. Since $\Psi$ is polynomial in $\tau$ with analytic radial coefficients, Proposition \ref{Prop.phase} gives
	\[|\Psi(\tau,x)|+|\partial_\tau\Psi(\tau,x)|+|\nabla_x^2\Psi(\tau,x)|\leq C(1+|\tau|)^M,\qquad \nabla_x\Psi(\tau,0)=0.\]
	Taylor expansion in $x$ and the bound $|\tau_Z-\tau_0|\leq CZ^2$ from \eqref{Eq.profile_Taylor_bounds} imply
	\[|\Psi(\tau_Z,X(Z)\omega)-\Psi(\tau_0,0)|\leq CZ^2(1+|\log(T-t)|)^M.\]
	Also $|\widehat\phi_0(Z)^{-1}-\widehat\phi_0(0)^{-1}|\leq CZ^2$. Thus the difference between the second terms of \eqref{Eq.physical_phase_definition} at $Z$ and at zero is bounded by $C(T-t)^{\beta-1}Z^2(1+|\log(T-t)|)^M$. Combining with \eqref{Eq.profile_Taylor_bounds} and dividing by $(\beta\ell/(2d))(T-t)^{1-\beta}Z^2$ proves \eqref{Eq.physical_phase_variation} for $Z>0$. The quotient extends continuously to $Z=0$ by evenness and smoothness.
\end{proof}

The identity in \eqref{Eq.physical_phase_variation} gives a phase change of order one when $Z$ has size $(T-t)^{(\beta-1)/2}$. We choose two phase values separated by $L$, at which $f$ takes the values $1$ and $-1$, and use \eqref{Eq.physical_phase_asymptotic} to control the error.

\begin{proof}[Proof of the cone estimates in Remark \ref{Rmk.cone}]
	Choose $K_*=(16Ld/(\beta\ell))^{1/2}$ and let \\ $M_R=\max_{[0,1]}\widehat\rho_0\geq1$. Since $\beta>1$, choose $t_1\in[t_2,T)$ sufficiently close to $T$ that, for every $t_1\leq t<T$,
	\begin{equation}\label{Eq.peak_threshold}
		\begin{gathered}
			K_*(T-t)^{(\beta-1)/2}\leq1,\\
			C_1\left(K_*^2(T-t)^{\beta-1}+(T-t)^{2\beta-2}(1+|\log(T-t)|)^M\right)\leq\frac12,\\
			C_RK_*^2(T-t)^{\beta-1}+M_RC_2(T-t)^{2\beta-2}(1+|\log(T-t)|)^M\leq\frac12.
		\end{gathered}
	\end{equation}
	Every expression on the left that must be small tends to zero as $t\uparrow T$, so a single $t_1$ suffices. By \eqref{Eq.physical_phase_variation},
	\[\Theta\bigl(t,K_*(T-t)^{(\beta-1)/2}\bigr)-\Theta(t,0)\geq\frac{\beta\ell}{4d}K_*^2=4L.\]
	Choose the integer $m(t)=\lceil\Theta(t,0)/(2L)\rceil$. Both $2Lm(t)$ and $2Lm(t)+L$ lie between the two endpoint phase values. The intermediate value theorem gives $Z_\pm(t)\in[0,K_*(T-t)^{(\beta-1)/2}]$ with
	\[\Theta(t,Z_+(t))=2Lm(t),\qquad \Theta(t,Z_-(t))=2Lm(t)+L.\]
	No monotonicity of $\Theta(t,\cdot)$ is required. Since $f(2Lm)=1$ and $f(2Lm+L)=-1$, define $y_\pm(t)=(T-t)Z_\pm(t)\omega$. Equations \eqref{Eq.selfsimilar}, \eqref{Eq.physical_phase_asymptotic} and \eqref{Eq.profile_Taylor_bounds} give
	\begin{align*}
		|y_\pm(t)|&\leq K_*(T-t)^{(\beta+1)/2},\\
		\left|(T-t)^{\frac{2\beta}{p-1}}u^{\mathrm c}(t,y_\pm(t))\mp1\right|&\leq C_RK_*^2(T-t)^{\beta-1}+M_RC_2(T-t)^{2\beta-2}(1+|\log(T-t)|)^M.
	\end{align*}
	The last bound tends to zero and is at most $1/2$ by \eqref{Eq.peak_threshold}. It proves \eqref{Eq.peaks} and the lower estimate in \eqref{Eq.cone_rate}. On the entire cone, \eqref{Eq.physical_phase_asymptotic} and $|f|\leq1$ give
	\[|u^{\mathrm c}(t,y)|\leq M_R(T-t)^{-\frac{2\beta}{p-1}}\left(1+C_2(T-t)^{2\beta-2}(1+|\log(T-t)|)^M\right)\leq2M_R(T-t)^{-\frac{2\beta}{p-1}}.\]
	This proves the upper estimate in \eqref{Eq.cone_rate}. The limits \eqref{Eq.beta_limits} follow from Proposition \ref{Prop.quantitative_background}.
\end{proof}

\section{Quantitative backgrounds and the large-exponent theorem}
\label{Sec.background_proof}

We prove $\chi_p\lesssim p^{-2}$ and construct Euler profiles with an outgoing lower bound independent of $p$ and a coefficient analytic radius $R_p\asymp p^{-3/2}$. The oscillator calculation proves \eqref{Eq.chi_p_bound} independently of the profile construction. The Euler argument starts at $\ell=1$, matches the origin and sonic solutions, and retains uniform control of the inverse coordinate used in \eqref{Eq.x_b}. Taking $\varrho_p$ to be a fixed small fraction of $R_p$ then gives $\chi_p/\varrho_p=O(p^{-1/2})$, as required by \eqref{Eq.analytic_margin}.

\subsection{A uniform one-period scattering estimate}
\label{Subsec.scattering_uniform}

We rescale one period of $h''+(pf^{p-1}+\xi^2)h=0$ to $\partial_s^2F+(Q_\zeta+\kappa^2)F=0$, with $\zeta=(p+1)^{-1}$. At $\zeta=0$, $Q_0(s)=2\operatorname{sech}^2s$ has the explicit solution \eqref{Eq.limiting_jost_solution}, which has no reflected wave. The solution $r_\zeta$ in \eqref{Eq.scattering_zero_solution} controls the additional degeneracy at $\kappa=0$. Uniform differentiation of the integral equation \eqref{Eq.jost_volterra} then yields the factor $\zeta\min\{\kappa^2,1\}$ in \eqref{Eq.reflection_numerator_bound}.

For the scattering calculation only, define
\begin{equation}\label{Eq.scattering_rescaling}
	\zeta:=(p+1)^{-1},\qquad \kappa:=\sqrt{\frac2{p+1}}\,\xi,\qquad s:=\sqrt{\frac{p+1}{2}}\,\theta.
\end{equation}
The oscillator maximum is at $s=0$. Define
\begin{equation}\label{Eq.scattering_coordinate}
	X_\zeta(w):=\int_0^w\operatorname{sech}^{2\zeta}t\,\mathrm dt,\qquad \Lambda_\zeta:=2X_\zeta(\infty)=\operatorname B(\zeta,1/2).
\end{equation}
The potential restricted to one period and extended by zero is
\begin{equation}\label{Eq.scattering_potential}
	Q_\zeta(X_\zeta(w)):=2(1-\zeta)\operatorname{sech}^{2-4\zeta}w\quad(w\in\R),\qquad Q_\zeta(s):=0\quad(|s|\geq\Lambda_\zeta/2).
\end{equation}
At $\zeta=0$ set $X_0(w):=w$, $\Lambda_0:=\infty$ and $Q_0(s):=2\operatorname{sech}^2s$. For $\zeta=(p+1)^{-1}$, the change $f^{p+1}=\operatorname{sech}^2w$ and \eqref{Eq.oscillator_energy} give \eqref{Eq.scattering_coordinate}--\eqref{Eq.scattering_potential}. In particular, $Q_\zeta(s)=2p(p+1)^{-1}f(\theta)^{p-1}$ inside the rescaled period and $L\xi=\Lambda_\zeta\kappa$. The Euler parameter $\varepsilon=\ell-1$ used later is different from $\zeta$.

\begin{lemma}\label{Lem.cell_potential}
	For $0\leq\zeta\leq1/16$, the functions $Q_\zeta$ are even and nonnegative. For every integer $m\geq0$ there is $C_m$ such that
	\begin{equation}\label{Eq.scattering_moments}
		\int_\R(1+|s|)^m\bigl(|Q_\zeta(s)|+|\partial_\zeta Q_\zeta(s)|\bigr)\,\mathrm ds\leq C_m.
	\end{equation}
	The map $\zeta\mapsto Q_\zeta$ is continuously differentiable in $L^1(\R,(1+|s|)^m\,\mathrm ds)$ for every integer $m\geq0$, including at zero. The function defined by
	\begin{equation}\label{Eq.scattering_zero_solution}
		r_\zeta(X_\zeta(w)):=\tanh w
	\end{equation}
	and extended as the constants $\pm1$ outside $[-\Lambda_\zeta/2,\Lambda_\zeta/2]$ solves $r_\zeta''+Q_\zeta r_\zeta=0$, satisfies $r_\zeta(0)=0$, and tends to one at positive infinity.
\end{lemma}

\begin{proof}
	Equation \eqref{Eq.oscillator_energy} and $f^{p+1}=\operatorname{sech}^2w$ give $\mathrm d\theta/\mathrm dw=\sqrt{2/(p+1)}\operatorname{sech}^{2/(p+1)}w$, proving the rescaling \eqref{Eq.scattering_coordinate}--\eqref{Eq.scattering_potential}. Differentiating \eqref{Eq.scattering_zero_solution} and using $X_\zeta'=\operatorname{sech}^{2\zeta}w$ gives
	\[r_\zeta'(X_\zeta(w))=\operatorname{sech}^{2-2\zeta}w,\qquad r_\zeta''(X_\zeta(w))=-2(1-\zeta)\tanh w\operatorname{sech}^{2-4\zeta}w.\]
	The first derivative vanishes at the endpoints, so extension by constants introduces no distributional boundary term. The equation and the normalization follow.
	
	To differentiate at fixed $s$, use $s=X_\zeta(w)$ from \eqref{Eq.scattering_coordinate}; implicit differentiation gives
	\begin{equation}\label{Eq.scattering_inverse_parameter}
		\partial_\zeta w=2\cosh^{2\zeta}w\int_0^w\log(\cosh t)\operatorname{sech}^{2\zeta}t\,\mathrm dt.
	\end{equation}
	For $w\geq0$, $|\partial_\zeta w|$ is at most $w^2\cosh^{2\zeta}w$. Differentiate the expression $2(1-\zeta)\operatorname{sech}^{2-4\zeta}w_\zeta(s)$ from \eqref{Eq.scattering_potential} at fixed $s$, using \eqref{Eq.scattering_inverse_parameter}. We obtain
	\begin{equation}\label{Eq.scattering_parameter_bound}
		|\partial_\zeta Q_\zeta(X_\zeta(w))|\leq C(1+w)^2\operatorname{sech}^{2-6\zeta}w.
	\end{equation}
	Since $|X_\zeta(w)|\leq|w|$ and $\,\mathrm ds=\operatorname{sech}^{2\zeta}w\,\mathrm dw$, the weighted integrals are bounded by integrals of a polynomial times $\ee^{-c|w|}$, with one fixed $c>0$. Near a finite endpoint $Q_\zeta$ vanishes to order $1/\zeta-2\geq14$. The vanishing of $Q_\zeta$ at both endpoints implies that $\partial_\zeta Q_\zeta$ is the zero extension of the interior parameter derivative.
	
	Fix $\zeta_0\in[0,1/16]$. If $\zeta_0>0$, the inverse coordinate $w_\zeta(s)$ and $\partial_\zeta w_\zeta(s)$ depend continuously on $\zeta$ on every compact subinterval strictly inside $(-\Lambda_{\zeta_0}/2,\Lambda_{\zeta_0}/2)$, for $\zeta$ sufficiently close to $\zeta_0$. Outside $[-\Lambda_{\zeta_0}/2,\Lambda_{\zeta_0}/2]$, both $Q_\zeta$ and $\partial_\zeta Q_\zeta$ vanish for $\zeta$ sufficiently close to $\zeta_0$. If $\zeta_0=0$, every fixed compact $s$ interval lies inside $(-\Lambda_\zeta/2,\Lambda_\zeta/2)$ for all sufficiently small positive $\zeta$, and the inverse-coordinate formulas converge there. Therefore $Q_\zeta$ and $\partial_\zeta Q_\zeta$ converge pointwise for almost every $s$ as $\zeta\to\zeta_0$.
	
	Equations \eqref{Eq.scattering_potential} and \eqref{Eq.scattering_parameter_bound} give a common bound on $Q_\zeta$ and $\partial_\zeta Q_\zeta$ on each compact $s$ interval: indeed, $(1+|w|)^2\operatorname{sech}^{2-6\zeta}w$ is uniformly bounded for $0\leq\zeta\leq1/16$. Dominated convergence therefore applies to the zero extensions on compact intervals. Moreover, $|s|>A$ implies $|w|>A$, and the weighted exponential tail estimates tend to zero uniformly as $A\to\infty$. Thus both $Q_\zeta$ and $\partial_\zeta Q_\zeta$ are continuous in each weighted $L^1$ space. The endpoint vanishing established above eliminates the boundary term when the zero extension is differentiated. Thus, writing $L^1_m:=L^1(\R,(1+|s|)^m\,\mathrm ds)$, we have
	\[\lim_{h\to0}\left\|\frac{Q_{\zeta+h}-Q_\zeta}{h}-\partial_\zeta Q_\zeta\right\|_{L^1_m}=0,\]
	with one-sided limits at the endpoints of $[0,1/16]$; the derivative is continuous in $L^1_m$.
\end{proof}

The bounds \eqref{Eq.scattering_moments} control the powers of $t-s$ obtained by differentiating the kernel $D_\kappa(t-s)$ in \eqref{Eq.jost_volterra}. We solve the integral equation for $m$ on $s\geq0$ and retain bounds uniform in $\kappa\in\R$ and $\zeta\in[0,1/16]$.

\begin{lemma}\label{Lem.jost_uniform}
	For $0\leq\zeta\leq1/16$ and real $\kappa$, there is a unique bounded solution on $s\geq0$ of
	\begin{equation}\label{Eq.jost_volterra}
		m(s,\kappa)=1-\int_s^\infty D_\kappa(t-s)Q_\zeta(t)m(t,\kappa)\,\mathrm dt,\qquad D_\kappa(a):=\int_0^a\ee^{2i\kappa v}\,\mathrm dv.
	\end{equation}
	For $a\in\{0,1\}$ and $j\in\{0,1,2\}$, the functions $\partial_\zeta^a\partial_\kappa^jm$ and $\partial_\zeta^a\partial_\kappa^j\partial_s m$ are uniformly bounded for $s\geq0$, $0\leq\zeta\leq1/16$ and $\kappa\in\R$. The function $F_+(s,\kappa):=\ee^{i\kappa s}m(s,\kappa)$ solves $\partial_s^2F_++(\kappa^2+Q_\zeta)F_+=0$ with $\ee^{-i\kappa s}F_+(s,\kappa)\to1$ and $\ee^{-i\kappa s}(\partial_sF_+-i\kappa F_+)(s,\kappa)\to0$ as $s\to\infty$.
\end{lemma}

\begin{proof}
	For a bounded function $g$ on $[0,\infty)$, define $(K_\zeta(\kappa)g)(s)=\int_s^\infty D_\kappa(t-s)Q_\zeta(t)g(t)\,\mathrm dt$. Then \eqref{Eq.jost_volterra} is $(I+K_\zeta(\kappa))m=1$. Since $|D_\kappa(t-s)|\leq t-s\leq t$, $K_\zeta(\kappa)^N$ has operator norm at most $M_1^N/N!$ on bounded functions, where $M_1:=\sup_\zeta\int_0^\infty tQ_\zeta(t)\,\mathrm dt<\infty$ by \eqref{Eq.scattering_moments}. The bound follows by integrating over the ordered simplex $s\leq t_1\leq\cdots\leq t_N$. Thus $(I+K_\zeta(\kappa))^{-1}=\sum_{N\geq0}(-K_\zeta(\kappa))^N$ converges in operator norm and is bounded by $\ee^{M_1}$. Therefore $m=(I+K_\zeta(\kappa))^{-1}1$ is the unique bounded solution and $\|m\|_{L^\infty}\leq\ee^{M_1}$.
	
	For $j\leq2$, one has $|\partial_\kappa^jD_\kappa(a)|\leq2^ja^{j+1}/(j+1)$. Write $K=K_\zeta(\kappa)$. For $a\in\{0,1\}$, $j\in\{0,1,2\}$ and $a+j>0$, differentiation of \eqref{Eq.jost_volterra} gives
	\begin{equation}\label{Eq.jost_parameter_equation}
		(I+K)\partial_\zeta^a\partial_\kappa^jm=-\sum_{\substack{0\leq b\leq a,\ 0\leq r\leq j\\
				b+r>0}}\binom ab\binom jr(\partial_\zeta^b\partial_\kappa^rK)\partial_\zeta^{a-b}\partial_\kappa^{j-r}m.
	\end{equation}
	More precisely,
	\begin{align*}
		\|\partial_\zeta^b\partial_\kappa^rK\|_{L^\infty\to L^\infty} &\leq \frac{2^r}{r+1}\int_0^\infty t^{r+1}|\partial_\zeta^bQ_\zeta(t)|\,\mathrm dt\leq C_r,\\
		\|\partial_\zeta^a\partial_\kappa^jm\|_{L^\infty} &\leq \ee^{M_1}\sum_{\substack{0\leq b\leq a,\ 0\leq r\leq j\\
				b+r>0}} \binom ab\binom jr C_r\|\partial_\zeta^{a-b}\partial_\kappa^{j-r}m\|_{L^\infty}.
	\end{align*}
	Induction on $a+j$ in \eqref{Eq.jost_parameter_equation}, starting with $\|m\|_{L^\infty}\leq\ee^{M_1}$, gives $\|\partial_\zeta^a\partial_\kappa^jm\|_{L^\infty}\leq C$ for $a\leq1$ and $j\leq2$. Difference quotients in $(I+K)m=1$ justify the differentiations in the bounded-function norm. Differentiation in $s$ gives
	\begin{equation}\label{Eq.jost_spatial_derivative}
		\partial_s m(s,\kappa)=\int_s^\infty \ee^{2i\kappa(t-s)}Q_\zeta(t)m(t,\kappa)\,\mathrm dt.
	\end{equation}
	Differentiating \eqref{Eq.jost_spatial_derivative} gives
	\begin{align*}
		&\|\partial_\zeta^a\partial_\kappa^j\partial_s m\|_{L^\infty} \leq \sum_{b=0}^a\sum_{r=0}^j\binom ab\binom jr 2^r\left(\int_0^\infty t^r|\partial_\zeta^bQ_\zeta(t)|\,\mathrm dt\right) \|\partial_\zeta^{a-b}\partial_\kappa^{j-r}m\|_{L^\infty}\leq C,\\
		&\qquad\qquad\qquad\partial_s^2m+2i\kappa \partial_s m=-Q_\zeta m,\qquad \partial_s^2F_++(\kappa^2+Q_\zeta)F_+=0.
	\end{align*}
	Finally,
	\[|m(s,\kappa)-1|+|\partial_s m(s,\kappa)|\leq C\int_s^\infty(1+t)Q_\zeta(t)\,\mathrm dt\longrightarrow0.\]
	Hence $\ee^{-i\kappa s}F_+(s,\kappa)\to1$ and $\ee^{-i\kappa s}(\partial_sF_+-i\kappa F_+)(s,\kappa)\to0$, also for $\zeta=0$.
\end{proof}

The solution $r_\zeta$ in \eqref{Eq.scattering_zero_solution} yields $m(0,0)=0$, while reality makes the numerator \eqref{Eq.reflection_numerator} even in $\kappa$. Thus the numerator vanishes to second order at zero frequency. The explicit formula \eqref{Eq.limiting_jost_solution} gives a further vanishing at $\zeta=0$; together, the cancellations produce the uniform factor $\zeta\min\{\kappa^2,1\}$.

\begin{lemma}[The uniform scattering bound]\label{Lem.small_scattering}
	There is $C_{\rm sc}>0$ independent of odd $p\geq15$ such that
	\[0<\chi_p\leq C_{\rm sc}(p+1)^{-2}.\]
\end{lemma}

\begin{proof}
	Set
	\begin{equation}\label{Eq.reflection_numerator}
		\mathcal N_\zeta(\kappa):=2\operatorname{Re}\{\overline{m(0,\kappa)}\partial_s m(0,\kappa)\}.
	\end{equation}
	Since the potential is real-valued, we have $m(s,-\kappa)=\overline{m(s,\kappa)}$, so $\mathcal N_\zeta$ is even in $\kappa$. At $\kappa=0$, the function $r_\zeta$ from \eqref{Eq.scattering_zero_solution} satisfies \eqref{Eq.jost_volterra}, as follows by integrating $r_\zeta''=-Q_\zeta r_\zeta$ twice from infinity with $r_\zeta(\infty)=1$ and $r_\zeta'(\infty)=0$. Uniqueness in Lemma \ref{Lem.jost_uniform} therefore gives $m(s,0)=r_\zeta(s)$. Thus $m(0,0)=0$ and $\mathcal N_\zeta(0)=\mathcal N_\zeta'(0)=0$. At $\zeta=0$, direct substitution gives
	\begin{equation}\label{Eq.limiting_jost_solution}
		F_+(s,\kappa)=\ee^{i\kappa s}\frac{\kappa+i\tanh s}{\kappa+i}.
	\end{equation}
	At $s=0$, \eqref{Eq.limiting_jost_solution} gives $m(0,\kappa)=\kappa/(\kappa+i)$ and $\partial_sm(0,\kappa)=i/(\kappa+i)$. Their product in \eqref{Eq.reflection_numerator} has real part zero, so $\mathcal N_0(\kappa)=0$ for all real $\kappa$.
	
	By Lemma \ref{Lem.jost_uniform},
	\[\sup_{0\leq u\leq1/16,\ \kappa\in\R} \left(|\partial_u\mathcal N_u(\kappa)|+|\partial_u\partial_\kappa^2\mathcal N_u(\kappa)|\right)\leq C.\]
	Since $\mathcal N_0\equiv0$ and $\mathcal N_u(0)=\partial_\kappa\mathcal N_u(0)=0$, Taylor's formula gives
	\begin{equation}\label{Eq.reflection_numerator_bound}
		\begin{aligned}
			\mathcal N_\zeta(\kappa) &=\int_0^\zeta\partial_u\mathcal N_u(\kappa)\,\mathrm du =\int_0^\zeta\int_0^\kappa(\kappa-v)\partial_u\partial_v^2\mathcal N_u(v)\,\mathrm dv\,\mathrm du,\\
			|\mathcal N_\zeta(\kappa)|&\leq C\zeta\min\{\kappa^2,1\}.
		\end{aligned}
	\end{equation}
	For $0<\zeta\leq1/16$ and $\kappa>0$, continue $F_+$ to all real $s$ by $F_+''+(\kappa^2+Q_\zeta)F_+=0$ and define $F_-(s,\kappa)=F_+(-s,\kappa)$. Evenness of $Q_\zeta$ makes $F_-$ a solution of the same equation. Outside the support of $Q_\zeta$, we have
\[F_+(s,\kappa)=\ee^{i\kappa s}\quad(s>\Lambda_\zeta/2),\qquad F_+(s,\kappa)=a(\kappa)\ee^{i\kappa s}+b(\kappa)\ee^{-i\kappa s}\quad(s<-\Lambda_\zeta/2).\]
For the Wronskian $W(g,h):=gh'-g'h$, evaluation to the left of the potential gives
\[W(\overline{F_-},F_+)=-2i\kappa b.\]
At zero, $F_-(0)=F_+(0)=m(0,\kappa)$ and $F_-'(0)=-F_+'(0)$, so
\begin{align*}
W(\overline{F_-},F_+)(0)
&=2\operatorname{Re}\bigl(\overline{F_+(0)}F_+'(0)\bigr)\\
&=2\operatorname{Re}\bigl(\overline{m(0,\kappa)}(i\kappa m(0,\kappa)+\partial_s m(0,\kappa))\bigr)=\mathcal N_\zeta(\kappa).
\end{align*}
Constancy of the Wronskian and \eqref{Eq.reflection_numerator_bound} therefore give
	\[|b(\kappa)|\leq C\zeta\min\{|\kappa|,|\kappa|^{-1}\}.\]
	We define $V_{\rm pw}:=2^{-1/2}\left(\begin{smallmatrix}1&1\\ \mathrm i&-\mathrm i\end{smallmatrix}\right)$. The off-diagonal entry of $V_{\rm pw}^*D_\xi^{1/2}\mathsf M D_\xi^{-1/2}V_{\rm pw}$ is
	\[\frac12\bigl(\mathsf M_{11}-\mathsf M_{22}-\mathrm i(\xi \mathsf M_{12}+\xi^{-1}\mathsf M_{21})\bigr).\]
	Since $\mathsf M_{11}=\mathsf M_{22}$, the modulus of the off-diagonal entry is $|b_{{\rm sc},p}(\xi)|$. Changing the reference points of the plane waves multiplies the off-diagonal entry by a unit-modulus factor. The scaling \eqref{Eq.scattering_rescaling} gives $\diag(\xi,1)(h,\partial_\theta h)^T=\sqrt{(p+1)/2}\diag(\kappa,1)(h,\partial_s h)^T$, so the balanced transfer matrices coincide. To check the orientation of the transfer, write a free solution as $A_+\ee^{i\kappa s}+A_-\ee^{-i\kappa s}$. The solution $F_+$ has right coefficients $(1,0)^T$ and left coefficients $(a,b)^T$; its conjugate has right coefficients $(0,1)^T$ and left coefficients $(\overline b,\overline a)^T$. The matrix sending right coefficients to left coefficients is therefore
\[\begin{pmatrix}a&\overline b\\b&\overline a\end{pmatrix},\qquad |a|^2-|b|^2=1.\]
The determinant identity follows by evaluating $W(\overline{F_+},F_+)=2i\kappa$ on both sides of the potential. Its inverse has off-diagonal entries $-\overline b,-b$. Thus reversing the direction of transfer does not change their moduli. Together with the scaling and the unit-modulus endpoint factors, this proves $|b(\kappa)|=|b_{{\rm sc},p}(\xi)|$.
	
	Finally $\Lambda_\zeta=\int_0^1t^{\zeta-1}(1-t)^{-1/2}\,\mathrm dt\geq1/\zeta$ and $L\xi=\Lambda_\zeta\kappa$. Since $\arsinh t\leq t$ for $t\geq0$, the bound $|b(\kappa)|\leq C\zeta\min\{\kappa,\kappa^{-1}\}$ for $\kappa>0$ gives
	\[\frac{\arsinh|b_{{\rm sc},p}(\xi)|}{L\xi}\leq\frac{C\zeta}{\Lambda_\zeta}\leq C\zeta^2.\]
	Taking the supremum proves the upper bound; positivity is Lemma \ref{Lem.zero_monodromy}.
\end{proof}

We now obtain Proposition \ref{Prop.scattering_growth} from the bound for $\chi_p$.

\begin{proof}[Proof of Proposition \ref{Prop.scattering_growth}]
	By Lemma \ref{Lem.small_scattering}, $0<\chi_p\leq C_{\rm sc}(p+1)^{-2}$ for every odd $p\geq15$, which proves \eqref{Eq.chi_p_bound}.
\end{proof}

\subsection{The sonic equation and the parameter curve}
\label{Subsec.sonic_parameters}

We rewrite \eqref{Eq.velocity} as the polynomial graph equation \eqref{Eq.sonic_graph_quant}, whose sonic point is $(z_{\rm s},U)=(0,e_*)$. Near $(\ell,\beta)=(1,k/2-1/3)$, the linearization $J_{\rm s}$ in \eqref{Eq.sonic_linearization_quant} has two distinct positive eigenvalues. We choose $\beta=\widehat\beta(\varepsilon,R)$ so that their ratio is $R$ near three; the reason for this choice is the cubic coefficient equation $(R-3)a_3=\varepsilon q_*(\varepsilon,R)$ in Lemma \ref{Lem.resonant_recurrence}.

For the rest of the section fix $k\in\{3,4\}$, set $d:=k+1$ and $\ell:=1+\varepsilon$. We denote the independent sonic coordinate by $z_{\rm s}$ and define
\begin{align}
		&g_*:=\frac{k}{\beta\ell}-1,\qquad e_*:=\ell g_*^2-1,\qquad A_*:=k+2-(k-2\ell)g_*,\qquad C_*:=2k+1-\ell,\label{Eq.sonic_parameter_coefficients}\\
		&\qquad\qquad\qquad\qquad\qquad F_*(z_{\rm s}):=-e_*-A_*z_{\rm s}+C_*z_{\rm s}^2.\nonumber
\end{align}
The change of variables and the sonic linearization below are proved in \cite[Lemma~3.4, Remark~3.5 and (4.1)]{Shao-Wei-Zhang}. The parameters denoted there by $(m,\gamma,\varepsilon,A,B)$ are $(\beta\ell,g_*,e_*,A_*,C_*)$ here. We retain the formulas needed for the parameter-uniform construction, without repeating their algebraic derivation.

\begin{lemma}\label{Lem.sonic_map_quant}
	Where $0<Z$, $0<v<1$ and $1-Zv>0$, set
	\begin{equation}\label{Eq.sonic_coordinate_map}
		z_{\rm s}:=\frac{(1+g_*v^2)Z-(1+g_*)v}{Z(1-v^2)},\qquad U:=\frac{(1+g_*)^2(1-Zv)^2}{Z^2(1-v^2)}.
	\end{equation}
	The inverse is
	\begin{equation}\label{Eq.sonic_inverse_quant}
		Z=\frac{(1+g_*)\sqrt{U+(1-z_{\rm s})^2}}{U+(1+g_*)(1-z_{\rm s})},\qquad v=\frac{1-z_{\rm s}}{\sqrt{U+(1-z_{\rm s})^2}}.
	\end{equation}
	The velocity equation is equivalent, wherever $z_{\rm s}$ is a regular coordinate, to
	\begin{equation}\label{Eq.sonic_graph_quant}
		\Delta_{\rm s}U'=\Delta_U,\qquad
		\begin{cases}
			\Delta_{\rm s}:=((k+1)z_{\rm s}-1)U+(z_{\rm s}-1)F_*(z_{\rm s}),\\
			\Delta_U:=2U[U+F_*(z_{\rm s})+kz_{\rm s}(1-z_{\rm s})].
		\end{cases}
	\end{equation}
	At $(z_{\rm s},U)=(0,e_*)$, the linearization of $(\Delta_{\rm s},\Delta_U)$ is the matrix
	\begin{equation}\label{Eq.sonic_linearization_quant}
		J_{\rm s}:=\begin{pmatrix}ke_*+A_*&-1\\
			2e_*(k-A_*)&2e_*\end{pmatrix},
	\end{equation}
	with trace $(k+2)e_*+A_*$ and determinant $2ke_*(e_*+1)$.
\end{lemma}

If $J_{\rm s}$ in \eqref{Eq.sonic_linearization_quant} has eigenvalues $\lambda_-,R\lambda_-$, then $(\operatorname{tr}J_{\rm s})^2/\det J_{\rm s}=(R+1)^2/R$. We solve this scalar equation for $\beta$ near $(\varepsilon,R)=(0,3)$, rather than differentiating the eigenvalues separately.

\begin{lemma}\label{Lem.beta_surface}
	There are $\epsilon_\beta,\delta_R>0$ and
	\[\widehat\beta\in\Hol(\mathbb D_{\epsilon_\beta}\times\{R\in\C:|R-3|<\delta_R\}),\]
	real-valued for real arguments, such that
	\begin{equation}\label{Eq.beta_parameter_expansion}
		\begin{aligned}
			\widehat\beta(0,R)&=k/2-1/R,\\
			\partial_\varepsilon\widehat\beta(0,R) &=-\frac{(Rk-2)(3R^2k-Rk+2R+2)}{8R^2k(R-1)}.
		\end{aligned}
	\end{equation}
	Set $\ell=1+\varepsilon$ and $\beta=\widehat\beta(\varepsilon,R)$ in \eqref{Eq.sonic_parameter_coefficients} and \eqref{Eq.sonic_linearization_quant}. For $\varepsilon,R\in\R$ with $|\varepsilon|<\epsilon_\beta$ and $|R-3|<\delta_R$, the eigenvalues $\lambda_-<\lambda_+$ of $J_{\rm s}$ satisfy
	\[0<\lambda_-<\lambda_+,\qquad \lambda_+/\lambda_-=R,\qquad \beta>1,\quad e_*>0,\quad k>\beta(\ell+\sqrt\ell).\]
\end{lemma}

\begin{proof}
	Substitution of the trace and determinant of \eqref{Eq.sonic_linearization_quant} into $(\operatorname{tr}J_{\rm s})^2/\det J_{\rm s}=(R+1)^2/R$ gives
	\begin{equation}\label{Eq.ratio_surface_quant}
		\frac{k[k+2-(\ell+1)\beta]^2}{2[k^2-2k\beta\ell+\beta^2\ell(\ell-1)]}=\frac{(R+1)^2}{R}.
	\end{equation}
	At $\varepsilon=0$, \eqref{Eq.ratio_surface_quant} holds for $\beta=k/2-1/R$. The $\beta$ derivative of the left-hand side at $\beta=k/2-1/R$ is $R^2-1$, which equals eight at $R=3$. The holomorphic implicit-function theorem therefore gives $\widehat\beta(\varepsilon,R)$. Differentiate \eqref{Eq.ratio_surface_quant} in $\varepsilon$ at zero and divide by $R^2-1$; we obtain $\partial_\varepsilon\widehat\beta(0,R)=-(Rk-2)(3R^2k-Rk+2R+2)/(8R^2k(R-1))$. The denominator $8R^2k(R-1)$ is nonzero near $R=3$.
	
	At $(0,3)$, $\beta_0=k/2-1/3>1$ and $g_0=(3k+2)/(3k-2)>1$, so $e_*>0$. The trace and determinant are positive. On \eqref{Eq.ratio_surface_quant}, the eigenvalues near $(0,3)$ are $\lambda_-=(\operatorname{tr}J_{\rm s})/(R+1)$ and $\lambda_+=R\lambda_-$, so their ratio is $R$ and they remain distinct and positive. Also $k-2\beta_0=2/3$. By continuity, choose a real rectangle in the holomorphic neighborhood on which $\beta\geq(1+\beta_0)/2$, $e_*\geq(g_0^2-1)/2$, the trace is positive, the eigenvalues remain distinct and positive, and $k-\beta(\ell+\sqrt\ell)\geq1/3$. The rectangle determines $\epsilon_\beta$ and $\delta_R$.
\end{proof}

\subsection{The limiting origin profile}
\label{Subsec.resonant_base}

At $\ell=1$ and $\beta_0=k/2-1/3$, the substitution $v_0=\Phi_0'/((\beta_0-1)\Phi_0+Z\Phi_0')$ reduces \eqref{Eq.velocity} to the linear equation \eqref{Eq.limiting_phase_ode}. Its hypergeometric solution determines the origin branch and proves $0<v_0(Z)<Z$ for $0<Z<1$. The expansion at $Z=1$ then gives the positive cubic coefficient $A_k$ in \eqref{Eq.base_cubic}, which will be the parameter used to match the perturbed sonic and origin branches.

We define
\begin{equation}\label{Eq.base_hypergeometric_parameters}
	\beta_0:=\frac k2-\frac13,\ \ B_0:=\beta_0-1,\ \ g_0:=\frac{3k+2}{3k-2},\quad a:=B_0/2,\quad b:=\beta_0/2,\quad c:=(k+1)/2.
\end{equation}
The constants $a,b,c$ are the hypergeometric parameters. We use the normalization
\[{}_2F_1(a,b;c;u):=\sum_{m\geq0}\frac{(a)_m(b)_m}{(c)_m m!}u^m,\qquad (a)_m:=\prod_{j=0}^{m-1}(a+j),\quad (a)_0:=1.\]

\begin{lemma}\label{Lem.base_origin}
	For $0\leq Z<1$, define
	\begin{equation}\label{Eq.base_origin_definition}
		\Phi_0(Z):=B_0^{-1}{}_2F_1(a,b;c;Z^2),\qquad v_0(Z):=\frac{\Phi_0'(Z)}{B_0\Phi_0(Z)+Z\Phi_0'(Z)}.
	\end{equation}
	Then $v_0$ solves \eqref{Eq.velocity} with $\ell=1$, $\beta=\beta_0$, is odd and analytic at zero, and satisfies $0<v_0(Z)<Z$ for $0<Z<1$. With $t:=1-Z^2$, we have
	\[v_0(Z)=g_0^{-1}-D_kt^{1/3}+\Order(t^{2/3}),\qquad D_k>0.\]
	The function $y\mapsto v_0(\sqrt{1-y^3})$ extends analytically to a neighborhood of zero.
\end{lemma}

\begin{proof}
	The hypergeometric equation and the substitution $u=Z^2$ give
	\begin{equation}\label{Eq.limiting_phase_ode}
		(1-Z^2)\Phi_0''+(k/Z-2\beta_0Z)\Phi_0'-B_0\beta_0\Phi_0=0.
	\end{equation}
	Substitution of \eqref{Eq.base_origin_definition} into \eqref{Eq.velocity}, using \eqref{Eq.limiting_phase_ode}, proves the velocity equation wherever the denominator $B_0\Phi_0+Z\Phi_0'$ is nonzero. The integral representation below proves positivity of that denominator. The power series of $\Phi_0$ at zero gives $v_0'(0)=\beta_0/(k+1)$ and odd analyticity.
	
	The parameters in \eqref{Eq.base_hypergeometric_parameters} satisfy $a,b,c-b>0$, so Euler's integral representation applies; see \cite[Sections 15.6 and 15.8]{DLMF}. Let $\mu_Z$ be the probability measure on $(0,1)$ with density proportional to $t_1^{b-1}(1-t_1)^{c-b-1}(1-Z^2t_1)^{-a}$. Differentiating Euler's integral and using $2a=B_0$ gives
	\[\frac{v_0(Z)}Z=\frac{\int_0^1\frac{t_1}{1-Z^2t_1}\,\mathrm d\mu_Z(t_1)}{\int_0^1\frac1{1-Z^2t_1}\,\mathrm d\mu_Z(t_1)}\in(0,1),\]
	where the strict inequalities follow from $0<t_1<1$ and positivity of $\mu_Z$. Also $\Phi_0>0$ and $\Phi_0'\geq0$, so $B_0\Phi_0+Z\Phi_0'>0$ on $[0,1)$.
	
	Since $c-a-b=4/3$, the connection formula at one gives
	\begin{equation}\label{Eq.base_connection_expansion}
		\Phi_0(Z)=P_0^*(t)+C_kt^{4/3}P_1^*(t),\qquad P_1^*(0)=1,
	\end{equation}
	with $P_0^*,P_1^*$ analytic near zero and
	\[\Phi_1:=\Phi_0(1)=\frac{\Gamma(c)\Gamma(4/3)}{B_0\Gamma(c-a)\Gamma(c-b)}>0,\qquad C_k:=\frac{\Gamma(c)\Gamma(-4/3)}{B_0\Gamma(a)\Gamma(b)}>0.\]
	The positivity of the second coefficient follows from $\Gamma(-4/3)>0$, by the functional equation of $\Gamma$. Gauss's formula for the value of the derivative gives $\Phi_0'(1)/\Phi_1=3B_0\beta_0/2$. Therefore $v_0(1)=3\beta_0/(2+3\beta_0)=g_0^{-1}$. Differentiating \eqref{Eq.base_connection_expansion} with $t=1-Z^2$ gives the first nonintegral term $-(8/3)C_kt^{1/3}$ in $\Phi_0'$. Substitution into the quotient \eqref{Eq.base_origin_definition} gives
	\begin{equation}\label{Eq.base_fractional_coefficient}
		D_k:=\frac{(8/3)C_k}{B_0\Phi_1(1+3\beta_0/2)^2}>0.
	\end{equation}
	Upon writing $t=y^3$, both $\Phi_0$ and $\Phi_0'$ are analytic in $y$; the denominator at zero is positive. Hence $y\mapsto v_0(\sqrt{1-y^3})$ is analytic near $y=0$ and $v_0=g_0^{-1}-D_ky+\Order(y^2)$.
\end{proof}

Under the coordinate map \eqref{Eq.sonic_coordinate_map}, the origin solution \eqref{Eq.base_origin_definition} approaches the line $U=g_0^2-1+2(1+g_0)z_{\rm s}$. The identity \eqref{Eq.lightcone_identity_sonic} converts $1-Z^2$ into the difference from this line and determines the coefficient of $z_{\rm s}^3$.

\begin{lemma}\label{Lem.base_cubic}
	Apply the coordinate map in Lemma \ref{Lem.sonic_map_quant}, with $g_*=g_0$, to $(Z,v_0(Z))$ from Lemma \ref{Lem.base_origin}. The resulting graph $U_0(z_{\rm s})$ is analytic near $z_{\rm s}=0$ and
	\begin{equation}\label{Eq.base_cubic}
		U_0(z_{\rm s})=g_0^2-1+2(1+g_0)z_{\rm s}+A_kz_{\rm s}^3+\Order(z_{\rm s}^4),\qquad A_k:=\frac{(g_0+1)^4}{(g_0-1)g_0^4D_k^3}>0.
	\end{equation}
\end{lemma}

\begin{proof}
	By Lemma \ref{Lem.base_origin}, $Z=\sqrt{1-y^3}$ and $v_0$ are analytic in $y$ at zero. The sonic coordinate formula has no vanishing denominator there. The derivative of $z_{\rm s}$ with respect to $y$ is
	\begin{equation}\label{Eq.base_sonic_coordinate_derivative}
		\left.\frac{\mathrm dz_{\rm s}}{\mathrm dy}\right|_{y=0}=\frac{g_0^2D_k}{g_0+1}>0.
	\end{equation}
	Indeed, $\left.\frac{\mathrm d}{\mathrm dy}\sqrt{1-y^3}\right|_{y=0}=0$, $\partial_vz_{\rm s}(1,g_0^{-1})=-g_0^2/(g_0+1)$ and $\left.\frac{\mathrm d}{\mathrm dy}v_0(\sqrt{1-y^3})\right|_{y=0}=-D_k$. The analytic inverse theorem makes $y$ an analytic function of $z_{\rm s}$. The formula $U_0=(1-z_{\rm s})^2(v_0^{-2}-1)$ is therefore analytic as well.
	
	For general $g$ the inverse map has the exact identity
	\begin{equation}\label{Eq.lightcone_identity_sonic}
		1-Z^2=\frac{U(U-L_g)}{[U+(1+g)(1-z_{\rm s})]^2},\qquad L_g:=g^2-1+2(1+g)z_{\rm s}.
	\end{equation}
	We obtain \eqref{Eq.lightcone_identity_sonic} by expanding the squared numerator of \eqref{Eq.sonic_inverse_quant}; the constant, linear and quadratic terms cancel leaving $U(U-L_g)$. At $(z_{\rm s},\varepsilon)=(0,0)$, the denominator in \eqref{Eq.lightcone_identity_sonic} is $g_0^2(g_0+1)^2$, and $U_0(0)=g_0^2-1$. Since $1-Z^2=y^3$, the difference $U_0-L_{g_0}$ has order exactly three. By \eqref{Eq.base_sonic_coordinate_derivative}, $y=((g_0+1)/(g_0^2D_k))z_{\rm s}+\Order(z_{\rm s}^2)$. Thus the cubic coefficients in \eqref{Eq.lightcone_identity_sonic} satisfy
	\[\frac{g_0-1}{g_0^2(g_0+1)}A_k=\left(\frac{g_0+1}{g_0^2D_k}\right)^3,\]
	which is the value of $A_k$ in \eqref{Eq.base_cubic}.
\end{proof}

\subsection{Convergent sonic curves and matching}
\label{Subsec.sonic_matching}

Substitution of $U=\sum_{m\geq0}a_mz_{\rm s}^m$ into \eqref{Eq.sonic_graph_quant} gives $(R-m)\delta_*a_m=E_m$ in \eqref{Eq.sonic_recurrence}. Near $R=3$, only the equation for $a_3$ has a vanishing denominator. We prescribe $a_3=a_{\rm s}$ near $A_k$ and solve $(R-3)a_{\rm s}=\varepsilon q_*(\varepsilon,R)$ for $R$. A majorant proves convergence of the remaining coefficients; the scalar mismatch \eqref{Eq.sonic_origin_mismatch} then determines $a_{\rm s}$ by matching the origin solution.

\begin{lemma}[The resonant coefficient equation]\label{Lem.resonant_recurrence}
	On the surface $\beta=\widehat\beta(\varepsilon,R)$ define
	\begin{equation}\label{Eq.sonic_series_initial_coefficients}
		\delta_*:=\frac{(k+2)e_*+A_*}{R+1},\qquad a_0:=e_*,\qquad a_1:=R\delta_*-2e_*.
	\end{equation}
	For a formal solution $U=\sum_{m\geq0}a_mz_{\rm s}^m$ of \eqref{Eq.sonic_graph_quant}, the coefficients satisfy
	\begin{equation}\label{Eq.sonic_recurrence}
		(R-m)\delta_*a_m=E_m\qquad(m\geq2),
	\end{equation}
	where
	\begin{equation}\label{Eq.sonic_recurrence_numerator}
		\begin{aligned}
			E_m:={}&-\frac{m+1}{2}\sum_{j=2}^{m-1}a_ja_{m+1-j}-\left(2-\frac{k+1}{2}m\right)\sum_{j=1}^{m-1}a_ja_{m-j}\\
			&-[(m-3)A_*+(m-1)C_*+2k]a_{m-1}+[(m-4)C_*+2k]a_{m-2}.
		\end{aligned}
	\end{equation}
	Empty sums are zero. The analytic coefficient $a_2$ vanishes at $\varepsilon=0$, and
	\[\partial_\varepsilon a_2(0,R)=\frac{R(Rk^2+2Rk+2R-2)}{2(R-2)(R-1)}.\]
	There is a holomorphic function $q_*$ near $(0,3)$ such that $E_3/\delta_*=\varepsilon q_*(\varepsilon,R)$ and
	\[q_*(0,R)=H_3(R):=\frac{R(Rk-2)(3Rk^2+4Rk+2R-2)}{4(R-2)(R-1)}>0\]
	for real $R$ in a neighborhood of three.
\end{lemma}

\begin{proof}
	Insert the power series into the two polynomials in \eqref{Eq.sonic_graph_quant}. The constant coefficient vanishes because $a_0=e_*$. The coefficient of $z_{\rm s}$ is $J_{\rm s}(1,a_1)^T=\delta_*(1,a_1)^T$, with $J_{\rm s}$ from \eqref{Eq.sonic_linearization_quant}. Thus \eqref{Eq.sonic_series_initial_coefficients} selects the smaller eigenvalue $\delta_*=\lambda_-$. In the coefficient of $z_{\rm s}^m$, collect the terms containing $a_m$ and use $a_1=R\delta_*-2e_*$ and $(k+2)e_*+A_*=(R+1)\delta_*$. The coefficient of $a_m$ in $\Delta_{\rm s}U'-\Delta_U$ is $(m-R)\delta_*$. The quadratic products in \eqref{Eq.sonic_graph_quant} give the two sums in \eqref{Eq.sonic_recurrence_numerator}, by symmetrizing $\sum j a_ja_{m+1-j}$ and $\sum j a_ja_{m-j}$. The remaining linear terms give the last line. Thus $(R-m)\delta_*a_m=E_m$ for $m\geq2$.
	
	To compute $\partial_\varepsilon a_2(0,R)$ and $q_*(0,R)$, we use
	\begin{equation}\label{Eq.sonic_quadratic_cubic_numerators}
		\begin{aligned}
			E_2&=(k-1)a_1^2+(A_*-C_*-2k)a_1+(2k-2C_*)e_*,\\
			E_3&=-2a_2^2+[(3k-1)a_1-2(C_*+k)]a_2+(2k-C_*)a_1.
		\end{aligned}
	\end{equation}
	At $\varepsilon=0$, $\ell=1$, $C_*=2k$, and $U=g_*^2-1+2(1+g_*)z_{\rm s}$ is an exact solution of \eqref{Eq.sonic_graph_quant}. Thus $a_2=E_3=0$. To differentiate $a_2$ and $E_3$ in $\varepsilon$, let $\beta_0(R):=k/2-1/R$, let $\beta_1(R)=\partial_\varepsilon\widehat\beta(0,R)$ from \eqref{Eq.beta_parameter_expansion}, and use the values and first derivatives
	\begin{equation}\label{Eq.sonic_parameter_derivatives}
		\begin{aligned}
			g_0&=k/\beta_0-1,&g_1&=-k(\beta_1+\beta_0)/\beta_0^2,\\
			e_0&=g_0^2-1,&e_1&=g_0^2+2g_0g_1,\\
			A_0^*&=k+2-(k-2)g_0,&A_1^*&=2g_0-(k-2)g_1,\qquad C_0^*=2k,\quad C_1^*=-1.
		\end{aligned}
	\end{equation}
	Define
	\begin{equation}\label{Eq.sonic_slope_variations}
		\delta_j:=\frac{(k+2)e_j+A_j^*}{R+1},\qquad u_j:=R\delta_j-2e_j=\left.\partial_\varepsilon^ja_1\right|_{\varepsilon=0},\qquad j=0,1.
	\end{equation}
	Since $E_2(0,R)=a_2(0,R)=0$, differentiating \eqref{Eq.sonic_quadratic_cubic_numerators} gives
	\begin{equation}\label{Eq.sonic_coefficient_variations}
		\begin{aligned}
			(\partial_\varepsilon E_2)|_0 &=2(k-1)u_0u_1+(A_1^*-C_1^*)u_0 +(A_0^*-C_0^*-2k)u_1\\
			&\quad-2C_1^*e_0+(2k-2C_0^*)e_1,\\
			(\partial_\varepsilon a_2)|_0 &=\frac{(\partial_\varepsilon E_2)|_0}{(R-2)\delta_0} =\frac{R(Rk^2+2Rk+2R-2)}{2(R-2)(R-1)},\\
			\left.\partial_\varepsilon\left(\frac{E_3}{\delta_*}\right)\right|_0 &=\frac{[(3k-1)u_0-2(C_0^*+k)](\partial_\varepsilon a_2)|_0-C_1^*u_0}{\delta_0} =H_3(R).
		\end{aligned}
	\end{equation}
	The denominators in \eqref{Eq.sonic_coefficient_variations} are products of $R$, $R-1$, $R-2$, $k$ and nonzero powers of $\beta_0$, so the differentiated expressions are holomorphic on a common neighborhood. Since $E_3(0,R)=0$, division by $\varepsilon$ is removable and defines $q_*$ holomorphically. The numerator $R(Rk-2)(3Rk^2+4Rk+2R-2)$ is positive for $k=3,4$ and $R$ near three.
\end{proof}

The equation for the cubic coefficient is the only resonant equation near $R=3$. We solve $(R-3)a_{\rm s}=\varepsilon q_*(\varepsilon,R)$ for $R$ and then estimate the remaining coefficients by a convergent majorant.

\begin{lemma}\label{Lem.sonic_family}
	There is a complex neighborhood of $(\varepsilon,a_{\rm s})=(0,A_k)$, a holomorphic function $\widehat R(\varepsilon,a_{\rm s})$, and a jointly holomorphic solution $U_{\rm s}(z_{\rm s};\varepsilon,a_{\rm s})$ on a fixed disk $|z_{\rm s}|<r_{\rm s}$, such that
	\[\widehat R(0,a_{\rm s})=3,\qquad \partial_\varepsilon\widehat R(0,a_{\rm s})=H_3(3)/a_{\rm s},\qquad a_3=a_{\rm s}.\]
	The parameters of $U_{\rm s}$ in \eqref{Eq.sonic_graph_quant} are $\ell=1+\varepsilon$, $R=\widehat R(\varepsilon,a_{\rm s})$ and $\beta=\widehat\beta(\varepsilon,R)$. At $(\varepsilon,a_{\rm s})=(0,A_k)$, we have $U_{\rm s}=U_0$ from Lemma \ref{Lem.base_cubic}.
\end{lemma}

\begin{proof}
	With $a_3=a_{\rm s}$, the $m=3$ equation in \eqref{Eq.sonic_recurrence} is
	\begin{equation}\label{Eq.sonic_cubic_compatibility}
		(R-3)a_{\rm s}-\varepsilon q_*(\varepsilon,R)=0.
	\end{equation}
	The $R$ derivative of $(R-3)a_{\rm s}-\varepsilon q_*(\varepsilon,R)$ at $(0,3,A_k)$ is $A_k>0$. The holomorphic implicit-function theorem gives $\widehat R(0,a_{\rm s})=3$ and $\partial_\varepsilon\widehat R(0,a_{\rm s})=q_*(0,3)/a_{\rm s}=H_3(3)/a_{\rm s}$. In particular $\partial_{a_{\rm s}}\widehat R(0,a_{\rm s})=0$. The coefficients $a_0,a_1,a_2$ are holomorphic in $(\varepsilon,a_{\rm s})$, $a_3=a_{\rm s}$, and for $m\geq4$ use \eqref{Eq.sonic_recurrence}.
	
	We verify convergence uniformly in the parameters. Choose a closed complex polydisk on which $|\widehat R-3|\leq1/4$, $|\delta_*|\geq c_0>0$, and the parameters and first four coefficients are bounded. Then $|(m-\widehat R)\delta_*|\geq c m$ for every $m\geq4$. The recurrence, absorbing terms containing $a_1$ into a fixed constant, gives
	\begin{equation}\label{Eq.sonic_majorant_recurrence}
		|a_m|\leq C\left(\sum_{j=2}^{m-1}|a_ja_{m+1-j}|+\sum_{j=2}^{m-2}|a_ja_{m-j}|+|a_{m-1}|+|a_{m-2}|\right).
	\end{equation}
	Take $M$ larger than the uniform bounds for $|a_2|,|a_3|$ and consider the scalar equation
	\begin{equation}\label{Eq.sonic_majorant_equation}
		Y(t)=M(1+t)+C(t+t^2)(Y(t)^2+Y(t)).
	\end{equation}
	The $Y$ derivative of $Y-M(1+t)-C(t+t^2)(Y^2+Y)$ equals one at $(0,M)$. The implicit-function theorem gives an analytic solution $Y(t)$ near zero. The coefficient recursion gives nonnegative Taylor coefficients for $Y$. For $b_j:=|a_{j+2}|$, the convolution sums in \eqref{Eq.sonic_majorant_recurrence} are dominated by the coefficients of $tY^2$ and $t^2Y^2$; the other terms are dominated by $tY$ and $t^2Y$. Induction therefore bounds each $b_j$ by the corresponding coefficient of $Y$. The radius of convergence of the scalar majorant \eqref{Eq.sonic_majorant_equation} is independent of $(\varepsilon,a_{\rm s})$ in the chosen polydisk. Hence $\sum_ma_mz_{\rm s}^m$ converges normally on one fixed smaller disk, uniformly in the parameters. The sum $U_{\rm s}$ is jointly holomorphic and solves the ODE by termwise differentiation.
	
	At $(0,A_k)$ the convergent series has the same first four coefficients as the analytic origin germ $U_0$. All later coefficients are uniquely determined, since $m\ne3$ for $m\geq4$. Thus the two germs agree.
\end{proof}

The family from Lemma \ref{Lem.sonic_family} still has the free coefficient $a_{\rm s}$. At a fixed small $z_*>0$, we require $U_{\rm s}(z_*;\varepsilon,a_{\rm s})=U_F(z_*;\varepsilon,\beta)$, where $U_F$ is obtained by continuing the regular origin solution. The derivative \eqref{Eq.sonic_matching_derivative} with respect to $a_{\rm s}$ is nonzero and determines the matched parameter.

\begin{lemma}\label{Lem.match_origin}
	There exist $\varepsilon_1>0$ and real-valued analytic functions $a_{\rm s}(\varepsilon)$, $R(\varepsilon)$ and $\beta(\varepsilon)$ for $|\varepsilon|<\varepsilon_1$ such that
	\[a_{\rm s}(0)=A_k,\qquad R(\varepsilon)=3+\frac{H_3(3)}{A_k}\varepsilon+\Order(\varepsilon^2),\qquad \beta(0)=k/2-1/3.\]
	The sonic curve $U_{\rm s}(z_{\rm s};\varepsilon,a_{\rm s}(\varepsilon))$ matches the solution issued from $Z=0$ at one fixed $z_{\rm s}=z_*>0$. On a fixed compact interval $0\leq Z\leq Z_*<1$ containing the matching point, the origin solution is jointly analytic in $Z$ and $\varepsilon$, is odd in $Z$, and converges to $v_0$.
\end{lemma}

\begin{proof}
	Write $v(Z)=Zu(Z^2)$ and $s:=Z^2$. Division of the velocity equation gives
	\begin{equation}\label{Eq.quantitative_origin_equation}
		2s u'(s)=\frac{(1-su^2)[\beta\ell(1-su^2)-ku(1-su)]}{(1-su)^2-\ell s(u-1)^2}-u=:\mathcal H(s,u;\varepsilon,\beta).
	\end{equation}
	The rational function is analytic near $s=0$ and $u=\beta\ell/(k+1)$, and $\mathcal H(0,u)=\beta\ell-(k+1)u$. Thus $\mathcal G:=[\mathcal H-\beta\ell+(k+1)u]/s$ extends analytically there. With $a_*:=(k+1)/2$, equation \eqref{Eq.quantitative_origin_equation} gives $\partial_s[s^{a_*}(u-\beta\ell/(k+1))]=\frac12s^{a_*}\mathcal G(s,u;\varepsilon,\beta)$. Integration from zero and the change of integration variable give
	\begin{equation}\label{Eq.quantitative_origin_integral}
		u(s)=\frac{\beta\ell}{k+1}+\frac s2\int_0^1t^{a_*}\mathcal G(st,u(st);\varepsilon,\beta)\,\mathrm dt.
	\end{equation}
	On a common parameter neighborhood choose a ball of holomorphic functions around the constant term. If the disk radius is $r_o$ and $L_o$ is a bound for $|\partial_u\mathcal G|$ there, choose $r_oL_o/(2(a_*+1))\leq1/2$, as well as the corresponding bound making the map preserve the ball. The contraction theorem gives a unique solution, holomorphic in $s$ and the parameters by uniform convergence of the fixed-point iterates. Analytic ODE continuation along any compact subinterval of $[0,1)$ is uniform near the base parameters, because the denominator on the base curve is $(1-Z^2)(1-v_0^2)>0$. The continuation gives the origin branch $v_F$ on the required compact interval.
	
	Choose $z_*\in(0,r_{\rm s}/2)$ on the base sonic curve. The corresponding coordinate $Z_*^0<1$ is a regular coordinate point: from \eqref{Eq.base_cubic} and \eqref{Eq.lightcone_identity_sonic}, $\mathrm dZ/\mathrm dz_{\rm s}<0$ for positive $z_{\rm s}$ in a fixed interval. The analytic inverse theorem therefore continues the crossing of $z_{\rm s}=z_*$ along $v_F$ in the parameters. Let $U_F(z_*;\varepsilon,\beta)$ denote the value of $U$ on $v_F$ at $z_{\rm s}=z_*$. Define the mismatch
	\begin{equation}\label{Eq.sonic_origin_mismatch}
		\mathcal M(\varepsilon,a_{\rm s}):=U_{\rm s}(z_*;\varepsilon,a_{\rm s})-U_F(z_*;\varepsilon,\widehat\beta(\varepsilon,\widehat R(\varepsilon,a_{\rm s}))).
	\end{equation}
	The function $\mathcal M$ is holomorphic and $\mathcal M(0,A_k)=0$. Since $\partial_{a_{\rm s}}\widehat R(0,a_{\rm s})=0$, the derivative of the origin branch with respect to $a_{\rm s}$ is zero at $\varepsilon=0$. The convergent sonic series gives
	\begin{equation}\label{Eq.sonic_matching_derivative}
		\partial_{a_{\rm s}}\mathcal M(0,A_k)=z_*^3+\Order(z_*^4).
	\end{equation}
	Choose $z_*>0$ so that the absolute value of the remainder is at most $z_*^3/2$. The implicit-function theorem yields $a_{\rm s}(\varepsilon)$, and then $R(\varepsilon):=\widehat R(\varepsilon,a_{\rm s}(\varepsilon))$ and $\beta(\varepsilon):=\widehat\beta(\varepsilon,R(\varepsilon))$. Differentiating \eqref{Eq.sonic_cubic_compatibility}, using $a_{\rm s}(0)=A_k$, gives $R'(0)=H_3(3)/A_k$; \eqref{Eq.beta_parameter_expansion} gives $\beta(0)=k/2-1/3$. At the matching point the graph equation is nonsingular. Equality of the value therefore identifies the two curves on an overlap by ODE uniqueness. The functions are real-valued for real parameters by conjugation and uniqueness.
\end{proof}

\subsection{The outgoing property and the quantitative complex analytic radius}
\label{Subsec.background_radius_proof}

We express the matched sonic curve in the normalized variables $X_{\rm s}=X\circ Z_{\rm s}$ and $W_{\rm s}=b\circ X_{\rm s}$. The bounds $-\partial_{z_{\rm s}}X_{\rm s},-\partial_{z_{\rm s}}W_{\rm s}\asymp\varepsilon+z_{\rm s}^2$ in \eqref{Eq.normalized_sonic_real_bounds} give $b'(X_{\rm s})=(\partial_{z_{\rm s}}W_{\rm s})/(\partial_{z_{\rm s}}X_{\rm s})\geq\gamma_*>0$. They also determine the inverse radius: a disk of radius comparable to $(\varepsilon+z_r^2)^{1/2}$ in $z_{\rm s}$ has an image containing a disk of radius comparable to $(\varepsilon+z_r^2)^{3/2}\geq\varepsilon^{3/2}$ in $X$. Lemmas \ref{Lem.complex_inverse_radius} and \ref{Lem.cartesian_uniform_background} convert this inverse estimate into the coefficient bounds \eqref{Eq.background_bounds}.

Use $a_{\rm s}(\varepsilon),\beta(\varepsilon)$ determined by \eqref{Eq.sonic_origin_mismatch} and Lemma \ref{Lem.match_origin}. Write $U(z_{\rm s},\varepsilon):=U_{\rm s}(z_{\rm s};\varepsilon,a_{\rm s}(\varepsilon))$ and use $Z_{\rm s},v_{\rm s}$ for the inverse coordinates of $U$. Define
\begin{equation}\label{Eq.normalized_sonic_functions}
	\begin{aligned}
		H_*&:=U+(1-z_{\rm s})^2,\qquad D_*:=U+(1+g_*)(1-z_{\rm s}),\\
		\mathcal M_*&:=\frac{UH_*}{(\beta-1)D_*^2},\qquad W_{\rm s}:=\frac{g_*+z_{\rm s}}{\sqrt{H_*}}.
	\end{aligned}
\end{equation}
Let $X_{\rm s}(z_{\rm s},\varepsilon)$ be the analytic primitive of $\partial_{z_{\rm s}}Z_{\rm s}/\mathcal M_*$ with the value of the normalized coordinate $X$ at the matching point. On the common real curve,
\begin{equation}\label{Eq.normalized_sonic_composition}
	X_{\rm s}=X\circ Z_{\rm s},\qquad W_{\rm s}=b\circ X_{\rm s}.
\end{equation}
Indeed $\mathcal M_*=(\mathrm dZ/\mathrm dX)$, by \eqref{Eq.x_b} and \eqref{Eq.sonic_inverse_quant}. The positive square root of $H_*$ extends to a fixed complex neighborhood of the base sonic point. In the formulas below, derivatives of $X_{\rm s},W_{\rm s},Z_{\rm s}$ are taken before evaluation, and an omitted $\varepsilon$ argument is held fixed.

\begin{lemma}\label{Lem.sonic_expansion_bounds}
	There are $r_*>0$, $\varepsilon_2>0$ and $c,C>0$ such that the functions $X_{\rm s},W_{\rm s}$ are holomorphic on $\{(z_{\rm s},\varepsilon)\in\C^2:|z_{\rm s}|<2r_*,\ |\varepsilon|<\varepsilon_2\}$, and
	\begin{equation}\label{Eq.normalized_sonic_expansions}
		\begin{aligned}
			\partial_{z_{\rm s}}X_{\rm s}&=-a_x\varepsilon-c_xz_{\rm s}^2+\Order(|\varepsilon|^2+|\varepsilon||z_{\rm s}|+|z_{\rm s}|^3),\\
			\partial_{z_{\rm s}}W_{\rm s}&=-a_w\varepsilon-c_wz_{\rm s}^2+\Order(|\varepsilon|^2+|\varepsilon||z_{\rm s}|+|z_{\rm s}|^3),
		\end{aligned}
	\end{equation}
	where
	\begin{equation}\label{Eq.normalized_sonic_constants}
		a_x:=\frac{(3k-2)(3k-8)}{16},\qquad a_w:=\frac{3k-2}{8},\qquad c_x:=\frac{3B_0A_k}{2g_0^2},\qquad c_w:=\frac{3A_k}{2g_0^2}.
	\end{equation}
	All four constants are positive. On the real set $|z_{\rm s}|\leq r_*$, $0<\varepsilon\leq\varepsilon_2$,
	\begin{equation}\label{Eq.normalized_sonic_real_bounds}
		c(\varepsilon+z_{\rm s}^2)\leq-\partial_{z_{\rm s}}X_{\rm s},-\partial_{z_{\rm s}}W_{\rm s}\leq C(\varepsilon+z_{\rm s}^2).
	\end{equation}
	For $z_{\rm s},\varepsilon\in\C$ with $|z_{\rm s}|<3r_*/2$ and $|\varepsilon|<\varepsilon_2/2$, one has $|\partial_{z_{\rm s}}^2X_{\rm s}|+|\partial_{z_{\rm s}}^2W_{\rm s}|\leq C(|\varepsilon|+|z_{\rm s}|)$.
\end{lemma}

\begin{proof}
	Differentiate \eqref{Eq.normalized_sonic_functions} and \eqref{Eq.sonic_inverse_quant}, using $\partial_{z_{\rm s}}X_{\rm s}=\mathcal M_*^{-1}\partial_{z_{\rm s}}Z_{\rm s}$. We obtain
	\begin{equation}\label{Eq.normalized_sonic_derivative_formulas}
		\begin{aligned}
			\partial_{z_{\rm s}}W_{\rm s}&=\frac{2H_*-(g_*+z_{\rm s})\partial_{z_{\rm s}}H_*}{2H_*^{3/2}},\\
			\partial_{z_{\rm s}}X_{\rm s}&=\frac{(1+g_*)(\beta-1)((\partial_{z_{\rm s}}H_*)D_*-2H_*\partial_{z_{\rm s}}D_*)}{2UH_*^{3/2}}.
		\end{aligned}
	\end{equation}
	The denominators $2H_*^{3/2}$ and $2UH_*^{3/2}$ are nonzero on a common neighborhood, since $U(0,0)=g_0^2-1>0$ and $H_*(0,0)=g_0^2>0$. At $\varepsilon=0$, substitute \eqref{Eq.base_cubic}. The linear curve $L_{g_0}$ gives zero for both derivatives, and the first correction $A_kz_{\rm s}^3$ gives $-c_xz_{\rm s}^2$ and $-c_wz_{\rm s}^2$, respectively.
	
	At $z_{\rm s}=0$, use $g_j,e_j$ from \eqref{Eq.sonic_parameter_derivatives} and $u_j=\left.\partial_\varepsilon^ja_1\right|_0$ from \eqref{Eq.sonic_slope_variations}, evaluated at $R=3$. Since
	\[H_*=e_*+1,\qquad \partial_{z_{\rm s}}H_*=a_1-2,\qquad D_*=e_*+1+g_*,\qquad \partial_{z_{\rm s}}D_*=a_1-1-g_*,\]
	we obtain
	\begin{equation}\label{Eq.normalized_sonic_parameter_derivatives}
		\begin{aligned}
			(\partial_\varepsilon\partial_{z_{\rm s}}W_{\rm s})(0,0) &=\frac{2e_1-g_1(u_0-2)-g_0u_1}{2g_0^3} =-\frac{3k-2}{8},\\
			J_1&:=u_1(e_0+1+g_0)+(u_0-2)(e_1+g_1)\\
			&\quad-2e_1(u_0-1-g_0)-2(e_0+1)(u_1-g_1),\\
			(\partial_\varepsilon\partial_{z_{\rm s}}X_{\rm s})(0,0) &=\frac{(1+g_0)(\beta_0-1)}{2e_0g_0^3}J_1 =-\frac{(3k-2)(3k-8)}{16}.
		\end{aligned}
	\end{equation}
	Before composition with $R(\varepsilon)$, both expressions in \eqref{Eq.normalized_sonic_derivative_formulas} vanish at $z_{\rm s}=\varepsilon=0$ for every nearby $R$, because $U=g_*^2-1+2(1+g_*)z_{\rm s}$ is the solution at $\varepsilon=0$. Their $R$-derivatives vanish there. The chain-rule term proportional to $R'(0)$ is therefore zero in \eqref{Eq.normalized_sonic_parameter_derivatives}. Moreover, the convergent family in Lemma \ref{Lem.sonic_family} satisfies $\partial_{a_{\rm s}}U_{\rm s}=z_{\rm s}^3+\Order(z_{\rm s}^4)$ at $\varepsilon=0$. Hence $a_{\rm s}'(0)$ does not enter the constant coefficient of either normalized derivative. Equations \eqref{Eq.base_cubic} and \eqref{Eq.normalized_sonic_parameter_derivatives} give the coefficients \eqref{Eq.normalized_sonic_constants} in \eqref{Eq.normalized_sonic_expansions}.
	
	Joint analyticity bounds the remainder by $C_0(|\varepsilon|^2+|\varepsilon z_{\rm s}|+|z_{\rm s}|^3)$ and the $z_{\rm s}$ derivative of the remainder by $C_1(|\varepsilon|+|z_{\rm s}|^2)$. Choose $r_*>0$ with $C_0r_*\leq\min\{a_x,a_w,c_x,c_w\}/8$, inside half the holomorphic disk. Then choose $\varepsilon_2\leq\min\{r_*^2,\varepsilon_1/2\}$ with $C_0\varepsilon_2\leq\min\{a_x,a_w\}/8$. For $|z_{\rm s}|\leq r_*$ and $0<\varepsilon\leq\varepsilon_2$, the remainder is at most $\tfrac12\min\{a_x,a_w,c_x,c_w\}(\varepsilon+z_{\rm s}^2)$. Thus $-\partial_{z_{\rm s}}X_{\rm s},-\partial_{z_{\rm s}}W_{\rm s}$ lie between $c(\varepsilon+z_{\rm s}^2)$ and $C(\varepsilon+z_{\rm s}^2)$. The differentiated remainder and the leading quadratic terms give the final complex derivative bound, after reducing the fixed complex neighborhood once.
\end{proof}

By \eqref{Eq.normalized_sonic_real_bounds}, $X_{\rm s}$ is strictly decreasing for $\varepsilon>0$, and the ratio $(\partial_{z_{\rm s}}W_{\rm s})/(\partial_{z_{\rm s}}X_{\rm s})$ is bounded above and below by positive constants. We continue the matched curve until $Z_{\rm s}=1$ and extend this bound for $b'$ to the origin interval.

\begin{lemma}\label{Lem.complete_background}
	There is $\varepsilon_3>0$ such that for every $0<\varepsilon<\varepsilon_3$ the matched origin curve extends through the sonic point up to $Z=1$ and defines
	\[v_\varepsilon\in C^\omega_{\rm o}([-1,1]),\qquad |v_\varepsilon(Z)|<1\quad(-1\leq Z\leq1),\qquad \beta(\varepsilon)>1.\]
	The profile $v_\varepsilon$ solves \eqref{Eq.velocity}. There are constants $\gamma_*,M_1>0$, independent of $\varepsilon$, such that
	\[\gamma_*\leq b'(X)\leq M_1\quad(0\leq X\leq X_1),\qquad \sym\mathrm D_x\bvec\geq\gamma_*I_d.\]
\end{lemma}

\begin{proof}
	The denominators and the square root in \eqref{Eq.sonic_inverse_quant} are positive on a fixed real neighborhood of $(z_{\rm s},\varepsilon)=(0,0)$. Lemma \ref{Lem.sonic_expansion_bounds} gives $\partial_{z_{\rm s}}Z_{\rm s}=\mathcal M_*\partial_{z_{\rm s}}X_{\rm s}<0$ for positive $\varepsilon$, so $Z_{\rm s}$ is a regular coordinate on the entire real sonic interval. Moreover
	\begin{equation}\label{Eq.profile_lightcone_crossing}
		U-L_{g_*}=\varepsilon g_*^2+A_kz_{\rm s}^3+\Order(\varepsilon|z_{\rm s}|+\varepsilon z_{\rm s}^2+\varepsilon|z_{\rm s}|^3+z_{\rm s}^4).
	\end{equation}
	In \eqref{Eq.profile_lightcone_crossing}, the constant coefficient of $U-L_{g_*}$ is $e_*-(g_*^2-1)=\varepsilon g_*^2$; the other coefficients follow from the analytic family and \eqref{Eq.base_cubic}. Set $c_*:=(g_0^2/A_k)^{1/3}$. At $z_{\rm s}=-(c_*/2)\varepsilon^{1/3}$ $U-L_{g_*}>0$, and at $z_{\rm s}=-2c_*\varepsilon^{1/3}$ we have $U-L_{g_*}<0$, for all $0<\varepsilon<\varepsilon_3$, after choosing $\varepsilon_3$ so that the errors divided by $\varepsilon$ are at most $g_0^2/4$. Both test points lie inside $(-r_*,0)$ by the same choice. Identity \eqref{Eq.lightcone_identity_sonic}, positivity of $U,D_*$, and strict monotonicity of $Z_{\rm s}$ imply that there is a unique point in $(-2c_*\varepsilon^{1/3},-(c_*/2)\varepsilon^{1/3})$ with $Z_{\rm s}=1$.
	
	Between $z_*$ and the point where $Z_{\rm s}=1$, the inverse formulas give $0<v_{\rm s}<1$. The origin branch has the same inequalities on $(0,Z_*]$ by Lemma \ref{Lem.match_origin} and the strict base inequalities, including the derivative at zero. The two branches agree on an overlap. The matched origin and sonic curves define a real-valued analytic solution up to and across $Z=1$; at the sonic point the multiplied ODE holds by continuity, and the inverse has nonzero derivative for $\varepsilon>0$. Odd extension at zero gives $v_\varepsilon\in C^\omega_{\rm o}([-1,1])$. Lemma \ref{Lem.beta_surface} gives $\beta(\varepsilon)>1$ and $k>\beta(\varepsilon)(\ell+\sqrt\ell)$.
	
	On the sonic interval, \eqref{Eq.normalized_sonic_composition} and \eqref{Eq.normalized_sonic_real_bounds} give $c/C\leq b'=(\partial_{z_{\rm s}}W_{\rm s})/(\partial_{z_{\rm s}}X_{\rm s})\leq C/c$. On the outer compact interval, we prove the outgoing property for the limiting normalized field and use continuity in $\varepsilon$. We define $W_0:=(Z-v_0)/(1-Zv_0)$. Substitution of the $\ell=1$ equation gives
	\begin{equation}\label{Eq.base_outgoing_equation}
		W_0'=\frac{Z[k+1-\beta_0+(\beta_0-1)W_0^2]-kW_0}{Z(1-Z^2)},\qquad W_0'(0)=1-\frac{\beta_0}{k+1}>0.
	\end{equation}
	If $W_0'$ had a first zero in $(0,1)$, differentiating \eqref{Eq.base_outgoing_equation} at that zero would give
	\[W_0''=\frac{k+1-\beta_0+(\beta_0-1)W_0^2}{Z(1-Z^2)}>0,\]
	contradicting the first zero of $W_0'$, which is positive to the left of the first zero. Thus $W_0'>0$. Since $X'(Z)>0$, the limiting $b'$ is positive on the outer compact interval, and analytic dependence preserves uniform positive lower and finite upper bounds. Taking the minimum and maximum of the two overlapping intervals gives $\gamma_*,M_1$. For $X>0$ and $e=x/|x|$,
	\begin{align*}
		\frac{b(X)}X&=\frac1X\int_0^Xb'(s)\,\mathrm ds\geq\gamma_*,\\
		\zeta^T \mathrm D_x\bvec(x)\zeta &=b'(X)|\zeta\cdot e|^2+\frac{b(X)}X\bigl(|\zeta|^2-|\zeta\cdot e|^2\bigr) \geq\gamma_*|\zeta|^2.
	\end{align*}
	The same inequality holds at $X=0$ by continuity. Also $b(X_1)=1$ by \eqref{Eq.x_b}, so $M_1^{-1}\leq X_1\leq\gamma_*^{-1}$. Thus the normalized real intervals have lengths bounded above and below independently of $\varepsilon$.
\end{proof}

The derivative lower bound in \eqref{Eq.normalized_sonic_real_bounds} degenerates as $\varepsilon\to0$. Instead of applying an inverse theorem without quantitative constants, we compare $X_{\rm s}(z)$ with its affine approximation at a real center $z_r$. The second-derivative bound in Lemma \ref{Lem.sonic_expansion_bounds} makes this comparison uniform on disks of radius proportional to $(\varepsilon+z_r^2)^{1/2}$.

\begin{lemma}\label{Lem.complex_inverse_radius}
	There are $c_*,C_*,\varepsilon_4>0$, independent of $\varepsilon\in(0,\varepsilon_4)$, such that, with
	\[\mathcal T_\varepsilon:=\{z\in\C:\operatorname{dist}(z,[0,X_1])<c_*\varepsilon^{3/2}\},\]
	the scalar function $b$ has an extension $\widetilde b_\varepsilon$ satisfying
	\[\widetilde b_\varepsilon\in\Hol(\mathcal T_\varepsilon),\qquad \widetilde b_\varepsilon(X)=b(X)\quad(0\leq X\leq X_1),\qquad \sup_{z\in\mathcal T_\varepsilon}|\widetilde b_\varepsilon'(z)|\leq C_*.\]
\end{lemma}

\begin{proof}
	First consider a real sonic point $z_r\in[-r_*,r_*]$ and define $h_r:=(\varepsilon+z_r^2)^{1/2}$. Lemma \ref{Lem.sonic_expansion_bounds} gives $|(\partial_{z_{\rm s}}X_{\rm s})(z_r)|\geq c h_r^2$. On $|z_{\rm s}-z_r|<\theta h_r$, the bound for $\partial_{z_{\rm s}}^2X_{\rm s}$ gives
	\begin{equation}\label{Eq.sonic_inverse_derivative_variation}
		|(\partial_{z_{\rm s}}X_{\rm s})(z_{\rm s})-(\partial_{z_{\rm s}}X_{\rm s})(z_r)|\leq C\theta h_r^2,
	\end{equation}
	provided $\theta\leq1/4$ and $0<\varepsilon_4\leq\min\{\varepsilon_3,\varepsilon_2/4,r_*^2/4\}$. Indeed $|z_{\rm s}|\leq r_*+\theta\sqrt{\varepsilon_4+r_*^2}<3r_*/2$ and $\varepsilon<\varepsilon_2/2$, so the complex derivative estimate in Lemma \ref{Lem.sonic_expansion_bounds} applies on every such disk. Choose $\theta>0$ with $C\theta\leq c/4$. Integrating \eqref{Eq.sonic_inverse_derivative_variation} along a segment in $|z_{\rm s}-z_r|<\theta h_r$ gives
	\[|X_{\rm s}(z_2)-X_{\rm s}(z_1)-(\partial_{z_{\rm s}}X_{\rm s})(z_r)(z_2-z_1)|\leq\tfrac14|(\partial_{z_{\rm s}}X_{\rm s})(z_r)||z_2-z_1|.\]
	Thus $X_{\rm s}$ is injective on $|z_{\rm s}-z_r|<\theta h_r$. Let $|y-X_{\rm s}(z_r)|<|(\partial_{z_{\rm s}}X_{\rm s})(z_r)|\theta h_r/4$. On $|z-z_r|=\theta h_r$,
	\begin{align*}
		&|X_{\rm s}(z)-y-(\partial_{z_{\rm s}}X_{\rm s})(z_r)(z-z_r)|\\
		&\qquad\leq \tfrac14|(\partial_{z_{\rm s}}X_{\rm s})(z_r)|\theta h_r+|y-X_{\rm s}(z_r)| <\tfrac12|(\partial_{z_{\rm s}}X_{\rm s})(z_r)(z-z_r)|.
	\end{align*}
	Rouch\'e's theorem therefore gives a unique preimage of $y$ in the disk. Since $|\partial_{z_{\rm s}}X_{\rm s}(z_r)|\geq ch_r^2$, the image contains a disk of radius $c\theta h_r^3/4\geq c\theta\varepsilon^{3/2}/4$. Denote a smaller fixed coefficient by $c'>0$. On the inverse image of $|y-X_{\rm s}(z_r)|<c'\varepsilon^{3/2}$,
	\[\left|\frac{\partial_{z_{\rm s}}W_{\rm s}}{\partial_{z_{\rm s}}X_{\rm s}}\right| \leq\frac{C(\varepsilon+(|z_r|+\theta h_r)^2)}{3ch_r^2/4} \leq\frac{C(1+2+2\theta^2)h_r^2}{3ch_r^2/4}\leq C.\]
	The composition $W_{\rm s}\circ X_{\rm s}^{-1}$ is a holomorphic extension of $b$ with bounded derivative on $|y-X_{\rm s}(z_r)|<c'\varepsilon^{3/2}$.
	
	On the remaining compact interval down to the origin, joint holomorphic dependence and a nonzero $X'(Z)$ give inverse disks of a fixed radius with the same type of derivative bound, uniformly in $\varepsilon$. A finite subcover suffices there. Shrink $c'$ so that all sonic and outer disks have at least radius $2c_*\varepsilon^{3/2}$. The disks have real centers. The intersection of two overlapping equal-radius disks is connected and contains a real interval on which both continuations equal $b$. The identity theorem therefore gives a single $\widetilde b\in\Hol(\{z:\dist(z,[0,X_1])<c_*\varepsilon^{3/2}\})$ satisfying $|\widetilde b'|\leq C_*$. The disks at $0$ and $X_1$ are obtained from the local analytic extensions across both endpoints.
\end{proof}

Lemma \ref{Lem.complex_inverse_radius} gives $\widetilde b_\varepsilon$ on a tube of width $c_*\varepsilon^{3/2}$ with $|\widetilde b_\varepsilon'|\leq C_*$. We use the odd expansion of $b$ at zero and $\bvec(x)=b(|x|)x/|x|$ away from zero to obtain a Cartesian tube. Cauchy estimates for $\mathrm D_x\bvec$ must have constants independent of $\varepsilon$; for the other coefficients in \eqref{Eq.background_coefficient_map}, only the common tube width is required to be uniform.

\begin{lemma}\label{Lem.cartesian_uniform_background}
	There are $\rho_*,M_*>0$ and $\varepsilon_0>0$, independent of $\varepsilon\in(0,\varepsilon_0)$ and common to $d=4,5$, such that, with $R_\varepsilon:=\rho_*\varepsilon^{3/2}$,
	\[\|\nabla^j\bvec\|_{L^\infty(\Omega)}\leq M_*(j-1)!R_\varepsilon^{1-j}\qquad(j\geq1).\]
	The coefficient map $\mathbf c$ in \eqref{Eq.background_coefficient_map} extends holomorphically to $\Omega_{R_\varepsilon}^{\C}$. For every component $a$ there is $C_{a,\varepsilon}<\infty$ such that $\|\na^ja\|_{L^\infty}\leq C_{a,\varepsilon}j!R_\varepsilon^{-j}$ for all $j\geq0$. The constants $C_{a,\varepsilon}$ may depend on $\varepsilon$.
\end{lemma}

\begin{proof}
	Near $X=0$, the origin construction and the inverse coordinate are jointly analytic on fixed disks; $b$ is odd there. Write $b(X)=X\mathfrak b_0(X^2)$. Then
	\[\bvec(x)=\mathfrak b_0\left(\sum_{i=1}^dx_i^2\right)x\]
	defines a holomorphic Cartesian map on a fixed complex neighborhood of the origin, with bounded derivative uniformly in $\varepsilon$. Away from that neighborhood, take the positive branch of $(\sum_i x_i^2)^{1/2}$ at every real point. Since the real radius is bounded below there, the derivative of $(\sum_i x_i^2)^{1/2}$ is uniformly bounded on fixed local disks. Lemma \ref{Lem.complex_inverse_radius} therefore gives a holomorphic extension of $\bvec=b(X)x/X$ on Cartesian balls of radius $c\varepsilon^{3/2}$. The formulas for $\mathrm D_x\bvec$ involve only $b'$, $b/X$ and bounded factors $x_i/X$, so that derivative is uniformly bounded. On overlaps, the continuations agree on a real open set and hence agree holomorphically, by the identity theorem applied in each complex coordinate. The local extensions define $\bvec$ on a common complex tube.
	
	Apply the polydisc Cauchy estimate to each entry of the bounded matrix $\mathrm D_x\bvec$, on a polydisc of coordinate radius $c\varepsilon^{3/2}/(2\sqrt d)$ contained in the available complex ball. Each component derivative of order $j-1$ is bounded by $C(j-1)!(2\sqrt d/(c\varepsilon^{3/2}))^{j-1}$. Taking the Euclidean norm over the ordered derivative indices introduces at most $d^{(j-1)/2}$; the fixed matrix indices contribute another constant. The total exponential factor is at most $(2d)^{j-1}$. We incorporate the factor $(2d)^{j-1}$ by choosing $R_\varepsilon:=\rho_*\varepsilon^{3/2}$ when $\rho_*\leq c/40$, uniformly for $d\in\{4,5\}$. Increase $M_*$ to absorb the remaining fixed factors. Hence $\|\nabla^j\bvec\|_{L^\infty}\leq M_*(j-1)!R_\varepsilon^{1-j}$ for every $j\geq1$.
	
	We check the other coefficients on a tube of the same scale. Along the sonic curve the logarithmic phase derivative is the exact formula
	\begin{equation}\label{Eq.quantitative_phase_amplitude}
		\frac{\mathrm d}{\mathrm dX}\log\widehat\phi_0=\frac{(1-z_{\rm s})\sqrt{H_*}}{D_*},\qquad \widehat\rho_0^{p-1}=\frac{(\beta-1)\widehat\phi_0^2}{\mathcal M_*}.
	\end{equation}
	The first identity in \eqref{Eq.quantitative_phase_amplitude} follows by multiplying $\widehat\phi_0'/\widehat\phi_0=(\beta-1)v/(1-Zv)$ from \eqref{Eq.profiles} by $\mathrm dZ/\mathrm dX=\mathcal M_*$ from \eqref{Eq.normalized_sonic_functions}; the second follows from \eqref{Eq.profiles}. The expressions $(1-z_{\rm s})\sqrt{H_*}/D_*$ and $(\beta-1)\widehat\phi_0^2/\mathcal M_*$ are holomorphic on the inverse disks, with nonvanishing denominators. Local integration and the positive real normalization give holomorphic nonzero continuations of the phase and amplitude, using the logarithm for the positive root. On the outer compact interval the same conclusion follows directly from \eqref{Eq.profiles}. The continuations agree on real overlaps.
	
	The equation defining $J_0$ is
	\[\frac{\mathrm d}{\mathrm dX}\log J_0=\frac{dc_0-b'-kb/X}{b}.\]
	Away from zero, $b$ is bounded below on each outer real compact interval and on the sonic interval; the bound for $\widetilde b_\varepsilon'$ preserves nonvanishing after one common reduction of the tube. At zero the numerator is even and vanishes to second order, whereas $b=c_0X+\Order(X^3)$ with $c_0$ bounded below. The quotient is thus odd analytic on a fixed disk. Integration gives holomorphic functions $J_0$ and $J_0^{-1}$. The formula for $V$ uses finitely many derivatives of the amplitude, $\mathrm d/\mathrm dZ=\mathcal M_*^{-1}\mathrm d/\mathrm dX$, and the nonzero denominators in \eqref{Eq.normalized_sonic_functions} and \eqref{Eq.quantitative_phase_amplitude}. The quotient $Z^{-1}\widehat\rho_0'$ is even analytic at zero. Thus $V$ extends as well. Radial scalar functions are converted to Cartesian functions by the even power series at zero and the square-root branches away from zero, as above. Let $\widehat R_\varepsilon=c\varepsilon^{3/2}$ be a common tube radius obtained for $\mathbf b,J_0,J_0^{-1},V,\widehat\rho_0\circ Z,\widehat\phi_0\circ Z$, and let $M_{a,\varepsilon}$ bound $a$ on $\Omega_{\widehat R_\varepsilon/2}^{\C}$. Such a bound is finite by compactness. At every real point, the polydisc with coordinate radius $\widehat R_\varepsilon/(4\sqrt d)$ lies in $\Omega_{\widehat R_\varepsilon/2}^{\C}$. Cauchy's formula and the ordered-tensor sum give
	\[\|\na^ja\|_{L^\infty}\leq M_{a,\varepsilon}j!\left(\frac{4d}{\widehat R_\varepsilon}\right)^j\qquad(j\geq0).\]
	Reduce the already chosen $\rho_*$ once, so that $R_\varepsilon\leq\widehat R_\varepsilon/20$ for both dimensions. Since $R_\varepsilon\leq\widehat R_\varepsilon/20$ and $d\leq5$, we obtain $\|\nabla^ja\|_{L^\infty}\leq M_{a,\varepsilon}j!R_\varepsilon^{-j}$. Reducing $R_\varepsilon$ preserves the estimate for $\nabla^j\mathbf b$.
\end{proof}

\subsection{Verification of the two assumptions}
\label{Subsec.verify_assumptions}

We now choose constants common to $d=4,5$ and set $\varepsilon=4/(p-1)$. Lemmas \ref{Lem.complete_background} and \ref{Lem.cartesian_uniform_background} give the profile bounds \eqref{Eq.background_bounds}; Lemma \ref{Lem.small_scattering} independently gives \eqref{Eq.chi_p_bound}. We finish by checking that the choices \eqref{Eq.quantitative_radii} satisfy the three inequalities in \eqref{Eq.analytic_margin}.

\begin{proof}[Proof of Proposition \ref{Prop.quantitative_background}]
	Apply Lemmas \ref{Lem.match_origin}--\ref{Lem.cartesian_uniform_background} for $k=3$ and $k=4$, and take the smaller parameter thresholds and analytic radius constants and the larger derivative constants for $k=3,4$. We obtain common $\varepsilon_0,\gamma_*,M_*,\rho_*>0$ and analytic functions $\beta_d(\varepsilon)$ with $\beta_d(0)=(d-1)/2-1/3$. For every $0<\varepsilon<\varepsilon_0$, the constructed $v_\varepsilon$ solves the velocity equation, extends oddly and analytically across the origin, obeys $|v_\varepsilon|<1$ on $[0,1]$, and has $\beta_d(\varepsilon)>1$. Thus for $\varepsilon=4/(p-1)$ the pair $(\beta_d(\varepsilon),v_\varepsilon)$ satisfies Assumption \ref{Ass.background}. Lemmas \ref{Lem.complete_background} and \ref{Lem.cartesian_uniform_background} give the uniform expansion and derivative estimates with $R_p=\rho_*(4/(p-1))^{3/2}$. Lemma \ref{Lem.small_scattering} gives $\chi_p\leq C_{\rm sc}(p+1)^{-2}$.
	
	Let $p_{\rm bg}$ be the least odd integer $P\geq15$ such that $4/(P-1)\leq\varepsilon_0/2$. The bounds of Lemmas \ref{Lem.complete_background} and \ref{Lem.cartesian_uniform_background} therefore hold for every odd $p\geq p_{\rm bg}$, proving Proposition \ref{Prop.quantitative_background}.
\end{proof}

For completeness, the choices \eqref{Eq.quantitative_radii} give the following direct verification of Assumption \ref{Ass.analytic_margin}. Let $\vartheta=\min\{1/16,\gamma_*/(32M_*)\}$, $\varrho_p=\vartheta R_p$ and $\chi=2C_{\rm sc}(p+1)^{-2}$. Then $\chi>\chi_p$, and
\[\mathcal L_v(4\varrho_p)\leq M_*\sum_{j\geq2}\frac{(4\vartheta)^{j-1}}j\leq4M_*\vartheta\leq\gamma_*/8.\]
Also
\[\frac\chi{\varrho_p}=\frac{C_{\rm sc}}{4\vartheta\rho_*}\frac{(p-1)^{3/2}}{(p+1)^2}\leq\frac{C_{\rm sc}}{4\vartheta\rho_*}(p-1)^{-1/2}.\]
The choice \eqref{Eq.p0_choice} gives $\chi/\varrho_p\leq\gamma_*/8$, and ensures $4/(p-1)<\varepsilon_0$. An odd integer satisfying \eqref{Eq.p0_choice} exists since $C_{\rm sc},\vartheta,\rho_*,\gamma_*>0$. Hence for every odd $p\geq p_0$, the constructed background and \eqref{Eq.quantitative_radii} satisfy
\[4\varrho_p<R_p,\qquad \gamma_v-\mathcal L_v(4\varrho_p)\geq7\gamma_*/8,\qquad \gamma_v-\mathcal L_v(\varrho_p)-\chi/\varrho_p\geq3\gamma_*/4.\]
Thus $(\beta_d(\varepsilon),v_\varepsilon,R_p,\chi,\varrho_p)$ satisfies Assumption \ref{Ass.analytic_margin} for every odd $p\geq p_0$, as stated in Corollary \ref{Cor.quantitative_margin}.

\appendix
\section{Euler profiles and the normalized equation}\label{App.analyticity}

In Appendix \ref{App.analyticity}, a pair $(\beta,v)$ satisfying Assumption \ref{Ass.background} is fixed. Subsections \ref{App.sonic}--\ref{App.regular_singular} prove $v\in C^\omega_{\rm o}([-1,1])$ and \eqref{Eq.profile_tube}. Subsection \ref{App.coefficients} proves the expansion bound for $\bvec$ and the extension \eqref{Eq.qualitative_coefficient_extension}. Subsection \ref{App.normalization} proves Lemma \ref{Lem.normalized_equation}.

\begin{proposition}\label{Prop.automatic_analyticity}
	Let $k\in\Z_{\geq1}$, $\ell>1$ and $\beta>1$. Suppose that $v\in C^\infty_{\rm o}([0,1])$ satisfies \eqref{Eq.velocity} and $|v(Z)|<1$ for every $Z\in[0,1]$. Then $k>\beta(\ell+\sqrt\ell)$ and the odd extension satisfies $v\in C^\omega_{\rm o}([-1,1])$. There exist $\eta>0$ and $\widetilde v\in\Hol(\mathcal U_\eta)$ such that
	\[\widetilde v|_{[-1,1]}=v,\qquad \widetilde v(-\zeta)=-\widetilde v(\zeta),\qquad \widetilde v(\overline\zeta)=\overline{\widetilde v(\zeta)} \quad(\zeta\in\mathcal U_\eta).\]
	There also exist constants $C_v,R_Z>0$ such that
	\begin{equation}\label{Eq.A_profile_Cauchy}
		\max_{0\leq Z\leq1}|v^{(j)}(Z)|\leq C_v j!R_Z^{-j}\qquad(j\in\Z_{\geq0}).
	\end{equation}
\end{proposition}

We first identify the singular points of the equation. We then prove a regular singular ODE lemma that applies at both the origin and the sonic point.

\subsection{The singular points of the profile equation}\label{App.sonic}

We recall the sonic geometry proved in \cite[Lemma~2.1 and Lemmas~A.2--A.3]{SWZ-NLW}. The arguments locating the sonic point and proving its simplicity use only the profile on $[0,1]$; they apply under Assumption \ref{Ass.background} without extending $v$ beyond $Z=1$. We keep the conclusions for reference in the analytic continuation argument.

\begin{lemma}\label{Lem.unique_sonic_point}
	Under the assumptions of Proposition \ref{Prop.automatic_analyticity}, one has $k>\beta(\ell+\sqrt\ell)$. Define
	\begin{equation}\label{Eq.A_sonic_location}
		Z_1:=\frac{k}{\sqrt\ell\,[k-\beta(\ell-1)]},\qquad v_1:=\frac{\beta\sqrt\ell}{k-\beta\ell},\qquad 0<v_1<Z_1<1.
	\end{equation}
	Then, with $\Delta_0(Z):=\Delta_Z(Z,v(Z))$,
	\begin{align*}
		&0<v_1=v(Z_1)<Z_1<1,\quad\{Z\in(0,1]:\Delta_0(Z)=0\}=\{Z_1\},\\
		&\Delta_0(Z)>0\quad(0<Z<Z_1),\qquad \Delta_0(Z)<0\quad(Z_1<Z\leq1),\\
		&\qquad\qquad\Delta_0'(0)=1,\qquad v'(0)=\beta\ell/(k+1).
	\end{align*}
\end{lemma}

The sonic location, uniqueness and signs follow from \cite[Lemma~2.1 and Lemmas~A.2--A.3]{SWZ-NLW}. The values at zero follow by inserting $v(0)=0$ into \eqref{Eq.velocity} and differentiating once. At the right endpoint, \eqref{Eq.velocity} also gives
\begin{equation}\label{Eq.A_endpoint}
\Delta_0(1)=(1-\ell)(1-v(1))^2<0.
\end{equation}

The simple-zero property in \cite[Lemma~A.3]{SWZ-NLW}, together with the signs in Lemma \ref{Lem.unique_sonic_point}, gives $\Delta_0'(Z_1)<0$. We collect the geometric conclusions used to reduce the profile equation to a regular singular equation.

\begin{lemma}\label{Lem.sonic_geometry}
	Under the assumptions of Proposition \ref{Prop.automatic_analyticity}, set $\Delta_0(Z):=\Delta_Z(Z,v(Z))$. Then $k>\beta(\ell+\sqrt\ell)$ and, with $Z_1,v_1$ defined in \eqref{Eq.A_sonic_location},
	\begin{equation}\label{Eq.A_sonic_properties}
		\begin{aligned}
			v(Z_1)&=v_1,\qquad \{Z\in(0,1]:\Delta_0(Z)=0\}=\{Z_1\},\\
			\Delta_0(Z)&>0\quad(0<Z<Z_1),\qquad \Delta_0(Z)<0\quad(Z_1<Z\leq1),\\
			\Delta_0'(0)&=1,\qquad v'(0)=\frac{\beta\ell}{k+1},\qquad \Delta_0'(Z_1)<0.
		\end{aligned}
	\end{equation}
\end{lemma}

\subsection{Analyticity at a regular singular point}\label{App.regular_singular}

At the zeros identified in Lemma \ref{Lem.sonic_geometry}, we will use the following regular singular ODE lemma.

\begin{lemma}\label{Lem.singular_integral}
	Let $\alpha,a_0,b_0>0$ and $K\in\Hol(\mathbb D_{a_0}\times\mathbb D_{b_0})$. There exists $\delta_1\in(0,a_0)$ such that, for each $0<\delta\leq\delta_1$, the equation
	\[tw'(t)+\alpha w(t)=tK(t,w(t))\]
	has a solution $w_\delta\in C(\overline{\mathbb D_\delta})\cap\Hol(\mathbb D_\delta)$ with $w_\delta(0)=0$ and $\|w_\delta\|_{L^\infty}\leq b_0/2$. If $J$ is any of $[-\delta,\delta]$, $[0,\delta]$ and $[-\delta,0]$, every $z\in C(J;\C)\cap C^1(J\setminus\{0\};\C)$ satisfying $tz'(t)+\alpha z(t)=tK(t,z(t))$ on $J\setminus\{0\}$ and $\sup_J|z|\leq b_0/2$ obeys $z(t)=w_\delta(t)$ for every $t\in J$.
\end{lemma}

\begin{proof}
	A bounded solution on a real interval from zero satisfies
	\begin{equation}\label{Eq.A_integral_fixedpoint}
		z(t)=t\int_0^1s^\alpha K(st,z(st))\,\mathrm ds.
	\end{equation}
	Indeed,
	\[\frac{\mathrm d}{\mathrm ds}\bigl(s^\alpha z(st)\bigr)=ts^\alpha K(st,z(st)),\qquad \lim_{s\downarrow0}s^\alpha z(st)=0.\]
	Integration proves \eqref{Eq.A_integral_fixedpoint} for either sign of $t$.
	
	We define $a_1:=a_0/2$, $b_1:=b_0/2$ and
	\[M_1:=\sup_{|t|\leq a_1,\,|u|\leq b_1}|K(t,u)|,\qquad L_1:=\sup_{|t|\leq a_1,\,|u|\leq b_1}|\partial_2K(t,u)|.\]
	Both constants are finite. Set
	\begin{equation}\label{Eq.A_contraction_radius}
		\delta_1:=\min\left\{\frac{a_1}{2},\frac{b_1(\alpha+1)}{2(M_1+1)},\frac{\alpha+1}{2(L_1+1)}\right\}>0.
	\end{equation}
	For $0<\delta\leq\delta_1$, work on the closed ball
	\[\mathscr B_\delta:=\{w\in C(\overline{\mathbb D_\delta})\cap\Hol(\mathbb D_\delta):\|w\|_{L^\infty}\leq b_1\}\]
	in the Banach space of continuous functions on the closed disc that are holomorphic on $\mathbb D_\delta$. Define
	\[(\mathscr T_\delta w)(t):=t\int_0^1s^\alpha K(st,w(st))\,\mathrm ds.\]
	Uniform convergence on compact subsets shows that $\mathscr T_\delta w$ is holomorphic, and
	\[\|\mathscr T_\delta w\|_{L^\infty}\leq\frac{\delta M_1}{\alpha+1}\leq\frac{b_1}{2},\qquad \|\mathscr T_\delta w_1-\mathscr T_\delta w_2\|_{L^\infty}\leq\frac{\delta L_1}{\alpha+1}\|w_1-w_2\|_{L^\infty}\leq\frac12\|w_1-w_2\|_{L^\infty}.\]
	Thus $\mathscr T_\delta$ has a unique fixed point $w_\delta\in\mathscr B_\delta$. With $k(t):=K(t,w_\delta(t))$, differentiation of the integral and integration by parts give
	\[tw_\delta'(t)+\alpha w_\delta(t)=t\int_0^1\bigl((\alpha+1)s^\alpha k(st)+ts^{\alpha+1}k'(st)\bigr)\,\mathrm ds=tk(t).\]
	At $t=0$, the equation follows from $w_\delta(0)=0$. On the closed ball $\{z\in C(J;\C):\|z\|_{L^\infty}\leq b_1\}$, the integral map has Lipschitz constant at most $1/2$. Every solution with $\sup_J|z|\leq b_1$ satisfies \eqref{Eq.A_integral_fixedpoint}; hence $z=w_\delta|_J$.
\end{proof}

We subtract a Taylor polynomial from a smooth solution of $th'=F(t,h)$. The remainder satisfies the equation in Lemma \ref{Lem.singular_integral}, whose uniqueness identifies the remainder with a holomorphic function.

\begin{lemma}\label{Lem.smooth_regular_singular}
	Let $a_0,b_0>0$ and $F\in C^\omega((-a_0,a_0)\times(-b_0,b_0))$, with $F(0,0)=0$. Let $0<a<a_0$, and let $I$ be one of $(-a,a)$, $[0,a)$ and $(-a,0]$. Suppose that $h\in C^\infty(I)$, $h(0)=0$, $h(I)\subset(-b_0,b_0)$, and
	\[th'(t)=F(t,h(t))\qquad(t\in I).\]
	Then there exist $\delta\in(0,a)$ and $h_{\rm an}\in\Hol(\mathbb D_\delta)$ such that
	\[h_{\rm an}(t)=h(t)\qquad(t\in I,\ |t|<\delta).\]
	
\end{lemma}

To prove Lemma \ref{Lem.smooth_regular_singular}, we first subtract a sufficiently long Taylor polynomial. The remainder equation then has the positive coefficient required in Lemma \ref{Lem.singular_integral}. Lemma \ref{Lem.singular_Taylor_reduction} gives the transformed equation and the vanishing order of the remainder.

\begin{lemma}\label{Lem.singular_Taylor_reduction}
	Let $F$ and $h$ satisfy the hypotheses of Lemma \ref{Lem.smooth_regular_singular}. We define
	\[\mu=\partial_2F(0,0),\qquad N=1+\max\{0,\lfloor\mu\rfloor\},\qquad P_N(t)=\sum_{j=1}^N\frac{h^{(j)}(0)}{j!}t^j.\]
	There are $r_0,b_N>0$ and $K\in\Hol(\mathbb D_{r_0}\times\mathbb D_{b_N})$ such that the smooth function $z(t)=(h(t)-P_N(t))/t^N$, extended by $z(0)=0$, satisfies
	\[tz'(t)+(N-\mu)z(t)=tK(t,z(t))\]
	whenever $t\in I$, $|t|<r_0$ and $|z(t)|<b_N$.
\end{lemma}

\begin{proof}
	We define
	\[\mu:=\partial_2F(0,0),\quad N:=1+\max\{0,\lfloor\mu\rfloor\},\quad \alpha:=N-\mu>0, \quad P_N(t):=\sum_{j=1}^N\frac{h^{(j)}(0)}{j!}t^j.\]
	The definition of $C^\omega$ gives $a_F\in(0,a_0)$, $b_F\in(0,b_0)$ and
	\[\widetilde F\in\Hol(\mathbb D_{a_F}\times\mathbb D_{b_F}),\qquad \widetilde F(t,u)=F(t,u)\quad(|t|<a_F,\ |u|<b_F,\ t,u\in\R).\]
	Choose $r_0>0$ such that
	\[r_0<\min\{a/2,a_F/2,1\},\qquad \sup_{|t|\leq r_0}|P_N(t)|\leq b_F/4, \qquad b_N:=b_F/4.\]
	The identity $P_N(0)=0$ ensures the existence of $r_0$. For $|t|<r_0$, $|w|<b_N$ and $0\leq\vartheta\leq1$, one has
	$|P_N(t)+\vartheta t^Nw|<b_F/2$.
	
	Taylor's formula on the real interval $I$ gives
	$h-P_N=\Order(t^{N+1})$ and $h'-P_N'=\Order(t^N)$. Thus, for real $t\in I$ with $|t|<r_0$,
	\begin{align*}
		\widetilde F(t,P_N(t))-tP_N'(t) =F(t,P_N(t))-F(t,h(t))+t(h'(t)-P_N'(t))=\Order(t^{N+1}).
	\end{align*}
	The function $\widetilde F(t,P_N(t))-tP_N'(t)$ is holomorphic on $\mathbb D_{r_0}$. The bound $O(t^{N+1})$ on a real half-interval implies that the Taylor coefficients of orders $0,\ldots,N$ vanish. Therefore there exists $r_N\in\Hol(\mathbb D_{r_0})$ such that
	\begin{equation}\label{Eq.A_Taylor_residual}
		\widetilde F(t,P_N(t))-tP_N'(t)=t^{N+1}r_N(t) \qquad(t\in\mathbb D_{r_0}).
	\end{equation}
	For $t\in I$, define $z(t):=(h(t)-P_N(t))/t^N$ when $t\ne0$ and $z(0):=0$. The integral remainder formula gives
	\[z(t)=\frac{t}{N!}\int_0^1(1-s)^Nh^{(N+1)}(st)\,\mathrm ds,\]
	so $z\in C^\infty(I)$. Define, on $\mathbb D_{r_0}\times\mathbb D_{b_N}$,
	\[A_N(t,w):=\int_0^1\partial_2\widetilde F(t,P_N(t)+\vartheta t^Nw)\,\mathrm d\vartheta, \qquad C_N(t,w):=\int_0^1\partial_1A_N(st,w)\,\mathrm ds.\]
	Differentiation under the integrals defining $A_N$ and $C_N$ is valid on compact subsets. Hence
	\[A_N,C_N\in\Hol(\mathbb D_{r_0}\times\mathbb D_{b_N}),\qquad A_N(0,w)=\mu,\qquad A_N(t,w)=\mu+tC_N(t,w).\]
	The identity
	\[\widetilde F(t,P_N(t)+t^Nw)-\widetilde F(t,P_N(t))=t^NwA_N(t,w)\]
	and \eqref{Eq.A_Taylor_residual} show that, wherever $|t|<r_0$ and $|z(t)|<b_N$,
	\begin{equation}\label{Eq.A_reduced_remainder}
		tz'(t)+\alpha z(t)=tK(t,z(t)),\qquad K(t,w):=r_N(t)+C_N(t,w)w.
	\end{equation}
	Here $K\in\Hol(\mathbb D_{r_0}\times\mathbb D_{b_N})$ and $K$ is real-valued for real arguments.
\end{proof}

We can now apply the analytic integral equation to the remainder. Since $z(0)=0$, the remainder lies in the uniqueness ball of Lemma \ref{Lem.singular_integral} on a sufficiently small interval.

\begin{proof}[Proof of Lemma \ref{Lem.smooth_regular_singular}]
	By Lemma \ref{Lem.singular_Taylor_reduction}, the remainder $z$ satisfies \eqref{Eq.A_reduced_remainder} with $\alpha=N-\mu>0$.
	Apply Lemma \ref{Lem.singular_integral} with the parameters $\alpha,r_0,b_N$ and the function $K$ in \eqref{Eq.A_reduced_remainder}. Lemma \ref{Lem.singular_integral} gives $\delta_1\in(0,r_0)$ and, for every $0<\delta\leq\delta_1$, an analytic solution $w_\delta$ on $\mathbb D_\delta$. Since $z(0)=0$, choose such a $\delta$ with
	\[\sup_{t\in I\cap[-\delta,\delta]}|z(t)|\leq b_N/2.\]
	Lemma \ref{Lem.singular_integral}, applied on $I\cap[-\delta,\delta]$, gives $z(t)=w_\delta(t)$ for every $t\in I$ with $|t|\leq\delta$. Therefore
	\[h_{\rm an}(t):=P_N(t)+t^Nw_\delta(t)\in\Hol(\mathbb D_\delta),\qquad h_{\rm an}(t)=h(t)\quad(t\in I,\ |t|<\delta).\qedhere\]
\end{proof}

We apply the regular singular lemma at the origin and at the sonic point. Away from $0$ and $Z_1$ the ordinary analytic ODE theorem applies directly.

\begin{lemma}\label{Lem.local_profile_extension}
	Under the assumptions of Proposition \ref{Prop.automatic_analyticity}, for each $Z_*\in\{0,Z_1\}$ there are $\delta_*>0$ and $\widetilde v_*\in\Hol(\mathbb D_{\delta_*})$ such that
	\[\widetilde v_*(t)=v(Z_*+t)\qquad(|t|<\delta_*,\ Z_*+t\in[0,1]).\]
\end{lemma}

\begin{proof}
	Fix $Z_*\in\{0,Z_1\}$ and define $v_*:=v(Z_*)$, $m_*:=v'(Z_*)$. Write
	\[t:=Z-Z_*,\qquad v(Z_*+t)=v_*+t\bigl(m_*+h(t)\bigr),\qquad h(0)=0.\]
	Taylor's formula gives $h\in C^\infty(I_*)$, where $I_*$ is an interval about zero when $Z_*=Z_1$ and a right half interval when $Z_*=0$. Define the polynomials
	\begin{align*}
		\mathcal A(t,u)&:=\int_0^1 (\partial_Z+(m_*+u)\partial_v)\Delta_Z (Z_*+st,v_*+st(m_*+u))\,\mathrm ds,\\
		\mathcal B(t,u)&:=\int_0^1 (\partial_Z+(m_*+u)\partial_v)\Delta_v (Z_*+st,v_*+st(m_*+u))\,\mathrm ds.
	\end{align*}
	Since $\Delta_Z(Z_*,v_*)=\Delta_v(Z_*,v_*)=0$, the fundamental theorem of calculus gives
	\begin{align*}
		\Delta_Z(Z_*+t,v_*+t(m_*+u))&=t\mathcal A(t,u),\\
		\Delta_v(Z_*+t,v_*+t(m_*+u))&=t\mathcal B(t,u).
	\end{align*}
	In particular, $\mathcal A,\mathcal B\in C^\omega(\R^2)$ and
	\[\mathcal A(0,0) =\partial_Z\Delta_Z(Z_*,v_*)+m_*\partial_v\Delta_Z(Z_*,v_*) =\Delta_0'(Z_*)\ne0\]
	by Lemma \ref{Lem.sonic_geometry}. Choose $a_*,b_*>0$ such that
	\[\sup_{|t|\leq a_*,\,|u|\leq b_*} |\mathcal A(t,u)-\mathcal A(0,0)| \leq\tfrac12|\mathcal A(0,0)|, \qquad U_*:=(-a_*,a_*)\times(-b_*,b_*).\]
	Then $\mathcal A^{-1}\in C^\omega(U_*)$. For $t\ne0$, substitution into \eqref{Eq.velocity} gives
	\[\mathcal A(t,h(t))\bigl(m_*+h(t)+th'(t)\bigr) =\mathcal B(t,h(t)).\]
	Taking $t\to0$ along $I_*$ yields $\mathcal B(0,0)=m_*\mathcal A(0,0)$. Therefore
	\begin{align*}
		th'(t)=F_*(t,h(t)),\quad F_*(t,u):=\frac{\mathcal B(t,u)}{\mathcal A(t,u)}-m_*-u,\quad F_*\in C^\omega(U_*),\quad F_*(0,0)=0.
	\end{align*}
	By $h(0)=0$ and continuity, there is $a\in(0,a_*)$ such that the appropriate interval $(-a,a)$ or $[0,a)$ is contained in $I_*$ and $|h(t)|<b_*$ there. Lemma \ref{Lem.smooth_regular_singular} gives $\delta_*>0$ and $h_*\in\Hol(\mathbb D_{\delta_*})$ agreeing with $h$ on $I_*\cap(-\delta_*,\delta_*)$. Hence the function
	\[\widetilde v_*(t):=v_*+t\bigl(m_*+h_*(t)\bigr),\qquad \widetilde v_*\in\Hol(\mathbb D_{\delta_*}),\]
	extends $v(Z_*+t)$ across $t=0$. At the origin, $\partial_2F_*(0,0)=-(k+1)$; at $Z_1$ the lemma applies for every value of $\partial_2F_*(0,0)$.
\end{proof}

The profile is analytic near both singular points and solves a regular analytic equation elsewhere. We patch the local extensions of $v$ and then apply Cauchy's formula on a smaller common complex neighborhood.

\begin{proof}[Proof of Proposition \ref{Prop.automatic_analyticity}]
	Lemma \ref{Lem.sonic_geometry} gives the parameter inequality and identifies the singular points. Lemma \ref{Lem.local_profile_extension} gives the analytic continuations at $Z=0,Z_1$.
	For $Z_*\in(0,1]\setminus\{Z_1\}$, choose a real open rectangle $U_{Z_*}$ containing $(Z_*,v(Z_*))$ on which $\Delta_Z\ne0$. Then
	\[\frac{\Delta_v}{\Delta_Z}\in C^\omega(U_{Z_*}).\]
	The local ODE existence theorem gives an open interval $J_{Z_*}\ni Z_*$ and $v_{Z_*}\in C^\omega(J_{Z_*})$ solving $v_{Z_*}'=\Delta_v(Z,v_{Z_*})/\Delta_Z(Z,v_{Z_*})$ with $v_{Z_*}(Z_*)=v(Z_*)$. Local uniqueness identifies $v_{Z_*}$ with $v$ on $J_{Z_*}\cap[0,1]$. At $Z_*=1$, \eqref{Eq.A_endpoint} gives the same conclusion across the endpoint. Combining $v_{Z_*}$ with the continuations at $0$ and $Z_1$ proves
	\[v\in C^\omega([0,1]).\]
	By \eqref{Eq.even_odd_spaces}, $v^{(2j)}(0)=0$ for every $j\geq0$. Thus the power series at zero contains only odd powers, and the odd extension satisfies
	\[v\in C^\omega_{\rm o}([-1,1]).\]
	At each real point of a real open interval containing $[-1,1]$, the local series defines a function in $\Hol(\mathbb D_r(Z_*))$ for some $r>0$. Here $\mathbb D_r(Z_*):=\{\zeta\in\C:|\zeta-Z_*|<r\}$. Choose the discs $\mathbb D_r(Z_*)$ so that the real diameters lie in the interval on which the extended profile $v$ is real analytic. Two overlapping discs centered on the real axis have connected intersection containing a real interval; both local holomorphic functions equal $v$ on the real interval. The identity theorem identifies the local functions on the whole intersection. Finitely many discs cover $[-1,1]$, and the union contains $\mathcal U_\eta$ for some $\eta>0$. Restriction of the patched function gives
	\[\widetilde v\in\Hol(\mathcal U_\eta),\qquad \widetilde v|_{[-1,1]}=v.\]
	On $[-1,1]$, we have $\widetilde v(-Z)=-\widetilde v(Z)$ and $\overline{\widetilde v(Z)}=\widetilde v(Z)$. Applying the identity theorem to the holomorphic functions $\widetilde v(-\zeta)$, $-\widetilde v(\zeta)$ and $\overline{\widetilde v(\overline\zeta)}$ gives $\widetilde v(-\zeta)=-\widetilde v(\zeta)$ and $\widetilde v(\overline\zeta)=\overline{\widetilde v(\zeta)}$ on $\mathcal U_\eta$.
	
	Finally, set
	\[R_Z:=\eta/4,\qquad C_v:=\sup_{\dist(\zeta,[-1,1])\leq\eta/2}|\widetilde v(\zeta)|<\infty.\]
	For every $Z\in[0,1]$, the closed complex disc of radius $R_Z$ about $Z$ lies in $\{\zeta:\operatorname{dist}(\zeta,[-1,1])\leq\eta/2\}$. Cauchy's integral formula gives \eqref{Eq.A_profile_Cauchy} for every $j\geq0$.
\end{proof}

\subsection{The outgoing property and analyticity of the normalized coefficients}\label{App.coefficients}

We now return to $k=d-1$, $\ell=1+4/(p-1)$ and the normalization in Subsection \ref{Subsec.leading}. Proposition \ref{Prop.automatic_geometry} proves $\sym \mathrm D_x\mathbf b\geq\gamma_vI_d$ and the common analytic extension of the normalized coefficients.

\begin{proposition}\label{Prop.automatic_geometry}
	Under Assumption \ref{Ass.background}, the maps in \eqref{Eq.x_b} satisfy
	\[X\in C^\omega_{\rm o}([-1,1]),\qquad Z\in C^\omega_{\rm o}([-X_1,X_1]),\qquad b\in C^\omega_{\rm o}([-X_1,X_1]),\]
	where $Z=Z(X)$ is the inverse of $X=X(Z)$ and $X'(Z)>0$ on $[-1,1]$. Moreover,
	\begin{equation}\label{Eq.A_strict_expansion}
		b(0)=0,\qquad b(X_1)=1,\qquad b'(X)>0\quad(0\leq X\leq X_1).
	\end{equation}
	Consequently, with $\gamma_v=\min_{0\leq X\leq X_1}b'(X)>0$, one has
	$\sym\mathrm D_x\bvec(x)\geq\gamma_v I_d$ on $\overline\Omega$.
	
	For $e_1:=(1,0,\ldots,0)\in\R^d$, we define $j_0(X):=J_0(Xe_1)$ and $V_{\rm r}(X):=V(Xe_1)$. Then
	\[\widehat\rho_0,\widehat\phi_0\in C^\omega_{\rm e}([-1,1]),\qquad j_0,j_0^{-1},V_{\rm r}\in C^\omega_{\rm e}([-X_1,X_1]).\]
	The map $\mathbf c$ in \eqref{Eq.background_coefficient_map} satisfies
	$\mathbf c\in C^\omega(\overline\Omega;\R^{d+5})$. There exist $R_v>0$ and
	\begin{equation}\label{Eq.A_Cartesian_tube}
		\widetilde{\mathbf c}\in\Hol(\Omega_{R_v}^{\C};\C^{d+5}),\qquad \widetilde{\mathbf c}|_{\overline\Omega}=\mathbf c,\qquad \Omega_{R_v}^{\C}:=\{\zeta\in\C^d:\dist(\zeta,\overline\Omega)<R_v\}.
	\end{equation}
	There is $M_v>0$ such that
	\[\norm{\na^j\bvec}_{L^\infty(\Omega)}\leq M_v(j-1)!R_v^{1-j}\qquad(j\geq1).\]
	Each scalar component $a$ of $\mathbf c$ satisfies
	$\norm{\na^ja}_{L^\infty(\Omega)}\leq C_a j!R_v^{-j}$ for some $C_a<\infty$ and every $j\geq0$. All constants may depend on the fixed profile.
\end{proposition}

We prove Proposition \ref{Prop.automatic_geometry} in two steps. First, Lemmas \ref{Lem.scalar_expansion_identity}--\ref{Lem.qualitative_expansion} establish the outgoing property of the normalized radial field. We then treat the scalar coefficients and the Cartesian extensions of $J_0,J_0^{-1},V,\widehat\rho_0\circ Z,\widehat\phi_0\circ Z$ in Lemmas \ref{Lem.qualitative_scalars} and \ref{Lem.qualitative_cartesian}. We start with the scalar quantity $(Z-v)/(1-Zv)$, whose derivative determines the sign of $b'$.

\begin{lemma}\label{Lem.scalar_expansion_identity}
	Under Assumption \ref{Ass.background}, we define $\mathfrak b(Z):=(Z-v(Z))/(1-Zv(Z))$. Then
	\[\mathfrak b\in C^\omega_{\rm o}([-1,1]),\qquad \mathfrak b'(0)>0,\quad \mathfrak b'(Z_1)>0,\quad \mathfrak b'(1)>0.\]
	For $0<Z<1$,
	\begin{equation}\label{Eq.A_bfrak_equation}
		\begin{aligned}
			Z(1-Z^2)(1-\ell\mathfrak b^2)\mathfrak b'&=(1-\mathfrak b^2)\bigl[Z A(\mathfrak b)-k\mathfrak b\bigr],\\
			A(a)&:=k+1-\beta\ell+\ell(\beta-1)a^2>0\qquad(a\in\R).
		\end{aligned}
	\end{equation}
	If $Z\in(0,1)\setminus\{Z_1\}$ and $\mathfrak b'(Z)=0$, then
	\begin{equation}\label{Eq.A_bfrak_critical}
		\mathfrak b''(Z)=\frac{(1-\mathfrak b(Z)^2)A(\mathfrak b(Z))}{Z(1-Z^2)(1-\ell\mathfrak b(Z)^2)}.
	\end{equation}
\end{lemma}

\begin{proof}
	Let
	\begin{equation}\label{Eq.A_bfrak}
		\mathfrak b(Z):=\frac{Z-v(Z)}{1-Zv(Z)}=b(X(Z)).
	\end{equation}
	The denominator is positive on $[-1,1]$. Proposition \ref{Prop.automatic_analyticity} gives $\mathfrak b\in C^\omega_{\rm o}([-1,1])$. One has
	\[1-\mathfrak b(Z)^2=\frac{(1-Z^2)(1-v(Z)^2)}{(1-Zv(Z))^2},\qquad \mathfrak b(0)=0,\qquad \mathfrak b(1)=1.\]
	In particular, $|\mathfrak b(Z)|<1$ for $0\leq Z<1$. Differentiating \eqref{Eq.A_bfrak} gives
	\[\mathfrak b'(Z)=\frac{1-v(Z)^2-(1-Z^2)v'(Z)}{(1-Zv(Z))^2}.\]
	Thus Lemma \ref{Lem.sonic_geometry} gives
	\[\mathfrak b'(0)=1-\frac{\beta\ell}{k+1}>0,\qquad \mathfrak b'(1)=\frac{1+v(1)}{1-v(1)}>0.\]
	At the sonic point, \eqref{Eq.A_sonic_location} gives $\mathfrak b(Z_1)=\ell^{-1/2}$. The factorization
	\[\Delta_0(Z)=Z(1-Zv(Z))^2(1-\ell\mathfrak b(Z)^2)\]
	therefore yields
	\[\Delta_0'(Z_1)=-2\ell Z_1(1-Z_1v_1)^2\mathfrak b(Z_1)\mathfrak b'(Z_1)<0,\]
	so $\mathfrak b'(Z_1)>0$.

	For $0<Z<1$, the inverse relation and the formula for $v'$ are
	\[v=\frac{Z-\mathfrak b}{1-Z\mathfrak b},\qquad v'=\frac{1-\mathfrak b^2-(1-Z^2)\mathfrak b'}{(1-Z\mathfrak b)^2}.\]
	Substitute the formulas for $v$ and $v'$ into \eqref{Eq.velocity}. The identities
	\[1-Zv=\frac{1-Z^2}{1-Z\mathfrak b},\qquad 1-v^2=\frac{(1-Z^2)(1-\mathfrak b^2)}{(1-Z\mathfrak b)^2}\]
	allow cancellation of $(1-Z^2)^2/(1-Z\mathfrak b)^4$, which is positive for $0<Z<1$. We obtain \eqref{Eq.A_bfrak_equation}. The inequality $A(a)>0$ follows from $k>\beta\ell$ and $\beta>1$.

	At any $Z\in(0,1)\setminus\{Z_1\}$ where $\mathfrak b'(Z)=0$, the right-hand side of \eqref{Eq.A_bfrak_equation} vanishes. Differentiating \eqref{Eq.A_bfrak_equation} at $Z$, all terms containing $\mathfrak b'$ disappear, which gives \eqref{Eq.A_bfrak_critical}. By Lemma \ref{Lem.sonic_geometry}, $Z(1-Z^2)(1-\ell\mathfrak b^2)$ is positive on $(0,Z_1)$ and negative on $(Z_1,1)$, whereas $(1-\mathfrak b^2)A(\mathfrak b)>0$ on $(0,1)$.
\end{proof}

We use the scalar identity to compare the velocity with the outgoing characteristic speed. The strict inequality below is the sign needed in the normalized expansion estimate.

\begin{lemma}\label{Lem.velocity_expansion}
	Under Assumption \ref{Ass.background}, the function $\mathfrak b(Z):=(Z-v(Z))/(1-Zv(Z))$ belongs to $C^\omega_{\rm o}([-1,1])$ and satisfies
	\[\mathfrak b(0)=0,\qquad \mathfrak b(1)=1,\qquad \mathfrak b'(Z)>0\quad(0\leq Z\leq1).\]
\end{lemma}

\begin{proof}
	By Lemma \ref{Lem.scalar_expansion_identity}, only positivity of $\mathfrak b'$ away from $0,Z_1,1$ remains to be proved.
	Suppose that $\mathfrak b'$ had a zero in $(0,Z_1)$. Since $\mathfrak b'(0)>0$, there is a first zero $Z_c$ of $\mathfrak b'$ in $(0,Z_1)$. For $h>0$ with $Z_c-h>0$, one has $\mathfrak b'(Z_c-h)>0=\mathfrak b'(Z_c)$, so the left derivative of $\mathfrak b'$ at $Z_c$ is nonpositive. The inequality contradicts $\mathfrak b''(Z_c)>0$ from \eqref{Eq.A_bfrak_critical}. Thus $\mathfrak b'>0$ on $[0,Z_1]$.

	If $\mathfrak b'$ had a zero in $(Z_1,1)$, let $Z_c$ be the first zero of $\mathfrak b'$ after $Z_1$. Equation \eqref{Eq.A_bfrak_critical} gives $\mathfrak b''(Z_c)<0$, so $\mathfrak b'<0$ immediately to the right of $Z_c$. Since $\mathfrak b'(1)>0$, the connected component of $\{Z:\mathfrak b'(Z)<0\}$ starting to the right of $Z_c$ has a right endpoint $Z_d\in(Z_c,1)$ with $\mathfrak b'(Z_d)=0$. The derivative of $\mathfrak b'$ from the left at $Z_d$ is nonnegative, contradicting \eqref{Eq.A_bfrak_critical} again. We have proved $\mathfrak b'(Z)>0$ on $[0,1]$.
\end{proof}

The velocity inequality implies positivity of both radial eigenvalues of $\mathrm D_x\bvec$. Compactness of the real interval then gives a single positive lower bound.

\begin{lemma}\label{Lem.qualitative_expansion}
	Let $(\beta,v)$ satisfy Assumption \ref{Ass.background}, and define $X$ and $b$ by \eqref{Eq.x_b}. Then
	\[X\in C^\omega_{\rm o}([-1,1]),\qquad X'(Z)>0\quad(-1\leq Z\leq1).\]
	With $X_1=X(1)$, the inverse $Z=Z(X)$ and the radial coefficient $b$ satisfy
	\begin{align*}
		Z,b\in C^\omega_{\rm o}([-X_1,X_1]),\quad b(0)=0,\quad b(X_1)=1,\quad b'(X)>0\quad(0\leq X\leq X_1).
	\end{align*}
	Consequently, $\gamma_v:=\min_{[0,X_1]}b'>0$ and
	\[\sym\mathrm D_x\bvec(x)\geq\gamma_vI_d\qquad(x\in\overline\Omega).\]
\end{lemma}

\begin{proof}
	Let $\mathfrak b$ be the function in Lemma \ref{Lem.velocity_expansion}. In particular, $\mathfrak b$ is analytic and $\mathfrak b'>0$ on $[0,1]$.
	The definition of $X$ gives
	\[X'(Z)=(\beta-1)\frac{1-v^2}{(1-Zv)^2}>0.\]
	Since $v\in C^\omega_{\rm o}([-1,1])$ and $1-Zv(Z)>0$ there, the right-hand side belongs to $C^\omega_{\rm e}([-1,1])$. Integration and $X(0)=0$ give $X\in C^\omega_{\rm o}([-1,1])$. By continuity on the compact interval, choose $a_X>1$ such that $X$ extends to $C^\omega_{\rm o}((-a_X,a_X))$ and $X'(Z)>0$ on $(-a_X,a_X)$. The inverse function theorem gives
	\[Z\in C^\omega\bigl((X(-a_X),X(a_X))\bigr),\qquad Z(-X)=-Z(X).\]
	Consequently $Z\in C^\omega_{\rm o}([-X_1,X_1])$, and composition yields
	\[b=\mathfrak b\circ Z\in C^\omega_{\rm o}([-X_1,X_1]),\qquad b'(X(Z))=\frac{\mathfrak b'(Z)}{X'(Z)}>0,\qquad b'(0)=\frac{k+1-\beta\ell}{(k+1)(\beta-1)}=c_0.\]
	Thus \eqref{Eq.A_strict_expansion} holds. Compactness gives $\gamma_v>0$, and $b(0)=0$ implies
	\[\frac{b(X)}X=\int_0^1b'(sX)\,\mathrm ds\geq\gamma_v \qquad(0<X\leq X_1).\]
	For $x\ne0$, we define $e_X:=x/X$, $X:=|x|$. Then
	\[\mathrm D_x\bvec=b'(X)e_X\otimes e_X +\frac{b(X)}X(I_d-e_X\otimes e_X).\]
	At $x=0$, $\mathrm D_x\mathbf b$ equals $c_0I_d$. Therefore, for $x\ne0$ and every $\zeta\in\R^d$,
	\[\zeta^T\sym\mathrm D_x\bvec(x)\zeta =b'(X)|\zeta\cdot e_X|^2+\frac{b(X)}X\bigl(|\zeta|^2-|\zeta\cdot e_X|^2\bigr) \geq\gamma_v|\zeta|^2.\qedhere\]
\end{proof}

We have proved expansion of the vector field. We next check the apparent singularities at the origin in the phase, amplitude, density and potential.

\begin{lemma}\label{Lem.qualitative_scalars}
	The scalar profiles satisfy
	\[\widehat\rho_0,\widehat\phi_0\in C^\omega_{\rm e}([-1,1]),\qquad j_0,j_0^{-1},V_{\rm r}\in C^\omega_{\rm e}([-X_1,X_1]).\]
	The profiles are defined by \eqref{Eq.profiles} and \eqref{Eq.normalized_coefficients}, with $j_0(X):=J_0(Xe_1)$ and $V_{\rm r}(X):=V(Xe_1)$. The apparent singularities at zero are removable, and the identity \eqref{Eq.density_identity} holds.
\end{lemma}

\begin{proof}
	The smooth evenness of $\widehat\rho_0,\widehat\phi_0$ is proved in \cite[Lemma~A.6]{SWZ-NLW}; the proof uses only the profile on the interval under consideration. Proposition \ref{Prop.automatic_analyticity} strengthens $v$ to an analytic odd function. Since $1-Zv(Z)>0$ and $1-v(Z)^2>0$ on $[-1,1]$, the integrals, exponentials and positive roots in \eqref{Eq.profiles} are analytic. Thus
\[\widehat\rho_0,\widehat\phi_0\in C^\omega_{\rm e}([-1,1]),\qquad \widehat\rho_0\circ Z,\widehat\phi_0\circ Z\in C^\omega_{\rm e}([-X_1,X_1]).\]
	To examine the density at zero, we define $c_3:=b'''(0)/6$. Then $b(X)=c_0X+c_3X^3+\Order(X^5)$, and
	\begin{equation}\label{Eq.A_density_cancellation}
		\frac{dc_0-b'(X)-kb(X)/X}{b(X)} =-\frac{(k+3)c_3}{c_0}X+\Order(X^3).
	\end{equation}
	More precisely, $b(X)/X\in C^\omega_{\rm e}([-X_1,X_1])$, with value $c_0>0$ at zero and no zero on $[-X_1,X_1]$. The numerator in \eqref{Eq.A_density_cancellation} belongs to $C^\omega_{\rm e}([-X_1,X_1])$ and vanishes to order at least two at zero. Hence the function
	\[\gamma_0(X):= \begin{cases} \displaystyle\frac{dc_0-b'(X)-kb(X)/X}{b(X)},&X\ne0,\\[4pt] 0,&X=0, \end{cases} \qquad \gamma_0\in C^\omega_{\rm o}([-X_1,X_1])\]
	is well defined. Formula \eqref{Eq.normalized_coefficients} becomes
	\[j_0(X)=\exp\left(\int_0^X\gamma_0(a)\,\mathrm da\right),\qquad j_0,j_0^{-1}\in C^\omega_{\rm e}([-X_1,X_1]),\qquad j_0>0.\]
	Differentiation gives, for $X\ne0$,
	\[\frac{j_0'}{j_0}=\gamma_0,\qquad \frac1{X^k}\frac{\mathrm d}{\mathrm dX}(X^kj_0b)=dc_0j_0.\]
	The identity $(X^kj_0b)'=dc_0X^kj_0$ extends to $X=0$ and proves \eqref{Eq.density_identity} in Cartesian variables.

	By \eqref{Eq.even_odd_spaces}, there exists $R_\square\in C^\omega([0,1])$ such that $\widehat\rho_0(Z)=R_\square(Z^2)$ for $|Z|\leq1$. Thus
	\[\frac{\widehat\rho_0'(Z)}Z=2R_\square'(Z^2) \in C^\omega_{\rm e}([-1,1]),\]
	with value $\widehat\rho_0''(0)$ at zero. All terms in the expression for $V$ in \eqref{Eq.normalized_coefficients} therefore belong to $C^\omega_{\rm e}([-1,1])$ before composition with $Z(X)$. Hence
	$V_{\rm r}\in C^\omega_{\rm e}([-X_1,X_1])$.
\end{proof}

The scalar profiles are analytic with the parity required at the origin. The even power series of the scalar profiles give Cartesian continuation near zero, and local square-root branches give the continuation away from zero.

\begin{lemma}\label{Lem.qualitative_cartesian}
	Let $(\beta,v)$ satisfy Assumption \ref{Ass.background}, and let $\mathbf c$ be the coefficient map in \eqref{Eq.background_coefficient_map}. There are $R_v,M_v>0$ and $\widetilde{\mathbf c}$ such that
	\[\widetilde{\mathbf c}\in\Hol(\Omega_{R_v}^{\C};\C^{d+5}),\qquad \widetilde{\mathbf c}|_{\overline\Omega}=\mathbf c,\]
	and
	\[\|\nabla^j\bvec\|_{L^\infty(\Omega)}\leq M_v(j-1)!R_v^{1-j}\qquad(j\geq1).\]
	For each scalar component $a$ of $\mathbf c$, there is $C_a<\infty$ such that
	\[\|\nabla^ja\|_{L^\infty(\Omega)}\leq C_a j!R_v^{-j}\qquad(j\geq0).\]
	The constants may depend on $(\beta,v)$.
\end{lemma}

\begin{proof}
	We now pass to Cartesian variables and construct a common complex domain. Set
	\[(a_1,a_2,a_3,a_4,a_5) :=(j_0,j_0^{-1},V_{\rm r},\widehat\rho_0\circ Z,\widehat\phi_0\circ Z).\]
	By \eqref{Eq.even_odd_spaces}, there exist $A_1,\ldots,A_5,B_0\in C^\omega([0,X_1^2])$ such that
	\[a_i(X)=A_i(X^2)\quad(1\leq i\leq5),\qquad b(X)=XB_0(X^2)\quad(|X|\leq X_1).\]
	In particular,
	\[\mathbf c(x)=\bigl(B_0(|x|^2)x,A_1(|x|^2),\ldots,A_5(|x|^2)\bigr)\]
	and $\mathbf c\in C^\omega(\overline\Omega;\R^{d+5})$. Apply the interval continuation argument in the proof of Proposition \ref{Prop.automatic_analyticity} to $A_1,\ldots,A_5,B_0$. There exist $\epsilon>0$ and
	\begin{align*}
		&\qquad\qquad\qquad\qquad\mathcal V_\epsilon:=\{w\in\C:\dist(w,[0,X_1^2])<\epsilon\},\\ 
		&\widetilde A_1,\ldots,\widetilde A_5,\widetilde B_0\in\Hol(\mathcal V_\epsilon),\quad\widetilde A_i|_{[0,X_1^2]}=A_i\quad(1\leq i\leq5),\quad \widetilde B_0|_{[0,X_1^2]}=B_0.
	\end{align*}
	We define
	\[\rho:=\min\left\{1,\frac{\epsilon}{2(2X_1+1)}\right\},\qquad \mathfrak q(\zeta):=\sum_{i=1}^d\zeta_i^2.\]
	For $\zeta\in\Omega_\rho^\C$, take $x\in\overline\Omega$ with $|\zeta-x|<\rho$. Then
	\[|\mathfrak q(\zeta)-|x|^2| \leq2|x||\zeta-x|+|\zeta-x|^2 <(2X_1+\rho)\rho\leq\epsilon/2.\]
	Hence $\mathfrak q(\zeta)\in\mathcal V_\epsilon$, and
	\[\widetilde{\mathbf c}(\zeta) :=\bigl(\widetilde B_0(\mathfrak q(\zeta))\zeta, \widetilde A_1(\mathfrak q(\zeta)),\ldots, \widetilde A_5(\mathfrak q(\zeta))\bigr) .\]
	Thus $\widetilde{\mathbf c}\in\Hol(\Omega_\rho^\C;\C^{d+5})$ is an extension of $\mathbf c$.

	Let $\widetilde\bvec$ be the first $d$ components of $\widetilde{\mathbf c}$. Define
	\begin{align*}
		r_0:=\frac{\rho}{4\sqrt d},\quad R_v:=\frac{r_0}{\sqrt d}=\frac{\rho}{4d},\quad M_v:=\max\left\{1, \sup_{\dist(\zeta,\overline\Omega)\leq\rho/2} \norm{\nabla\widetilde\bvec(\zeta)}_{\mathrm{HS}}\right\}<\infty.
	\end{align*}
	The complex polydisc of component radius $r_0$ centered at any $x\in\overline\Omega$ is contained in $\Omega_{\rho/2}^\C$. The Banach-valued Cauchy estimate applied to $\nabla\widetilde\bvec$ bounds each ordered derivative of length $j-1$ in Hilbert--Schmidt norm by $M_v(j-1)!r_0^{-(j-1)}$. There are $d^{j-1}$ such ordered derivatives, so
	\[\norm{\na^j\bvec(x)} \leq d^{(j-1)/2}M_v(j-1)!r_0^{-(j-1)} =M_v(j-1)!R_v^{1-j}\qquad(j\geq1).\]
	For a scalar component $a$ of $\mathbf c$, let $\widetilde a$ be the corresponding component of $\widetilde{\mathbf c}$ and define
	\[C_a:=\sup_{\dist(\zeta,\overline\Omega)\leq\rho/2}|\widetilde a(\zeta)|<\infty.\]
	The same argument gives
	$\norm{\na^ja(x)}\leq d^{j/2}C_a j!r_0^{-j}=C_a j!R_v^{-j}$ for $j\geq0$. Restricting $\widetilde{\mathbf c}$ to $\Omega_{R_v}^\C$ proves \eqref{Eq.A_Cartesian_tube}.
\end{proof}

The scalar expansion bound and the Cartesian continuation are now proved. We choose one common complex neighborhood for the finite list of coefficients and obtain the derivative bounds by Cauchy's formula.

\begin{proof}[Proof of Proposition \ref{Prop.automatic_geometry}]
	Combine Lemmas \ref{Lem.qualitative_expansion}--\ref{Lem.qualitative_cartesian}.
\end{proof}

\subsection{Proof of Lemma \ref{Lem.normalized_equation}}\label{App.normalization}

Proposition \ref{Prop.automatic_geometry} gives the regularity of the coefficients in the normalized equation. We now prove Lemma \ref{Lem.normalized_equation} by direct computation. We first expand the weighted wave operator acting on $z^\sharp$ in the variables $(\sigma,X)$, and then compute the coefficient of each derivative of $z$.

\begin{proof}[Proof of Lemma \ref{Lem.normalized_equation}]
	Let $I\subset(0,\infty)$ be an open interval and let $z\in C^\infty(I\times\overline\Omega;\R)$ be radial in $x$. We use $z^\sharp$ as defined in Lemma \ref{Lem.normalized_equation}. For $e_1:=(1,0,\ldots,0)\in\R^d$, define $z_{\rm r}(\sigma,X):=z(\sigma,Xe_1)$ for $(\sigma,X)\in I\times[0,X_1]$. We first work at points with $t<T$, $0<r<T-t$, and $\sigma(t,r)\in I$, and write $Z=r/(T-t)$ and $X=X(Z)$. When taking $t$- and $r$-derivatives, we identify $z^\sharp$ with its radial profile $(t,r)\mapsto z^\sharp(t,re_1)=z_{\rm r}(\sigma(t,r),X(Z))$.
	
	Since $\Box=-\partial_t^2+\Delta_y$ and $k=d-1$, the product rule gives
	\begin{align*}
		-\rho_0^{-p-1}\partial^\mu(\rho_0^2\partial_\mu z^\sharp)&=\rho_0^{-p-1}\left[\partial_t(\rho_0^2\partial_tz^\sharp)-r^{-k}\partial_r(r^k\rho_0^2\partial_rz^\sharp)\right]\\
		&=\rho_0^{1-p}\left[\partial_t^2z^\sharp-\partial_r^2z^\sharp-\frac{k}{r}\partial_rz^\sharp+2\frac{\partial_t\rho_0}{\rho_0}\partial_tz^\sharp-2\frac{\partial_r\rho_0}{\rho_0}\partial_rz^\sharp\right].
	\end{align*}
	The chain rule gives $\partial_tz^\sharp=(\partial_\sigma z_{\rm r})(\partial_t\sigma)+(\partial_Xz_{\rm r})(\partial_tX)$ and $\partial_rz^\sharp=(\partial_\sigma z_{\rm r})(\partial_r\sigma)+(\partial_Xz_{\rm r})(\partial_rX)$. Differentiating once more and collecting the coefficients, we obtain
	\begin{equation}\label{Eq.normalization_chain_rule}
		\begin{aligned}
			&-\rho_0^{-p-1}\partial^\mu(\rho_0^2\partial_\mu z^\sharp)=\rho_0^{1-p}((\partial_t\sigma)^2-(\partial_r\sigma)^2)\partial_\sigma^2z_{\rm r}\\
			&\qquad\qquad+2\rho_0^{1-p}((\partial_t\sigma)(\partial_tX)-(\partial_r\sigma)(\partial_rX))\partial_\sigma\partial_Xz_{\rm r}+\rho_0^{1-p}((\partial_tX)^2-(\partial_rX)^2)\partial_X^2z_{\rm r}\\
			&\qquad\qquad+\rho_0^{1-p}\left[\partial_t^2\sigma-\partial_r^2\sigma-\frac{k}{r}\partial_r\sigma+2\frac{\partial_t\rho_0}{\rho_0}\partial_t\sigma-2\frac{\partial_r\rho_0}{\rho_0}\partial_r\sigma\right]\partial_\sigma z_{\rm r}\\
			&\qquad\qquad+\rho_0^{1-p}\left[\partial_t^2X-\partial_r^2X-\frac{k}{r}\partial_rX+2\frac{\partial_t\rho_0}{\rho_0}\partial_tX-2\frac{\partial_r\rho_0}{\rho_0}\partial_rX\right]\partial_Xz_{\rm r}.
		\end{aligned}
	\end{equation}
	Here all derivatives of $z_{\rm r}$ are evaluated at $(\sigma(t,r),X(Z))$. The symbols $\partial_tX,\partial_rX,\partial_t^2X,\partial_r^2X$ denote derivatives of the composition $(t,r)\mapsto X(r/(T-t))$, whereas $X'(Z),X''(Z)$ below denote derivatives of $Z\mapsto X(Z)$. Primes on $\widehat\phi_0,\widehat\rho_0$ denote $Z$-derivatives.
	
	We now compute the five coefficients in \eqref{Eq.normalization_chain_rule}. Since $\sigma=\phi_0$, the first equation in \eqref{Eq.leading_system} gives $\rho_0^{1-p}((\partial_t\sigma)^2-(\partial_r\sigma)^2)=1$. The second equation in \eqref{Eq.leading_system} gives
	\begin{align*}
		\rho_0^{1-p}\left[\partial_t^2\sigma-\partial_r^2\sigma-\frac{k}{r}\partial_r\sigma+2\frac{\partial_t\rho_0}{\rho_0}\partial_t\sigma-2\frac{\partial_r\rho_0}{\rho_0}\partial_r\sigma\right] =-\rho_0^{-p-1}\partial^\mu(\rho_0^2\partial_\mu\sigma)=0.
	\end{align*}
	Thus the coefficient of $\partial_\sigma z_{\rm r}$ vanishes.
	
	To compute the other second-order coefficients, we use \eqref{Eq.phase_derivatives} and $Z=r/(T-t)$ to obtain
	\begin{equation}\label{Eq.coordinate_derivatives}
		\partial_t\sigma=\frac{(\beta-1)\sigma}{(T-t)(1-Zv(Z))},\quad \partial_r\sigma=v(Z)\partial_t\sigma,\quad \partial_tX=\frac{ZX'(Z)}{T-t},\quad \partial_rX=\frac{X'(Z)}{T-t}.
	\end{equation}
	The definitions in \eqref{Eq.profiles} and \eqref{Eq.x_b} give
	\begin{equation}\label{Eq.normalization_direct_profiles}
		\begin{aligned}
			X'(Z)&=(\beta-1)\frac{1-v(Z)^2}{(1-Zv(Z))^2},\qquad \widehat\rho_0(Z)^{p-1}=(\beta-1)\widehat\phi_0(Z)^2X'(Z),\\
			\frac{\widehat\phi_0'(Z)}{\widehat\phi_0(Z)}&=\frac{(\beta-1)v(Z)}{1-Zv(Z)},\qquad (1-Z^2)X'(Z)=(\beta-1)(1-b(X)^2).
		\end{aligned}
	\end{equation}
	In particular, \eqref{Eq.selfsimilar} implies 
	\[\sigma^2\rho_0^{1-p}=(T-t)^2\frac{\widehat\phi_0(Z)^2}{\widehat\rho_0(Z)^{p-1}}=\frac{(T-t)^2}{(\beta-1)X'(Z)}.\]
	Hence
	\begin{align*}
		\rho_0^{1-p}((\partial_t\sigma)(\partial_tX)-(\partial_r\sigma)(\partial_rX))&=\frac{(\beta-1)\sigma\rho_0^{1-p}X'(Z)}{(T-t)^2}\frac{Z-v(Z)}{1-Zv(Z)}=\frac{b(X)}{\sigma},\\
		\rho_0^{1-p}((\partial_tX)^2-(\partial_rX)^2)&=\frac{\rho_0^{1-p}(Z^2-1)(X'(Z))^2}{(T-t)^2}=\frac{b(X)^2-1}{\sigma^2}.
	\end{align*}
	
	It remains to compute the coefficient of $\partial_Xz_{\rm r}$. Differentiating \eqref{Eq.coordinate_derivatives} and \eqref{Eq.selfsimilar}, we obtain
	\begin{align*}
		\partial_t^2X&=\frac{Z^2X''(Z)+2ZX'(Z)}{(T-t)^2},\qquad \partial_r^2X=\frac{X''(Z)}{(T-t)^2},\\
		\frac{\partial_t\rho_0}{\rho_0}&=\frac1{T-t}\left(\frac{2\beta}{p-1}+Z\frac{\widehat\rho_0'(Z)}{\widehat\rho_0(Z)}\right),\qquad \frac{\partial_r\rho_0}{\rho_0}=\frac1{T-t}\frac{\widehat\rho_0'(Z)}{\widehat\rho_0(Z)}.
	\end{align*}
	Using $r=Z(T-t)$ and the identity for $\sigma^2\rho_0^{1-p}$ above, we find
	\begin{align*}
		&\rho_0^{1-p}\left[\partial_t^2X-\partial_r^2X-\frac{k}{r}\partial_rX+2\frac{\partial_t\rho_0}{\rho_0}\partial_tX-2\frac{\partial_r\rho_0}{\rho_0}\partial_rX\right]\\
		=&\,\frac{\rho_0^{1-p}X'(Z)}{(T-t)^2}\left[(Z^2-1)\left(\frac{X''(Z)}{X'(Z)}+2\frac{\widehat\rho_0'(Z)}{\widehat\rho_0(Z)}\right)+2\left(1+\frac{2\beta}{p-1}\right)Z-\frac{k}{Z}\right]\\
		=&\,\frac1{(\beta-1)\sigma^2}\left[(Z^2-1)\left(\frac{X''(Z)}{X'(Z)}+2\frac{\widehat\rho_0'(Z)}{\widehat\rho_0(Z)}\right)+2\left(1+\frac{2\beta}{p-1}\right)Z-\frac{k}{Z}\right].
	\end{align*}
	The expression in brackets depends only on $Z$. Since $X:(0,1)\to(0,X_1)$ is a diffeomorphism, we define $a:(0,X_1)\to\R$ by
	\begin{equation}\label{Eq.normalization_direct_coefficient}
		a(X(Z)):=\frac1{\beta-1}\biggl[(Z^2-1)\left(\frac{X''(Z)}{X'(Z)}+2\frac{\widehat\rho_0'(Z)}{\widehat\rho_0(Z)}\right)+2\left(1+\frac{2\beta}{p-1}\right)Z-\frac{k}{Z}\biggr]
	\end{equation}
	for $Z\in(0,1)$. Thus the coefficient of $\partial_Xz_{\rm r}$ is $a(X)/\sigma^2$. Substituting the five coefficients into \eqref{Eq.normalization_chain_rule}, we obtain
	\begin{equation}\label{Eq.normalization_direct_operator}
		-\rho_0^{-p-1}\partial^\mu(\rho_0^2\partial_\mu z^\sharp)=\partial_\sigma^2z_{\rm r}+\frac{2b(X)}{\sigma}\partial_\sigma\partial_Xz_{\rm r}+\frac{b(X)^2-1}{\sigma^2}\partial_X^2z_{\rm r}+\frac{a(X)}{\sigma^2}\partial_Xz_{\rm r}.
	\end{equation}
	
	We next express $a$ in terms of $b$ and $J_0$. Primes on $b$ denote $X$-derivatives. Differentiating the last identity in \eqref{Eq.normalization_direct_profiles} with respect to $Z$, and using the other identities there, we obtain
	\begin{align*}
		(Z^2-1)\frac{X''(Z)}{X'(Z)}&=2(\beta-1)b(X)b'(X)-2Z,\\
		(Z^2-1)\frac{\widehat\phi_0'(Z)}{\widehat\phi_0(Z)}&=(\beta-1)(b(X)-Z),\\
		2\frac{\widehat\rho_0'(Z)}{\widehat\rho_0(Z)}&=\frac{2}{p-1}\left(2\frac{\widehat\phi_0'(Z)}{\widehat\phi_0(Z)}+\frac{X''(Z)}{X'(Z)}\right).
	\end{align*}
	Substituting these identities into \eqref{Eq.normalization_direct_coefficient}, and recalling that $\ell=1+4/(p-1)$, we find
	\[a(X)=(\ell+1)b(X)b'(X)+(\ell-1)b(X)-\frac{k}{(\beta-1)Z}.\]
	On the other hand, \eqref{Eq.A_bfrak_equation}, $\mathfrak b(Z)=b(X(Z))$ and $(1-Z^2)X'(Z)=(\beta-1)(1-b(X)^2)$ imply
	\[(1-\ell b(X)^2)b'(X)=dc_0+\ell b(X)^2-\frac{kb(X)}{(\beta-1)Z},\]
	where we have used $(\beta-1)dc_0=d-\beta\ell$. The factor $1-b(X)^2$ is positive for $0<X<X_1$ by Proposition \ref{Prop.automatic_geometry}. No division by $1-\ell b(X)^2$ is used, so this identity holds also at $X=X(Z_1)$. Let $j_0(X)=J_0(Xe_1)$ be the radial profile defined in Proposition \ref{Prop.automatic_geometry}, and let $j_0'$ denote its $X$-derivative. Since $b(X)>0$ for $0<X<X_1$, it follows that
	\begin{align*}
		a(X)&=b(X)b'(X)-b(X)-\frac{dc_0-b'(X)}{b(X)}\\
		&=b(X)b'(X)-b(X)-\frac{k}{X}-\frac{j_0'(X)}{j_0(X)},
	\end{align*}
	where the last equality follows from \eqref{Eq.normalized_coefficients}. The same definition gives
	\[b'(X)+\left(\frac{k}{X}+\frac{j_0'(X)}{j_0(X)}\right)b(X)=dc_0,\qquad 0<X<X_1.\]
	Together with the Cartesian smoothness from Proposition \ref{Prop.automatic_geometry} and $b(X_1)=1$, this also proves \eqref{Eq.density_identity} on $\overline\Omega$.
	
	For $(\sigma,X)\in I\times(0,X_1)$, the operators in \eqref{Eq.normalized_operators} satisfy
	\begin{align*}
		(D^2z)(\sigma,Xe_1)&=\partial_\sigma^2z_{\rm r}+\frac{2b(X)}{\sigma}\partial_\sigma\partial_Xz_{\rm r}+\frac{b(X)^2}{\sigma^2}\partial_X^2z_{\rm r}+\frac{b(X)b'(X)-b(X)}{\sigma^2}\partial_Xz_{\rm r},\\
		(A_0z)(\sigma,Xe_1)&=-\partial_X^2z_{\rm r}-\left(\frac{k}{X}+\frac{j_0'(X)}{j_0(X)}\right)\partial_Xz_{\rm r},
	\end{align*}
	where the derivatives of $z_{\rm r}$ on the right are evaluated at $(\sigma,X)$. Comparing these identities with \eqref{Eq.normalization_direct_operator} proves \eqref{Eq.physical_P} for $0<r<T-t$.
	
	It remains to compute the potential term. By \eqref{Eq.selfsimilar},
	\begin{align*}
		\Box\rho_0=(T-t)^{-\frac{2\beta}{p-1}-2}\biggl[&(1-Z^2)\widehat\rho_0''(Z)+\left(\frac{k}{Z}-2\left(1+\frac{2\beta}{p-1}\right)Z\right)\widehat\rho_0'(Z)\\
		&-\frac{2\beta}{p-1}\left(1+\frac{2\beta}{p-1}\right)\widehat\rho_0(Z)\biggr].
	\end{align*}
	Thus \eqref{Eq.selfsimilar} and the definition of $V$ in \eqref{Eq.normalized_coefficients} give 
	\[\rho_0^{-p}\Box\rho_0=(T-t)^{2\beta-2}\widehat\phi_0(Z)^{-2}V(x(t,y))=\sigma^{-2}V(x(t,y)).\]
	The product rule and \eqref{Eq.physical_P} now imply
	\begin{align*}
		\Box(\rho_0z^\sharp)-(\rho_0z^\sharp)^p&=z^\sharp\Box\rho_0+\rho_0^{-1}\partial^\mu(\rho_0^2\partial_\mu z^\sharp)-\rho_0^p(z^\sharp)^p\\
		&=-\rho_0^p\bigl(\cP z-\sigma^{-2}Vz+z^p\bigr)^\sharp=-\rho_0^p\bigl(\cN(z)\bigr)^\sharp.
	\end{align*}
	This proves the first identity in \eqref{Eq.exact_reconstruction} for $0<r<T-t$. Since $B$ is independent of $\sigma$, we have $D^2z=\partial_\sigma^2z+2\sigma^{-1}B\partial_\sigma z+\sigma^{-2}(B^2-B)z$, which gives the second identity in \eqref{Eq.exact_reconstruction}.
	
	We now include the origin and the lateral boundary. By Proposition \ref{Prop.automatic_geometry}, $\widehat\phi_0,\widehat\rho_0$ and $X(Z)/Z$ belong to $C^\omega_{\rm e}([-1,1])$, where the last quotient has value $\beta-1$ at zero. Thus
	\[x(t,y)=\frac1{T-t}\frac{X(Z)}Z\,y,\qquad \sigma(t,y)=(T-t)^{1-\beta}\widehat\phi_0(Z),\qquad Z=\frac{|y|}{T-t},\]
	and \eqref{Eq.even_odd_spaces} show that $x,\sigma,\rho_0$ are smooth in $(t,y)$ at $y=0$ and up to $|y|=T-t$ for $t<T$. Since $z\in C^\infty(I\times\overline\Omega;\R)$ and the coefficients $\bvec,J_0,J_0^{-1},V$ are smooth on $\overline\Omega$, both sides of \eqref{Eq.physical_P} and \eqref{Eq.exact_reconstruction} are smooth there. These identities therefore extend to $r=0$ and $r=T-t$ by continuity.
	
	Finally, for $(\sigma,x)\in I\times\overline\Omega$, the inverse coordinates are
	\begin{align*}
		t(\sigma,x)=T-\left(\frac{\widehat\phi_0(Z(|x|))}{\sigma}\right)^{\frac1{\beta-1}},\quad y(\sigma,x)=\left(\frac{\widehat\phi_0(Z(|x|))}{\sigma}\right)^{\frac1{\beta-1}}\frac{Z(|x|)}{|x|}\,x.
	\end{align*}
	Here $Z(|x|)/|x|$ has value $(\beta-1)^{-1}$ at $x=0$. The odd analyticity of $Z$ and the even analyticity and positivity of $\widehat\phi_0$ imply that these formulas are smooth at $x=0$ and up to $|x|=X_1$. Direct substitution shows that they are inverse to the coordinates in \eqref{Eq.sigma_z} and take $I\times\overline\Omega$ onto $\{(t,y):t<T,\ |y|\leq T-t,\ \sigma(t,y)\in I\}$. Since $\rho_0>0$, the first identity in \eqref{Eq.exact_reconstruction} proves the equivalence of \eqref{Eq.wave_eq} and \eqref{Eq.N}. The proof is complete.
\end{proof}

\section{Properties of the analytic function spaces}\label{App.analytic_spaces}

We collect here the properties of $\mathcal B_\varrho$ used in the evolution and approximation arguments. Throughout this appendix, $\Omega=B_{X_1}(0)\subset\R^d$, the indices $s_d,s_*$ are given by \eqref{Eq.dimension_indices}, and $E$ is a finite-dimensional real or complex Euclidean space. We use \eqref{Eq.analytic_norm} and the full ordered derivative convention. For holomorphic extension, a real Euclidean space is complexified with the same norm.

\begin{lemma}\label{Lem.analytic_derivatives}
	For $0<\varrho<\varrho'$, every integer $q\geq0$ and every $h\in\mathcal B_{\varrho'}(\Omega;E)$,
	\[\|\nabla^qh\|_{\mathcal B_\varrho}\leq\frac{q!}{(\varrho'-\varrho)^q}\|h\|_{\mathcal B_{\varrho'}}.\]
\end{lemma}

\begin{proof}
	Let $0<\varrho<\varrho'$. For $N, q\in\Z_{\geq 0}$, we define $M=N+q$ and $t=\varrho/\varrho'\in(0,1)$. The binomial theorem gives
	\[\binom Mq t^{M-q}(1-t)^q\leq\sum_{j=0}^M\binom Mj t^{M-j}(1-t)^j=1.\]
	Consequently,
	\begin{align*}
		\|\nabla^qh\|_{\mathcal B_\varrho}^2 &=\sum_{M= q}^\infty(1+M-q)^{2s_*}\frac{\varrho^{2(M-q)}}{((M-q)!)^2}\|\nabla^Mh\|_{L^2}^2\\
		&\leq\frac{(q!)^2}{(\varrho'-\varrho)^{2q}}\sum_{M=q}^\infty (1+M)^{2s_*}\left[\binom{M}{q}t^{M-q}(1-t)^q\right]^2\frac{(\varrho')^{2M}}{(M!)^2}\|\nabla^Mh\|_{L^2}^2\\
		&\leq \frac{(q!)^2}{(\varrho'-\varrho)^{2q}}\sum_{M\geq q}(1+M)^{2s_*}\frac{(\varrho')^{2M}}{(M!)^2}\|\nabla^Mh\|_{L^2}^2.
	\end{align*}
	Taking square roots gives \eqref{Eq.analytic_radius_recovery}.
\end{proof}

\begin{lemma}[Multiplier bounds at a smaller analytic radius]\label{Lem.analytic_to_multiplier}
	There is a constant $C>0$, depending only on $\Omega$ and $d$, such that, for every $0<r<r'$, every integer $c\geq0$ and every $u\in\mathcal B_{r'}(\Omega;E)$,
	\begin{equation}\label{Eq.analytic_to_multiplier}
		\mathfrak M_r(\nabla^cu)\leq Cr'\left(\sum_{q=0}^{s_d}\frac{(c+q)!}{(r'-r)^{c+q+1}}\right)\|u\|_{\mathcal B_{r'}}.
	\end{equation}
	The multiplier norm uses the Euclidean norm on $E$ and on the ordered derivative tensors. The constant $C$ is independent of $r,r',c,u$ and $\dim E$.
\end{lemma}

\begin{proof}
	Let $0<r<r'$, let $c\in\Z_{\geq0}$, and let $u\in\mathcal B_{r'}(\Omega;E)$. Since $s_d>d/2$, Sobolev embedding gives
	\[\|\nabla^{N+c}u\|_{L^\infty}\leq C\sum_{q=0}^{s_d}\|\nabla^{N+c+q}u\|_{L^2},\qquad N\in\Z_{\geq0}.\]
	The constant is independent of $N,c$ and $\dim E$: apply the scalar Sobolev inequality to each component of the ordered derivative tensor and then sum the squares. By \eqref{Eq.analytic_norm},
	\[\|\nabla^{N+c+q}u\|_{L^2}\leq\frac{(N+c+q)!}{(1+N+c+q)^{s_*}(r')^{N+c+q}}\|u\|_{\mathcal B_{r'}}.\]
	Therefore,
	\begin{align*}
		\mathfrak M_r(\nabla^cu)&\leq C\sum_{q=0}^{s_d}\sum_{N=0}^\infty(1+N)^{s_*}\frac{r^N}{N!}\|\nabla^{N+c+q}u\|_{L^2}\\
		&\leq C\|u\|_{\mathcal B_{r'}}\sum_{q=0}^{s_d}(r')^{-c-q}\sum_{N=0}^\infty\left(\frac r{r'}\right)^N\frac{(1+N)^{s_*}(N+c+q)!}{N!(1+N+c+q)^{s_*}}\\
		&\leq C\|u\|_{\mathcal B_{r'}}\sum_{q=0}^{s_d}(r')^{-c-q}\sum_{N=0}^\infty\left(\frac r{r'}\right)^N\frac{(N+c+q)!}{N!}.
	\end{align*}
	The binomial series gives
	\begin{align*}
		\sum_{N=0}^\infty\frac{(N+c+q)!}{N!}\left(\frac r{r'}\right)^N=(c+q)!\sum_{N=0}^\infty\binom{N+c+q}{c+q}\left(\frac r{r'}\right)^N=\frac{(c+q)!}{(1-r/r')^{c+q+1}}.
	\end{align*}
	Consequently,
	\begin{align*}
		\mathfrak M_r(\nabla^cu)\leq C\|u\|_{\mathcal B_{r'}}\sum_{q=0}^{s_d}\frac{(c+q)!(r')^{-c-q}}{(1-r/r')^{c+q+1}}=Cr'\left(\sum_{q=0}^{s_d}\frac{(c+q)!}{(r'-r)^{c+q+1}}\right)\|u\|_{\mathcal B_{r'}}.
	\end{align*}
	This proves \eqref{Eq.analytic_to_multiplier}.
\end{proof}

\begin{lemma}\label{Lem.analytic_algebra}
	For $\varrho>0$ and scalar functions $u,v\in\mathcal B_\varrho$,
	\[\|uv\|_{\mathcal B_\varrho}\leq C_\varrho\|u\|_{\mathcal B_\varrho}\|v\|_{\mathcal B_\varrho}.\]
	The constant $C_\varrho$ is locally uniform for $\varrho\in(0,\infty)$. If $a:\Omega\to\mathcal L(E)$ is a smooth coefficient with $\mathfrak M_\varrho(a)<\infty$, then, for every $U\in\mathcal B_\varrho(\Omega;E)$,
	\[\|aU\|_{\mathcal B_\varrho}\leq\mathfrak M_\varrho(a)\|U\|_{\mathcal B_\varrho}.\]
	A scalar coefficient is identified with the corresponding multiple of the identity on $E$. Moreover, for a smooth scalar function $a$ and a smooth $E$-valued function $b$ with finite multiplier norms,
	\begin{equation}\label{Eq.multiplier_algebra}
		\mathfrak M_\varrho(ab)\leq\mathfrak M_\varrho(a)\mathfrak M_\varrho(b).
	\end{equation}
	Inequality \eqref{Eq.multiplier_algebra} also holds for smooth matrix-valued functions $a,b$, using the Euclidean operator norm and the ordered matrix derivative convention.
\end{lemma}

\begin{proof}
	For the product estimate, fix $u,v\in\mathcal B_\varrho$. In the Leibniz expansion, we place the factor with fewer derivatives in $L^\infty$ and the other factor in $L^2$. Thus, for every $N\geq0$,
	\begin{align*}
		\|\nabla^N(uv)\|_{L^2} \leq\sum_{0\leq r\leq N/2}\binom Nr\|\nabla^ru\|_{L^\infty}\|\nabla^{N-r}v\|_{L^2}+\sum_{N/2<r\leq N}\binom Nr\|\nabla^ru\|_{L^2}\|\nabla^{N-r}v\|_{L^\infty}.
	\end{align*}
	Define the nonnegative sequences
	\begin{align*}
		u_N&:=(1+N)^{s_*}\frac{\varrho^N}{N!}\|\nabla^Nu\|_{L^2},\qquad \widetilde u_N:=\frac{\varrho^N}{N!}\|\nabla^Nu\|_{L^\infty},\\
		v_N&:=(1+N)^{s_*}\frac{\varrho^N}{N!}\|\nabla^Nv\|_{L^2},\qquad \widetilde v_N:=\frac{\varrho^N}{N!}\|\nabla^Nv\|_{L^\infty}.
	\end{align*}
	Then $\|(u_N)_{N\geq0}\|_{\ell^2}=\|u\|_{\mathcal B_\varrho}$ and $\|(v_N)_{N\geq0}\|_{\ell^2}=\|v\|_{\mathcal B_\varrho}$. For $0\leq r\leq N/2$, we have $1+N\leq2(1+N-r)$, and hence
	\begin{align*}
		\frac{(1+N)^{s_*}\varrho^N}{N!}\binom Nr\|\nabla^ru\|_{L^\infty}\|\nabla^{N-r}v\|_{L^2}=\left(\frac{1+N}{1+N-r}\right)^{s_*}\widetilde u_rv_{N-r} \leq2^{s_*}\widetilde u_rv_{N-r}.
	\end{align*}
	For $r>N/2$, the same argument applies with $u$ and $v$ interchanged and $r$ replaced by $N-r$. As a consequence, we have
	\[\frac{(1+N)^{s_*}\varrho^N}{N!}\|\nabla^N(uv)\|_{L^2}\leq2^{s_*}\sum_{r=0}^N\bigl(\widetilde u_rv_{N-r}+\widetilde v_ru_{N-r}\bigr).\]
	
	We next show that $(\widetilde u_r)_{r\geq0}$ and $(\widetilde v_r)_{r\geq0}$ belong to $\ell^1$. Since $s_d>d/2$, Sobolev embedding applied to the Euclidean tensor $\nabla^ru$ gives
	\[\|\nabla^ru\|_{L^\infty}\leq C\sum_{q=0}^{s_d}\|\nabla^{r+q}u\|_{L^2}.\]
	The constant is independent of $r$: one applies the scalar Sobolev inequality to each component and sums the squares. Therefore,
	\begin{align}\label{Eq.L^infty_est} \widetilde u_r\leq C\sum_{q=0}^{s_d}\varrho^{-q}\frac{(r+q)!}{r!(1+r+q)^{s_*}}u_{r+q}\leq\frac{C_\varrho}{1+r}\sum_{q=0}^{s_d}u_{r+q},
	\end{align}
	where we used $s_*=s_d+1$ and
	\[\frac{(r+q)!}{r!(1+r+q)^{s_*}}\leq(1+r+q)^{q-s_*}\leq\frac1{1+r},\qquad 0\leq q\leq s_d.\]
	By Cauchy--Schwarz,
	\begin{align*}
		\sum_{r=0}^\infty\widetilde u_r &\leq C_\varrho\sum_{q=0}^{s_d}\sum_{r=0}^\infty\frac{u_{r+q}}{1+r}\leq C_\varrho\sum_{q=0}^{s_d} \left(\sum_{r=0}^\infty\frac1{(1+r)^2}\right)^{1/2} \left(\sum_{r=0}^\infty u_{r+q}^2\right)^{1/2}\leq C_\varrho\|u\|_{\mathcal B_\varrho}.
	\end{align*}
	The same argument gives $\|(\widetilde v_r)_{r\geq 0}\|_{\ell^1}\leq C_\varrho\|v\|_{\mathcal B_\varrho}$.
	
	Young's inequality finally implies that
	\begin{align*}
		\|uv\|_{\mathcal B_\varrho}&=\left\|\left(\frac{(1+N)^{s_*}\varrho^N}{N!}\|\nabla^N(uv)\|_{L^2}\right)_{N\geq 0}\right\|_{\ell^2}\\
		&\leq 2^{s_*}\left(\left\|(\widetilde u_r)_{r\geq 0}*(v_r)_{r\geq 0}\right\|_{\ell^2}+\left\|(\widetilde v_r)_{r\geq 0}*(u_r)_{r\geq 0}\right\|_{\ell^2}\right)\\
		&\leq 2^{s_*}\left(\left\|(\widetilde u_r)_{r\geq 0}\right\|_{\ell^1}\left\|(v_r)_{r\geq 0}\right\|_{\ell^2}+\left\|(\widetilde v_r)_{r\geq 0}\right\|_{\ell^1}\left\|(u_r)_{r\geq 0}\right\|_{\ell^2}\right)\\
		&\leq C\|u\|_{\mathcal B_\varrho}\|v\|_{\mathcal B_\varrho}.
	\end{align*}
	The constants are locally uniform for $\varrho\in(0,\infty)$.

	Next we prove the multiplier estimate. For $N\geq0$, the ordered Leibniz formula and $(1+N)\leq(1+r)(1+N-r)$ give
	\begin{align*}
		&\frac{(1+N)^{s_*}\varrho^N}{N!}\|\nabla^N(aU)\|_{L^2}\leq\sum_{r=0}^N\frac{(1+N)^{s_*}\varrho^N}{r!(N-r)!}\|\nabla^ra\|_{L^\infty}\|\nabla^{N-r}U\|_{L^2}\\
		&\leq\sum_{r=0}^N\left((1+r)^{s_*}\frac{\varrho^r}{r!}\|\nabla^ra\|_{L^\infty}\right)\left((1+N-r)^{s_*}\frac{\varrho^{N-r}}{(N-r)!}\|\nabla^{N-r}U\|_{L^2}\right).
	\end{align*}
	The first sequence in the last convolution has $\ell^1$ norm $\mathfrak M_\varrho(a)$, and the second has $\ell^2$ norm $\|U\|_{\mathcal B_\varrho}$. Young's inequality proves
	\[\|aU\|_{\mathcal B_\varrho}\leq\mathfrak M_\varrho(a)\|U\|_{\mathcal B_\varrho}.\]
	For a matrix-valued coefficient, the same calculation uses the Euclidean operator norm of each component derivative of $a$, as specified in the notation conventions. 
	
	Finally, let $a$ be scalar-valued and let $b$ take values in $E$. The ordered Leibniz formula and the Euclidean tensor norm give
	\[\|\nabla^N(ab)\|_{L^\infty}\leq\sum_{j=0}^N\binom Nj\|\nabla^ja\|_{L^\infty}\|\nabla^{N-j}b\|_{L^\infty},\qquad N\in\Z_{\geq0}.\]
	Since $1+N\leq(1+j)(1+N-j)$, multiplying by the analytic weights and rearranging the nonnegative double series yield
	\begin{align*}
		\mathfrak M_\varrho(ab)&\leq\sum_{N=0}^\infty\sum_{j=0}^N\frac{(1+N)^{s_*}\varrho^N}{j!(N-j)!}\|\nabla^ja\|_{L^\infty}\|\nabla^{N-j}b\|_{L^\infty}\\
		&\leq\sum_{j=0}^\infty\sum_{k=0}^\infty\left((1+j)^{s_*}\frac{\varrho^j}{j!}\|\nabla^ja\|_{L^\infty}\right)\left((1+k)^{s_*}\frac{\varrho^k}{k!}\|\nabla^kb\|_{L^\infty}\right)=\mathfrak M_\varrho(a)\mathfrak M_\varrho(b).
	\end{align*}
	This proves \eqref{Eq.multiplier_algebra}. The proof remains exactly the same for matrix-valued $a,b$.
\end{proof}

\begin{lemma}\label{Lem.analytic_compactness}
	For $0<\varrho<\varrho'$, the inclusion $\mathcal B_{\varrho'}(\Omega;E)\hookrightarrow\mathcal B_\varrho(\Omega;E)$ is compact, and
	\[\|h\|_{\mathcal B_\varrho}\leq\|h\|_{\mathcal B_{\varrho'}}\qquad(h\in\mathcal B_{\varrho'}).\]
\end{lemma}

\begin{proof}
	The continuous inclusion follows directly from the weights in \eqref{Eq.analytic_norm}. Let $(h_n)$ be a sequence satisfying $\sup_{n\geq1}\|h_n\|_{\mathcal B_{\varrho'}}\leq M$. For every fixed $N\geq0$, the definition gives
	\begin{align*}
		\|\nabla^Nh_n\|_{H^1(\Omega)}^2 &=\|\nabla^Nh_n\|_{L^2}^2+\|\nabla^{N+1}h_n\|_{L^2}^2\\
		&\leq M^2\left[ \frac{(N!)^2}{(1+N)^{2s_*}(\varrho')^{2N}} +\frac{((N+1)!)^2}{(2+N)^{2s_*}(\varrho')^{2N+2}} \right].
	\end{align*}
	Rellich's theorem, applied to each fixed finite-dimensional derivative tensor, and a diagonal extraction give a subsequence for which $\nabla^Nh_n$ converges in $L^2$ for every $N\geq0$. We retain the same notation for this subsequence. For every integer $K\geq0$,
	\begin{align*}
		\sum_{N>K}(1+N)^{2s_*}\frac{\varrho^{2N}}{(N!)^2} \|\nabla^N(h_n-h_j)\|_{L^2}^2\leq \left(\frac{\varrho}{\varrho'}\right)^{2(K+1)} \|h_n-h_j\|_{\mathcal B_{\varrho'}}^2\leq4M^2\left(\frac{\varrho}{\varrho'}\right)^{2(K+1)}.
	\end{align*}
	The finitely many terms with $N\leq K$ tend to zero as $n,j\to\infty$. Hence
	\[\limsup_{n,j\to\infty}\|h_n-h_j\|_{\mathcal B_\varrho}^2\leq4M^2\left(\frac{\varrho}{\varrho'}\right)^{2(K+1)}.\]
	Letting $K\to\infty$ shows that $(h_n)$ is Cauchy in $\mathcal B_\varrho$. Completeness of $\mathcal B_\varrho$ gives convergence.
\end{proof}

\begin{lemma}\label{Lem.analytic_density}
	For $\varrho>0$, $h\in\mathcal B_\varrho(\Omega;E)$ and $0<a<1$, define $h_a(x):=h(ax)$. Then
	\[h_a\in\mathcal B_{\varrho/a},\qquad\|h_a\|_{\mathcal B_{\varrho/a}}\leq a^{-d/2}\|h\|_{\mathcal B_\varrho},\qquad\lim_{a\uparrow1}\|h_a-h\|_{\mathcal B_\varrho}=0.\]
	In particular,
	\[\overline{\bigcup_{\varrho'>\varrho}\mathcal B_{\varrho'}(\Omega;E)}^{\,\mathcal B_\varrho(\Omega;E)}=\mathcal B_\varrho(\Omega;E).\]
\end{lemma}

\begin{proof}
	Sobolev embedding at every finite order gives continuous boundary values for all derivatives of $h$. Since $a\Omega\subset\Omega$ for $0<a<1$, we have
	\begin{align*}
		\nabla^Nh_a(x)=a^N(\nabla^Nh)(ax),\quad \|\nabla^Nh_a\|_{L^2(\Omega)}^2=a^{2N-d}\int_{a\Omega}|\nabla^Nh(y)|^2\,\mathrm dy\leq a^{2N-d}\|\nabla^Nh\|_{L^2(\Omega)}^2.
	\end{align*}
	Consequently,
	\begin{align*}
		\|h_a\|_{\mathcal B_{\varrho/a}}^2 &=a^{-d}\sum_{N=0}^\infty(1+N)^{2s_*} \frac{\varrho^{2N}}{(N!)^2} \int_{a\Omega}|\nabla^Nh(y)|^2\,\mathrm dy\leq a^{-d}\|h\|_{\mathcal B_\varrho}^2.
	\end{align*}
	This proves $h_a\in\mathcal B_{\varrho/a}$ and $\|h_a\|_{\mathcal B_{\varrho/a}}\leq a^{-d/2}\|h\|_{\mathcal B_\varrho}$.
	
	To prove convergence in $\mathcal B_\varrho$, we first fix $N$. The fundamental theorem of calculus gives
	\begin{align*}
		\nabla^Nh_a(x)-\nabla^Nh(x)=(a^N-1)\nabla^Nh(x)-a^N\int_a^1\sum_{l=1}^d x_l(\partial_l\nabla^Nh)(tx)\,\mathrm dt.
	\end{align*}
	For $1/2\leq a<1$, the inequalities $1-a^N\leq N(1-a)$ and $|x|\leq X_1$ imply
	\begin{align*}
		\|\nabla^N(h_a-h)\|_{L^2} &\leq N(1-a)\|\nabla^Nh\|_{L^2} +X_1\int_a^1t^{-d/2}\|\nabla^{N+1}h\|_{L^2}\,\mathrm dt\\
		&\leq(1-a)\left[ N\|\nabla^Nh\|_{L^2} +2^{d/2}X_1\|\nabla^{N+1}h\|_{L^2} \right].
	\end{align*}
	Thus $\nabla^Nh_a\to\nabla^Nh$ in $L^2$ for each fixed $N$. On the other hand, the dilation estimate gives the uniform bound
	\[\|\nabla^N(h_a-h)\|_{L^2}^2\leq2(a^{2N-d}+1)\|\nabla^Nh\|_{L^2}^2\leq2(2^d+1)\|\nabla^Nh\|_{L^2}^2.\]
	Therefore, for every $K\geq0$,
	\begin{align*}
		\sum_{N>K}(1+N)^{2s_*}\frac{\varrho^{2N}}{(N!)^2} \|\nabla^N(h_a-h)\|_{L^2}^2\leq2(2^d+1) \sum_{N>K}(1+N)^{2s_*}\frac{\varrho^{2N}}{(N!)^2} \|\nabla^Nh\|_{L^2}^2.
	\end{align*}
	The right-hand side tends to zero as $K\to\infty$, uniformly for $1/2\leq a<1$. Combining this tail estimate with convergence at the finitely many levels $N\leq K$ proves $\lim_{a\uparrow1}\|h_a-h\|_{\mathcal B_\varrho}=0$. Since $\varrho/a>\varrho$, the density follows.
\end{proof}

\begin{lemma}\label{Lem.analytic_extension}
	Let $\varrho>0$, and let $E$ be complex. Every $h\in\mathcal B_\varrho$ has a unique holomorphic extension $\widetilde h\in\Hol(\Omega_\varrho^\C;E)$ agreeing with its smooth representative on $\overline\Omega$, where
	\[\Omega_\varrho^\C:=\{\zeta\in\C^d:\dist(\zeta,\overline\Omega)<\varrho\}.\]
	For every integer $m\geq0$ and $0<r<\varrho$,
	\begin{equation}\label{Eq.derivative_holomorphic_extension}
		\sup_{\zeta\in\Omega_r^\C}|\nabla_\zeta^m\widetilde h(\zeta)| \leq C_{\Omega,\varrho}\frac{m!\varrho}{(\varrho-r)^{m+1}}\|h\|_{\mathcal B_\varrho}.
	\end{equation}
	Here $\nabla_\zeta^m$ denotes the full ordered tensor of holomorphic derivatives, and $\nabla_\zeta^0\widetilde h=\widetilde h$. The constant $C_{\Omega,\varrho}$ is independent of $m,r,h$ and $\dim E$, and is locally uniform for $\varrho\in(0,\infty)$.
\end{lemma}

\begin{proof}
	We first obtain a derivative bound with the same analytic radius $\varrho$. Similar to \eqref{Eq.L^infty_est}, we have 
	\begin{align}\label{Eq.L^infty_est2}
		\frac{\varrho^N}{N!}\|\nabla^Nh\|_{L^\infty}\leq C_\Omega\sum_{j=0}^{s_d}\varrho^{-j}\frac{(N+j)!}{N!(1+N+j)^{s_*}}\|h\|_{\mathcal B_\varrho}\leq\frac{C_{\Omega,\varrho}}{1+N}\|h\|_{\mathcal B_\varrho}.
	\end{align}
	The constant $C_{\Omega,\varrho}$ is independent of $N$ and $\dim E$, and can be chosen locally uniformly for $\varrho>0$. For $x\in\overline\Omega$ and $w\in\C^d$, we define the complex multilinear contraction
	\[\nabla^Nh(x)[w^{\otimes N}]:=\sum_{i_1,\ldots,i_N=1}^d(\partial_{i_1}\cdots\partial_{i_N}h)(x)w_{i_1}\cdots w_{i_N},\]
	with value $h(x)$ when $N=0$. No complex conjugation is taken in this contraction. Cauchy-Schwarz's inequality gives
	\begin{align*}
		|\nabla^Nh(x)[w^{\otimes N}]|\leq|\nabla^Nh(x)|\left(\sum_{i_1,\ldots,i_N=1}^d|w_{i_1}\cdots w_{i_N}|^2\right)^{1/2}=|\nabla^Nh(x)|\,|w|^N.
	\end{align*}
	For each $x\in\overline\Omega$, define
	\[H_x(\zeta):=\sum_{N=0}^\infty\frac1{N!}\nabla^Nh(x)[(\zeta-x)^{\otimes N}],\qquad |\zeta-x|<\varrho, \quad \zeta\in \C^d.\]
	For every $0<r<\varrho$, we have
	\begin{align*}
		\sum_{N=0}^\infty\sup_{|\zeta-x|\leq r}\frac{|\nabla^Nh(x)[(\zeta-x)^{\otimes N}]|}{N!}\leq C_{\Omega,\varrho}\|h\|_{\mathcal B_\varrho}\sum_{N=0}^\infty\frac{(r/\varrho)^N}{1+N}\leq\frac{C_{\Omega,\varrho}}{1-r/\varrho}\|h\|_{\mathcal B_\varrho}.
	\end{align*}
	The series therefore converges locally uniformly and defines an $E$-valued holomorphic function on the complex ball $|\zeta-x|<\varrho$.
	
	We verify that $H_x$ agrees with $h$ on $\Omega\cap\{y:|y-x|<\varrho\}$. For such a real point $y$, convexity of $\Omega$ allows Taylor's formula along the segment from $x$ to $y$, with boundary values understood by continuity. For every integer $K\geq0$,
	\begin{align*}
		h(y)-\sum_{N=0}^K\frac1{N!}\nabla^Nh(x)[(y-x)^{\otimes N}]=\frac1{K!}\int_0^1(1-t)^K \nabla^{K+1}h(x+t(y-x))[(y-x)^{\otimes(K+1)}]\,\mathrm dt.
	\end{align*}
	The derivative bound above yields
	\begin{align*}
		\left|h(y)-\sum_{N=0}^K\frac1{N!}\nabla^Nh(x)[(y-x)^{\otimes N}]\right|\leq\frac{|y-x|^{K+1}}{(K+1)!}\|\nabla^{K+1}h\|_{L^\infty}\leq C_{\Omega,\varrho}\|h\|_{\mathcal B_\varrho}\left(\frac{|y-x|}{\varrho}\right)^{K+1},
	\end{align*}
	which tends to $0$ as $K\to\infty$. Hence $H_x(y)=h(y)$.
	
	It remains to check agreement on overlapping complex balls. Let $x,y\in\overline\Omega$, and suppose that the balls $|\zeta-x|<\varrho$ and $|\zeta-y|<\varrho$ intersect. Then
	\[|x-y|<2\varrho,\qquad \left|\frac{x+y}{2}-x\right|=\left|\frac{x+y}{2}-y\right|=\frac{|x-y|}{2}<\varrho.\]
	The midpoint belongs to $\overline\Omega$. A real point of $\Omega$ sufficiently close to the midpoint therefore has a real open neighborhood contained in both balls and in $\Omega$. On this neighborhood, $H_x-H_y=0$. All real derivatives of $H_x-H_y$ vanish at any point of that neighborhood; since the difference is holomorphic, its complex Taylor coefficients vanish there as well. Thus $H_x=H_y$ on a complex neighborhood of that point. The intersection of the two complex balls is convex and hence connected, so the identity theorem gives agreement throughout the intersection.
	
	The functions $H_x$ consequently define a holomorphic function $\widetilde h$ on
	\[\bigcup_{x\in\overline\Omega}\{\zeta\in\C^d:|\zeta-x|<\varrho\}=\Omega_\varrho^\C.\]
	It agrees with $h$ on $\Omega$ and, by continuity, on $\overline\Omega$. The uniqueness of such $\widetilde h$ follows directly from the identity theorem on the connected domain $\Omega_\varrho^\C$.
	
	We now estimate the full ordered tensor of holomorphic derivatives. Fix an integer $m\geq0$ and an ordered index $J=(j_1,\ldots,j_m)\in\{1,\ldots,d\}^m$. We write $\partial_{\zeta,J}:=\partial_{\zeta_{j_1}}\cdots\partial_{\zeta_{j_m}}$ and $w:=\zeta-x$. The empty index gives the identity operator. For $k\geq m$, differentiating the homogeneous polynomial of degree $k$ gives
	\begin{align*}
		&\partial_{\zeta,J}\left(\frac1{k!}\sum_{i_1,\ldots,i_k=1}^d(\partial_{i_1}\cdots\partial_{i_k}h)(x)w_{i_1}\cdots w_{i_k}\right)\\
		&\quad=\frac1{(k-m)!}\sum_{i_1,\ldots,i_{k-m}=1}^d(\partial_{j_1}\cdots\partial_{j_m}\partial_{i_1}\cdots\partial_{i_{k-m}}h)(x)w_{i_1}\cdots w_{i_{k-m}}.
	\end{align*}
	The derivative is zero when $k<m$. Indeed, the $m$ derivatives must act on distinct factors of the monomial. There are $k(k-1)\cdots(k-m+1)=k!/(k-m)!$ ordered choices, and symmetry of the derivatives of $h$ gives the same remaining sum for each choice. This calculation remains valid when some entries of $J$ coincide.

	For $N\geq0$, let $I=(i_1,\ldots,i_N)\in\{1,\ldots,d\}^N$, write $w_I:=w_{i_1}\cdots w_{i_N}$, and let $\partial_{(J,I)}$ denote the real derivative with the concatenated ordered index $(J,I)$. The empty product is $1$. Cauchy--Schwarz in $I$, followed by summation over $J$, yields
	\begin{align*}
		&\sum_{J\in\{1,\ldots,d\}^m}\left|\sum_{I\in\{1,\ldots,d\}^N}(\partial_{(J,I)}h)(x)w_I\right|^2\\
		&\quad\leq\sum_{J\in\{1,\ldots,d\}^m}\left(\sum_{I\in\{1,\ldots,d\}^N}|(\partial_{(J,I)}h)(x)|^2\right)\left(\sum_{I\in\{1,\ldots,d\}^N}|w_I|^2\right)=|\nabla^{N+m}h(x)|^2|w|^{2N}.
	\end{align*}
	The last equality uses the full ordered derivative convention and the identity $\sum_I|w_I|^2=(\sum_{i=1}^d|w_i|^2)^N$. 

	For $0\leq t<1$, differentiating the geometric series $m$ times and reindexing gives
	\[\sum_{N=0}^\infty\frac{(N+m)!}{N!}t^N=\frac{\mathrm d^m}{\mathrm dt^m}\left(\sum_{k=0}^\infty t^k\right)=\frac{m!}{(1-t)^{m+1}}.\]
	The derivative estimate \eqref{Eq.L^infty_est2} therefore implies, for every $0<r<\varrho$,
	\begin{align*}
		\sum_{N=0}^\infty\frac{r^N}{N!}|\nabla^{N+m}h(x)|\leq C_{\Omega,\varrho}\|h\|_{\mathcal B_\varrho}\varrho^{-m}\sum_{N=0}^\infty\frac{(N+m)!}{N!(1+N+m)}\left(\frac r\varrho\right)^N\leq\frac{C_{\Omega,\varrho}m!\varrho}{(\varrho-r)^{m+1}}\|h\|_{\mathcal B_\varrho}.
	\end{align*}
	Hence the series obtained by differentiating $m$ times converges uniformly on $|\zeta-x|\leq r$, including the full ordered tensor norm. Applying this conclusion at every derivative order up to $m$ justifies termwise differentiation of the Taylor series and gives
	\[\partial_{\zeta,J}H_x(\zeta)=\sum_{N=0}^\infty\frac1{N!}\sum_{I\in\{1,\ldots,d\}^N}(\partial_{(J,I)}h)(x)(\zeta-x)_I,\qquad |\zeta-x|<\varrho.\]
	The triangle inequality in the ordered tensor norm and the contraction estimate above now yield
	\begin{align*}
		|\nabla_\zeta^mH_x(\zeta)|\leq\sum_{N=0}^\infty\frac{|\zeta-x|^N}{N!}|\nabla^{N+m}h(x)|\leq\frac{C_{\Omega,\varrho}m!\varrho}{(\varrho-|\zeta-x|)^{m+1}}\|h\|_{\mathcal B_\varrho}.
	\end{align*}
	For every $\zeta\in\Omega_r^\C$, choose $x\in\overline\Omega$ with $|\zeta-x|<r$. Since $\widetilde h=H_x$ on $|\zeta-x|<\varrho$, the preceding bound gives \eqref{Eq.derivative_holomorphic_extension}. The constant $C_{\Omega,\varrho}$ is the one in the real derivative estimate and is independent of $m,r,h$ and $\dim E$.
\end{proof}

\section*{Acknowledgments}
F. Shao is supported by China National Postdoctoral Program for Innovative Talents under Grant BX2026005. Z. Zhang is partially supported by National Natural Science Foundation of China under Grant 12288101.

\section*{Disclosure of AI assistance}
The authors began investigating the problem studied in this paper in May 2025. The proof framework was developed by the authors, who also carried out the mathematical calculations underlying the proofs themselves. The initial drafts of the manuscript were written by the authors. During drafting and revision, the authors used AI tools, including GPT-5.6 and GPT-6, to varying extents. The scope and extent of the assistance differed between sections as the available tools evolved over the course of the project. The assistance was used mainly to improve grammar and exposition, clarify and supplement selected proof details, and identify possible errors. The authors take full responsibility for the accuracy, originality, and integrity of the entire paper.

\raggedbottom

\end{document}